\documentclass[letterpaper]{article} 
\usepackage{aaai23}  
\usepackage{times}  
\usepackage{helvet}  
\usepackage{courier}  
\usepackage[hyphens]{url}  
\usepackage{graphicx} 
\usepackage{natbib}  
\usepackage{caption} 
\usepackage{algorithm}
\usepackage{algorithmic}

\usepackage{newfloat}
\usepackage{listings}
\DeclareCaptionStyle{ruled}{labelfont=normalfont,labelsep=colon,strut=off} 
\floatstyle{ruled}
\newfloat{listing}{tb}{lst}{}
\floatname{listing}{Listing}
\usepackage{microtype}
\usepackage{booktabs} 

\usepackage{amsmath}
\usepackage{amssymb}
\usepackage{mathtools}
\usepackage{amsthm}

\usepackage[capitalize,noabbrev]{cleveref}

\theoremstyle{plain}
\newtheorem{theorem}{Theorem}[section]
\newtheorem{proposition}[theorem]{Proposition}
\newtheorem{lemma}[theorem]{Lemma}

\theoremstyle{definition}

\newtheorem{definition}[theorem]{Definition}
\newtheorem{assumption}[theorem]{Assumption}

\theoremstyle{remark}

\theoremstyle{example}

\usepackage[textsize=tiny]{todonotes}

\def\beq{\begin{eqnarray}}
\def\eeq{\end{eqnarray}}
\def\noi{\noindent}
\def\nn{\nonumber}
\def\la{\langle}
\def\ra{\rangle}
\def\bal#1\eal{\begin{align}#1\end{align}}

\def\etas{\boldsymbol{\eta}}
\def\x{\mathbf{x}}
\def\y{\mathbf{y}}
\def\z{\mathbf{z}}
\def\g{\mathbf{g}}

\def\d{\mathbf{d}}
\def\v{\mathbf{v}}
\def\u{\mathbf{u}}

\def\p{\mathbf{p}}

\def\C{\mathbf{C}}
\def\U{\mathbf{U}}
\def\G{\mathbf{G}}
\def\Q{\mathbf{Q}}
\def\O{\mathbf{O}}
\def\I{\mathbf{I}}
\def\A{\mathbf{A}}
\def\X{\mathbf{X}}

\def\lam{\boldsymbol{\lambda}}

\def\E{\mathbb{E}}

\def\QQ{\mathcal{Q}}
\def\F{{F}}

\def\K{\mathcal{K}}

\def\c{\mathbf{c}}
\def\ei{{e}_i}
\def\eit{{e}_{i^t}}
\usepackage{subcaption}
\def\objimghei{0.10\textheight}

\def\conenn{\textcolor[rgb]{0.9961,0,0}}
\newcommand{\cone}[1]{\textbf{\conenn{#1}}}
\def\ctwo{\textcolor[rgb]{0,0.7,0}}
\def\cthree{\textcolor[rgb]{0,0.5820,0.9}}
\definecolor{mygray}{gray}{.9}
\definecolor{mypink}{rgb}{.99,.91,.95}

\definecolor{mycyana}{RGB}{153,204,153}
\definecolor{mycyanb}{RGB}{255,153,102}
\definecolor{mycyanc}{RGB}{255,204,51}
\definecolor{mycyand}{RGB}{204,204,102}

\title{Coordinate Descent Methods for DC Minimization: \\Optimality Conditions and Global Convergence}

\author {
  Ganzhao Yuan
}
\affiliations{
    Peng Cheng Laboratory, China\\
    yuangzh@pcl.ac.cn
}

\usepackage{bibentry}

\begin{document}

\maketitle

\begin{abstract}
Difference-of-Convex (DC) minimization, referring to the problem of minimizing the difference of two convex functions, has been found rich applications in statistical learning and studied extensively for decades. However, existing methods are primarily based on multi-stage convex relaxation, only leading to weak optimality of critical points. This paper proposes a coordinate descent method for minimizing a class of DC functions based on sequential nonconvex approximation. Our approach iteratively solves a nonconvex one-dimensional subproblem globally, and it is guaranteed to converge to a coordinate-wise stationary point. We prove that this new optimality condition is always stronger than the standard critical point condition and directional point condition under a mild \textit{locally bounded nonconvexity assumption}. For comparisons, we also include a naive variant of coordinate descent methods based on sequential convex approximation in our study. When the objective function satisfies a \textit{globally bounded nonconvexity assumption} and \textit{Luo-Tseng error bound assumption}, coordinate descent methods achieve \textit{Q-linear} convergence rate. Also, for many applications of interest, we show that the nonconvex one-dimensional subproblem can be computed exactly and efficiently using a breakpoint searching method. Finally, we have conducted extensive experiments on several statistical learning tasks to show the superiority of our approach.


\end{abstract}

\section{Introduction}
This paper mainly focuses on the following DC minimization problem (`$\triangleq$' means define):
\beq \label{eq:main}
\bar{\x} \in \arg \min_{\x\in\mathbb{R}^n}~ {F}(\x)\triangleq f(\x) + h(\x)  - g(\x).~
\eeq
\noi Throughout this paper, we make the following assumptions on Problem (\ref{eq:main}). \textit{\textbf{(i)}} $f(\cdot)$ is convex and continuously differentiable, and its gradient is coordinate-wise Lipschitz continuous with constant $\c_i\geq0$ that \cite{nesterov2012efficiency,beck2013convergence}:
\beq \label{eq:f:Lipschitz}
f(\x + \eta \ei  ) \leq f(\x) +  \la  \nabla f(\x),~ \eta \ei ) + \frac{\c_i}{2} \|\eta \ei \|_2^2
\eeq
\noi $\forall \x,\eta, i=1,...,n$. Here $\c \in \mathbb{R}^n$, and $\ei \in \mathbb{R}^n$ is an indicator vector with one on the $i$-th entry and zero everywhere else. \textit{\textbf{(ii)}} $h(\cdot)$ is convex and coordinate-wise separable with $h(\x) = \sum_{i=1}^n h_i(\x_i)$. Typical examples of $h(\x)$ include the bound constrained function and the $\ell_1$ norm function. \textit{\textbf{(iii)}} $g(\cdot)$ is convex and its associated proximal operator:
\beq \label{prox:operator}
\min_{\eta \in \mathbb{R}} p(\eta) \triangleq \frac{a}{2}\eta^2 + b \eta + h_i(\x + \eta \ei ) -g(\x + \eta \ei ), \label{eq:prox}
\eeq
\noi can be computed exactly and efficiently for given $a \in \mathbb{R}_+$, $b\in \mathbb{R}$ and $i\in\{1,...,n\}$. We remark that $g(\cdot)$ is neither necessarily differentiable nor coordinate-wise separable, and typical examples of $g(\x)$ are the $\ell_p$ norm function $g(\x)=\|\A\x\|_p$ with $p=\{1,2,\infty\}$, the RELU function $g(\x) = \|\max(0,\A\x)\|_1$, and the top-$s$ norm function $g(\x)=\sum_{i=1}^s |\x_{[i]}|$. Here $\A \in \mathbb{R}^{m\times n}$ is an arbitrary given matrix and $\x_{[i]}$ denotes the $i$th largest component of $\x$ in magnitude. \textit{\textbf{(iv)}} $\F(\x)$ only takes finite values.







\textbf{DC programming.} DC Programming/minimization is an extension of convex maximization over a convex set \cite{tao1997convex,ThiD18a}. It is closely related to the concave-convex procedure and alternating minimization in the literature. The class of DC functions is very broad, and it includes many important classes of nonconvex functions, such as twice continuously differentiable function on compact convex set and multivariate polynomial functions \cite{AhmadiH18}. DC programs have been mainly considered in global optimization and some algorithms have been proposed to find global solutions to such problem \cite{horst1999dc,horst2013global}. Recent developments on DC programming primarily focus on designing local solution methods for some specific DC programming problems. For example, proximal bundle DC methods \cite{joki2017proximal}, double bundle DC methods \cite{joki2018double}, inertial proximal methods \cite{mainge2008convergence}, and enhanced proximal methods \cite{LuZ19DC} have been proposed. DC programming has been applied to solve a variety of statistical learning tasks, such as sparse PCA \cite{SriperumbudurTL07,beck2021dual}, variable selection \cite{GotohTT18,GongZLHY13}, single source localization \cite{BeckH20SIOPT}, positive-unlabeled learning \cite{kiryo2017positive,XuQLJY19}, and deep Boltzmann machines \cite{NitandaS17}.


\textbf{Coordinate descent methods.} Coordinate Descent (CD) is a popular method for solving large-scale optimization problems. Advantages of this method are that compared with the full gradient descent method, it enjoys faster convergence \cite{tseng2009coordinate,xu2013block}, avoids tricky parameters tuning, and allows for easy parallelization \cite{liu2015asynchronous}. It has been well studied for convex optimization such as Lasso \cite{tseng2009coordinate}, support vector machines \cite{hsieh2008dual}, nonnegative matrix factorization \cite{hsieh2011fast}, and the PageRank problem \cite{nesterov2012efficiency}. Its convergence and worst-case complexity are well investigated for different coordinate selection rules such as cyclic rule \cite{beck2013convergence}, greedy rule \cite{hsieh2011fast}, and random rule \cite{lu2015complexity,RichtarikTakac2014}. It has been extended to solve many nonconvex problems such as penalized regression \cite{breheny2011coordinate,deng2020efficiency}, eigenvalue complementarity problem \cite{patrascu2015efficient}, $\ell_0$ norm minimization \cite{BeckE13,yuan2020block}, resource allocation problem \cite{necoara2013random}, leading eigenvector computation \cite{li2019coordinatewise}, and sparse phase retrieval \cite{shechtman2014gespar}.

\textbf{Iterative majorization minimization.} Iterative majorization / upper-bound minimization is becoming a standard principle in developing nonlinear optimization algorithms. Many surrogate functions such as Lipschitz gradient surrogate, proximal gradient surrogate, DC programming surrogate, variational surrogate, saddle point surrogate, Jensen surrogate, quadratic surrogate, cubic surrogate have been considered, see \cite{Mairal13,RazaviyaynHL13}. Recent work extends this principle to the coordinate update, incremental update, and stochastic update settings. However, all the previous methods are mainly based on multiple-stage convex relaxation, only leading to weak optimality of critical points. In contrast, our method makes good use of sequential nonconvex approximation to find stronger stationary points. Thanks to the coordinate update strategy, we can solve the one-dimensional nonconvex subproblem \textit{globally} by using a novel exhaustive breakpoint searching method even when $g(\cdot)$ is \textit{nonseparable} and \textit{non-differentiable}.

\textbf{Theory for nonconvex optimization.} We pay specific attention to two contrasting approaches on the theory for nonconvex optimization. \textit{\textbf{(i)}} Strong optimality. The first approach is to achieve stronger optimality guarantees for nonconvex problems. For smooth optimization, canonical gradient methods only converge to a first-order stationary point, recent works aim at finding a second-order stationary point \cite{Jin0NKJ17}. For cardinality minimization, the work of \cite{BeckE13,yuan2020block} introduces a new optimality condition of (block) coordinate stationary point which is stronger than that of the Lipschitz stationary point \cite{YuanLZ17}. \textit{\textbf{(ii)}} Strong convergence. The second approach is to provide convergence analysis for nonconvex problems. The work of \cite{Jin0NKJ17} establishes a global convergence rate for nonconvex matrix factorization using a regularity condition. The work of \cite{AttouchBRS10} establishes the convergence rate for general nonsmooth problems by imposing Kurdyka-\L{}ojasiewicz inequality assumption of the objective function. The work of \cite{Dong2021,YueZS19} establish linear convergence rates under the \textit{Luo-Tseng error bound assumption}. Inspired by these works, we prove that the proposed CD method has strong optimality guarantees and convergence guarantees.


\textbf{Contributions.} The contributions of this paper are as follows: \textit{\textbf{(i)}} We propose a new CD method for minimizing DC functions based on sequential nonconvex approximation (See Section \ref{sect:alg}). \textit{\textbf{(ii)}} We prove that our method converge to a coordinate-wise stationary point, which is always stronger than the optimality of standard critical points and directional points when the objective function satisfies a \textit{locally bounded nonconvexity assumption}. When the objective function satisfies a \textit{globally bounded nonconvexity assumption} and \textit{Luo-Tseng error bound assumption}, CD methods achieve \emph{Q-linear} convergence rate (See Section \ref{sect:theory}). \textit{\textbf{(iii)}} We show that, for many applications of interest, the one-dimensional subproblem can be computed exactly and efficiently using a breakpoint searching method (See Section \ref{sect:breakpoint}). \textit{\textbf{(iv)}} We have conducted extensive experiments on some statistical learning tasks to show the superiority of our approach (See Section \ref{sect:exp}). \textit{\textbf{(v)}} We also provide several important discussions of the proposed method (See Section \ref{sect:disc:ext} in the Appendix).



\textbf{Notations.} Vectors are denoted by boldface lowercase letters, and matrices by boldface uppercase letters. The Euclidean inner product between $\x$ and $\y$ is denoted by $\la \x,\y \ra$ or $\x^T\y$. We denote $\|\x\| = \|\x\|_2 = \sqrt{\la \x,\x\ra}$. $\x_i$ denotes the $i$-th element of the vector $\x$. $\E[\cdot]$ represents the expectation of a random variable. $\odot$ and $\div$ denote the element-wise multiplication and division between two vectors, respectively. For any extended real-valued function $h:\mathbb{R}^{n}\rightarrow(-\infty,+\infty]$, the set of all subgradients of $h$ at $\x$ is defined as $\partial h(\x) =\{ \g \in \mathbb{R}^n: h(\y)\geq h(\x)+\la \g,\y-\x \ra\}$, the conjugate of $h(\x)$ is defined as $h^*(\x)\triangleq \max_{\y}\{ \la \y,\x \ra - h(\y)\}$, and $(\partial h(\x))_i$ denotes the subgradient of $h(\x)$ at $\x$ for the $i$-th componnet. $\text{diag}(\c)$ is a diagonal matrix with $\c$ as the main diagonal entries. We define $\|\d\|_{\c}^2 =  \sum_{i} \c_i\d_i^2$. $\text{sign}(\cdot)$ is the signum function. $\I$ is the identity matrix of suitable size. The directional derivative of $F(\cdot)$ at a point $\x$ in its domain along a direction $\d$ is defined as: $\F'(\x;\d) \triangleq \lim_{t \downarrow  0} \frac{1}{t} (\F(\x + t \d) - \F(\x))$. $\text{dist}(\Omega,\Omega') \triangleq \inf_{\v \in \Omega,\v' \in \Omega'}\| \v - \v'\|$ denotes the distance between two sets.


\section{Motivating Applications} \label{sect:app}
A number of statistical learning models can be formulated as Problem (\ref{eq:main}), which we present some instances below.


\noi $\bullet$ \textbf{Application I: $\ell_p$ Norm Generalized Eigenvalue Problem}. Given arbitrary data matrices $\G \in \mathbb{R}^{m\times n}$ and $\Q \in \mathbb{R}^{n\times n}$ with $\Q \succ \mathbf{0}$, it aims at solving the following problem:
\beq \label{eq:pca0}
\bar{\v} \in \arg \max_{\v} \|\G\v\|_p,~s.t.~\v^T\Q\v = 1.
\eeq with $p\geq 1$. Using the Lagrangian dual, we have the following equivalent unconstrained problem:
\beq \label{eq:Lpnorm:PCA}
\bar{\x} \in \arg \min_{\x}~ \tfrac{\alpha}{2} \x^T \Q\x-\|\G\x\|_p,
\eeq
\noi for any given $\alpha>0$. The optimal solution to Problem (\ref{eq:pca0}) can be recovered as $\bar{\v}=\pm \bar{\x} \cdot (\bar{\x}^T\Q\bar{\x})^{-\frac{1}{2}}$. Refer to Section \ref{sect:disc:ref} in the appendix for a detailed discussion.

\noi$\bullet$ \textbf{Application II: Approximate Sparse/Binary Optimization}. Given a channel matrix $\G \in \mathbb{R}^{m\times n}$, a structured signal $\x$ is transmitted through a communication channel, and received
as $\y = \G \x + \v$, where $\v$ is the Gaussian noise. If $\x$ has $s$-sparse or binary structure, one can recover $\x$ by solving the following optimization problem \cite{GotohTT18,Forney72}:
\beq
\textstyle\min_{\x}~\tfrac{1}{2}\|\G\x - \y\|_2^2,~s.t.~\|\x\|_0 \leq s, \nn \\
\textstyle\text{or}~~~\min_{\x}~\tfrac{1}{2}\|\G\x - \y\|_2^2,~s.t.~\x \in \{-1+1\}^n. \nn
\eeq
\noi Here, $\|\cdot\|_{0}$ is the number of non-zero components. Using the equivalent variational reformulation of the $\ell_0$ (pseudo) norm $\|\x\|_0\leq s \Leftrightarrow \|\x\|_1 = \sum_{i=1}^s | \x_{[i]}|$ and the binary constraint $\{-1,+1\}^n \Leftrightarrow \{\x|-\mathbf{1}\leq\x\leq \mathbf{1},\|\x\|_2^2 = n\}$, one can solve the following approximate sparse/binary optimization problem \cite{GotohTT18,YuanG17}: 
\beq
\textstyle \min_{\x}~\tfrac{1}{2}\|\G\x - \y\|_2^2 + \rho ( \|\x\|_1 -  \sum_{i=1}^s | \x_{[i]}| )\label{eq:sparse:opt} \\
\textstyle \min_{\|\x\|_{\infty} \leq 1}~\tfrac{1}{2}\|\G\x - \y\|_2^2 + \rho (\sqrt{n}-\|\x\|). \label{eq:binary:opt}
\eeq


\noi $\bullet$ \textbf{Application III: Generalized Linear Regression}. Given a sensing matrix $\G \in \mathbb{R}^{m\times n}$ and measurements $\y\in \mathbb{R}^{m}$, it deals with the problem of recovering a signal $\x$ by solving $\bar{\x} = \arg \min_{\x\in\mathbb{R}^{n}} \tfrac{1}{2}\|\sigma(\G\x)-\y\|_2^2$. When $\sigma(\z) = \max(0,\z)$ or $\sigma(\z) = |\z|$, this problem reduces to the one-hidden-layer ReLU networks \cite{zhang2019learning} or the amplitude-base phase retrieval problem \cite{CandesLS15}. When $\y\geq \mathbf{0}$, we have the following equivalent DC program:
\beq
\min_{\x\in\mathbb{R}^n}~  \tfrac{1}{2}\|\sigma(\G\x)\|_2^2 - \la \mathbf{1}, \sigma(\text{diag}(\y)\G)\x)\ra + \tfrac{1}{2}\|\y\|_2^2. \label{eq:oneNN}
\eeq

\section{Related Work} \label{sect:related}



We now present some related DC minimization algorithms. 

\textit{\textbf{(i)}} Multi-Stage Convex Relaxation (MSCR)\cite{zhang2010analysis,bi2014exact}. It solves Problem (\ref{eq:main}) by generating a sequence $\{\x^{t}\}$ as:
\beq \label{eq:dc:mscr}
\x^{t+1} \in \arg \min_{\x}~f(\x) + h(\x) - \la \x - \x^t,~ \g^t \ra
\eeq
\noi where $\g^t \in \partial g(\x^t)$. Note that Problem (\ref{eq:dc:mscr}) is convex and can be solved via standard proximal gradient method. The computational cost of MSCR could be expensive for large-scale problems, since it is $K$ times that of solving Problem (\ref{eq:dc:mscr}) with $K$ being the number of outer iterations.

\textit{\textbf{(ii)}} Proximal DC algorithm (PDCA) \cite{GotohTT18}. To alleviate the computational issue of solving Problem (\ref{eq:dc:mscr}), PDCA exploits the structure of $f(\cdot)$ and solves Problem (\ref{eq:main}) by generating a sequence $\{\x^t\}$ as:
\beq \label{eq:conv:sub2}
\x^{t+1} = \arg \min_{\x}~\QQ(\x,\x^t)+ h(\x) - \la \x - \x^t,~ \g^t \ra \nn
\eeq
\noi where $\QQ(\x,\x^t) \triangleq f(\x^t) + \la  \nabla f(\x^t),~\x - \x^t \ra + \frac{L}{2} \|\x - \x^t\|_2^2$, and $L$ is the Lipschitz constant of $\nabla f(\cdot)$.

\textit{\textbf{(iii)}} Toland's duality method \cite{toland1979duality,beck2021dual}. Assuming $g(\x)$ has the following structure $g(\x) = \bar{g}(\A\x) = \max_{\y}\{ \la \A\x,\y \ra - \bar{g}^*(\y) \}$. This approach rewrites Problem (\ref{eq:main}) as the following equivalent problem using the conjugate of $g(\x)$: $\min_{\x} \min_{\y}~ f(\x) + h(\x) - \la \y,\A\x \ra + \bar{g}^*(\y)$. Exchanging the order of minimization yields the equivalent problem: $\min_{\y}\min_{\x} f(\x) + h(\x) - \la \y,\A\x \ra + \bar{g}^*(\y)$. The set of minimizers of the inner problem with respect to $\x$ is $ \partial h^*(\A^T\y) + \nabla f^*(\A^T\y)$, and the minimal value is $-f^*(\A^T\y) - h^*(\A^T\y)+ \bar{g}^*(\y)$.
We have the Toland-dual problem which is also a DC program:
\beq \label{eq:toland:dual}
\min_{\y} \bar{g}^*(\y) - f^*(\A^T\y)  -h^*(\A^T\y)
\eeq
\noi This method is only applicable when the minimization problem with respect to $\x$ is simple so that it has an analytical solution. Toland's duality method could be useful if one of the subproblems is easier to solve than the other.

\textit{\textbf{(iv)}} Subgradient descent method (SubGrad). It uses the iteration $\x^{t+1} = \mathcal{P}(\x^t - \eta^t \g^t)$, where $\g^t \in \partial \F(\x^t)$, $\eta^t$ is the step size, and $\mathcal{P}$ is the projection operation on some convex set. This method has received much attention recently due to its simplicity \cite{zhang2019learning,DavisDMP18,DavisG19,li2021weakly}.



\section{Coordinate Descent Methods for DC Minimization} \label{sect:alg}

This section presents a new Coordinate Descent (CD) method for solving Problem (\ref{algo:main}), which is based on Sequential NonConvex Approximation (SNCA). For comparisons, we also include a naive variant of CD methods based on Sequential Convex Approximation (SCA) in our study. These two methods are denoted as \textit{\textbf{CD-SNCA}} and \textit{\textbf{CD-SCA}}, respectively.

Coordinate descent is an iterative algorithm that sequentially minimizes the objective function along coordinate directions. In the $t$-th iteration, we minimize $\F(\cdot)$ with respect to the $i^t$ variable while keeping the remaining $(n-1)$ variables $\{\x^t_j\}_{j\neq i^t}$ fixed. This is equivalent to performing the following one-dimensional search along the $i^t$-th coordinate:
\beq \label{eq:pure:cd}
\bar{\eta}^t \in \arg \min_{\eta \in \mathbb{R}}~f(\x^t + \eta \eit) + h(\x^t + \eta \eit) -g(\x^t + \eta \eit).\nn
\eeq
\noi Then $\x^t$ is updated via: $\x^{t+1} = \x^t + \bar{\eta}^t \eit$. However, the one-dimensional problem above could be still hard to solve when $f(\cdot)$ and/or $g(\cdot)$ is complicated. One can consider replacing $f(\cdot)$ and $g(\cdot)$ with their majorization function:
\bal
&f(\x^t + \eta \eit) \leq \mathcal{S}_{i^t}(\x^t,\eta)  \nn \\
\text{with}~&\mathcal{S}_{i}(\x,\eta) \triangleq f(\x) +  \la  \nabla f(\x),~ \eta \ei \ra + \frac{\c_{i}}{2} \eta^2, \label{eq:upbound:1}  \\
&-g(\x^t + \eta \eit) \leq \mathcal{G}_{i^t}(\x^t,\eta) \nn \\
\text{with}~&\mathcal{G}_{i}(\x,\eta) \triangleq -g(\x) -  \la  \partial g(\x),~ (\x + \eta \ei) - \x  \ra. \label{eq:upbound:2}
\eal

\begin{algorithm}[!t]
\caption{ {\bf Coordinate Descent Methods for Minimizing DC functions using \textit{\textbf{SNCA}} or \textit{\textbf{SCA}} strategy.} }
\begin{algorithmic}
  \STATE Input: an initial feasible solution $\x^0$, $\theta>0$. Set $t=0$.
  \WHILE{not converge}
  \STATE (\textbf{S1}) Use some strategy to find a coordinate $i^t \in \{1,...,n\}$ for the $t$-th iteration.
  \STATE (\textbf{S2}) Solve the following nonconvex or convex subproblem globally and exactly.

$\bullet$ Option I: Sequential NonConvex Approximation (\textit{\textbf{SNCA}}) strategy.
   \beq \label{eq:subprob:nonconvex}
    \bar{\eta}^t \in \bar{\mathcal{M}}_{i^t}(\x^t)\triangleq  \arg \min_{\eta}~\mathcal{M}_{i^t}(\x^t,\eta) \\
 \text{with}~\mathcal{M}_{i}(\x,\eta)\triangleq  \mathcal{S}_{i}(\x,\eta)+ h_{i}(\x + \eta \ei)  \nn\\
	- g(\x + \eta \ei) +   \tfrac{\theta}{2} \| (\x+\eta \ei) - \x \|_2^2\nn
\eeq

  $\bullet$ Option II: Sequential Convex Approximation (\textit{\textbf{SCA}}) strategy.
  \beq \label{eq:subprob:convex}
    \bar{\eta}^t \in \bar{\mathcal{P}}_{i^t}(\x^t)\triangleq \arg \min_{\eta}~ \mathcal{P}_{i^t}(\x^t,\eta) \\
    \mathcal{P}_{i}(\x,\eta) \triangleq  \mathcal{S}_{i}(\x,\eta)  + h_{i}(\x + \eta \ei) \nn\\
	 + \mathcal{G}_{i}(\x,\eta)+   \tfrac{\theta}{2} \| (\x+\eta \ei) - \x\|_2^2\nn
\eeq

\STATE (\textbf{S3}) $\x^{t+1}  = \x^t  + \bar{\eta}^t \cdot \eit~~~ (\Leftrightarrow \x_{i^t}^{t+1}  = \x_{i^t}^t  + \bar{\eta}^t )$
  \STATE (\textbf{S4}) Increment $t$ by 1.
  \ENDWHILE
\end{algorithmic}\label{algo:main}
\end{algorithm}
\noi $\blacktriangleright$ \textbf{Choosing the Majorization Function}
\begin{enumerate}
  \item \textbf{Sequential NonConvex Approximation Strategy}. If we replace $f(\x^t + \eta \eit)$ with its upper bound $\mathcal{S}_{i^t}(\x^t,\eta)$ as in (\ref{eq:upbound:1}) while keep the remaining two terms unchanged, we have the resulting subproblem as in (\ref{eq:subprob:nonconvex}), which is a nonconvex problem. It reduces to the proximal operator computation as in (\ref{prox:operator}) with $a=\c_{i^t}+\theta$ and $b=\nabla_{i^t} f(\x^t)$. Setting the subgradient with respect to $\eta$ of the objective function in (\ref{eq:subprob:nonconvex}) to zero, we have the following \emph{necessary but not sufficient} optimality condition for (\ref{eq:subprob:nonconvex}):
\beq \label{eq:optimal:1}
  0 \in [\nabla f(\x^t)  + \partial h(\x^{t+1}) - \partial g(\x^{t+1})]_{i^t} + (\c_{i^t}+\theta) \bar{\eta}^t.\nn
\eeq

  \item \textbf{Sequential Convex Approximation Strategy}. If we replace $f(\x^t + \eta \eit)$ and $-g(\x^t + \eta \eit)$ with their respective upper bounds $\mathcal{S}_{i^t}(\x^t,\eta)$ and $\mathcal{G}_{i^t}(\x^t,\eta)$ as in (\ref{eq:upbound:1}) and (\ref{eq:upbound:2}), while keep the term $h(\x^t + \eta e_{i^t})$ unchanged, we have the resulting subproblem as in (\ref{eq:subprob:convex}), which is a convex problem. We have the following \emph{necessary and sufficient} optimality condition for (\ref{eq:subprob:convex}):
\beq \label{eq:optimal:2}
  0 \in [\nabla f(\x^t)  + \partial h(\x^{t+1}) - \partial g(\x^{t}) ]_{i^t}+ (\c_{i^t}+\theta) \bar{\eta}^t.\nn
\eeq

\end{enumerate}
\noi $\blacktriangleright$\textbf{Selecting the Coordinate to Update}

There are several fashions to decide which coordinate to update in the literature \cite{tseng2009coordinate}. {\textbf{(i) \textit{Random rule}}.} $i^t$ is randomly selected from $\{1,...,n\}$ with equal probability. {\textbf{(ii) \textit{Cyclic rule}}.} $i^t$ takes all coordinates in cyclic order $1 \rightarrow 2 \rightarrow ... \rightarrow n \rightarrow 1$. {\textbf{(iii) \textit{Greedy rule}}.} Assume that $\nabla f(\x)$ is Lipschitz continuous with constant $L$. The index $i^t$ is chosen as $i^t = \arg \max_{j}|\d^t_j|$
where $\d^t = \arg \min_{\d}  ~h(\x^t + \d) + \frac{ L}{2} \|\d\|_2^2 +  \la \nabla f(\x^t) - \partial g(\x^t)),\d\ra$. Note that $\d^t=\mathbf{0}$ implies that $\x^t$ is a critical point. 


\noi We summarize \textit{\textbf{CD-SNCA}} and \textit{\textbf{CD-SCA}} in Algorithm \ref{algo:main}.

\noi \textbf{Remarks}. \textit{\textbf{(i)}} We use a proximal term for the subproblems in (\ref{eq:subprob:nonconvex}) and (\ref{eq:subprob:convex}) with $\theta$ being the proximal point parameter. This is to guarantee sufficient descent condition and global convergence for Algorithm \ref{algo:main}. As can be seen in Theorem \ref{the:global:conv} and Theorem \ref{the:rate:ncvx}, the parameter $\theta$ is critical for \textit{\textbf{CD-SNCA}}. \textit{\textbf{(ii)}} Problem (\ref{eq:subprob:nonconvex}) can be viewed as \textit{globally} solving the following nonconvex problem which has a bilinear structure: $(\bar{\eta}^t,\y) = \arg \min_{\eta,\y}~\mathcal{S}_{i^t}(\x^t,\eta) +  \frac{\theta}{2} \eta^2 + h(\x^t + \eta \eit) - \la \y,\x^t + \eta \eit\ra  + g^*(\y)$. \textit{\textbf{(iii)}} While we apply CD to the primal, one may apply to the dual as in Problem (\ref{eq:toland:dual}). \textit{\textbf{(iv)}} The nonconvex majorization function used in \textit{\textbf{CD-SNCA}} is always a lower bound of the convex majorization function used in \textit{\textbf{CD-SCA}}, i.e., $\mathcal{M}_{i}(\x,\eta) \leq \mathcal{P}_{i}(\x,\eta),~\forall i, \x, \eta$.

\section{Theoretical Analysis}
\label{sect:theory}

This section provides a novel optimality analysis and a novel convergence analysis for Algorithm \ref{algo:main}. Due to space limit, all proofs are placed in Section \ref{sect:proof} in the appendix.

We introduce the following useful definition.

\begin{definition} \label{def:1}
(\textbf{Globally or Locally Bounded Nonconvexity}) A function $z(\x)$ is called to be globally $\rho$-bounded nonconvex if: $\forall \x,\y,~z(\x)   \leq z(\y)  + \la \x - \y,~ \partial z(\x) \ra + \frac{\rho}{2}\|\x - \y\|_2^2$ with $\rho<+\infty$. In particular, $z(\x)$ is locally $\rho$-bounded nonconvex if $\x$ is restricted to some point $\ddot{\x}$ with $\x = \ddot{\x}$.
\end{definition}

%

%
\noi \textbf{Remarks}. \textit{\textbf{(i)}} Globally $\rho$-bounded nonconvexity of $z(\x)$ is equivalent to $z(\x) + \frac{\rho}{2}\|\x\|_2^2$ is convex, and this notation is also referred as \textit{semi-convex}, \textit{approximate convex}, or \textit{weakly-convex} in the literature (cf. \cite{BohmW21,DavisDMP18,li2021weakly}). \textit{\textbf{(ii)}} Many nonconvex functions in the robust statistics literature are \textit{globally} $\rho$-bounded nonconvex, examples of which includes the \textit{minimax concave penalty}, the \textit{fractional penalty}, the \textit{smoothly clipped absolute deviation}, and the \textit{Cauchy loss} (c.f. \cite{BohmW21}). \textit{\textbf{(iii)}} Any globally $\rho$-bounded nonconvex function $z(\x)$ can be rewritten as a DC function that $z(\x) = \frac{\rho}{2}\|\x\|^2 - g(\x)$, where $g(\x)=\frac{\rho}{2}\|\x\|^2 - z(\x)$ is convex and $(-g(\x))$ is globally $(2\rho)$-bounded nonconvex.

Globally bounded nonconvexity could be a strong definition, one may use a weaker definition of locally bounded nonconvexity instead. The following lemma shows that some nonconvex functions are locally bounded nonconvex.

\begin{lemma} \label{lemma:bounde:nonconvex}
The function $z(\x) \triangleq -\|\x\|_p$ with $p\in [1,\infty)$ is concave and locally $\rho$-bounded nonconvex with $\rho<+\infty$.
\end{lemma}

\noi \textbf{Remarks}. By Lemma \ref{lemma:bounde:nonconvex}, we have that the functions $z(\x) = -\|\G\x\|_p$ in (\ref{eq:Lpnorm:PCA}) and $z(\x) = - \rho\|\x\|$ in (\ref{eq:binary:opt}) are locally $\rho$-bounded nonconvex. Using similar strategies, one can conclude that the functions $z(\x)=-  \sum_{i=1}^s | \x_{[i]}| $ and $z(\x) = - \la \mathbf{1}, \sigma(\text{diag}(\y)\G)\x)\ra  $ as in (\ref{eq:sparse:opt}) and (\ref{eq:oneNN}) are locally $\rho$-bounded nonconvex.


We assume that the random-coordinate selection rule is used. After $t$ iterations, Algorithm \ref{algo:main} generates a random output $\x^t$, which depends on the observed realization of the random variable: $\xi^{t-1} \triangleq \{i^0,~i^1,...,i^{t-1}\}$. 

We now develop the following technical lemma that will be used to analyze Algorithm \ref{algo:main} subsequently.

\begin{lemma} \label{lemma:tech}
\noi For any $\x \in \mathbb{R}^n$, $\d\in \mathbb{R}^n$, $\bar{\c}\in \mathbb{R}^n$, we define $h'(\x) \triangleq \sum_{i=1}^n h(\x + \d_i \ei) $, $g'(\x) \triangleq   \sum_{i=1}^n g(\x + \d_i \ei)$, and $f'(\x) \triangleq \sum_{i=1}^n f(\x + \d_i \ei)$. We have:
{
\beq
\textstyle \sum_{i=1}^n \| \x+\d_i \ei\|_{\bar{\c}}^2  = \|\x+\d\|_{\bar{\c}}^2 + (n-1) \|\x\|_{\bar{\c}}^2 \label{eq:tech:xx} \\
h'(\x) =   h(\x+\d) + (n- 1) h(\x)   \label{eq:tech:hh}  \\
f'(\x)  \leq  f(\x) + \la \nabla f(\x),\d  \ra  + \tfrac{1}{2}\|\d\|_{\c}^2 + (n-1)f(\x) \label{eq:tech:ff}\\
-g'(\x)  \leq   - g(\x)   - \la  \partial g(\x),\d\ra -(n-1)g(\x) \label{eq:tech:gg0} 
\eeq}
\end{lemma}

\subsection{Optimality Analysis}

We now provide an optimality analysis of our method. Since the coordinate-wise optimality condition is novel in this paper, we clarify its relations with existing optimality conditions formally.


\begin{definition}
(Critical Point) A solution $\check{\x}$ is called a critical point if \cite{toland1979duality}: $0 \in \nabla f(\check{\x}) + \partial h(\check{\x})  - \partial g(\check{\x})$.
\end{definition}

\noi \textbf{Remarks}. \textit{\textbf{(i)}} The expression above is equivalent to $(f(\check{\x}) + \partial h(\check{\x})) \cap \partial g(\check{\x}) \neq \emptyset $. The sub-differential is always non-empty on convex functions; that is why we assume that $\F(\cdot)$ can be repressed as the difference of two convex functions. \textit{\textbf{(ii)}} Existing methods such as MSCR, PDCA, and  SubGrad as shown in Section (\ref{sect:related}) are only guaranteed to find critical points of Problem (\ref{eq:main}).

\begin{definition}
(Directional Point) A solution $\grave{\x}$ is called a directional point if \cite{PangRA17}: $\F'(\grave{\x};\y-\grave{\x}) \geq 0,~\forall \y \in \text{dom}(\F) \triangleq \{\x: |\F(\x)|< + \infty \} $.

\end{definition}


\noi \textbf{Remarks}. The work of \cite{PangRA17} characterizes different types of stationary points, and proposes an enhanced DC algorithm that subsequently converges to a directional point. However, they only consider the case $g(\x)  = \max_{i \in I} g_i(\x)$ where each $g_i(\x)$ is continuously differentiable and convex and $I$ is a finite index set.


%
%

\begin{definition} \label{def:cwsp}
(Coordinate-Wise Stationary Point) A solution $\ddot{\x}$ is called a coordinate-wise stationary point if the following holds: $0 \in \arg \min_{\eta} \mathcal{M}_{i}(\ddot{\x},\eta)$ for all $i=1,...,n$, where $\mathcal{M}_{i}(\x,\eta)\triangleq f(\x) +  \la  \nabla f(\x),~ \eta \ei \ra + \frac{\c_{i}}{2} \eta^2 + h_{i}(\x + \eta \ei)  - g(\x + \eta \ei) +   \tfrac{\theta}{2}\eta^2$, and $\theta\geq0$ is a constant.

\end{definition}
\noi \textbf{Remarks}. \textbf{(i)} Coordinate-wise stationary point states that if we minimize the majorization function $\mathcal{M}_i(\x,\eta)$, we can not improve the objective function value for $\mathcal{M}_i(\x,\eta)$ for all $i \in \{1,...,n\}$. \textbf{(ii)} For any coordinate-wise stationary point $\ddot{\x}$, we have the following necessary but not sufficient condition: $\forall i\in \{1,...,n\}, ~0 \in \partial \mathcal{M}_{i}(\ddot{\x},\eta) \triangleq  (\c_i+\theta) \eta + [\nabla f(\ddot{\x}) + \partial h(\ddot{\x} +\eta \ei ) - \partial h(\ddot{\x} +\eta \ei)]_i$ with $\eta=0$, which coincides with the critical point condition. Therefore, any coordinate-wise stationary point is a critical point.

The following lemma reveals a \textit{quadratic growth condition} for any coordinate-wise stationary point.
\begin{lemma} \label{lemma:cw:point}
Let $\ddot{\x}$ be any coordinate-wise stationary point. Assume that $z(\x)\triangleq -g(\x)$ is locally $\rho$-bounded nonconvex at the point $\ddot{\x}$. We have: $\forall \d,~\F(\ddot{\x})  -\F(\ddot{\x}+\d) \leq \frac{1}{2} \|\d\|^2_{(\c+\theta+\rho)}$.
\end{lemma}

\noi \textbf{Remarks}. Recall that a solution $\dot{\x}$ is said to be a local minima if $\F({\dot{\x}}) \leq \F({\dot{\x}}+\d)$ for a sufficiently small constant $\delta$ that $\|\d\|\leq \delta$. The coordinate-wise optimality condition does not have any restriction on $\d$ with $\|\d\|\leq +\infty$. Thus, neither the optimality condition of coordinate-wise stationary point nor that of the local minima is stronger than the other.

We use $\check{\x}$, $\grave{\x}$, $\ddot{\x}$, and $\bar{\x}$ to denote any critical point, directional point, coordinate-wise stationary point, and optimal point, respectively. The following theorem establishes the relations between different types of stationary points list above.

\begin{theorem} \label{the:optimality}
\textbf{(Optimality Hierarchy between the Optimality Conditions).} Assume that the assumption made in Lemma \ref{lemma:cw:point} holds, we have: $\{\bar{\x}\} \overset{\textbf{(a)}}{ \subseteq} \{\ddot{\x}\} \overset{\textbf{(b)}}{ \subseteq} \{\grave{\x}\}\overset{\textbf{(c)}}{ \subseteq} \{\check{\x}\}$.

\end{theorem}
\noi \textbf{Remarks}. \textit{\textbf{(i)}} The coordinate-wise optimality condition is stronger than the critical point
condition \cite{GotohTT18,zhang2010analysis,bi2014exact} and the directional point condition \cite{PangRA17} when the function $(-g(\x))$ is locally $\rho$-bounded nonconvex. \textit{\textbf{(ii)}} Our optimality analysis can be also applied to the equivalent dual problem which is also a DC program as in (\ref{eq:toland:dual}). \textit{\textbf{(iii)}} We explain the optimality of coordinate-wise stationary point is stronger than that of previous definitions using the following one-dimensional example: $\min_{x} (x-1)^2 - 4|x|$. This problem contains three critical points $\{-1,0,3\}$, two directional points / local minima $\{-1,3\}$, and a unique coordinate-wise stationary point $\{3\}$. This unique coordinate-wise stationary point can be found using a clever breakpoint searching method (discussed later in Section \ref{sect:breakpoint}). 


\subsection{Convergence Analysis}

We provide a convergence analysis for \textit{\textbf{CD-SNCA}} and \textit{\textbf{CD-SCA}}. First, we define the approximate critical point and approximate coordinate-wise stationary point as follows.

\begin{definition}

(Approximate Critical Point) Given any constant $\epsilon>0$, a point $\check{\x}$ is called a $\epsilon$-approximate critical point if: $\text{dist}(\nabla f(\check{\x}) , \partial g(\check{\x}) - \partial h(\check{\x}))^2 \leq \epsilon$.

\end{definition}

\begin{definition}

(Approximate Coordinate-Wise Stationary Point) Given any constant $\epsilon>0$, a point $\ddot{\x}$ is called a $\epsilon$-approximate coordinate-wise stationary point if: $\frac{1}{n}\sum_{i=1}^n \text{dist}(0,\arg \min_{\eta} \mathcal{M}_{i}(\ddot{\x},\eta)   )^2 \leq \epsilon$, where $\mathcal{M}_{i}(\x,\eta)$ is defined in Definition \ref{def:cwsp}.

\end{definition}

\begin{theorem} \label{the:global:conv}
We have the following results. \textbf{(a)} For \textbf{\textit{CD-SNCA}}, it holds that $\F(\x^{t+1}) - \F(\x^{t}) \leq - \tfrac{\theta}{2} \|\x^{t+1}-\x^{t}\|^2$. Algorithm \ref{algo:main} finds an $\epsilon$-approximate {\textbf{coordinate-wise stationary point}} of Problem (\ref{eq:main}) in at most $T$ iterations in the sense of expectation, where $T \leq \lceil \frac{ 2n (\F(\x^0) - \F(\bar{\x})) }{\theta \epsilon} \rceil = \mathcal{O}(\epsilon^{-1})$. \textbf{(b)} For \textbf{\textit{CD-SCA}}, it holds that $\F(\x^{t+1}) - \F(\x^t) \leq - \tfrac{\beta}{2} \|\x^{t+1}-\x^{t}\|^2$ with $\beta \triangleq \min(\c) + 2 \theta$. Algorithm \ref{algo:main} finds an $\epsilon$-approximate {\textbf{critical point}} of Problem (\ref{eq:main}) in at most $T$ iterations in the sense of expectation, where $T \leq \lceil \frac{ 2n (\F(\x^0) - \F(\bar{\x})) }{\beta \epsilon} \rceil = \mathcal{O}(\epsilon^{-1})$.

\end{theorem}

\noi \textbf{Remarks}. While existing methods only find critical points or directional points of Problem (\ref{eq:main}), \textit{\textbf{CD-SNCA}} is guaranteed to find a coordinate-wise stationary point which has stronger optimality guarantees (See Theorem \ref{the:optimality}).

To achieve stronger convergence result for Algorithm \ref{algo:main}, we make the following \textit{Luo-Tseng error bound assumption}, which has been extensively used in all aspects of mathematical optimization (cf. \cite{Dong2021,YueZS19}).

\begin{assumption} \label{ass:luo}
(\textit{Luo-Tseng Error Bound} \cite{luo1993error,tseng2009coordinate}) We define a residual function as $\mathcal{R}(\x)= \frac{1}{n} \sum_{i=1}^n| \text{dist}(0,\bar{\mathcal{M}}_{i}(\x))|$ or $\mathcal{R}(\x)= \frac{1}{n} \sum_{i=1}^n|\text{dist}(0,\bar{\mathcal{P}}_{i}(\x))|$, where $\bar{\mathcal{M}}_{i}(\x)$ and $\bar{\mathcal{P}}_{i}(\x)$ are respectively defined in (\ref{eq:subprob:nonconvex}) and (\ref{eq:subprob:convex}). For any $\varsigma \geq \min_{\x} F(\x)$, there exist scalars $\delta>0$ and $ \varrho >0$ such that:
\beq
\forall \x,~\text{dist}(\x,\mathcal{X}) \leq \delta \mathcal{R}(\x),~\text{whenever}~F(\x) \leq \varsigma, \mathcal{R}(\x) \leq \varrho.\nn
\eeq
\noi Here, $\mathcal{X}$ is the set of stationary points satisfying $\mathcal{R}(\x)=0$.


\end{assumption}

%
%

We have the following theorems regarding to the convergence rate of \textit{\textbf{CD-SNCA}} and \textit{\textbf{CD-SCA}}.

\begin{theorem}\label{the:rate:ncvx}
\textbf{(Convergence Rate for \textit{\textbf{CD-SNCA}})}. Let $\ddot{\x}$ be any coordinate-wise stationary point. We define $\ddot{q}^t \triangleq F(\x^t) - F(\ddot{\x})$, $\ddot{r}^t  \triangleq \frac{1}{2}\|\x^t - \ddot{\x}\|_{\bar{\c}}^2,~\bar{\c} \triangleq \c + \theta$, $\bar{\rho} = \frac{\rho}{ \min(\bar{\c})}$, $\gamma \triangleq 1+\frac{\rho}{\theta}$, and $\varpi \triangleq 1-\bar{\rho}$. Assume that $z(\x) \triangleq -g(\x)$ is globally $\rho$-bounded non-convex. \textbf{(a)} We have $\varpi \E[\ddot{r}^{t+1}] + \gamma \E[\ddot{q}^{t+1}] \leq (\varpi+ \frac{\bar{\rho} }{n}) \ddot{r}^{t} + (\gamma -  \frac{1}{n}) \ddot{q}^{t}$. \textbf{(b)} If $\theta$ is sufficiently large such that $\varpi\geq 0$, $\mathcal{M}_{i^t}(\x^t,\eta)$ in (\ref{eq:subprob:nonconvex}) is convex w.r.t. $\eta$ for all $t$, and it holds that: $\E[\ddot{q}^{t+1}] \leq (\frac{\kappa_1 - \frac{1}{n}}{\kappa_1  })^{t+1} \ddot{q}^0$, where $\kappa_0 \triangleq \max(\bar{\c}) \frac{\delta^2}{\theta}$ and $\kappa_1\triangleq n\kappa_0 (\varpi + \frac{\bar{\rho} }{n} ) + \gamma$.

\end{theorem}

%
%
%
%

\begin{theorem} \label{the:rate:cvx}
\textbf{(Convergence Rate for \textit{\textbf{CD-SCA}})}. Let $\check{\x}$ be any critical point. We define $\check{q}^t \triangleq F(\x^t) - F(\check{\x})$, $\check{r}^t  \triangleq \frac{1}{2}\|\x^t - \check{\x}\|_{\bar{\c}}^2,~\bar{\c} \triangleq \c + \theta$, and $\bar{\rho} = \frac{\rho}{ \min(\bar{\c})}$. Assume that $z(\x) \triangleq -g(\x)$ is globally $\rho$-bounded non-convex. \textbf{(a)} We have $\E[\check{r}^{t+1}] + \E[\check{q}^{t+1}]  \leq (1+ \frac{\bar{\rho}}{n}) \check{r}^t     +  (1- \frac{1}{n})\check{q}^t$. \textbf{(b)} It holds that: $\E[\check{q}^{t+1}] \leq (\frac{\kappa_2 - \frac{1}{n}}{\kappa_2  })^{t+1} \check{q}^0$, where $\kappa_0 \triangleq \max(\bar{\c}) \frac{\delta^2}{\theta}$ and $\kappa_2 =n\kappa_0 ( 1 + \frac{\bar{\rho} }{n} ) + 1 $.
\end{theorem}

%
%
%
%
%

\noi \textbf{Remarks}. \textbf{(i)} Under the \textit{Luo-Tseng error bound assumption}, \textit{\textbf{CD-SNCA}} (or \textit{\textbf{CD-SCA}}) converges to the coordinate-wise stationary point (or critical point) Q-linearly. \textbf{(ii)} Note that the convergence rate $\kappa_1$ of \textit{\textbf{CD-SNCA}} and $\kappa_2$ of \textit{\textbf{CD-SCA}} depend on the same coefficients $\kappa_0$. When $n$ is large, the terms $n\kappa_0 (\varpi + \frac{\bar{\rho} }{n} )$ and $n\kappa_0 (1 + \frac{\bar{\rho} }{n} )$ respectively dominate the value of $\kappa_1$ and $\kappa_2$. If we choose $0\leq \varpi<1$ for \textit{\textbf{CD-SNCA}}, we have $\kappa_1 \ll \kappa_2$. Thus, the convergence rate of \textit{\textbf{CD-SNCA}} could be much faster than that of \textit{\textbf{CD-SCA}} for high-dimensional problems.



\section{A Breakpoint Searching Method for Proximal Operator Computation}
\label{sect:breakpoint}
This section presents a new breakpoint searching method to solve Problem (\ref{eq:prox}) exactly and efficiently for different $h(\cdot)$ and $g(\cdot)$. This method first identifies all the possible critical points / breakpoints $\Theta$ for $\min_{\eta \in \mathbb{R}} p(\eta)$ as in Problem (\ref{eq:prox}), and then picks the solution that leads to the lowest value as the optimal solution. We denote $\A\in \mathbb{R}^{m\times n}$ be an arbitrary matrix, and define $\g=\A\ei\in\mathbb{R}^m, \d=\A\x\in\mathbb{R}^m$.

\subsection{When $g(\y) =  \|\A \y\|_1$ and $h_i(\cdot) \triangleq 0$} \label{sect:L1norm}

Consider the problem: $\min_\eta \tfrac{a}{2}\eta^2 + b \eta  -  \|\A(\x+\eta e_i)\|_1$. It can be rewritten as: $\min_\eta p(\eta) \triangleq \frac{a}{2}\eta^2 + b \eta -  \|\g \eta + \d\|_1$. Setting the gradient of $p(\cdot)$ to zero yields: $0 = a \eta + b - \la \text{sign}(\eta\g+\d),\g \ra  = a \eta + b - \la \text{sign}(\eta+\d \div |\g|),\g \ra $, where we use: $\forall \rho>0,\text{sign}(\x) = \text{sign}(\rho \x)$. We assume $\g_i \neq 0$. If this does not hold and there exists $\g_j = 0$ for some $j$, then $\{\g_j,\d_j\}$ can be removed since it does not affect the minimizer of the problem. We define $\z \triangleq \{ +\tfrac{\d_1}{\g_1},-\tfrac{\d_1}{\g_1},...,+\tfrac{\d_m}{\g_m},-\tfrac{\d_m}{\g_m}  \} \in \mathbb{R}^{2m \times 1}$, and assume $\z$ has been sorted in ascending order. The domain $p(\eta)$ can be divided into $2m+1$ intervals: $(-\infty,\z_1)$, $(\z_1,\z_2)$,..., and $(\z_{2m},+\infty)$. There are $2m + 1$ breakpoints $\etas \in \mathbb{R}^{(2m+1)\times 1}$. In each interval, the sign of $(\eta + \d \div |\g|)$ can be determined. Thus, the $i$-th breakpoints for the $i$-th interval can be computed as $\etas_i = ( \la \text{sign}(\eta+\d \div |\g|),\g \ra - b) / {a}$. Therefore, Problem (\ref{eq:prox}) contains $2m+1$ breakpoints $\Theta=\{\etas_1,\etas_2,...,\etas_{(2m+1)}\}$ for this example.

\subsection{When $g(\y) =  \sum_{i=1}^s |\y_{[i]}|$ and $h_i(\y) \triangleq |\y_i|$} \label{eq:proximal:topk}

Consider the problem: $\min_{\eta}~\frac{a}{2}\eta^2 + b\eta +  |\x_i+\eta | - \sum_{i=1}^s | (\x+\eta\ei)_{[i]}|$. Since the variable $\eta$ only affects the value of $\x_i$, we consider two cases for $\x_i+\eta$. \textit{\textbf{(i)}} $\x_i+\eta$ belongs to the top-$s$ subset. This problem reduces to $\min_{\eta}~\frac{a}{2}\eta^2 + b\eta$, which contains one unique breakpoint: $\{-b/a\}$. \textit{\textbf{(ii)}} $\x_i+\eta$ does not belong to the top-$s$ subset. This problem reduces to $\min_{\eta}~\frac{a}{2}\eta^2 + bt +  |\x_i+\eta |$, which contains three breakpoints $\{-\x_i,~ (-1-b)/a,~(1-b)/a\}$. Therefore, Problem (\ref{eq:prox}) contains $4$ breakpoints $\Theta =\{-b/a,-\x_i,~ (-1-b)/a,~(1-b)/a\}$ for this example.

When we have found the breakpoint set $\Theta$, we pick the solution that results in the lowest value as the global optimal solution $\bar{\eta}$, i.e., $\bar{\eta} = \arg \min_{\eta} p(\eta),~s.t.~\eta \in \Theta$. Note that the coordinate-wise separable function $h_i(\cdot)$ does not bring much difficulty for solving Problem (\ref{eq:prox}).

\if
We let $\bar{p}(\eta) \triangleq \frac{a}{2}\eta^2 + b \eta -g(\x + \eta \ei )$. (\textbf{i}) For the bound constrained function with $h_i(z) =  {\tiny \left\{
        \begin{array}{cc}
          0, & lb \leq z \leq ub \\
          \infty, & else
        \end{array}
      \right.}$, we have: $\min_{\eta}~ \bar{p}(\eta),s.t.~ lb \leq \x_i + \eta \leq ub$. One can consider the breakpoint set $\Theta \triangleq \{lb,\max(lb,\min(ub,\varsigma)),ub\}$ with $\varsigma \triangleq \arg \min_{\eta}\bar{p}(\eta)$.  (\textbf{ii}) For the $\ell_1$ norm function with $h_i(z) = {\lambda} |z|$ for some ${\lambda}>0$, we have: $\min_{\eta}~ \bar{p}(\eta) + \lambda   | \x_i + \eta|$. One can consider the breakpoint set $\Theta \triangleq \{[\arg \min_{\eta} \bar{p}(\eta) + {\lambda} ( \x_i + \eta)],[\arg \min_{\eta} \bar{p}(\eta) - {\lambda} ( \x_i + \eta)],-\x_i\}$.


\fi

\section{Experiments}
\label{sect:exp}
This section demonstrates the effectiveness and efficiency of Algorithm \ref{algo:main} on two statistical learning tasks, namely the $\ell_p$ norm generalized eigenvalue problem and the approximate sparse optimization problem. For more experiments, please refer to Section \ref{sect:more:exp} in the Appendix.

\subsection{Experimental Settings}
We consider the following four types of data sets for the sensing/channel matrix $\G\in \mathbb{R}^{m\times n}$. \textit{\textbf{(i)}} `randn-m-n': $\G = \text{randn}(m,n)$. \textit{\textbf{(ii)}} `e2006-m-n': $\G = \X$. \textit{\textbf{(iii)}} `randn-m-n-C': $\G = \mathcal{N}(\text{randn}(m,n))$. \textit{\textbf{(iv)}} `e2006-m-n-C': $\G = \mathcal{N}(\X)$. Here, $\text{randn}(m,n)$ is a function that returns a standard Gaussian random matrix of size $m\times n$. $\X$ is generated by sampling from the original real-world data set `e2006-tfidf'. $\mathcal{N}(\G)$ is defined as: $[\mathcal{N}(\G)]_{I} = 100\cdot\G_I,[\mathcal{N}(\G)]_{\bar{I}} = \G_{\bar{I}}$, where $I$ is a random subset of $\{1,...,m n\}$, $\bar{I} = \{1,...,mn\} \setminus I$, and $|I| = 0.1\cdot mn$. The last two types of data sets are designed to verify the robustness of the algorithms.

All methods are implemented in MATLAB on an Intel 2.6 GHz CPU with 32 GB RAM. Only our breakpoint searching procedure is developed in C and wrapped into the MATLAB code, since it requires elementwise loops that are less efficient in native MATLAB. We keep a record of the relative changes of the objective by $\z_t = [{F}(\x^t)-{F}(\x^{t+1})]/{F}(\x^t)$, and let all algorithms run up to $T$ seconds and stop them at iteration $t$ if $\text{mean}([{\z}_{t-\text{min}(t,\upsilon)+1},\z_{t-min(t,\upsilon)+2},...,\z_t]) \leq \epsilon$. The default value $(\theta,\epsilon,\upsilon,T)=(10^{-6},10^{-10},500,60)$ is used. All methods are executed 10 times and the average performance is reported. Some Matlab code can be found in the supplemental material.



\subsection{$\ell_p$ Norm Generalized Eigenvalue Problem}

We consider Problem (\ref{eq:pca0}) with $p=1$ and $\Q=\textbf{I}$. We have the following problem: $\min_{\x}~\frac{\alpha}{2}\|\x\|_2^2 - \|\G\x\|_1$. It is consistent with Problem (\ref{eq:main}) with $f(\x) \triangleq \frac{\alpha}{2}\|\x\|_2^2$, $h(\x)\triangleq0$, and $g(\x) \triangleq \|\G\x\|_1$. The subgradient of $g(\x)$ at $\x^t$ can be computed as $\g^t \triangleq \G^T \text{sign}(\G\x^t)$. $\nabla f(\x)$ is $L$-Lipschitz with $L=1$ and coordinate-wise Lipschitz with $\c=\mathbf{1}$. We set $\alpha=1$.

We compare with the following methods. \textit{\textbf{(i)}} Multi-Stage Convex Relaxation (MSCR). It generates the new iterate using: $\x^{t+1} = \arg \min_{\x} f(\x) - \la \x-\x^t,\g^t\ra$. \textit{\textbf{(ii)}} Toland's dual method (T-DUAL). It rewrite the problem as: $\min_{-\textbf{1}\leq \y\leq \textbf{1}}\min_{\x}~f(\x) - \la \G\x,\y\ra$. Setting the gradient of $\x$ to zero, we have: $\alpha \x - \G^T\y = \mathbf{0}$, leading to the following dual problem: $\min_{-\textbf{1}\leq \y\leq \textbf{1}} -\frac{1}{2\alpha} \y^T \G\G^T\y$. Toland's dual method uses the iteration: $\y^{t+1} = \text{sign}(\G\G^T\y^t)$, and recovers the primal solution via $\x = \frac{1}{\alpha} \G^T\y$. Note that the method in \cite{kim2019simple} is essentially the Toland's duality method and they consider a constrained problem: $\min_{\|\x\| = 1} -\|\G\x\|_1$. \textit{\textbf{(iii)}} Subgradient method (SubGrad). It generates the new iterate via: $\x^{t+1} = \x^t - \frac{0.1}{t} \cdot ( \nabla f(\x^t) - \g^t)$. \textit{\textbf{(iv)}} \textit{\textbf{CD-SCA}} solves a convex problem: $\bar{\eta}^t = \arg \min_{\eta} \frac{\c_i+\theta}{2}\eta^2 + (\nabla_{i^t}f(\x^t)- \g^t_{i^t} )\eta$ and update $\x^t$ via $\x_{i^t}^{t+1}  = \x_{i^t}^t  + \bar{\eta}^t$. \textit{\textbf{(v)}} \textit{\textbf{CD-SNCA}} computes the nonconvex proximal operator of $\ell_1$ norm (see Section \ref{sect:L1norm}) as: $\bar{\eta}^t = \arg \min_{\eta} \frac{\c_i+\theta}{2}\eta^2 + \nabla_{i^t}f(\x^t) \eta - \|\G(\x+\eta \ei)\|_1$ and update $\x^t$ via $\x_{i^t}^{t+1}  = \x_{i^t}^t  + \bar{\eta}^t$.

\begin{table}[!t]
\begin{center}
\scalebox{0.48}{\begin{tabular}{|c|c|c|c|c|c|c|c|c|}
\hline
            & MSCR  & PDCA & T-DUAL  & \textit{\textbf{CD-SCA}} & \textit{\textbf{CD-SNCA}} \\
\hline
randn-256-1024 & \cthree{-1.329 $\pm$ 0.038} & \cthree{-1.329 $\pm$ 0.038} & \cthree{-1.329 $\pm$ 0.038} & \ctwo{-1.426 $\pm$ 0.056} & \cone{-1.447 $\pm$ 0.053}  \\
randn-256-2048 & \cthree{-1.132 $\pm$ 0.021} & \cthree{-1.132 $\pm$ 0.021} & \cthree{-1.132 $\pm$ 0.021} & \ctwo{-1.192 $\pm$ 0.019} & \cone{-1.202 $\pm$ 0.016}  \\
randn-1024-256 & \cthree{-5.751 $\pm$ 0.163} & \cthree{-5.751 $\pm$ 0.163} & -5.664 $\pm$ 0.173 & \ctwo{-5.755 $\pm$ 0.108} & \cone{-5.817 $\pm$ 0.129}  \\
randn-2048-256 & \cthree{-9.364 $\pm$ 0.183} & \cthree{-9.364 $\pm$ 0.183} & -9.161 $\pm$ 0.101 & \ctwo{-9.405 $\pm$ 0.182} & \cone{-9.408 $\pm$ 0.164}  \\
e2006-256-1024 & \ctwo{-28.031 $\pm$ 37.894} & \ctwo{-28.031 $\pm$ 37.894} & \cthree{-27.996 $\pm$ 37.912} & -27.880 $\pm$ 37.980 & \cone{-28.167 $\pm$ 37.826}  \\
e2006-256-2048 & \ctwo{-22.282 $\pm$ 24.007} & \ctwo{-22.282 $\pm$ 24.007} & \ctwo{-22.282 $\pm$ 24.007} & \cthree{-22.113 $\pm$ 23.941} & \cone{-22.448 $\pm$ 23.908}  \\
e2006-1024-256 & \ctwo{-43.516 $\pm$ 77.232} & \ctwo{-43.516 $\pm$ 77.232} & \cthree{-43.364 $\pm$ 77.265} & -43.283 $\pm$ 77.297 & \cone{-44.269 $\pm$ 76.977}  \\
e2006-2048-256 & \ctwo{-44.705 $\pm$ 47.806} & \ctwo{-44.705 $\pm$ 47.806} & \ctwo{-44.705 $\pm$ 47.806} & \cthree{-44.633 $\pm$ 47.789} & \cone{-45.176 $\pm$ 47.493}  \\
randn-256-1024-C & \cthree{-1.332 $\pm$ 0.019} & \cthree{-1.332 $\pm$ 0.019} & \cthree{-1.332 $\pm$ 0.019} & \ctwo{-1.417 $\pm$ 0.027} & \cone{-1.444 $\pm$ 0.029}  \\
randn-256-2048-C & \cthree{-1.161 $\pm$ 0.024} & \cthree{-1.161 $\pm$ 0.024} & \cthree{-1.161 $\pm$ 0.024} & \ctwo{-1.212 $\pm$ 0.022} & \cone{-1.219 $\pm$ 0.023}  \\
randn-1024-256-C & \cthree{-5.650 $\pm$ 0.141} & \cthree{-5.650 $\pm$ 0.141} & -5.591 $\pm$ 0.145 & \ctwo{-5.716 $\pm$ 0.159} & \cone{-5.808 $\pm$ 0.134}  \\
randn-2048-256-C & \cthree{-9.236 $\pm$ 0.125} & \cthree{-9.236 $\pm$ 0.125} & -9.067 $\pm$ 0.137 & \ctwo{-9.243 $\pm$ 0.145} & \cone{-9.377 $\pm$ 0.233}  \\
e2006-256-1024-C & \ctwo{-4.841 $\pm$ 6.410} & \ctwo{-4.841 $\pm$ 6.410} & \cthree{-4.840 $\pm$ 6.410} & -4.837 $\pm$ 6.411 & \cone{-5.027 $\pm$ 6.363}  \\
e2006-256-2048-C & \ctwo{-4.297 $\pm$ 2.825} & \ctwo{-4.297 $\pm$ 2.825} & \ctwo{-4.297 $\pm$ 2.823} & \cthree{-4.259 $\pm$ 2.827} & \cone{-4.394 $\pm$ 2.814}  \\
e2006-1024-256-C & \cthree{-6.469 $\pm$ 3.663} & \cthree{-6.469 $\pm$ 3.663} & \cthree{-6.469 $\pm$ 3.663} & \ctwo{-6.470 $\pm$ 3.663} & \cone{-6.881 $\pm$ 3.987}  \\
e2006-2048-256-C & \ctwo{-31.291 $\pm$ 60.597} & \ctwo{-31.291 $\pm$ 60.597} & \ctwo{-31.291 $\pm$ 60.597} & \cthree{-31.284 $\pm$ 60.599} & \cone{-32.026 $\pm$ 60.393}  \\
\hline
\end{tabular}}

\caption{Comparisons of objective values of all the methods for solving the $\ell_1$ norm PCA problem. The $1^{st}$, $2^{nd}$, and $3^{rd}$ best results are colored with \cone{red}, \ctwo{green} and \cthree{blue}, respectively.}
\label{tab:acc:l1pca}
\end{center}

\vspace{-13pt}

\begin{center}
\scalebox{0.55}{\begin{tabular}{|c|c|c|c|c|c|c|c|c|}
\hline
& MSCR  & PDCA & SubGrad  & \textit{\textbf{CD-SCA}} & \textit{\textbf{CD-SNCA}} \\
\hline
randn-256-1024 & \ctwo{0.090 $\pm$ 0.017} & \ctwo{0.090 $\pm$ 0.016} & 0.775 $\pm$ 0.040 & \cthree{0.092 $\pm$ 0.018} & \cone{0.034 $\pm$ 0.004}  \\
randn-256-2048 & \ctwo{0.052 $\pm$ 0.009} & \ctwo{0.052 $\pm$ 0.010} & 1.485 $\pm$ 0.030 & \cthree{0.061 $\pm$ 0.012} & \cone{0.027 $\pm$ 0.002}  \\
randn-1024-256 & 1.887 $\pm$ 0.353 & \cthree{1.884 $\pm$ 0.352} & 2.215 $\pm$ 0.379 & \ctwo{1.881 $\pm$ 0.337} & \cone{1.681 $\pm$ 0.346}  \\
randn-2048-256 & 3.795 $\pm$ 0.518 & \cthree{3.794 $\pm$ 0.518} & 4.127 $\pm$ 0.525 & \ctwo{3.772 $\pm$ 0.522} & \cone{3.578 $\pm$ 0.484}  \\
e2006-256-1024 & \ctwo{0.217 $\pm$ 0.553} & \ctwo{0.217 $\pm$ 0.553} & 0.597 $\pm$ 0.391 & \cthree{0.218 $\pm$ 0.556} & \cone{0.087 $\pm$ 0.212}  \\
e2006-256-2048 & \ctwo{0.050 $\pm$ 0.068} & \ctwo{0.050 $\pm$ 0.068} & \cthree{0.837 $\pm$ 0.209} & \ctwo{0.050 $\pm$ 0.068} & \cone{0.025 $\pm$ 0.032}  \\
e2006-1024-256 & \ctwo{3.078 $\pm$ 2.928} & \ctwo{3.078 $\pm$ 2.928} & 3.112 $\pm$ 2.844 & \cthree{3.097 $\pm$ 2.960} & \cone{2.697 $\pm$ 2.545}  \\
e2006-2048-256 & \ctwo{1.799 $\pm$ 1.453} & \ctwo{1.799 $\pm$ 1.453} & 1.918 $\pm$ 1.518 & \cthree{1.805 $\pm$ 1.456} & \cone{1.688 $\pm$ 1.398}  \\
randn-256-1024-C & \cthree{0.086 $\pm$ 0.012} & 0.087 $\pm$ 0.012 & 0.775 $\pm$ 0.038 & \ctwo{0.083 $\pm$ 0.011} & \cone{0.033 $\pm$ 0.002}  \\
randn-256-2048-C & \ctwo{0.043 $\pm$ 0.006} & \cthree{0.044 $\pm$ 0.006} & 1.472 $\pm$ 0.027 & 0.051 $\pm$ 0.009 & \cone{0.026 $\pm$ 0.001}  \\
randn-1024-256-C & \cthree{1.997 $\pm$ 0.250} & 1.998 $\pm$ 0.250 & 2.351 $\pm$ 0.297 & \ctwo{1.979 $\pm$ 0.265} & \cone{1.781 $\pm$ 0.244}  \\
randn-2048-256-C & \cthree{3.618 $\pm$ 0.681} & \ctwo{3.617 $\pm$ 0.682} & 3.965 $\pm$ 0.717 & 3.619 $\pm$ 0.679 & \cone{3.420 $\pm$ 0.673}  \\
e2006-256-1024-C  & \cthree{0.031 $\pm$ 0.031} & \cthree{0.031 $\pm$ 0.031} & 0.339 $\pm$ 0.073 & \ctwo{0.030 $\pm$ 0.028} & \cone{0.015 $\pm$ 0.014}  \\
e2006-256-2048-C  & \cthree{0.217 $\pm$ 0.575} & \cthree{0.217 $\pm$ 0.575} & 0.596 $\pm$ 0.418 & \ctwo{0.215 $\pm$ 0.568} & \cone{0.071 $\pm$ 0.176}  \\
e2006-1024-256-C  & \ctwo{3.789 $\pm$ 4.206} & \cthree{3.798 $\pm$ 4.213} & 3.955 $\pm$ 4.363 & 3.851 $\pm$ 4.339 & \cone{3.398 $\pm$ 3.855}  \\
e2006-2048-256-C  & \cthree{4.480 $\pm$ 6.916} & 4.482 $\pm$ 6.918 & 4.710 $\pm$ 7.292 & \ctwo{4.461 $\pm$ 6.844} & \cone{4.200 $\pm$ 6.608}  \\
\hline
\end{tabular}}
\caption{Comparisons of objective values of all the methods for solving the approximate sparse optimization problem. The $1^{st}$, $2^{nd}$, and $3^{rd}$ best results are colored with \cone{red}, \ctwo{green} and \cthree{blue}, respectively.}
\label{tab:sparse:opt}
\end{center}

\vspace{-18pt}
\end{table}
\begin{figure} [!t]
\centering
      \begin{subfigure}{0.25\textwidth}\includegraphics[height=\objimghei,width=\textwidth]{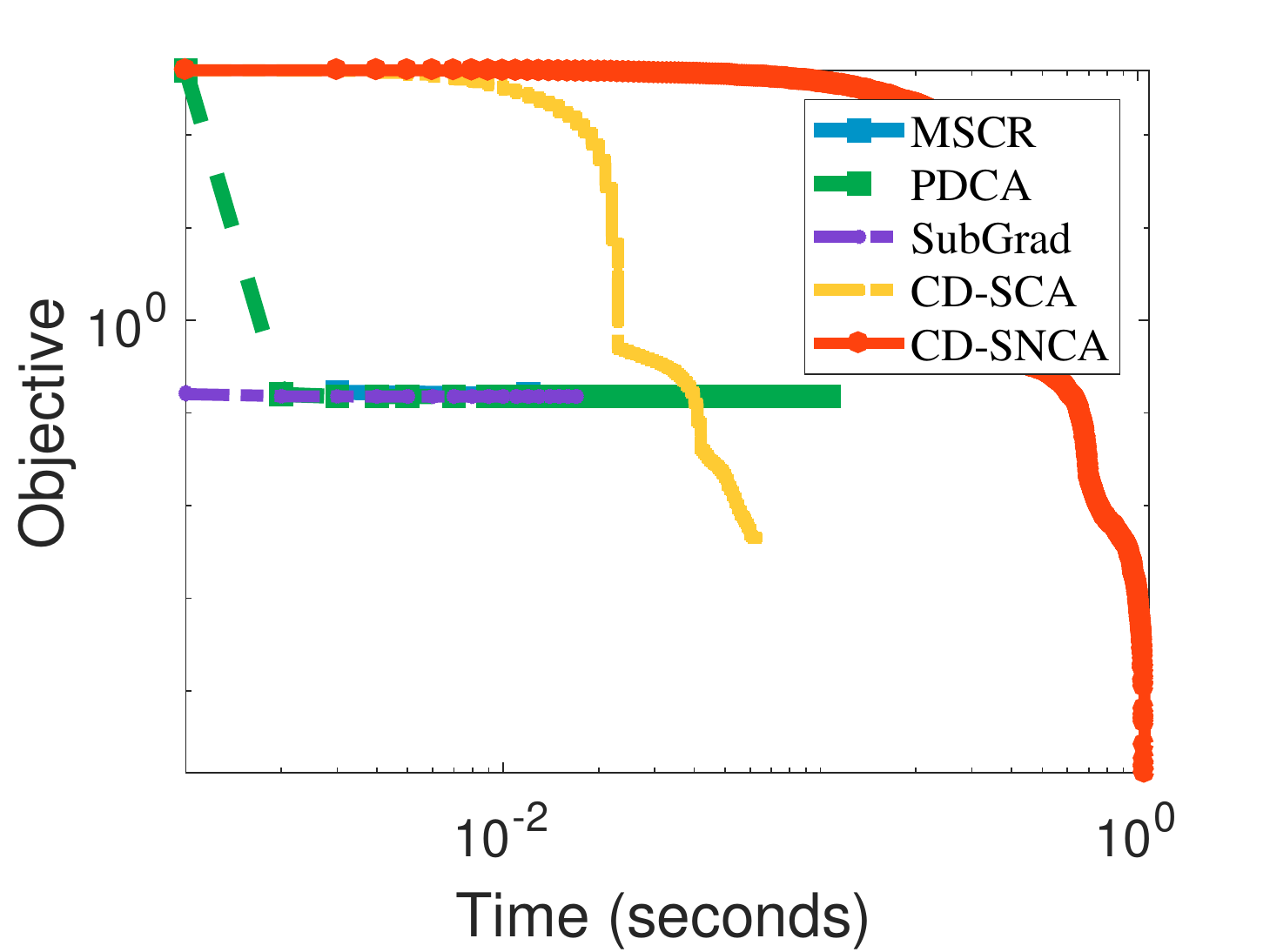}\vspace{-6pt} \caption{\scriptsize randn-256-1024  }\end{subfigure}\begin{subfigure}{0.25\textwidth}\includegraphics[height=\objimghei,width=\textwidth]{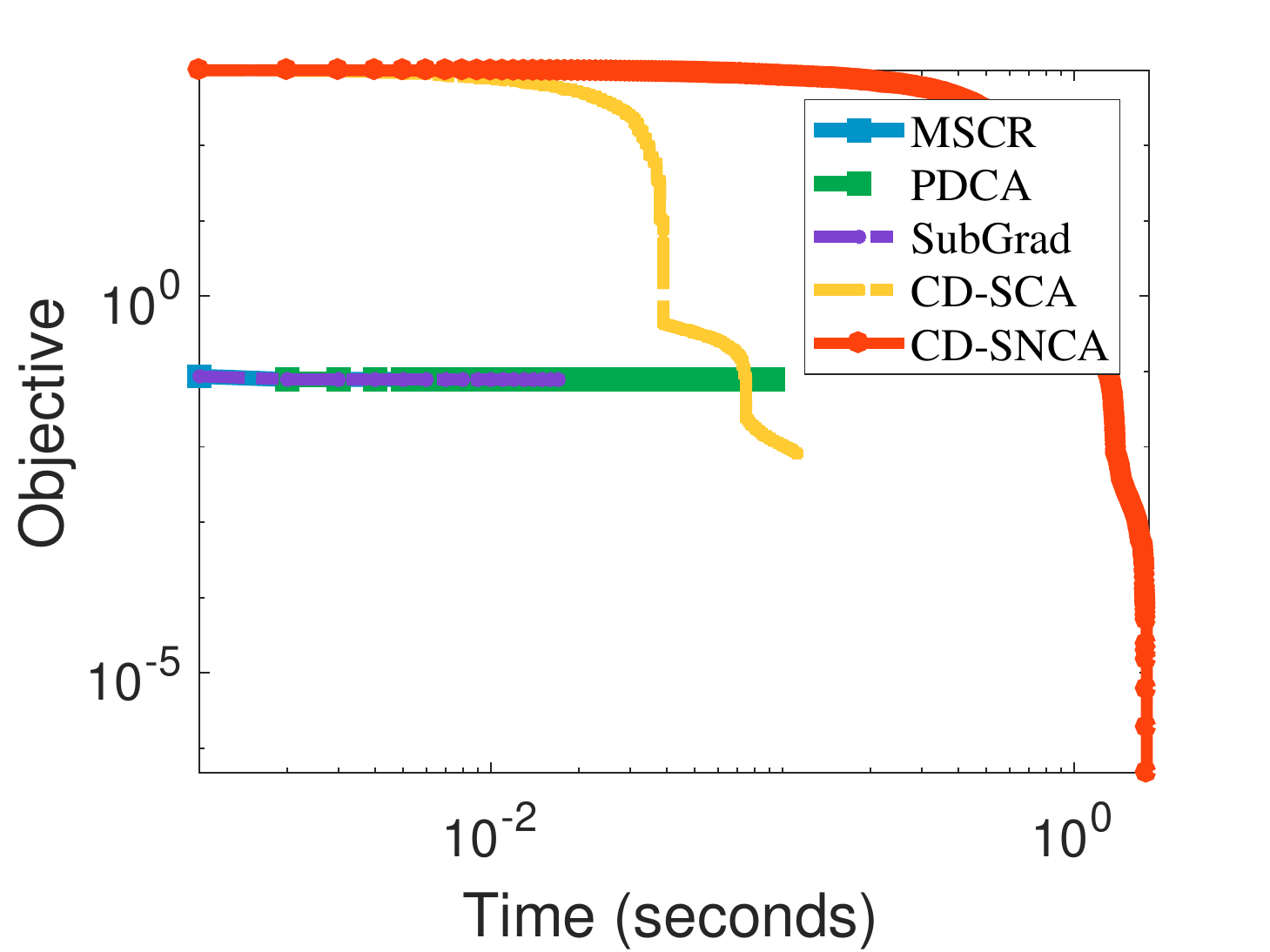}\vspace{-6pt} \caption{ \scriptsize randn-256-2048 } \end{subfigure}\\

      \begin{subfigure}{0.25\textwidth}\includegraphics[height=\objimghei,width=\textwidth]{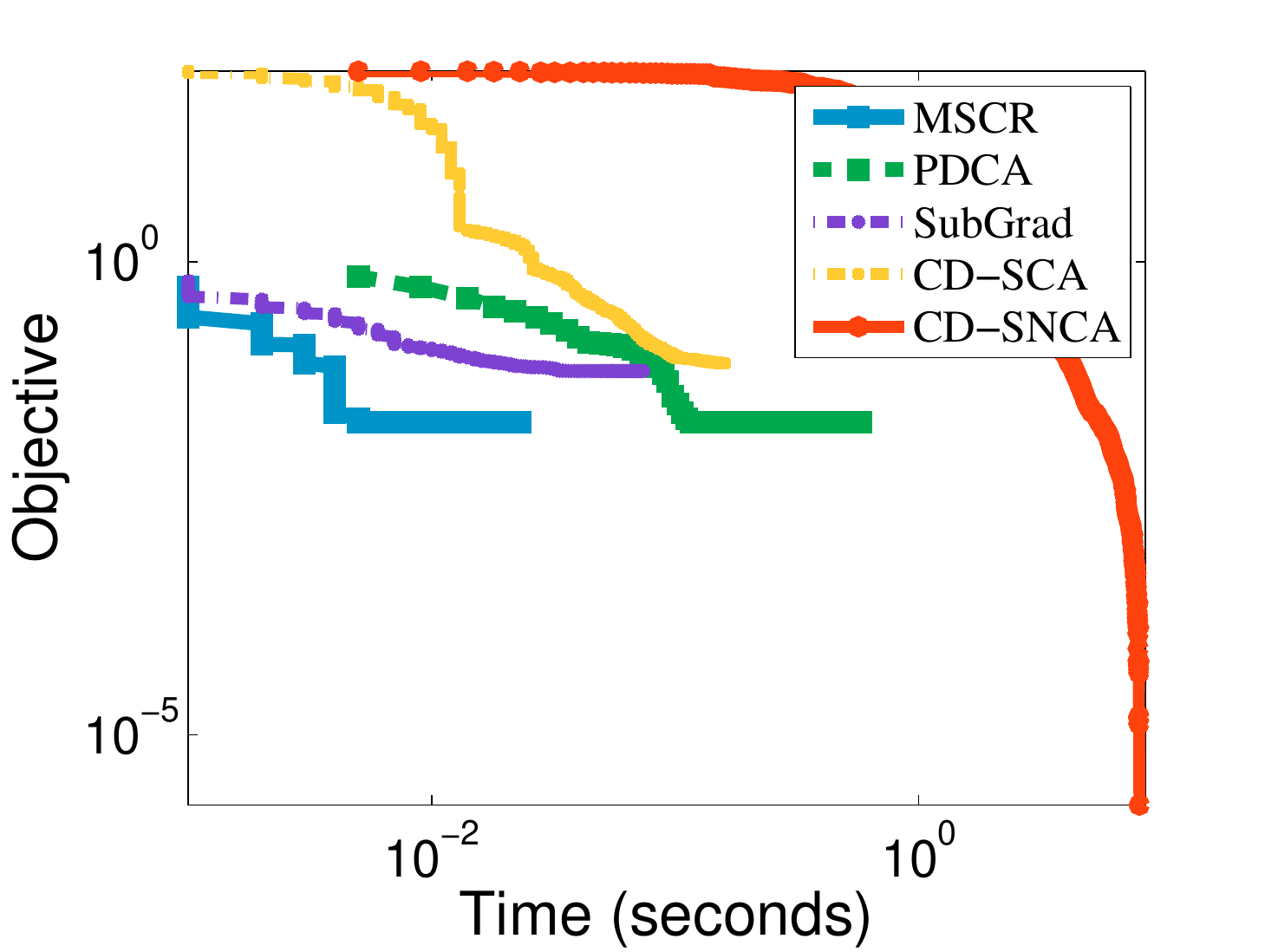}\vspace{-6pt} \caption{\scriptsize randn-1024-256}\end{subfigure}\begin{subfigure}{0.25\textwidth}\includegraphics[height=\objimghei,width=\textwidth]{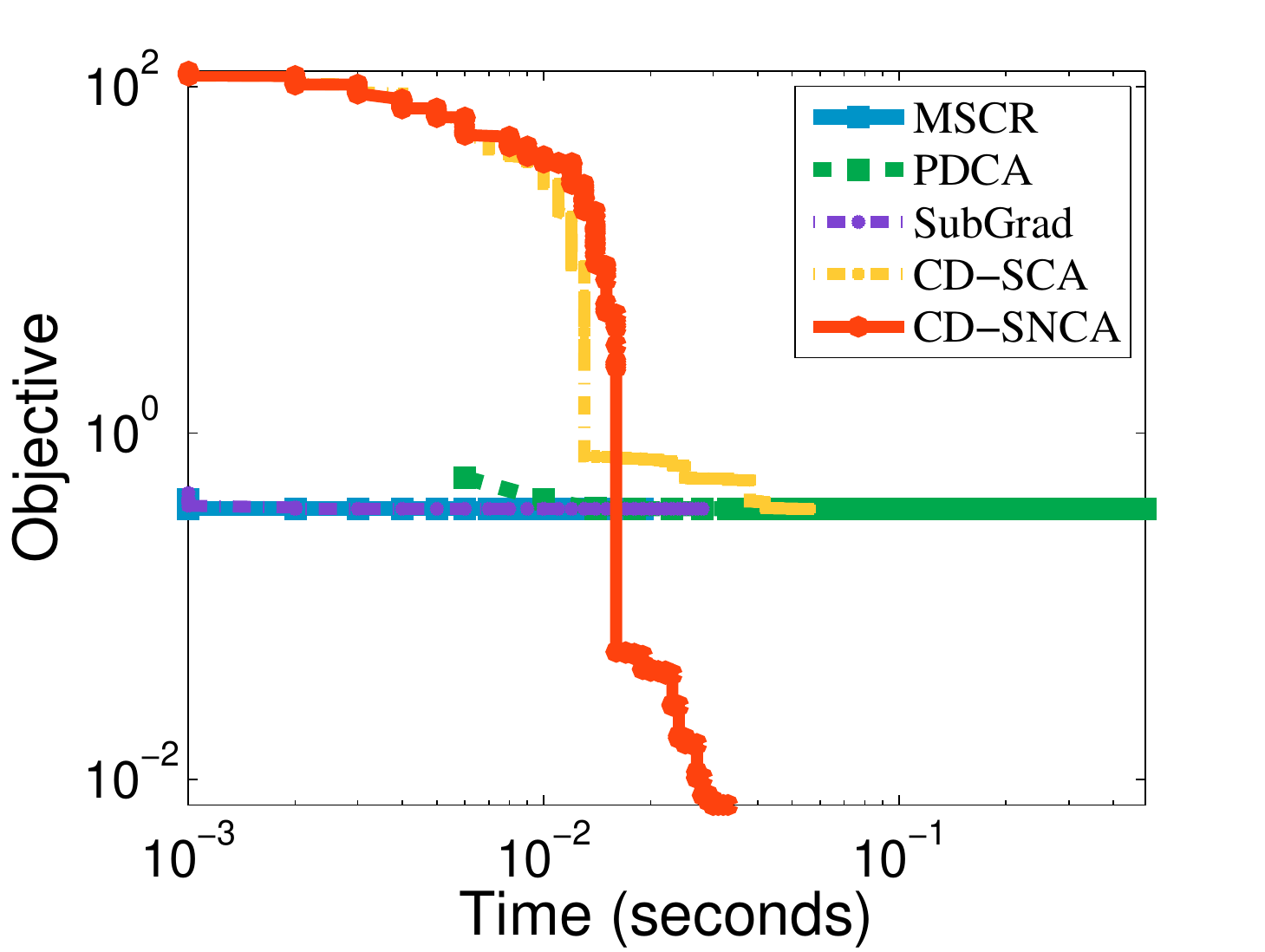}\vspace{-6pt} \caption{\scriptsize randn-2048-256 }\end{subfigure}
\caption{The convergence curve of the compared methods for solving the $\ell_p$ norm generalized eigenvalue problem on different data sets.}
\label{exp:cpu:1}
\vspace{-18pt}
\end{figure}


As can be seen from Table \ref{tab:acc:l1pca}, the proposed method \textit{\textbf{CD-SNCA}} consistently gives the best performance. Such results are not surprising since \textit{\textbf{CD-SNCA}} is guaranteed to find stronger stationary points than the other methods (while \textit{\textbf{CD-SNCA}} finds a coordinate-wise stationary point, all the other methods only find critical points).



\subsection{Approximate Sparse Optimization}

We consider solving Problem (\ref{eq:sparse:opt}). To generate the original signal $\ddot{\x}$ of $s$-sparse structure, we randomly select a support set $S$ with $|S|=200$ and set $\ddot{\x}_{\{1,...,n\} \setminus S} = \mathbf{0},~\ddot{\x}_{S} =\text{randn}(|S|,1)$. The observation vector is generated via $\y = \A\ddot{\x} + \text{randn}(m,1) \times 0.1\times\|\A\ddot{\x}\|$. This problem is consistent with Problem (\ref{eq:main}) with $f(\x) \triangleq \frac{1}{2}\|\G\x-\y\|_2^2$, $h(\x) \triangleq \rho \|\x\|_1$, and $g(\x) \triangleq \rho\sum_{i=1}^s | \x^t_{[i]}|$. $\nabla f(\x)$ is $L$-Lipschitz with $L=\|\G\|_2^2$ and coordinate-wise Lipschitz with $\c_i=(\G^T\G)_{ii},~\forall i$. The subgradient of $g(\x)$ at $\x^t$ can be computed as: $\g^t = \rho\cdot\arg \max_{\y}\la \y,\x^t \ra,s.t.\|\y\|_{\infty}\leq 1,\|\y\|_1\leq k$. We set $\rho=1$.


We compare with the following methods. \textit{\textit{\textbf{(i)}}} Multi-Stage Convex Relaxation (MSCR). It generate a sequence $\{\x^t\}$ as: $\x^{t+1} = \arg \min_{\x} \frac{1}{2}\|\G\x-\y\|_2^2 + \rho \|\x\|_1 -  \la \x-\x^t,\g^t \ra $. \textit{\textbf{(ii)}} Proximal DC algorithm (PDCA). It generates the new iterate using: $\x^{t+1} = \arg \min_{\x} \frac{L}{2}\|\x-\x^t\|_2^2 + \la \x - \x^t, \nabla f(\x)\ra +  \rho \|\x\|_1 -  \la \x-\x^t,\g^t \ra $. \textit{\textbf{(iii)}} Subgradient method (SubGrad). It uses the following iteration: $\x^{t+1} = \x^t  - \frac{0.1}{t} \cdot (\nabla f(\x) + \rho \text{sign}(\x^t) -  \g^t )$. \textit{\textbf{(iv)}} \textit{\textbf{CD-SCA}} solves a convex problem: $\bar{\eta}^t = \arg \min_{\eta} 0.5 (\c_{i^t}+\theta)\eta^2 + \rho|\x^t_{i^t} + \eta| + [\nabla f(\x^t) -\g^t]_{i^t} \cdot \eta$      and update $\x^t$ via $\x_{i^t}^{t+1}  = \x_{i^t}^t  + \bar{\eta}^t$. \textit{\textbf{(v)}} \textit{\textbf{CD-SNCA}} computes the nonconvex proximal operator of the top-$s$ norm function (see Section \ref{eq:proximal:topk}) as: $\bar{\eta}^t = \arg \min_{\eta} \frac{\c_i+\theta}{2}\eta^2 +  \nabla_{i^t} f(\x^t)  \eta + \rho | \x_{i^t}^t+\eta  | - \rho  \sum_{i=1}^s | (\x^t + \eta\ei)_{[i]}|$ and update $\x^t$ via $\x_{i^t}^{t+1}  = \x_{i^t}^t  + \bar{\eta}^t$.

As can be seen from Table \ref{tab:sparse:opt}, \textit{\textbf{CD-SNCA}} consistently gives the best performance. 

\vspace{-10pt}

\subsection{Computational Efficiency}

Figure \ref{exp:cpu:1} shows the convergence curve for solving the $\ell_p$ norm generalized eigenvalue problem. All methods take about 30 seconds to converge. \textit{\textbf{CD-SNCA}} generally takes a little more time to converge than the other methods. However, we argue that the computational time is acceptable and pays off as \textit{\textbf{CD-SNCA}} generally achieves higher accuracy.
%


\bibliographystyle{aaai23}
\bibliography{my}

\onecolumn
\clearpage
{
\huge
Appendix
}

\noi The appendix is organized as follows. \\
Section \ref{sect:proof} presents the mathematical proofs for the theoretical analysis. \\
Section \ref{sect:more:break} shows more examples of the breakpoint searching methods for proximal operator computation.\\
Section \ref{sect:more:exp} demonstrates some more experiments.\\
Section \ref{sect:disc:ext} provides some discussions of our methods.

\appendix

\section{Mathematical Proofs}
\label{sect:proof}

\subsection{Proof for Lemma \ref{lemma:bounde:nonconvex}}

\begin{proof}

Recall that the function $\tilde{z}(\x)\triangleq \|\x\|_p$ is convex when $p\geq 1$, and its subgradient w.r.t. $\x$ can be computed as $\partial \tilde{z}(\x) = \|\x\|_p^{1-p} \text{sign}(\x) \odot |\x|^{p-1}$. Therefore, the function $z(\x) = - \|\x\|_p$ with $p\geq 1$ is concave, and $\partial z(\x) = -  \|\x\|_p^{1-p} \text{sign}(\x) \odot |\x|^{p-1}$.

As the two reference points are different with $\x\neq \y$, we assume that there exists a constant $\epsilon>0$ satisfying $\|\x-\y\|\geq\epsilon$. We consider two cases for $p\geq1$ and derive the following results.

\noi \textbf{(a)} When $p\geq 2$, we have:
\bal
~& z(\x) - z(\y) - \la \x-\y, \partial z(\x)\ra \nn\\
\overset{(a)}{ = } ~& - \|\x\|_p + \|\y\|_p + \la \x-\y,   \|\x\|_p^{1-p} \text{sign}(\x) \odot |\x|^{p-1}\ra   \nn\\
\overset{(b)}{ \leq } ~& -\|\x\|_p  + \|\y\|_p  +   \|\x\|_p^{1-p}  \|\y-\x\|\| \text{sign}(\x) \odot |\x|^{p-1}\|   \nn\\
\overset{(c)}{ \leq } ~&   \|\y-\x\|_p +   \|\y-\x\|\|\x\|_p^{1-p}   \|\x\|_p^{p-1}   \nn\\
=~ &\textstyle\|\x - \y\|_p  + \|\x-\y\|_2  \nn\\
\overset{(d)}{ \leq } ~& 2 \|\x-\y\|_2 \nn\\
=~&\frac{4}{\epsilon} \cdot \frac{1}{2}\|\x-\y\|_2^2, \label{eq:global:rho}
\eal
\noi where step $(a)$ uses $z(\x)=-\|\x\|_p$ and $\partial z(\x) = -  \|\x\|_p^{1-p} \text{sign}(\x) \odot |\x|^{p-1}$; step $(b)$ uses the Cauchy-Schwarz inequality; step $(c)$ uses triangle inequality and the fact that $\||\x|^{p-1}\|_2\leq \|\x\|_p^{p-1}$ when $p\geq 2$; step $(d)$ uses $\|\x-\y\|_p\leq \|\x-\y\|$ for all $p\geq 2$.

\noi \textbf{(b)} When $1\leq p< 2$, we have:
\bal
~& z(\x) - z(\y) - \la \x-\y, \partial z(\x)\ra \nn\\
\overset{}{ = } ~& - \|\x\|_p + \|\y\|_p + \la \x-\y,   \|\x\|_p^{1-p} \text{sign}(\x) \odot |\x|^{p-1}\ra   \nn\\
\overset{}{ \leq } ~& -\|\x\|_p  + \|\y\|_p  +   \|\x\|_p^{1-p}  \|\y-\x\|\| \text{sign}(\x) \odot |\x|^{p-1}\|   \nn\\
\overset{(a)}{ \leq } ~&   \|\y-\x\|_p +  \|\x-\y\|\|\x\|_p^{1-p}   \|\x\|_p^{p-1} \cdot n^{1/p} \nn\\
\overset{(b)}{ \leq } ~&  \|\x-\y\| \cdot n^{1/p} + \|\x-\y\| \cdot n^{1/p} \nn\\
\overset{}{ = } ~& 2 \|\x-\y\| \cdot n^{1/p} \nn\\
\leq~&  \frac{4}{\epsilon} \cdot n^{1/p} \cdot \frac{1}{2}\|\x-\y\|_2^2 ,  \label{eq:global:rho:2}
\eal
\noi where step $(a)$ uses $\||\x|^{p-1}\|_2\leq n^{1/p} \|\x\|_p^{p-1}$ for all $1\leq p< 2$; step $(b)$ uses $\|\y-\x\|_p\leq \|\y-\x\| \cdot n^{1/p}$.

Combining the two inequalities as in (\ref{eq:global:rho}) and (\ref{eq:global:rho:2}), we conclude that there exists $\rho<+\infty$ such that $z(\x) - z(\y) - \la \x-\y, \partial z(\x)\ra \leq \frac{\rho}{2}\|\x-\y\|_2^2 $ with $\rho= {\tiny \{
                               \begin{array}{ll}
                                 {4}/{\epsilon},  & \hbox{$p\geq 2$;} \\
                                n^{1/p} \cdot 4 /\epsilon , & \hbox{ $1\leq p\leq 1$.}
                               \end{array}}$. In other words, $z(\x)= - \|\x\|_p$ is locally $\rho$-bounded nonconvex.

\end{proof}

\subsection{Proof for Lemma \ref{lemma:tech}}

\begin{proof}

\noi \textbf{(a)} For any $\x\in \mathbb{R}^n$, $\d\in \mathbb{R}^n$, and $\bar{\c}\in \mathbb{R}^n$, we derive the following equalities:
\beq
\frac{1}{n}\sum_{i=1}^n\|\x + \d_i \ei\|_{\bar{\c}}^2 &=& \frac{1}{n}\|\d\|_{\bar{\c}}^2 + \frac{2}{n} \la \x, {\bar{\c}} \odot \d \ra +  \|\x\|_{\bar{\c}}^2 \nn\\
&=& \frac{1}{n}\|\d\|_{\bar{\c}}^2  + \frac{2}{n} \la \x, \d \odot {\bar{\c}} \ra + (\frac{1}{n}\|\x\|_{\bar{\c}}^2 -  \frac{1}{n}\|\x\|_{\bar{\c}}^2) + \|\x\|_{\bar{\c}}^2 \nn\\
&=& \frac{1}{n}\|\d + \x\|_{\bar{\c}}^2  + (1-\frac{1}{n})\|\x\|_{\bar{\c}}^2. \nn
\eeq

\noi \textbf{(b)} The proof for this equality is almost the same as Lemma 1 in \cite{lu2015complexity}. For completeness, we include the proof here. We have the following results:
\beq
\frac{1}{n} \sum_{i=1}^n h(\x + \d_i \ei) & = & \frac{1}{n}\sum_{i=1}^n \left( h_i(\x_i + \d_i) + \sum_{j\neq i} h_j(\x_j) \right) \nn\\
&=& \frac{1}{n}\sum_{i=1}^n \left( h_i(\x_i + \d_i) \right) + \frac{1}{n}\sum_{i=1}^n \sum_{j\neq i} h_j(\x_j) \nn\\
&=& \frac{1}{n} h(\x+\d) + \frac{n-1}{n} h(\x).\nn
\eeq

\noi \textbf{(c)} We obtain the following results:
\beq
&& \frac{1}{n} \sum_{i=1}^n f(\x + \d_{i} \ei) \nn\\
&\overset{(a)}{ \leq } &\frac{1}{n} \left( \sum_{i=1}^n f(\x) + \la \nabla f(\x), \d_i \ei \ra + \frac{1}{2}\|\d\|_{\c}^2  \right) \nn\\
&\overset{(b)}{ = } & f(\x) + \frac{1}{n}  [\la \nabla f(\x),\d  \ra  + \frac{1}{2}\|\d\|_{\c}^2] \nn\\
&\overset{}{ = } & (1-\frac{1}{n})f(\x) + \frac{1}{n}  [f(\x) + \la \nabla f(\x),\d  \ra  + \frac{1}{2}\|\d\|_{\c}^2], \nn
\eeq
\noi where step $(a)$ uses the coordinate-wise Lipschitz continuity of $\nabla f(\x)$ as in (\ref{eq:f:Lipschitz}); step $(b)$ uses $ \sum_{i=1}^n \la \nabla f(\x^t),~\d_{i} \ei \ra =  \la \nabla f(\x^t),~\d \ra$ and $\sum_{i=1}^n \c_{i} \d_{i}^2 = \|\d\|_{\c}^2$.

\noi \textbf{(d)} We have the following inequalities:
\beq
 - \frac{1}{n}\sum_{i=1}^n g(\x + \d_i \ei)  &\overset{(a)}{ \leq } & \frac{1}{n} \sum_{i=1}^n \left(  -g(\x) - \la \partial g(\x),~\d_i \ei \ra \right),  \nn\\
&\overset{(b)}{ = }&  -g(\x) - \frac{1}{n} \la  \partial g(\x),~\d\ra , \label{eq:continou:gg}
\eeq
\noi where step $(a)$ uses the fact $g(\x)$ is convex that $\forall \x,\y,~-g(\y) \leq -g(\x) - \la \partial g(\x),~\y-\x\ra$; step $(b)$ uses $\sum_{i=1}^n \la \y,\d_i \ei\ra = \la \y,\d\ra$.

\end{proof}

\subsection{Proof for Lemma \ref{lemma:cw:point}}

First, since $z(\x)  \triangleq  -g(\x)$ is locally $\rho$-bounded nonconvex at the point $\ddot{\x}$, we have:
 \beq
-g(\ddot{\x})   &\leq& -g(\y)  - \la \ddot{\x} - \y,~ \partial g(\ddot{\x}) \ra + \frac{\rho}{2}\|\ddot{\x} - \y\|_2^2,~\forall \y.\nn
\eeq
\noi Applying the inequality above with $\y = \ddot{\x} + \d$ for any $\d\in\mathbb{R}^n$, we obtain:
\beq
\forall \d,~-g(\ddot{\x}) &\overset{}{ \leq } & -g(\ddot{\x} + \d)  - \la \ddot{\x} - (\ddot{\x} + \d),~ \partial g(\ddot{\x}) \ra + \frac{\rho}{2}\|\ddot{\x} - (\ddot{\x} + \d)\|_2^2\nn\\
&\overset{}{ = } & -g(\ddot{\x} + \d)  + \la \d,~ \partial g(\ddot{\x}) \ra + \frac{\rho}{2}\|\d\|_2^2\nn\\
&\overset{(a)}{ \leq } & -g(\ddot{\x} + \d)  + \sum_{i=1}^n g(\ddot{\x} + \d_i \ei)  - n g(\ddot{\x})  + \frac{\rho}{2}\|\d\|_2^2, \label{eq:dfdfa11}
\eeq
\noi where step $(a)$ uses claim \textbf{(d)} in Lemma \ref{lemma:tech} that $-\sum_{i=1}^n g(\ddot{\x} + \d_i \ei)  \leq   - g(\x)   - \la  \partial g(\x),\d\ra -(n-1)g(\x)$.

Second, by the optimality of $\ddot{\x}$, we obtain:
\beq
 h(\ddot{\x} ) - g(\ddot{\x}  ) \leq   \la \d_i \ei,~\nabla f(\ddot{\x}) \ra + \frac{\c_i+\theta}{2} \d_i^2 + h(\ddot{\x} + \d_i \ei) - g(\ddot{\x} + \d_i \ei),~\forall \d_i. \nn
\nn
\eeq
\noi Summing the inequality above over $i=1,...,n$, we have:
\beq
\forall \d,~0 &\leq&  n g(\ddot{\x})  - n h(\ddot{\x}) + \frac{1}{2} \|\d  \|^2_{(\c+\theta)}  + \la \d,~ \nabla f(\ddot{\x}) \ra + \sum_{i=1}^n h(\ddot{\x} + \d_i \ei) -\sum_{i=1}^n g(\ddot{\x} + \d_i \ei) \nn\\
&\overset{(a)}{ \leq } & \frac{1}{2} \|\d  \|^2_{(\c+\theta)} + \la \d,~ \nabla f(\ddot{\x}) \ra - h(\ddot{\x}) + h(\ddot{\x}+\d) + g(\ddot{\x})- g(\ddot{\x}+\d) + \frac{\rho}{2}\|\d\|_2^2 \nn\\
&\overset{(b)}{ \leq } & \frac{1}{2} \|\d  \|^2_{(\c+\theta)} + f(\ddot{\x}+ \d) - f(\ddot{\x}) - h(\ddot{\x}) + h(\ddot{\x}+\d) + g(\ddot{\x})- g(\ddot{\x}+\d) + \frac{\rho}{2}\|\d\|_2^2 \nn\\
&\overset{(c)}{ = } &   \frac{1}{2} \|\d  \|^2_{(\c+\theta)} + F(\ddot{\x}+ \d) - F(\ddot{\x}) + \frac{\rho}{2}\|\d\|_2^2, \nn
\eeq
\noi where step $(a)$ uses (\ref{eq:tech:hh}) in Lemma \ref{lemma:tech} and (\ref{eq:dfdfa11}); step $(b)$ uses the convexity of $f(\cdot)$ that:
\beq \label{eq:strong:conve:2}
\forall \d,~\la \nabla f(\ddot{\x}),  (\ddot{\x}+ \d) - \ddot{\x}\ra \leq f(\ddot{\x}+ \d) - f(\ddot{\x}); \nn
\eeq
\noi step $(c)$ uses the definition of $\F(\x) = f(\x) +h(\x) - g(\x)$. Rearranging terms, we obtain:
\beq
\forall \d,~F(\ddot{\x})   \leq F(\ddot{\x}+ \d) +  \frac{1}{2} \|\d\|^2_{(\c+\theta+\rho)}. \nn
\eeq

\subsection{Proof for Theorem \ref{the:optimality}}

\begin{proof}

\noi \textbf{(a)} We show that any optimal point $\bar{\x}$ is a coordinate-wise stationary point $\ddot{\x}$, i.e., $\{\bar{\x}\} \subseteq \{\ddot{\x}\}$. By the optimality of $\bar{\x}$, we have:
\beq
f(\bar{\x}) +   h(\bar{\x})  -   g(\bar{\x})  \leq f(\x) +   h(\x)  -   g(\x),~\forall \x .\nn
\eeq
\noi Letting $\x = \bar{\x} + \d_i \ei$, we have:
\beq \label{eq:dfsaddfd}
f(\bar{\x}) +   h(\bar{\x})  -   g(\bar{\x})  &\leq& f(\bar{\x} + \d_i \ei) +   h(\bar{\x} + \d_i \ei)  -   g(\bar{\x} + \d_i \ei),~\forall \d_i,\forall i \nn \\
&\overset{(a)}{ \leq } & f(\bar{\x}) +  \la  \nabla_i f(\bar{\x}),~ \d_i \ei ) + \frac{\c_i}{2} \d_i^2+   h(\bar{\x} + \d_i \ei)  -   g(\bar{\x} + \d_i \ei),~\forall \d_i,\forall i ,
\eeq
\noi where step $(a)$ uses the coordinate-wise Lipschitz continuity of $\nabla f(\cdot)$ that:
\beq \label{eq:f:Lipschitz:111}
f(\bar{\x} + \d_i \ei  ) \leq f(\bar{\x}) +  \la  \nabla_i f(\bar{\x}),~ \d_i \ei ) + \frac{\c_i}{2} \d_i^2,~\forall \d_i\nn.
\eeq

\noi We denote $\bar{\eta}_i$ as the minimizer of the following problem:
\beq
\forall i,~\bar{\eta}_i \in \arg \min_{\eta}~\mathcal{M}_{i}(\bar{\x},\eta).\nn
\eeq
\noi Rearranging terms for (\ref{eq:dfsaddfd}) and using the fact that $\theta\geq 0$, we have:
\beq
  && h(\bar{\x})  -   g(\bar{\x})  \leq h(\bar{\x} + \d_i \ei)  -   g(\bar{\x} + \d_i \ei)   +  \la  \nabla_i f(\bar{\x}),~ \d_i \ei ) + \frac{\c_i+\theta}{2} \d_i^2,~\forall \d_i , \forall i\nn\\
   &\overset{(a)}{ \Rightarrow } &  h(\bar{\x})  -   g(\bar{\x})  \leq h(\bar{\x} + \bar{\eta}_i \ei)  -   g(\bar{\x} + \d_i \ei)   +  \la  \nabla_i f(\bar{\x}),~ \bar{\eta}_i \ei ) + \frac{\c_i+\theta}{2} (\bar{\eta}_i)^2 , \forall i\nn\\
   &\overset{(b)}{ \Rightarrow } &  \mathcal{M}_i(\bar{\x},0) \leq \min_{\eta}~\mathcal{M}_{i}(\bar{\x},\eta), \forall i,\nn\\
   &\overset{}{ \Rightarrow } &  0 \in \bar{\mathcal{M}}_i(\bar{\x}),\forall i.\nn
\eeq
\noi where step $(a)$ uses the choice $\d_i =\bar{\eta}_i$ for all $i$; step $(b)$ uses the fact that $\forall i,~\mathcal{M}_i(\bar{\x},0)=h(\bar{\x})-g(\bar{\x})$ and the definition of $\min_{\eta}~\mathcal{M}_{i}(\bar{\x},\eta)$. Therefore, any optimal point $\bar{\x}$ is also a coordinate-wise stationary point $\ddot{\x}$.

\textbf{(b)} We show that any coordinate-wise stationary point $\ddot{\x}$ is a directional point $\grave{\x}$, i.e., $\{\ddot{\x}\} \subseteq \{\grave{\x}\}$. Applying the inequality in Lemma \ref{lemma:cw:point} with $\d = t (\y-\ddot{\x})$, we directly obtain the following results:
\beq
\lim_{t \downarrow  0}~ \frac{\F(\ddot{\x} + t (\y-\ddot{\x})) - \F(\ddot{\x}) }{t} &\geq & \lim_{t \downarrow  0}~\frac{- \frac{1}{2}\|t (\y-\ddot{\x})\|_{(\c+\theta+\rho)}^2}{ t} \nn\\
&\overset{(a)}{ = } &  \lim_{t \downarrow  0}~  \frac{-t^2 \frac{1}{2}\| \y-\ddot{\x}\|_{(\c+\theta+\rho)}^2}{t} =  0, \nn
\eeq
\noi where step $(a)$ uses the boundedness of $\rho$ that $\rho<+\infty$. Therefore, any coordinate-wise stationary point $\ddot{\x}$ is also a directional point $\grave{\x}$.

\textbf{(c)} We show that any directional point $\grave{\x}$ is a critical point $\check{\x}$, i.e., $\{\grave{\x}\} \subseteq \{\check{\x}\}$. Noticing $f(\x)$, $h(\x)$, and $g(\x)$ are convex, we have:
\beq
f(\z)\leq  f(\grave{\x}) - \la \grave{\x} - \z,~ \nabla f(\z) \ra,~\forall \z,\nn \\
h(\z) \leq  h(\grave{\x}) - \la \grave{\x} - \z,~ \partial h(\z) \ra,~\forall \z,\nn\\
-g(\z)   \leq -g(\grave{\x})  - \la \z - \grave{\x},~ \partial  g(\grave{\x}) \ra,~\forall \z.\nn
\eeq
\noi Adding these three inequalities together, we obtain:
\beq\label{eq:opt:ooo}
\F(\z) - \F(\grave{\x}) \leq   \la  \z - \grave{\x}, - \partial  g(\grave{\x}) +   \nabla f(\z) + \partial h(\z) \ra,~\forall \z.
\eeq

\noi We derive the following inequalities:
\beq
\forall \y \in \text{dom}(\F),~ 0 &\leq&~ \lim_{t \downarrow  0} \frac{ \F(\grave{\x} + t (\y-\grave{\x}))   -   \F(\grave{\x})  }{t} \nn\\
 &\overset{(a)}{ \leq } &~ \lim_{t \downarrow  0}  \frac{ \la  (\grave{\x} + t (\y-\grave{\x})) - \grave{\x}, - \partial  g(\grave{\x}) +   \nabla f(\grave{\x} + t (\y-\grave{\x})) + \partial h(\grave{\x} + t (\y-\grave{\x})) \ra }{t},  \nn\\
 &\overset{}{ = } &~ \lim_{t \downarrow  0}  \la \y-\grave{\x}, - \partial  g(\grave{\x}) +   \nabla f(\grave{\x} + t (\y-\grave{\x})) + \partial h(\grave{\x} + t (\y-\grave{\x})) \ra,  \nn\\
&\overset{(b)}{ = } &~ \lim_{t \downarrow  0}   \la  \y-\grave{\x}    , \partial \F(\grave{\x}) \ra , \nn
\eeq
\noi where step $(a)$ uses (\ref{eq:opt:ooo}) with $\z = \grave{\x} + t (\y-\grave{\x})$; step $(b)$ uses $\grave{\x} + t (\y-\grave{\x}) = \grave{\x}$ as $t \downarrow0$. Noticing the inequality above holds for all $\y$ only when $0 \in \partial \F(\grave{\x})$, we conclude that any directional point $\grave{\x}$ is also a critical point $\check{\x}$.

%
%

\end{proof}

\subsection{Proof for Theorem \ref{the:global:conv}}

\begin{proof}

\textbf{(a)} We now focus on \textit{\textbf{CD-SNCA}}. Since $\bar{\eta}^t$ is the global optimal solution to Problem (\ref{eq:subprob:nonconvex}), we have:
\beq
&& f(\x^t) + \la \bar{\eta}^t \eit,~\nabla f(\x^t) \ra + \frac{\c_{{i^t}}+\theta}{2} (\bar{\eta}^t)^2 + h(\x^t + \bar{\eta}^t \eit) - g(\x^t + \bar{\eta}^t \eit) \nn\\
& \leq & f(\x^t) + \la \eta \eit,~\nabla f(\x^t) \ra + \frac{\c_{{i^t}}+\theta}{2} \eta^2 + h(\x^t + \eta \eit) - g(\x^t + \eta \eit), \forall \eta.\nn
\eeq
\noi Letting $\eta=0$ and using the fact that $\x^{t+1}  = \x^t  + \bar{\eta}^t \cdot \eit$, we obtain:
\beq \label{eq:nonconvex:1d}
f(\x^t) + \la \x^{t+1}-\x^{t} ,~\nabla f(\x^t) \ra + \frac{1}{2} \|\x^{t+1}-\x^{t}\|_{\c+\theta}^2 + h(\x^{t+1}) - g(\x^{t+1}) \leq  \F(\x^t)
\eeq

\noi We derive the following results:
\beq
&&\F(\x^{t+1}) - \F(\x^{t})\nn\\
&\overset{(a)}{ \leq } & \F(\x^{t+1})  - f(\x^t) - \la \x^{t+1}-\x^{t}, \nabla f(\x^t) \ra - \tfrac{1}{2} \|\x^{t+1}-\x^{t}\|_{\c+\theta}^2 - h(\x^{t+1}) + g(\x^{t+1}),\nn\\
&\overset{(b)}{ = } & f(\x^{t+1})  - f(\x^t) - \la \x^{t+1}-\x^{t},~\nabla f(\x^t) \ra - \tfrac{1}{2} \|\x^{t+1}-\x^{t}\|_{\c+\theta}^2 ,\nn\\
&\overset{(c)}{ \leq } & - \frac{\theta}{2} \|\x^{t} - \x^{t+1}\|_2^2 ,\nn
\eeq
\noi where step $(a)$ uses (\ref{eq:nonconvex:1d}); step $(b)$ uses the definition $\F(\x^{t+1}) = f(\x^{t+1}) + h(\x^{t+1}) - g(\x^{t+1})$; step $(c)$ uses the coordinate-wise Lipschitz continuity of $\nabla f(\cdot)$.

\noi Taking the expectation for the inequality above, we obtain a lower bound on the expected progress made by each iteration for \textit{\textbf{CD-SNCA}}:
\beq
\E[\F(\x^{t+1})] - \F(\x^{t}) \leq  -\E[\tfrac{\theta}{2n} \|\x^{t+1}-\x^{t}\|^2]\nn
\eeq

\noi Summing up the inequality above over $t=0,1,...,T-1$, we have:
\beq
&&\E[\frac{\theta}{2 } \sum_{t=0}^{T-1}\| \x^{t+1}-\x^{t} \|_2^2] \leq n\E[\F(\x^0) - \F(\x^{T+1})] \leq n\E[\F(\x^0) - \F(\bar{\x})] \nn
\eeq
\noi As a result, there exists an index $\bar{t}$ with $0\leq \bar{t}\leq T-1$ such that:
\beq
\E[\| \x^{\bar{t}+1}-\x^{\bar{t}} \|_2^2 ] \leq \frac{2 n (\F(\x^0) - \F(\bar{\x}))}{\theta T}.
\eeq
\noi Furthermore, for any $t$, we have:
\beq
\E[\| \x^{\bar{t}+1}-\x^{\bar{t}} \|_2^2 ] = \frac{1}{n} \sum_{i=1}^n (  \bar{\mathcal{M}}_{i}(\x^{\bar{t}})  )^2.
\eeq
\noi Combining the two inequalities above, we have the following result:
\beq
\frac{1}{n} \sum_{i=1}^n (  \bar{\mathcal{M}}_{i}(\x^{\bar{t}})  )^2 \leq \frac{2 n (\F(\x^0) - \F(\bar{\x}))}{\theta T}.
\eeq
Therefore, we conclude that \textit{\textbf{CD-SNCA}} finds an $\epsilon$-approximate coordinate-wise stationary point in at most $T$ iterations in the sense of expectation, where
\beq
T \leq \lceil\frac{2 n (\F(\x^0) - \F(\bar{\x}))}{\theta \epsilon}\rceil = \mathcal{O}(\epsilon^{-1}).\nn
\eeq


\textbf{(b)} We now focus on \textit{\textbf{CD-SCA}}. Since $\bar{\eta}^t$ is the global optimal solution to Problem (\ref{eq:subprob:convex}), we have:
\beq \label{eq:optimadafdfdf}
0 \in [\nabla f(\x^t)  + \partial h(\x^{t+1}) - \partial g(\x^{t+1})]_{i^t} + (\c_{i^t}+\theta) \bar{\eta}^t.
\eeq
\noi Using the coordinate-wise Lipschitz continuity of $\nabla f(\cdot)$, we obtain:
\beq\label{eq:lip:f:conv}
f(\x^{t+1}) \leq f(\x^{t}) + \la \x^{t+1} - \x^t,~\nabla f(\x^t) \ra + \frac{1}{2}\|\x^{t+1} - \x^t\|_{\c}^2
\eeq
\noi Since both $h(\cdot)$ and $g(\cdot)$ are convex, we have:
\beq
h(\x^{t+1})   \leq h(\x^{t}) - \la \x^{t} - \x^{t+1},~\nabla h(\x^{t+1}) \ra, \label{eq:cvx:h} \\
-g(\x^{t+1}) \leq -g(\x^{t}) - \la \partial g(\x^{t}),\x^{t+1}-\x^{t}\ra. \label{eq:cav:p}
\eeq
\noi Adding these three inequalities in (\ref{eq:lip:f:conv}), (\ref{eq:cvx:h}), and (\ref{eq:cav:p}) together, we have:
\beq
&& \F(\x^{t+1}) -\F(\x^t) \nn\\
&\leq&  \la \x^{t+1}-\x^{t},~\nabla f(\x^t)  + \partial h(\x^{t+1}) - \partial g(\x^t) \ra + \tfrac{1 }{2}\|\x^{t+1}-\x^{t}\|_{\c}^2 \nn\\
&\overset{(a)}{ = } &  \la \bar{\eta}^t \eit,~\nabla f(\x^t)  + \partial h(\x^{t+1}) - \partial g(\x^t) \ra + \frac{\c_{i^t}  }{2}\|\bar{\eta}^t \eit\|_2^2 \nn\\
& \overset{}{ = } &  \bar{\eta}^t (\nabla f(\x^t)  + \partial h(\x^{t+1}) - \partial g(\x^t) )_{i^t} + \frac{\c_i}{2}(\bar{\eta}^t)^2 \nn\\
&\overset{(b)}{ = } &  - \frac{\c_{}  + 2 \theta}{2}(\bar{\eta}^t)^2 \nn \\
&\overset{(c)}{ \leq  } &  - \frac{ \min(\c) + 2 \theta}{2} \|\x^{t+1}-\x^{t}\|^2, \nn
\eeq
\noi where step $(a)$ uses the fact that $\x^{t+1} - \x^t =  \bar{\eta}^t \eit$; step $(b)$ uses (\ref{eq:optimadafdfdf}); step $(c)$ uses $(\bar{\eta}^t)^2=\|\x^{t+1}-\x^{t}\|^2$.

Using similar strategies as in deriving the results for \textit{\textbf{CD-SNCA}}, we conclude that Algorithm \ref{algo:main} finds an $\epsilon$-approximate critical point of Problem (\ref{eq:main}) in at most $T$ iterations in the sense of expectation, where $T \leq \lceil \frac{ 2n (\F(\x^0) - \F(\bar{\x})) }{\beta \epsilon} \rceil = \mathcal{O}(\epsilon^{-1})$.

\end{proof}

\subsection{Proof for Theorem \ref{the:rate:ncvx}}

\begin{proof}

Let $\ddot{\x}$ be any coordinate-wise stationary point. First, the optimality condition for the nonconvex subproblem as in (\ref{eq:subprob:nonconvex}) can be written as:
\beq \label{eq:nonconvex:cd:optimality}
0 \in \nabla_{i^t} f(\x^t) +  \bar{\c}_{i^t}  \bar{\eta}^t + \partial_{i^t} h(\x^{t+1}) - \partial_{i^t} g(\x^{t+1}).
\eeq

\noi Second, for any $\x^t$, $\x^{t+1}$, and $\ddot{\x}$, since $i^t$ is chosen uniformly and randomly, we have:
\beq \label{eq:tech:xx:00:1}
\E[\|\x^{t+1} - \ddot{\x}\|_{\bar{\c}}^2 ] = \frac{1}{n} \sum_{i=1}^n \|\x^t + (\x^{t+1}-\x^t)_i \ei - \ddot{\x}\|_{\bar{\c}}^2.
\eeq

\noi Applying the inequality in (\ref{eq:tech:xx}) with $\x = \x^t-\ddot{\x}$ and $\d = \x^{t+1}-\x^t$, we have:
\beq\label{eq:tech:xx:00:2}
\frac{1}{n}\|\x^{t+1}-\ddot{\x}\|_{\bar{\c}}^2 + \frac{1}{n}(n-1) \|\x^t - \ddot{\x}\|_{\bar{\c}}^2 = \frac{1}{n}\sum_{i=1}^n \|(\x^t -\ddot{\x}) + (\x^{t+1}-\x^t)_i\ei \|_{\bar{\c}}^2.
\eeq
\noi Combining the two inequalities in (\ref{eq:tech:xx:00:1}) and (\ref{eq:tech:xx:00:2}), we have:
\beq
\E[\|\x^{t+1} - \ddot{\x}\|_{\bar{\c}}^2 ] = \frac{1}{n}\|\x^{t+1}-\ddot{\x}\|_{\bar{\c}}^2 + \frac{1}{n}(n-1) \|\x^t - \ddot{\x}\|_{\bar{\c}}^2. \label{eq:tech:xx:2ddd}
\eeq

\noi Third, since $z(\x) \triangleq -g(\x)$ is globally $\rho$-bounded nonconvex, we have $$\forall \x,\y,~-g({\x})   \leq -g(\y)  - \la {\x} - \y,~ \partial g({\x}) \ra + \frac{\rho}{2}\|{\x} - \y\|_2^2.$$
\noi Applying this inequality with $\y = \x + \d$, we have:
\beq\label{eq:dadsafddf}
\forall \x,\forall \d,-g({\x}) + g({\x} + \d)  - \frac{\rho}{2}\|\d\|_2^2  &\leq&    \la   \d,~ \partial g({\x}) \ra   \nn\\
&\overset{(a)}{ \leq }& \sum_{i=1}^n g({\x} + \d_i \ei) - g({\x})  -(n-1)g({\x}) \nn\\
&= & \sum_{i=1}^n g({\x} + \d_i \ei)   - n g({\x}) ,
\eeq
\noi where step $(a)$ uses (\ref{eq:tech:gg0}) in Lemma \ref{lemma:tech}.

\noi We apply (\ref{eq:tech:hh}), (\ref{eq:tech:ff}), and (\ref{eq:dadsafddf}) with $\x=\x^t$ and $\d=\x^{t+1}-\x^t$, and obtain the following inequalities:
\beq
\E[ h(\x^{t+1})] &= &  \frac{1}{n}h(\x^{t+1}) + (1- \frac{1}{n}) h(\x^t)   \label{eq:111:hhh}  \\
\E[ f(\x^{t+1})]  &\leq&  \frac{1}{n}f(\x^{t}) + \frac{1}{n}\la \nabla f(\x^t),\x^{t+1}-\x^t  \ra  + \frac{1}{2n}\|\x^{t+1}-\x^{t}\|_{\c}^2 + (1-\frac{1}{n})f(\x^t) \label{eq:111:fff}\\
-\E[ g(\x^{t+1})]  &\leq&     - \frac{1}{n} g(\x^{t+1}) + \tfrac{\rho}{2n}\|\x^{t+1}-\x^t\|_2^2 - (1-\frac{1}{n}) g(\x^t). \label{eq:111:ggg}
\eeq

\textbf{(a)} We derive the following results:
\beq \label{eq:nonconvex:conv:0}
&&\E[\frac{1}{2}\|\x^{t+1}-\ddot{\x}\|_{\bar{\c}}^2] -   \frac{1}{2}\|\x^{t}-\ddot{\x}\|_{\bar{\c}}^2  \nn\\
&\overset{(a)}{ = }&  \E[\la \x^{t+1}-\ddot{\x},\bar{\c} \odot (\x^{t+1}- \x^t)  \ra ] - \E[\frac{1}{2}\|\x^{t+1}- \x^t\|_{\bar{\c}}^2 ] \nn\\
&\overset{(b)}{ = }& \E[\la \x^{t+1}-\ddot{\x}, -(\nabla_{i^t} f(\x^t) +  \partial_{i^t} h(\x^{t+1}) - \partial_{i^t} g(\x^{t+1}) )\cdot \eit  \ra ]- \E[\frac{1}{2}\|\x^{t+1}- \x^t\|_{\bar{\c}}^2 ]\nn\\
&\overset{(c)}{ = }& \frac{1}{n} \la \ddot{\x}- \x^{t+1}, \nabla f(\x^t) +  \partial h(\x^{t+1}) - \partial g(\x^{t+1}) \ra -  \frac{1}{2n}\| \x^{t+1}- \x^t \|_{\bar{\c}}^2,
\eeq
\noi where step $(a)$ uses uses the Pythagoras relation that: $\forall \bar{\c},\x,\y, \z,\frac{1}{2}\|\y-\z\|_{\bar{\c}}^2 - \frac{1}{2}\|\x-\z\|_{\bar{\c}}^2 = \la \y-\z,\bar{\c}\odot (\y - \x) \ra- \frac{1}{2}\|\x-\y\|_{\bar{\c}}^2$; step $(b)$ uses the optimality condition in (\ref{eq:nonconvex:cd:optimality}); step $(c)$ uses the fact that $\E[\la \x_{i^t} \eit, \y\ra] = \frac{1}{n} \sum_{j=1}^n  \x_j \y_j = \frac{1}{n}\la \x,\y \ra$.

\noi We now bound the term $\frac{1}{n}\la \ddot{\x}- \x^{t+1},  - \partial g(\x^{t+1}) \ra$ in (\ref{eq:nonconvex:conv:0}) by the following inequalities:
\beq \label{eq:nonconvex:conv:1}
&& \frac{1}{n}\la \ddot{\x}- \x^{t+1},   - \partial g(\x^{t+1}) \ra \nn\\
& \overset{(a)}{ \leq } &   - \frac{1}{n} g(\ddot{\x}) + \frac{1}{n} g(\x^{t+1}) +   \frac{{\rho}}{2n}\|\ddot{\x}-\x^{t+1}\|_{2}^2\nn\\
& \overset{(b)}{ \leq } &   - \frac{1}{n} g(\ddot{\x}) + \frac{1}{n} g(\x^{t+1}) +   \frac{\bar{\rho}}{2n}\|\ddot{\x}-\x^{t+1}\|_{\bar{\c}}^2\nn\\
& \overset{(c)}{ = } &   - \frac{1}{n} g(\ddot{\x}) + \frac{1}{n} g(\x^{t+1}) +  \frac{\bar{\rho}}{2} \left( \E[\|\ddot{\x}-\x^{t+1}\|_{\bar{\c}}^2] -(1-\frac{1}{n})\|\ddot{\x}-\x^{t}\|_{\bar{\c}}^2 \right)\nn\\
& \overset{(d)}{ \leq } &   - \frac{1}{n} g(\ddot{\x}) + \E[g(\x^{t+1})] + \frac{\rho}{2n}\|\x^{t+1}-\x^t\|_2^2 - (1-\frac{1}{n})g(\x^t) +  \frac{\bar{\rho}}{2} \left( \E[\|\ddot{\x}-\x^{t+1}\|_{\bar{\c}}^2] -(1-\frac{1}{n})\|\ddot{\x}-\x^{t}\|_{\bar{\c}}^2 \right),~~~~~~~
\eeq
\noi where step $(a)$ uses the globally $\rho$-bounded nonconvexity of $-g(\cdot)$; step $(b)$ uses the fact that $\rho \|\v\|_2^2 \leq \rho \cdot \frac{1}{\min(\bar{\c})}\|\v\|_{\bar{\c}}^2 = \bar{\rho}\|\v\|_{\bar{\c}}^2$ for all $\v$; step $(c)$ uses (\ref{eq:tech:xx:2ddd}); step $(d)$ uses (\ref{eq:111:ggg}).

\noi We now bound the term $\frac{1}{n} \la \ddot{\x}- \x^{t+1}, \nabla f(\x^t) +  \partial h(\x^{t+1})\ra$ in (\ref{eq:nonconvex:conv:0}) by the following inequalities:
\beq \label{eq:nonconvex:conv:2}
&&\frac{1}{n} \la \ddot{\x}- \x^{t+1}, \nabla f(\x^t) +  \partial h(\x^{t+1})\ra \nn\\
& = & \left( \frac{1}{n} \la \ddot{\x} - \x^{t}, \nabla f(\x^t) \ra  +   \frac{1}{n} \la \x^{t} - \x^{t+1}, \nabla f(\x^t) \ra \right) + \frac{1}{n} \la \ddot{\x} - \x^{t+1}, \partial h(\x^{t+1})\ra  \nn\\
& \overset{(a)}{ \leq } & \frac{1}{n} \left(f(\ddot{\x}) - f(\x^{t})\right)  - \frac{1}{n} \la \x^{t+1}-\x^{t}  , \nabla f(\x^t) \ra + \frac{1}{n} \left(  h(\ddot{\x}) - h(\x^{t+1}) \right) \nn\\
& \overset{(b)}{ \leq } & \frac{1}{n} f(\ddot{\x}) -\E[ f(\x^{t+1})]  + \frac{1}{2n}\|\x^{t+1}-\x^{t}\|_{\c}^2 + (1-\frac{1}{n})f(\x^t) + \frac{1}{n} \left(  h(\ddot{\x}) - h(\x^{t+1}) \right) \nn\\
& \overset{(c)}{ \leq } & \frac{1}{n} f(\ddot{\x})   -\E[ f(\x^{t+1})]  + \tfrac{1}{2n}\|\x^{t+1}-\x^{t}\|_{\c}^2 + (1-\frac{1}{n})f(\x^t) + \frac{1}{n}  h(\ddot{\x})   -\E[ h(\x^{t+1})]  + (1- \frac{1}{n}) h(\x^t),
\eeq
\noi where step $(a)$ uses the convexity of $f(\cdot)$ and $h(\cdot)$ that:
\beq
\la   \ddot{\x} - \x^t ,\nabla f(\x^t) \ra &\leq& f(\ddot{\x})    - f(\x^{t})  ,\nn \\
\la  \ddot{\x} - \x^{t+1}, \partial h(\x^{t+1}) \ra &\leq &h(\ddot{\x})    - h(\x^{t+1}) ; \nn
\eeq
\noi step $(b)$ uses (\ref{eq:111:fff}); step $(c)$ uses (\ref{eq:111:hhh}).

\noi Combining (\ref{eq:nonconvex:conv:0}), (\ref{eq:nonconvex:conv:1}), and (\ref{eq:nonconvex:conv:2}) together, and using the fact that $F(\x) = f(\x)+h(\x)-g(\x)$, we obtain:
\beq \label{eq:nonconvex:conv:3}
&& \E[\frac{1-\bar{\rho}}{2}\|\x^{t+1}-\ddot{\x}\|_{\bar{\c}}^2 ]  -  \frac{1 - \bar{\rho} + \tfrac{\bar{\rho}}{n} }{2}\|\x^{t}-\ddot{\x}\|_{\bar{\c}}^2   \nn\\
& \leq &   \frac{1}{n} ( F(\ddot{\x}) - F(\x^t))- \E[F(\x^{t+1})] + F(\x^t)    + \E[\frac{\rho}{2n} \|\x^{t+1}-\x^{t}\|_2^2] \nn\\
& \overset{(a)}{ \leq } &   \frac{1}{n} ( F(\ddot{\x}) - F(\x^t))- \E[F(\x^{t+1})] + F(\x^t)    + \frac{\rho}{\theta} \E[F(\x^{t}) - F(\x^{t+1})] \nn\\
& \overset{(b)}{ = } &   - \frac{1}{n}  \ddot{q}^t + (1+ \frac{\rho}{\theta}) ( \ddot{q}^t - \E[\ddot{q}^{t+1}]),
\eeq
\noi where step $(a)$ uses the sufficient decrease condition that $\E[\frac{1}{2n}\|\x^{t+1}-\x^{t}\|_2^2]\leq \frac{1}{\theta} \E[F(\x^t)-F(\x^{t+1})]$; step $(b)$ uses the definition of $\ddot{q}^t \triangleq F(\x^t)- F(\ddot{\x})$ and the fact that $F(\x^t)-F(\x^{t+1})=\ddot{q}^t-\ddot{q}^{t+1}$. Using the definitions that $\ddot{r}^{t+1}\triangleq \frac{1}{2}\|\x^{t}-\ddot{\x}\|_{\bar{\c}}^2$, $\varpi  \triangleq 1-\bar{\rho}$, and $\gamma \triangleq (1+\frac{\rho}{\theta})$, we rewrite (\ref{eq:nonconvex:conv:3}) as:
\beq \label{eq:nonconvex:conv:4}
&&\E[{(1-\bar{\rho})}\ddot{r}^{t+1}]   - ({1 -\bar{\rho}}+\tfrac{\bar{\rho}}{n}) \ddot{r}^{t}   \leq    - \frac{1}{n}  \ddot{q}^t + {(1+ \frac{\rho}{\theta})} ( \ddot{q}^t -\E[\ddot{q}^{t+1}]) \nn\\
&\Leftrightarrow&\varpi \E[\ddot{r}^{t+1}] + \gamma \E[\ddot{q}^{t+1}]    \leq (\varpi +\frac{\bar{\rho} }{n} ) \ddot{r}^{t} + ( \gamma- \frac{1}{n})\ddot{q}^t.
\eeq

\textbf{(b)} We now discuss the situation when $\varpi\geq 0$. We notice that the function $\mathcal{S}_{i^t}(\x,\eta)+ h_{i^t}(\x + \eta \eit) +   \tfrac{\theta}{2} \| (\x+\eta \eit) - \x \|_2^2$ is ($\min(\c) + \theta$)-strongly convex w.r.t. $\eta$ and the term $- g(\x + \eta \eit)$ is globally $\rho$-bounded nonconvex w.r.t. $\eta$ for all $t$. Therefore, $\mathcal{M}_{i^t}(\x^t,\eta)$ in (\ref{eq:subprob:nonconvex}) is convex if:
\beq
\min(\c) + \theta - \rho\geq 0 \quad \Leftrightarrow \quad \varpi\geq 0.\nn
\eeq

\noi We now discuss the case when $F(\cdot)$ satisfies the Luo-Tseng error bound assumption. We bound the term $\ddot{r}^t$ using the following inequalities:
\beq\label{eq:conv:rate:tsengluo:1}
\ddot{r}^t\nn & \triangleq &  \max(\bar{\c}) \frac{1}{2}\|\x^t - \ddot{\x}\|_{2}^2 \nn\\
&\overset{(a)}{ \leq } & \max(\bar{\c}) \frac{1}{2}\frac{\delta^2}{n^2} (   \sum_{i=1}^n |\bar{\mathcal{M}}_i(\x^t)|)^2 \nn\\
&\overset{(b)}{ \leq } & \max(\bar{\c}) \frac{1}{2}\frac{\delta^2}{n^2} n \cdot (   \sum_{i=1}^n |\bar{\mathcal{M}}_i(\x^t)|^2) \nn\\
&\overset{(c)}{ \leq } & \max(\bar{\c}) \frac{1}{2}\frac{\delta^2}{n^2} n \cdot \left(n \E[\|\x^{t+1}-\x^t\|^2_2]\right)  \nn\\
&\overset{}{ = } &  \max(\bar{\c})  \frac{\delta^2 }{\theta} \cdot \frac{\theta}{2} \E[\|\x^{t+1}-\x^t\|^2_2] \nn\\
&\overset{(d)}{ \leq } &  \max(\bar{\c})\delta^2 \frac{n}{\theta}(F(\x^t)-\E[F(\x^{t+1})]) \nn\\
& = &  \max(\bar{\c})\delta^2 \frac{n}{\theta}(\ddot{q}^t-\E[\ddot{q}^{t+1}]),\nn\\
& \overset{(e)}{ = } &  \kappa_0 n (\ddot{q}^t-\E[\ddot{q}^{t+1}]),
\eeq
\noi where step $(a)$ uses Assumption \ref{ass:luo} that $\|\x^t - \ddot{\x}\|_{2}^2\leq \delta^2 ( \frac{1}{n} \sum_{i=1}^n| \text{dist}(0,\bar{\mathcal{M}}_{i}(\x))|)^2$ for any coordinate-wise stationary point $\ddot{\x}$; step $(b)$ uses the fact that $\forall \x,~\|\x\|_1^2\leq n \|\x\|_2^2$; step $(c)$ uses the fact that $\E[\|\x^{t+1}-\x^t\|^2_2] = \E[\|(\x^t + \bar{\mathcal{M}}_{i^t}(\x^t)\eit )-\x^t\|^2_2] =  \E[|\bar{\mathcal{M}}_{i^t}(\x^t)|^2] = \frac{1}{n}\sum_{i=1}^n |\bar{\mathcal{M}}_i(\x^t)|^2$; step $(d)$ uses the sufficient decrease condition that $\E[\frac{1}{2n}\|\x^{t+1}-\x^{t}\|_2^2]\leq \frac{1}{\theta} \E[F(\x^t)-F(\x^{t+1})]$; step $(e)$ uses the definition of $\kappa_0 \triangleq \max(\bar{\c}) \frac{ \delta^2}{\theta} $. Since $\varpi\geq0$, we have form (\ref{eq:nonconvex:conv:4}):
\beq
&&\gamma \E[\ddot{q}^{t+1}]    \leq \left(\varpi + \frac{\bar{\rho} }{n} \right)\kappa_0n(\ddot{q}^t-\E[\ddot{q}^{t+1}]) + (\gamma - \frac{1}{n} ) \ddot{q}^{t} \nn\\
&\Rightarrow &\underbrace{ ( \gamma + \kappa_0n ( \varpi + \frac{\bar{\rho} }{n} )  )}_{\triangleq \kappa_1} \E[\ddot{q}^{t+1}] \leq ( \underbrace{\gamma+ \kappa_0 n( \varpi + \frac{\bar{\rho} }{n} )}_{ = \kappa_1 } - \frac{1}{n})\ddot{q}^t\nn\\
&\Rightarrow & \E[\ddot{q}^{t+1}] \leq  \frac{\kappa_1 - \frac{1}{n}}{\kappa_1} \ddot{q}^t\nn\\
&\Rightarrow & \E[\ddot{q}^{t+1}] \leq  (\frac{\kappa_1 - \frac{1}{n}}{\kappa_1})^{t+1} \ddot{q}^0.\nn
\eeq

\noi Thus, we finish the proof of this theorem.

\if
We now discuss the case when $r^t \leq \varepsilon$ with $\varepsilon$ being sufficiently small such that $\varepsilon \leq \frac{2}{\max(\bar{\c})} (\frac{\mu \chi}{ \bar{\rho}  })^2$. We first obtain the following inequalities:
\beq
r^t \triangleq \frac{1}{2}\|\x^{t} - \ddot{\x}\|_{\bar{\c}}^2 \overset{(a)}{ \leq } \frac{\max(\bar{\c})}{2}\|\x^{t} - \ddot{\x}\|_{2}^2 \overset{(b)}{ \leq } \frac{\max(\bar{\c})}{2} \cdot ( \frac{1}{\mu} (F(\x^t) - F(\ddot{\x}))  )^2 \overset{(c)}{ = } \frac{\max(\bar{\c})}{2} \cdot ( \frac{q^t}{\mu}  )^2, \label{eq:nonconvex:conv:5}
\eeq
\noi where step $(a)$ uses the fact that $\tfrac{1}{\max(\bar{\c})} \|\x\|_{\bar{\c}}^2 \leq \|\x\|_2^2 \leq \tfrac{1}{\min(\bar{\c})} \|\x\|_{\bar{\c}}^2,\forall \x$; step $(b)$ uses the sharpness of $F(\cdot)$; step $(c)$ uses the definition that $q^t = F(\x^t) - F(\ddot{\x})$.

We define
\beq\label{eq:nonconvex:conv:6}
A  \triangleq \frac{\chi}{ \bar{\rho} },~\kappa \triangleq \frac{ \frac{1}{n} - \frac{\bar{\rho} A}{n}}{ A \varpi + \gamma }
\eeq
\noi and obtain:
\beq\label{eq:nonconvex:conv:7}
r^t = \sqrt{r^t} \cdot \sqrt{r^t} \overset{(a)}{ \leq } \frac{\sqrt{2}}{\sqrt{\max(\bar{\c})}} \frac{\mu \chi}{ \bar{\rho}} \sqrt{r^t} \overset{(b)}{ \leq }  \frac{ \chi}{ \bar{\rho} } q^t = A q^t
\eeq
\noi where step $(a)$ uses the condition that $r^t\leq \frac{2}{\max(\bar{\c})}  (\frac{\mu \chi}{ \bar{\rho}  })^2$; step $(b)$ uses (\ref{eq:nonconvex:conv:5}).

\noi Based on (\ref{eq:nonconvex:conv:4}), we derive the following inequalities:
\beq\label{eq:rec:conv:rate:nonconvex}
\varpi \E[r^{t+1}] + \gamma \E[q^{t+1}]   & \leq& (\varpi r^{t} + \gamma q^{t}) + \frac{\bar{\rho} }{n} r^t - \frac{q^t}{n} \nn\\
&\overset{(a)}{ \leq } & (\varpi r^{t} + \gamma q^{t}) + \frac{\bar{\rho} }{n} A q^t - \frac{q^t}{n} \nn\\
&\overset{(b)}{ \leq } & (\varpi r^{t} + \gamma q^{t}) - \kappa (\varpi Aq^{t} + \gamma q^{t}) \nn\\
&\overset{(c)}{ \leq } & (\varpi r^{t} + \gamma q^{t}) - \kappa (\varpi r^{t} + \gamma q^{t}) \nn\\
&\overset{}{ = } & (1-\kappa)(\varpi r^{t} + \gamma q^{t})
\eeq
\noi where step $(a)$ uses $r^t \leq A q^t$ as shown in (\ref{eq:nonconvex:conv:7}); step $(b)$ uses the definition of $\kappa$ in (\ref{eq:nonconvex:conv:6}); step $(c)$ uses (\ref{eq:nonconvex:conv:7}) again. Note that the parameter $\kappa$ in (\ref{eq:nonconvex:conv:6}) can be simplified as:
\beq
\kappa = \frac{ \frac{1}{n} -  \frac{\chi}{n}}{ \frac{\chi}{ \bar{\rho} } \varpi + \gamma } = \frac{ 1 -  \chi}{ n(\gamma+ \frac{1-\bar{\rho}}{ \bar{\rho} } \chi  ) } \nn
\eeq

\noi Solving the recursive formulation as in (\ref{eq:rec:conv:rate}), we have:
\beq
\E[\varpi r^{t}  + \gamma q^{t}] &\leq&  (1-\kappa)^t \cdot ( \varpi r^{0} + \gamma q^0).\nn
\eeq

\noi Using the fact that $q^t\geq 0$, we obtain the following result:
\beq
\E[r^{t}]        &\leq&  (1-\kappa)^t \cdot (  r^{0} + \tfrac{\gamma}{\varpi}  q^0), \nn
\eeq
\fi

\end{proof}

\subsection{Proof for Theorem \ref{the:rate:cvx}}

\begin{proof}

Let $\check{\x}$ be any coordinate-wise stationary point. First, the optimality condition for the convex subproblem as in (\ref{eq:subprob:convex}) can be written as:
\beq \label{eq:convex:cd:optimality}
0 \in \nabla_{i^t} f(\x^t) +  \bar{\c}_{i^t}  \bar{\eta}^t + \partial_{i^t} h(\x^{t+1}) - \partial_{i^t} g(\x^{t}).
\eeq

\noi Second, we apply (\ref{eq:tech:hh}), (\ref{eq:tech:ff}), and (\ref{eq:tech:gg0}) in Lemma \ref{lemma:tech} with $\x=\x^t$ and $\d=\x^{t+1}-\x^t$, and obtain the following inequalities:
\beq
\E[ h(\x^{t+1})] &= &  \frac{1}{n}h(\x^{t+1}) + (1- \frac{1}{n}) h(\x^t)   \label{eq:tech:hh:3}  \\
\E[ f(\x^{t+1})]  &\leq&  \frac{1}{n}f(\x^{t}) + \frac{1}{n}\la \nabla f(\x^t),\x^{t+1}-\x^t  \ra  + \tfrac{1}{2n}\|\x^{t+1}-\x^{t}\|_{\c}^2 + (1-\frac{1}{n})f(\x^t) \label{eq:tech:ff:3}\\
-\E[ g(\x^{t+1})]  &\leq&     - \frac{1}{n} \la\partial g(\x^t),\x^{t+1}-\x^{t}  \ra - g(\x^t).  \label{eq:tech:gg:3}
\eeq

\textbf{(a)} We derive the following results:
\beq \label{eq:convex:conv:0}
&&\E[\frac{1}{2}\|\x^{t+1}-\check{\x}\|_{\bar{\c}}^2] -   \frac{1}{2}\|\x^{t}-\check{\x}\|_{\bar{\c}}^2  \nn\\
&\overset{(a)}{ = }&  \E[\la \x^{t+1}-\check{\x},\bar{\c} \odot (\x^{t+1}- \x^t)  \ra ] - \E[\frac{1}{2}\|\x^{t+1}- \x^t\|_{\bar{\c}}^2 ] \nn\\
&\overset{(b)}{ = }& \E[\la \x^{t+1}-\check{\x}, -(\nabla_{i^t} f(\x^t) +  \partial_{i^t} h(\x^{t+1}) - \partial_{i^t} g(\x^{t}) )\cdot \eit  \ra ]- \E[\frac{1}{2}\|\x^{t+1}- \x^t\|_{\bar{\c}}^2 ]\nn\\
&\overset{(c)}{ = }& \frac{1}{n} \la \check{\x}- \x^{t+1}, \nabla f(\x^t) +  \partial h(\x^{t+1}) - \partial g(\x^{t}) \ra -  \frac{1}{2n}\| \x^{t+1}- \x^t \|_{\bar{\c}}^2,
\eeq
\noi where step $(a)$ uses uses the Pythagoras relation that: $\forall \bar{\c},\x,\y, \z,\frac{1}{2}\|\y-\z\|_{\bar{\c}}^2 - \frac{1}{2}\|\x-\z\|_{\bar{\c}}^2 = \la \y-\z,\bar{\c}\odot (\y - \x) \ra- \frac{1}{2}\|\x-\y\|_{\bar{\c}}^2$; step $(b)$ uses the optimality condition in (\ref{eq:convex:cd:optimality}); step $(c)$ uses the fact that $\E[\la \x_{i^t} \eit, \y\ra] = \frac{1}{n} \sum_{j=1}^n  \x_j \y_j = \frac{1}{n}\la \x,\y \ra$.

\noi We now bound the term $\frac{1}{n}\la \check{\x}- \x^{t+1},   - \partial g(\x^{t}) \ra$ in (\ref{eq:convex:conv:0}) by the following inequalities:
\beq \label{eq:convex:conv:1}
&& \frac{1}{n}\la \check{\x}- \x^{t+1},   - \partial g(\x^{t}) \ra \nn\\
& \overset{}{ = }& \frac{1}{n}\la \check{\x}- \x^{t},   - \partial g(\x^{t}) \ra + \frac{1}{n}\la  \x^t - \x^{t+1},   - \partial g(\x^{t}) \ra \nn\\
& \overset{(a)}{ \leq } &   - \frac{1}{n} g(\check{\x}) + \frac{1}{n} g(\x^{t}) +   \frac{\rho}{2n}\|\check{\x}-\x^{t}\|_2^2 + \frac{1}{n}\la  \x^{t+1} - \x^{t},    \partial g(\x^{t}) \nn\\ & \overset{(b)}{ \leq } &   - \frac{1}{n} g(\check{\x}) + \frac{1}{n} g(\x^{t}) +   \frac{\rho}{2n}\|\check{\x}-\x^{t}\|_2^2 - g(\x^t) + \E[ g(\x^{t+1})],
 \eeq
\noi where step $(a)$ uses the globally $\rho$-bounded nonconvexity of $-g(\cdot)$; step $(b)$ uses (\ref{eq:tech:gg:3}).

\noi We now bound the term $\frac{1}{n} \la \check{\x}- \x^{t+1}, \nabla f(\x^t) +  \partial h(\x^{t+1})\ra$ in (\ref{eq:convex:conv:0}) by the following inequalities:
\beq \label{eq:convex:conv:2}
&&\frac{1}{n} \la \check{\x}- \x^{t+1}, \nabla f(\x^t) +  \partial h(\x^{t+1})\ra \nn\\
& = & \frac{1}{n} \la \check{\x} - \x^{t}, \nabla f(\x^t) \ra + \frac{1}{n} \la \check{\x} - \x^{t+1}, \partial h(\x^{t+1})\ra  +   \frac{1}{n} \la \x^{t} - \x^{t+1}, \nabla f(\x^t) \ra   \nn\\
& \overset{(a)}{ \leq } & \frac{1}{n} [f(\check{\x}) - f(\x^{t}) + h(\check{\x}) - h(\x^{t+1}) ]  - \frac{1}{n} \la \x^{t+1}-\x^{t}  , \nabla f(\x^t) \ra  \nn\\
& \overset{(b)}{ \leq } & \frac{1}{n} [f(\check{\x}) + h(\check{\x}) - h(\x^{t+1}) ] -\E[ f(\x^{t+1})]  + \tfrac{1}{2n}\|\x^{t+1}-\x^{t}\|_{\c}^2 + (1-\frac{1}{n})f(\x^t) \nn\\
& \overset{(c)}{ \leq } & \frac{1}{n} [f(\check{\x}) + h(\check{\x})   ] -\E[ h(\x^{t+1})]  + (1- \frac{1}{n}) h(\x^t)  -\E[ f(\x^{t+1})]  + \tfrac{1}{2n}\|\x^{t+1}-\x^{t}\|_{\c}^2 + (1-\frac{1}{n})f(\x^t),
\eeq
\noi where step $(a)$ uses the convexity of $f(\cdot)$ and $h(\cdot)$ that:
\beq
\la   \check{\x} - \x^t ,\nabla f(\x^t) \ra &\leq& f(\check{\x})    - f(\x^{t})  ,\nn \\
\la  \check{\x} - \x^{t+1}, \partial h(\x^{t+1}) \ra &\leq &h(\check{\x})    - h(\x^{t+1}) ; \nn
\eeq
\noi step $(b)$ uses (\ref{eq:tech:ff:3}); step $(c)$ uses (\ref{eq:tech:hh:3}).

\noi Combining (\ref{eq:convex:conv:0}), (\ref{eq:convex:conv:1}), (\ref{eq:convex:conv:2}), and using the fact that $F(\x) = f(\x)+h(\x)-g(\x)$, we obtain:
\beq \label{eq:convex:conv:3}
&& \E[\frac{1}{2}\|\x^{t+1}-\check{\x}\|_{\bar{\c}}^2 ]  - \frac{1}{2}\|\x^{t}-\check{\x}\|_{\bar{\c}}^2  \nn\\
& \overset{}{ \leq } & \frac{\rho}{2n}\|\x^{t}-\check{\x}\|_{2}^2 + \frac{1}{n}F(\check{\x}) - \E[ F(\x^{t+1})] + (1-\frac{1}{n})F(\x^t)\nn\\
& \overset{(a)}{ \leq } & \frac{\bar{\rho}}{2n}\|\x^{t}-\check{\x}\|_{\bar{\c}}^2 + \frac{1}{n}F(\check{\x}) - \E[ F(\x^{t+1})] + (1-\frac{1}{n})F(\x^t),
\eeq
\noi where step $(a)$ uses $\|\x\|_2^2 \leq \tfrac{1}{\min(\bar{\c})} \|\x\|_{\bar{\c}}^2,\forall \x$. The inequality in (\ref{eq:convex:conv:3}) can be rewritten as:
\beq \label{eq:convex:conv:4}
\E[\check{r}^{t+1}] + \E[\check{q}^{t+1}]  \leq \check{r}^t  + \frac{\bar{\rho}}{n} \check{r}^{t}  - \frac{1}{n} \check{q}^t + \check{q}^t.
\eeq

\textbf{(b)} We now discuss the case when $F(\cdot)$ satisfies the Luo-Tseng error bound assumption. We first bound the term $\check{r}^t$:
\beq
\check{r}^t & \triangleq &  \max(\bar{\c}) \frac{1}{2}\|\x^t - \check{\x}\|_{2}^2 \nn\\
&\overset{(a)}{ \leq } & \max(\bar{\c}) \frac{1}{2}\frac{\delta^2}{n^2} (   \sum_{i=1}^n |\bar{\mathcal{P}}_i(\x^t)|)^2 \nn\\
&\overset{(b)}{ \leq } &   \kappa_0 n (\check{q}^t-\E[\check{q}^{t+1}]) ,\label{eq:conv:rate:tsengluo:2}\nn
\eeq
\noi where step $(a)$ uses Assumption \ref{ass:luo} with the residual function defining as $\mathcal{R}(\x)= \frac{1}{n} \sum_{i=1}^n|\bar{\mathcal{P}}_{i}(\x)|$; step $(b)$ uses the same strategy as in deriving the results in (\ref{eq:conv:rate:tsengluo:1}). Finally, we have form (\ref{eq:convex:conv:4}):
\beq
&&\E[\check{q}^{t+1}]    \leq \left(1 + \frac{\bar{\rho} }{n} \right)  \kappa_0 n (\check{q}^t-\E[\check{q}^{t+1}]) + (1 - \frac{1}{n} ) \check{q}^{t} \nn\\
&\Rightarrow & \underbrace {(1 +  (1 + \frac{\bar{\rho} }{n} )  \kappa_0 n )}_{ \triangleq \kappa_2}\E[\check{q}^{t+1}]    \leq  ( \underbrace{(1 + \frac{\bar{\rho} }{n} )  \kappa_0 n + 1}_{=\kappa_2} - \frac{1}{n} ) \check{q}^{t} \nn\\
&\Rightarrow & \E[\check{q}^{t+1}] \leq  \frac{\kappa_2 - \frac{1}{n}}{\kappa_2} \check{q}^t\nn\\
&\Rightarrow & \E[\check{q}^{t+1}] \leq  (\frac{\kappa_2 - \frac{1}{n}}{\kappa_2})^{t+1} \check{q}^0.\nn
\eeq

\noi Thus, we finish the proof of this theorem.

\if

We now discuss the case when $r^t \leq \varepsilon$ with $\varepsilon$ being sufficiently small such that $\varepsilon \leq \frac{2}{\max(\bar{\c})} (\frac{\mu \chi}{ \bar{\rho}  })^2$. We first obtain the following inequality:
\beq
r^t \leq \frac{\max(\bar{\c})}{2} \cdot ( \frac{q^t}{\mu}  )^2, \label{eq:convex:conv:5}
\eeq
\noi which is the same as (\ref{eq:nonconvex:conv:5}).

We define
\beq\label{eq:convex:conv:6}
A \triangleq \frac{\chi}{\bar{\rho}},~\kappa \triangleq\frac{ \frac{1}{n} - \frac{\bar{\rho} A}{n}}{ A  + 1 }
\eeq
\noi and obtain:
\beq\label{eq:convex:conv:7}
r^t = \sqrt{r^t} \cdot \sqrt{r^t} \overset{(a)}{ \leq } \frac{\sqrt{2}}{\sqrt{\max(\bar{\c})}} \frac{\mu \chi}{ \bar{\rho}} \sqrt{r^t} \overset{(b)}{ \leq }  \frac{ \chi}{ \bar{\rho} } q^t = A q^t
\eeq
\noi where step $(a)$ uses the condition that $r^t\leq \frac{2}{\max(\bar{\c})}  (\frac{\mu \chi}{ \bar{\rho}  })^2$; step $(b)$ uses (\ref{eq:convex:conv:5}).

\noi Based on (\ref{eq:convex:conv:4}), we derive the following inequalities:
\beq\label{eq:rec:conv:rate}
\E[r^{t+1}] +  \E[q^{t+1}]   & \leq& ( r^{t} +  q^{t}) + \frac{\bar{\rho} }{n} r^t - \frac{q^t}{n} \nn\\
&\overset{(a)}{ \leq } & ( r^{t} +  q^{t}) + \frac{\bar{\rho} }{n} A q^t - \frac{q^t}{n} \nn\\
&\overset{(b)}{ \leq } & ( r^{t} +  q^{t}) - \kappa ( Aq^{t} +  q^{t}) \nn\\
&\overset{(c)}{ \leq } & ( r^{t} +  q^{t}) - \kappa ( r^{t} +  q^{t}) \nn\\
&\overset{}{ = } & (1-\kappa)( r^{t} +  q^{t})
\eeq
\noi where step $(a)$ uses $r^t \leq A q^t$ as shown in (\ref{eq:convex:conv:7}); step $(b)$ uses the definition of $\kappa$ in (\ref{eq:convex:conv:6}); step $(c)$ uses (\ref{eq:convex:conv:7}) again. Note that the parameter $\kappa$ in (\ref{eq:convex:conv:6}) can be simplified as:
\beq
\kappa = \frac{ \frac{1}{n} - \frac{\bar{\rho} A}{n}}{ A  + 1 } = \frac{ 1 - \chi}{n(1+ \frac{\chi}{\bar{\rho}}  ) } \nn
\eeq

\noi Solving the recursive formulation as in (\ref{eq:rec:conv:rate}), we have:
\beq
\E[ r^{t}  +  q^{t}] &\leq&  (1-\kappa)^t \cdot (  r^{0} +  q^0).\nn
\eeq

\noi Using the fact that $q^{t}\geq0$, we obtain the following result:
\beq
\E[r^{t}]        &\leq&  (1-\kappa)^t \cdot (  r^{0} +   q^0), \nn
\eeq

\fi

\end{proof}

\section{More Examples of the Breakpoint Searching Method for Proximal Operator Computation}
\label{sect:more:break}
\subsection{When $g(\y) =  \|\A\y\|_{\infty}$ and $h_i(\cdot) \triangleq 0$}
Consider the problem: $\min_{\eta} \tfrac{a}{2}\eta^2 + b\eta-\|\A(\x+\eta e_i)\|_{\infty}$. It can be rewritten as: $\min_{\eta} \frac{a}{2}  \eta^2 + b \eta -  \|\g \eta + \d\|_{\infty}$. It is equivalent to $\min_\eta p(\eta) \triangleq \frac{a}{2} \eta^2 + b \eta -  \max_{i=1}^{2m}(\bar{\g}_i \eta + \bar{\d}_i)$ with $\bar{\g}=[\g_1,\g_2,...,\g_m,-\g_1,-\g_2,...,-\g_m]$ and $\bar{\d}=[\d_1,\d_2,...,\d_m,-\d_1,-\d_2,...,-\d_m]$. Setting the gradient of $p(\cdot)$ to zero yields: $a\eta + b + \bar{\g}_i=0$ with $i=1,2,...,(2m)$. We have $\etas = (-b-\bar{\g})/a$. Therefore, Problem (\ref{eq:prox}) contains $2m$ breakpoints $\Theta =\{\etas_1,\etas_2,...,\etas_{2m}\}$ for this example.

\subsection{When $g(\y) =  \|\max(0,\A\y)\|_1$ and $h_i(\cdot) \triangleq 0$} \label{eq:proximal:lmax0}
Consider the problem: $\min_\eta \tfrac{a}{2}\eta^2 + b \eta -  \|\max(0,\A(\x+\eta\ei))\|_1$. Using the fact that $\max(0,a) = \tfrac{1}{2}(a+|a|)$, we have the following equivalent problem: $\min_\eta \tfrac{a}{2}\eta^2 + b \eta  -  \frac{1}{2}  \la \mathbf{1}, \A\ei\ra \eta -  \frac{1}{2}\|\A(\x+\eta\ei)\|_1$. Therefore, the proximal operator of $g(\x) =  \|\max(0,\A\x)\|_1$ can be transformed to the proximal operator of $g(\x) =  \|\A\x\|_1$ with suitable parameters.

\subsection{When $g(\y) =  \|\A\y\|_2$ and $h_i(\cdot) \triangleq 0$} \label{eq:proximal:l2}

Consider the problem: $\min_\eta \tfrac{a}{2}\eta^2 + b \eta  -  \|\A(\x+\eta e_i)\|_p$. It can be rewritten as: $\min_\eta p(\eta)\triangleq \frac{a}{2} \eta^2 + b \eta -  \|\g \eta + \d\|_p$. Setting the gradient of $p(\cdot)$ to zero yields: $0 = a\eta + b -  \|\g \eta - \d\|_p^{1-p} \la \g, \text{sign}(\g \eta + \d) \odot |\g \eta +\d|^{p-1}\ra$. We only focus on $p=2$. We obtain: $a\eta + b = \frac{  \la \g, \g \eta + \d\ra }{  \|\g \eta - \d\| }  \Leftrightarrow \|\g \eta - \d\| ( a \eta + b) =  {  \la \g, \g \eta + \d\ra } \Leftrightarrow \|\g \eta - \d\|_2^2 (a \eta + b)^2 =   ( \la \g,\g \eta + \d\ra)^2 $. Solving this quartic equation we obtain all of its real roots $\{\etas_1,\etas_2,...,\etas_c\}$ with $1\leq c\leq 4$. Therefore, Problem (\ref{eq:prox}) at most contains $4$ breakpoints $\Theta =\{\etas_1,\etas_2,...,\etas_c\}$ for this example.

\section{More Experiments}
\label{sect:more:exp}

In this section, we present the experiment results of the approximate binary optimization problem and the generalized linear regression problem.

\subsection{Approximate Binary Optimization}

We consider Problem (\ref{eq:binary:opt}). We generate the observation vector via $\y = \max(0,\A\ddot{\x} + \text{randn}(m,1) \times 0.1\times\|\A\ddot{\x}\|)$ with $\ddot{\x}=\text{randn}(d,1)$. This problem is consistent with $f(\x)\triangleq \frac{1}{2}\|\G\x-\y\|_2^2,~g(\x) \triangleq - \rho\|\x\|$, and $h(\x) = \sum_{i}^n h_i(\x_i)$ with $h_i(z) \triangleq I_{[-1,1]}(z)$ where $I_{[-1,1]}(z)$ denotes an indicator function on the box constraint ($h_i(z) =0$ if $-1 \leq z\leq 1$, $+\infty$ otherwise). We notice that $\nabla f(\cdot)$ is $L$-Lipschitz with constant $L=\|\G\|_2^2$ and coordinate-wise Lipschitz with constant $\c=\text{diag}(\G^T\G)$. The subgradient of $g(\x)$ at $\x^t$ can be computed as: $\partial g(\x^t) = - \frac{\rho\x^t}{\|\x^t\|}\triangleq \g^t$. We set $\rho=5$.


We compare with the following methods. \textit{\textbf{(i)}} Multi-Stage Convex Relaxation (MSCR). It solves the following problem: $\x^{t+1} = \arg \min_{\x} \frac{1}{2}\|\G\x-\y\|_2^2 -   \la \x-\x^t,~\g^t\ra,s.t.\|\x\|_{\infty}\leq1$. This is essentially equivalent to the alternating minimization method in \cite{YuanG17}. \textit{\textbf{(ii)}} Proximal DC algorithm (PDCA). It considers the following iteration: $\x^{t+1} = \arg \min_{\x}\frac{L}{2}\|\x-\x^t\|_2^2 + \la \x - \x^t,~ \nabla f(\x^t)\ra  -  \la \x-\x^t,~\g^t\ra$. \textit{\textbf{(iii)}} Subgradient method (SubGrad). It uses the following iteration: $\x^{t+1} = \mathcal{P}_{\Omega}(\x^t  - \frac{0.1}{t} \cdot (\nabla f(\x) -  \g^t ))$, where $\mathcal{P}_{\Omega}(\x) \triangleq \max(-1,\min(\x,1))$ is the projection operation on the convex set $\Omega \triangleq \{\x| \|\x\|_{\infty} \leq 1\}$. \textit{\textbf{(iv)}} \textit{\textbf{CD-SCA}} solves a convex problem $\bar{\eta}^t = \arg \min_{\eta} 0.5 (\c_{i^t}+\theta)\eta^2   + [\nabla f(\x^t) - \g^t]_{i^t}  \eta,~s.t.~-1\leq \x^t_{i^t} +\eta\leq 1$ and update $\x^t$ via $\x_{i^t}^{t+1}  = \x_{i^t}^t  + \bar{\eta}^t$. \textit{\textbf{(v)}} \textit{\textbf{CD-SNCA}} computes the nonconvex proximal operator of $\ell_2$ norm (see Section \ref{eq:proximal:l2}) as $\bar{\eta}^t = \arg \min_{\eta} \frac{\c_i+\theta}{2}\eta^2 +  \nabla_{i^t} f(\x^t)  \eta  - \rho\| \x^t + \eta\ei\|,s.t.-1\leq \x_{i^t}^t+\eta\leq 1$ and update $\x^t$ via $\x_{i^t}^{t+1}  = \x_{i^t}^t  + \bar{\eta}^t$.

As can be seen from Table \ref{tab:binary:opt}, the proposed method \textit{\textbf{CD-SNCA}} consistently gives the best performance. This is due to the fact that \textit{\textbf{CD-SNCA}} finds stronger stationary points than the other methods. Such results consolidate our previous conclusions.


\begin{table}[!th]
\begin{center}
\scalebox{0.45}{\begin{tabular}{|c|c|c|c|c|c|c|c|c|}
\hline
            & MSCR  &PDCA & SubGrad  & CD-SCA & CD-SNCA \\
\hline
randn-256-1024 & \cthree{1.336 $\pm$ 0.108} & \cthree{1.336 $\pm$ 0.108} & \ctwo{1.280 $\pm$ 0.098} & 1.540 $\pm$ 0.236 & \cone{0.046 $\pm$ 0.010}  \\
randn-256-2048 & \cthree{1.359 $\pm$ 0.207} & \cthree{1.359 $\pm$ 0.207} & \ctwo{1.305 $\pm$ 0.199} & 1.503 $\pm$ 0.242 & \cone{0.021 $\pm$ 0.004}  \\
randn-1024-256 & \cthree{2.275 $\pm$ 0.096} & \cthree{2.275 $\pm$ 0.096} & \ctwo{2.268 $\pm$ 0.092} & 2.380 $\pm$ 0.180 & \cone{1.203 $\pm$ 0.043}  \\
randn-2048-256 & \cthree{3.569 $\pm$ 0.144} & \cthree{3.569 $\pm$ 0.144} & \ctwo{3.561 $\pm$ 0.143} & 3.614 $\pm$ 0.162 & \cone{2.492 $\pm$ 0.084}  \\
e2006-256-1024 & 1.069 $\pm$ 0.313 & 1.069 $\pm$ 0.313 & \ctwo{0.605 $\pm$ 0.167} & \cthree{0.809 $\pm$ 0.222} & \cone{0.291 $\pm$ 0.025}  \\
e2006-256-2048 & 0.936 $\pm$ 0.265 & 0.936 $\pm$ 0.265 & \ctwo{0.640 $\pm$ 0.187} & \cthree{0.798 $\pm$ 0.255} & \cone{0.263 $\pm$ 0.028}  \\
e2006-1024-256 & 2.245 $\pm$ 0.534 & 2.245 $\pm$ 0.534 & \ctwo{1.670 $\pm$ 0.198} & \cthree{1.780 $\pm$ 0.238} & \cone{1.266 $\pm$ 0.057}  \\
e2006-2048-256 & 3.507 $\pm$ 0.529 & 3.507 $\pm$ 0.529 & \ctwo{3.053 $\pm$ 0.250} & \cthree{3.307 $\pm$ 0.396} & \cone{2.532 $\pm$ 0.191}  \\
randn-256-1024-C & \cthree{1.357 $\pm$ 0.130} & \cthree{1.357 $\pm$ 0.130} & \ctwo{1.302 $\pm$ 0.134} & 1.586 $\pm$ 0.192 & \cone{0.051 $\pm$ 0.012}  \\
randn-256-2048-C & \cthree{1.260 $\pm$ 0.126} & 1.261 $\pm$ 0.126 & \ctwo{1.202 $\pm$ 0.122} & 1.444 $\pm$ 0.099 & \cone{0.019 $\pm$ 0.003}  \\
randn-1024-256-C & \cthree{2.254 $\pm$ 0.097} & \cthree{2.254 $\pm$ 0.097} & \ctwo{2.226 $\pm$ 0.084} & 2.315 $\pm$ 0.154 & \cone{1.175 $\pm$ 0.045}  \\
randn-2048-256-C & \cthree{3.531 $\pm$ 0.159} & \cthree{3.531 $\pm$ 0.159} & \ctwo{3.520 $\pm$ 0.150} & 3.544 $\pm$ 0.184 & \cone{2.445 $\pm$ 0.082}  \\
e2006-256-1024-C & 1.281 $\pm$ 0.628 & 1.323 $\pm$ 0.684 & \ctwo{0.473 $\pm$ 0.128} & \cthree{0.671 $\pm$ 0.257} & \cone{0.302 $\pm$ 0.043}  \\
e2006-256-2048-C & 1.254 $\pm$ 0.535 & 1.254 $\pm$ 0.535 & \ctwo{0.577 $\pm$ 0.144} & \cthree{0.717 $\pm$ 0.218} & \cone{0.287 $\pm$ 0.029}  \\
e2006-1024-256-C & 2.308 $\pm$ 0.640 & 2.308 $\pm$ 0.640 & \ctwo{1.570 $\pm$ 0.237} & \cthree{1.837 $\pm$ 0.322} & \cone{1.303 $\pm$ 0.060}  \\
e2006-2048-256-C & 3.429 $\pm$ 0.687 & 3.429 $\pm$ 0.687 & \ctwo{2.693 $\pm$ 0.335} & \cthree{2.790 $\pm$ 0.287} & \cone{2.431 $\pm$ 0.150}  \\
\hline
\end{tabular}}
\caption{Comparisons of objective values of all the methods for solving the approximate binary optimization problem. The $1^{st}$, $2^{nd}$, and $3^{rd}$ best results are colored with \cone{red}, \ctwo{green} and \cthree{blue}, respectively.}
\label{tab:binary:opt}
\end{center}
\end{table}

\begin{table}[!th]
\begin{center}
\scalebox{0.45}{\begin{tabular}{|c|c|c|c|c|c|c|c|c|}
\hline
              & MSCR &  PDCA & SubGrad & CD-SCA & CD-SNCA \\
\hline
randn-256-1024 & \ctwo{0.046 $\pm$ 0.019} & \ctwo{0.046 $\pm$ 0.019} & \cthree{0.077 $\pm$ 0.017} & \cone{0.039 $\pm$ 0.018} & \cone{0.039 $\pm$ 0.019}  \\
randn-256-2048 & \cthree{0.023 $\pm$ 0.008} & \ctwo{0.022 $\pm$ 0.007} & 0.060 $\pm$ 0.006 & \cone{0.021 $\pm$ 0.007} & \cone{0.021 $\pm$ 0.007}  \\
randn-1024-256 & 0.480 $\pm$ 0.063 & \cthree{0.473 $\pm$ 0.057} & 0.771 $\pm$ 0.089 & \ctwo{0.464 $\pm$ 0.059} & \cone{0.461 $\pm$ 0.060}  \\
randn-2048-256 & 1.335 $\pm$ 0.205 & \cthree{1.330 $\pm$ 0.205} & 1.810 $\pm$ 0.262 & \ctwo{1.329 $\pm$ 0.197} & \cone{1.325 $\pm$ 0.197}  \\
e2006-256-1024 & \ctwo{0.046 $\pm$ 0.093} & \cthree{0.047 $\pm$ 0.105} & 0.050 $\pm$ 0.088 & \cthree{0.047 $\pm$ 0.100} & \cone{0.045 $\pm$ 0.097}  \\
e2006-256-2048 & \ctwo{0.022 $\pm$ 0.009} & \cthree{0.025 $\pm$ 0.012} & 0.036 $\pm$ 0.040 & 0.029 $\pm$ 0.039 & \cone{0.020 $\pm$ 0.020}  \\
e2006-1024-256 & \ctwo{0.922 $\pm$ 0.754} & \cthree{0.925 $\pm$ 0.758} & 0.941 $\pm$ 0.792 & \cthree{0.925 $\pm$ 0.757} & \cone{0.858 $\pm$ 0.717}  \\
e2006-2048-256 & \cthree{1.031 $\pm$ 0.835} & 1.035 $\pm$ 0.838 & 1.075 $\pm$ 0.867 & \ctwo{1.024 $\pm$ 0.827} & \cone{1.010 $\pm$ 0.817}  \\
randn-256-1024-C & \cthree{0.036 $\pm$ 0.012} & \cthree{0.036 $\pm$ 0.012} & 0.069 $\pm$ 0.014 & \ctwo{0.031 $\pm$ 0.012} & \cone{0.030 $\pm$ 0.010}  \\
randn-256-2048-C & \cthree{0.019 $\pm$ 0.003} & \ctwo{0.018 $\pm$ 0.003} & 0.058 $\pm$ 0.004 & \cone{0.016 $\pm$ 0.003} & \cone{0.016 $\pm$ 0.003}  \\
randn-1024-256-C & \ctwo{0.462 $\pm$ 0.089} & 0.465 $\pm$ 0.092 & 0.768 $\pm$ 0.127 & \cthree{0.463 $\pm$ 0.088} & \cone{0.457 $\pm$ 0.092}  \\
randn-2048-256-C & \ctwo{1.155 $\pm$ 0.159} & \cthree{1.157 $\pm$ 0.165} & 1.570 $\pm$ 0.238 & 1.161 $\pm$ 0.168 & \cone{1.147 $\pm$ 0.160}  \\
e2006-256-1024-C & \ctwo{0.023 $\pm$ 0.020} & \cthree{0.025 $\pm$ 0.023} & 0.032 $\pm$ 0.026 & 0.031 $\pm$ 0.038 & \cone{0.019 $\pm$ 0.018}  \\
e2006-256-2048-C & \ctwo{0.034 $\pm$ 0.029} & 0.037 $\pm$ 0.034 & \cthree{0.036 $\pm$ 0.066} & \ctwo{0.034 $\pm$ 0.052} & \cone{0.025 $\pm$ 0.043}  \\
e2006-1024-256-C & \cthree{1.772 $\pm$ 2.200} & 1.788 $\pm$ 2.200 & 1.797 $\pm$ 2.294 & \ctwo{1.768 $\pm$ 2.195} & \cone{1.702 $\pm$ 2.162}  \\
e2006-2048-256-C & \cthree{1.474 $\pm$ 1.247} & 1.486 $\pm$ 1.249 & 1.520 $\pm$ 1.278 & \ctwo{1.446 $\pm$ 1.233} & \cone{1.431 $\pm$ 1.224}  \\
\hline
\end{tabular}}
\caption{Comparisons of objective values of all the methods for solving the generalized linear regression problem. The $1^{st}$, $2^{nd}$, and $3^{rd}$ best results are colored with \cone{red}, \ctwo{green} and \cthree{blue}, respectively.}
\label{tab:regression}

\end{center}
\end{table}
\subsection{Generalized Linear Regression}

We consider Problem (\ref{eq:oneNN}). We have the following optimization problem: $\min_{\x}\tfrac{1}{2}\|\max(0,\G\x)-\y\|_2^2 $. We generate the observation vector via $\y = \max(0,\A\ddot{\x} + \text{randn}(m,1) \times 0.1\times\|\A\ddot{\x}\|)$ with $\ddot{\x}=\text{randn}(d,1)$. This problem is consistent with Problem (\ref{eq:main}) with $f(\x)  \triangleq \frac{1}{2}\|\max(0,\G\x)\|_2^2$ and $g(\x) \triangleq \|\max(0,\A\x)\|_1$ with $\A=\text{diag}(\y)\G$. We notice that $\nabla f(\cdot)$ is $L$-Lipschitz with $L=\|\G\|_2^2$ and coordinate-wise Lipschitz with $\c=\text{diag}(\G^T\G)$. The subgradient of $g(\x)$ at $\x^t$ can be computed as: $\partial g(\x^t) = \A^T \max(0,\A\x^t) \triangleq \g^t$.

We compare with the following methods. \textit{\textbf{(i)}} Multi-Stage Convex Relaxation (MSCR). It solves the following problem: $\x^{t+1} = \arg \min_{\x} f(\x) -\la \x-\x^t,~\g^t\ra$. \textit{\textbf{(ii)}} Proximal DC algorithm (PDCA). It considers the following iteration: $\x^{t+1} = \arg \min_{\x} \la \x - \x^t,~ \nabla f(\x)\ra + \frac{L}{2}\|\x-\x^t\|_2^2  -\la \x-\x^t,~\g^t\ra$. \textit{\textbf{(iii)}} Subgradient method (SubGrad). It uses the following iteration: $\x^{t+1} = \x^t  - \frac{0.1}{t} \cdot (\nabla f(\x) -  \g^t )$. \textit{\textbf{(iv)}} \textit{\textbf{CD-SCA}} solves a convex problem $\bar{\eta}^t = \arg \min_{\eta} 0.5 (\c_{i^t}+\theta)\eta^2   + [\nabla f(\x^t) -\rho \g^t]_{i^t} \cdot \eta$ with  and update $\x^t$ via $\x_{i^t}^{t+1}  = \x_{i^t}^t  + \bar{\eta}^t$. \textit{\textbf{(v)}} \textit{\textbf{CD-SNCA}} computes the nonconvex proximal operator (see Section \ref{eq:proximal:lmax0}) as $\bar{\eta}^t = \arg \min_{\eta} \frac{\c_i+\theta}{2}\eta^2 +  \nabla_{i^t} f(\x^t) \eta  - \|\max(0,\A (\x^t + \eta\ei))\|_1$ and update $\x^t$ via $\x_{i^t}^{t+1}  = \x_{i^t}^t  + \bar{\eta}^t$.

As can be seen from Table \ref{tab:regression}, the proposed method \textit{\textbf{CD-SNCA}} consistently outperforms the other methods.

\subsection{More Experiments on Computational Efficiency}

\begin{figure*} [!t]
\centering

      \begin{subfigure}{0.25\textwidth}\includegraphics[height=\objimghei,width=\textwidth]{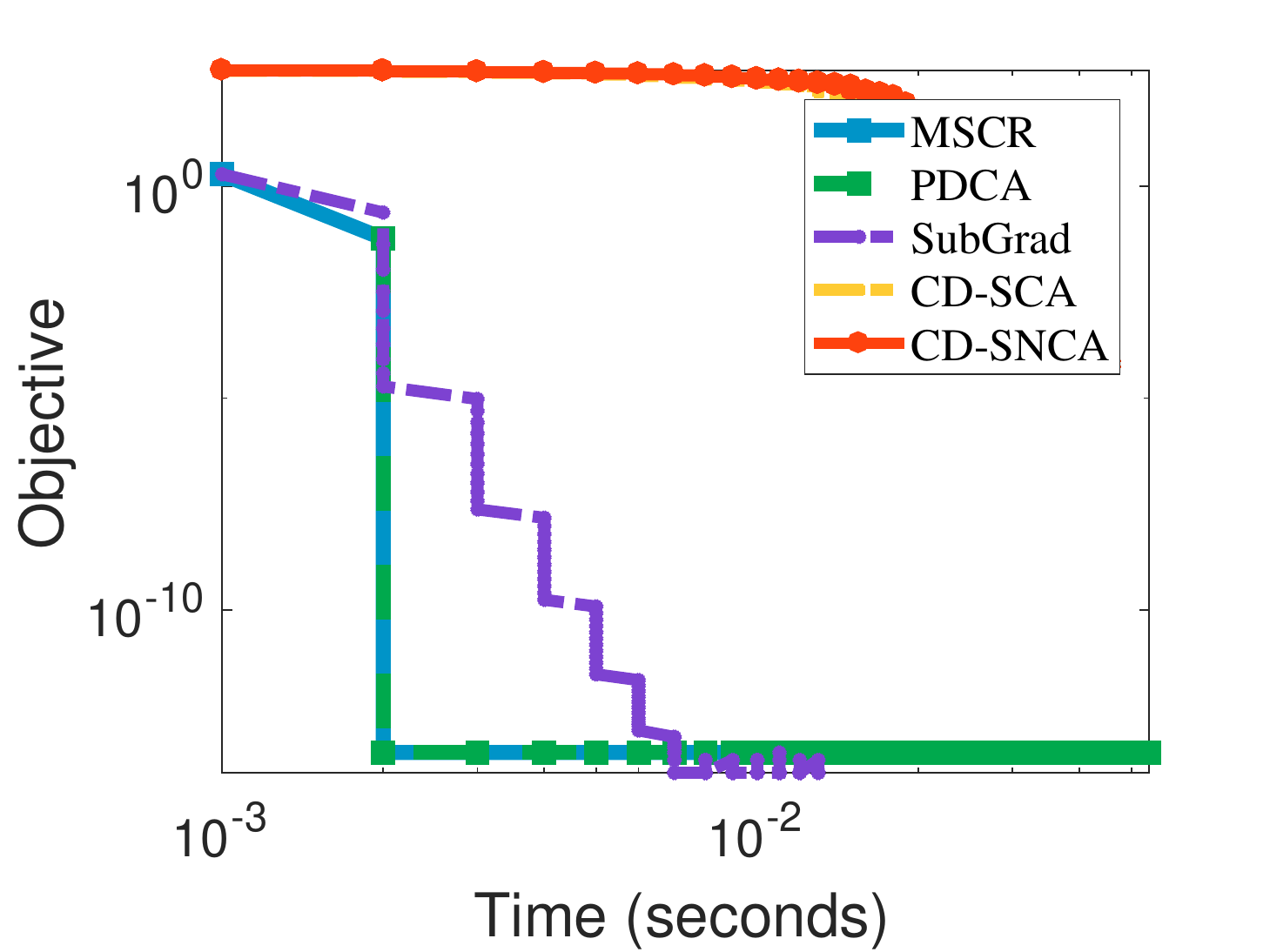}\vspace{-6pt} \caption{\scriptsize e2006-256-1024  }\end{subfigure}~~\begin{subfigure}{0.25\textwidth}\includegraphics[height=\objimghei,width=\textwidth]{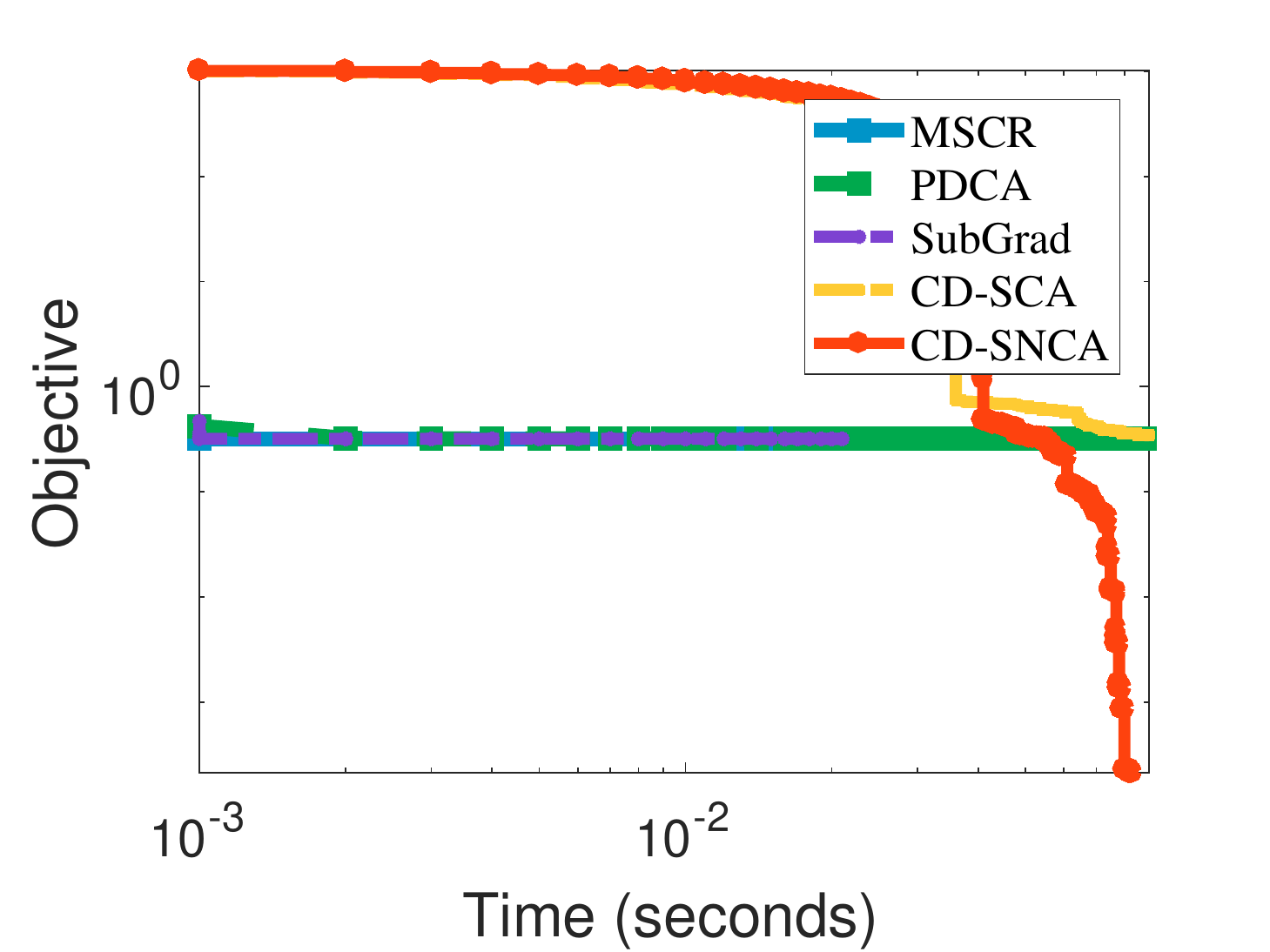}\vspace{-6pt} \caption{ \scriptsize e2006-256-2048 } \end{subfigure}~~\begin{subfigure}{0.25\textwidth}\includegraphics[height=\objimghei,width=\textwidth]{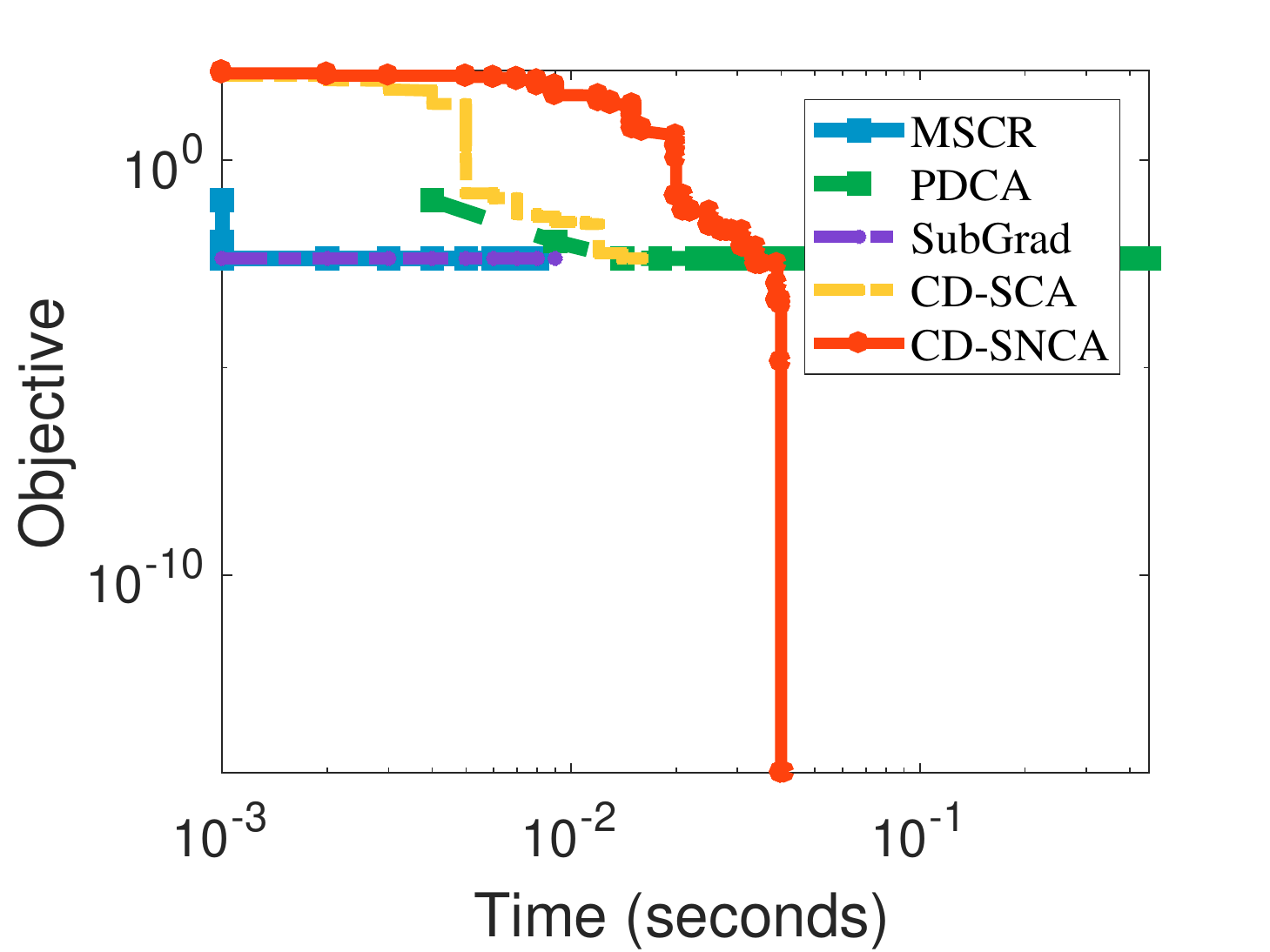}\vspace{-6pt} \caption{\scriptsize e2006-1024-256}\end{subfigure}~~\begin{subfigure}{0.25\textwidth}\includegraphics[height=\objimghei,width=\textwidth]{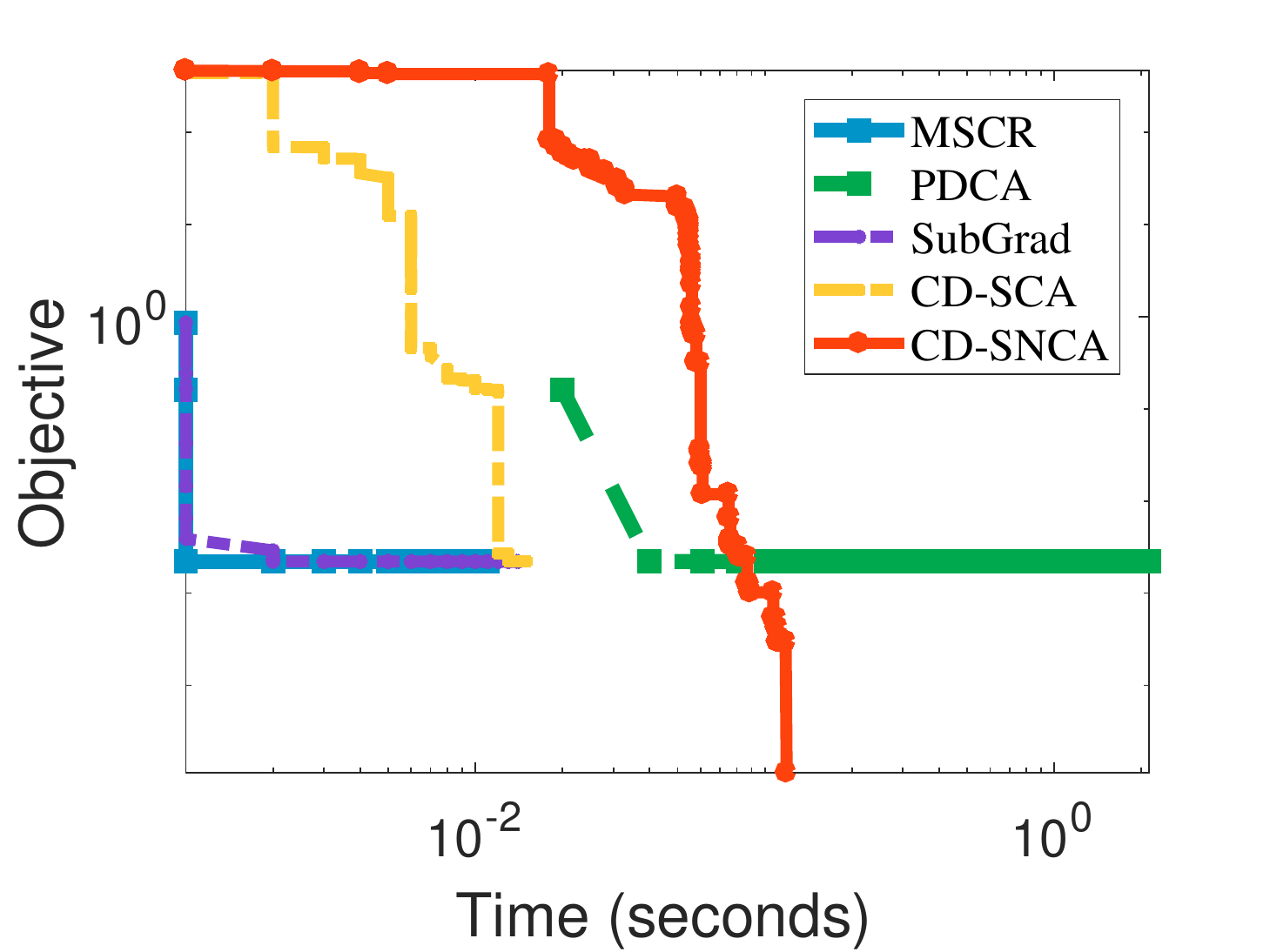}\vspace{-6pt} \caption{\scriptsize e2006-2048-256  }\end{subfigure}\\

      \begin{subfigure}{0.25\textwidth}\includegraphics[height=\objimghei,width=\textwidth]{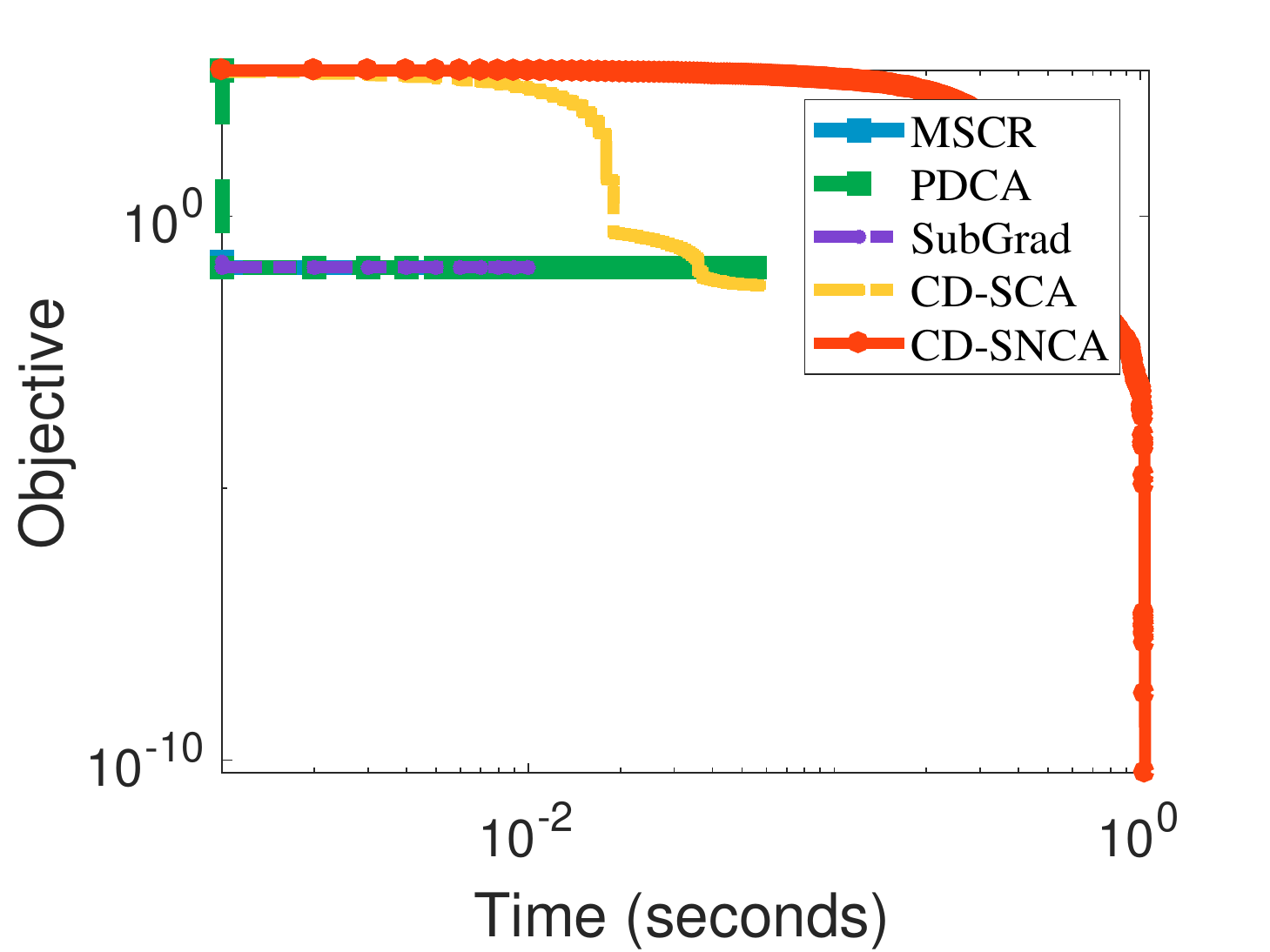}\vspace{-6pt} \caption{\scriptsize randn-256-1024-C  }\end{subfigure}~~\begin{subfigure}{0.25\textwidth}\includegraphics[height=\objimghei,width=\textwidth]{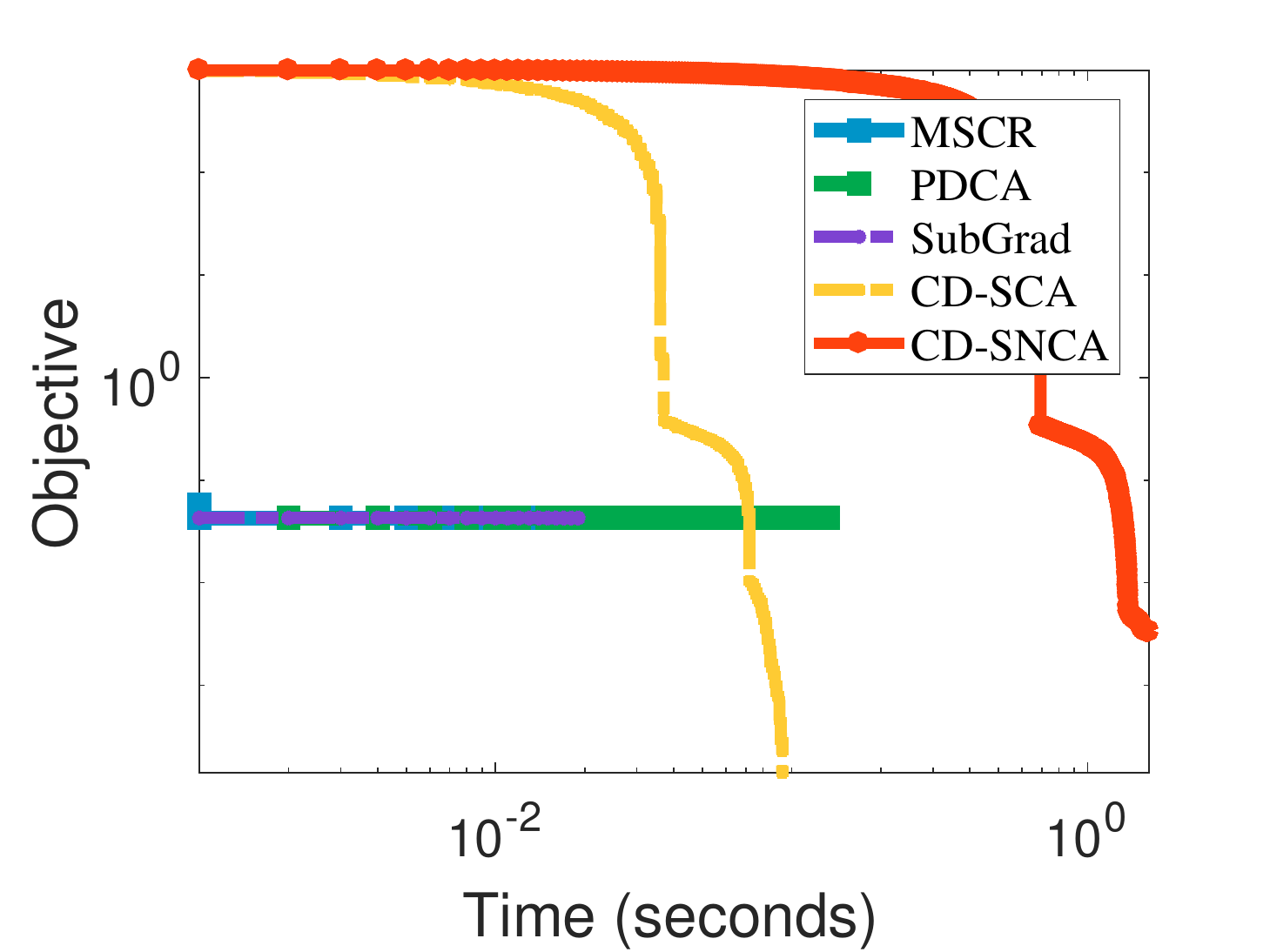}\vspace{-6pt} \caption{ \scriptsize randn-256-2048-C } \end{subfigure}~~\begin{subfigure}{0.25\textwidth}\includegraphics[height=\objimghei,width=\textwidth]{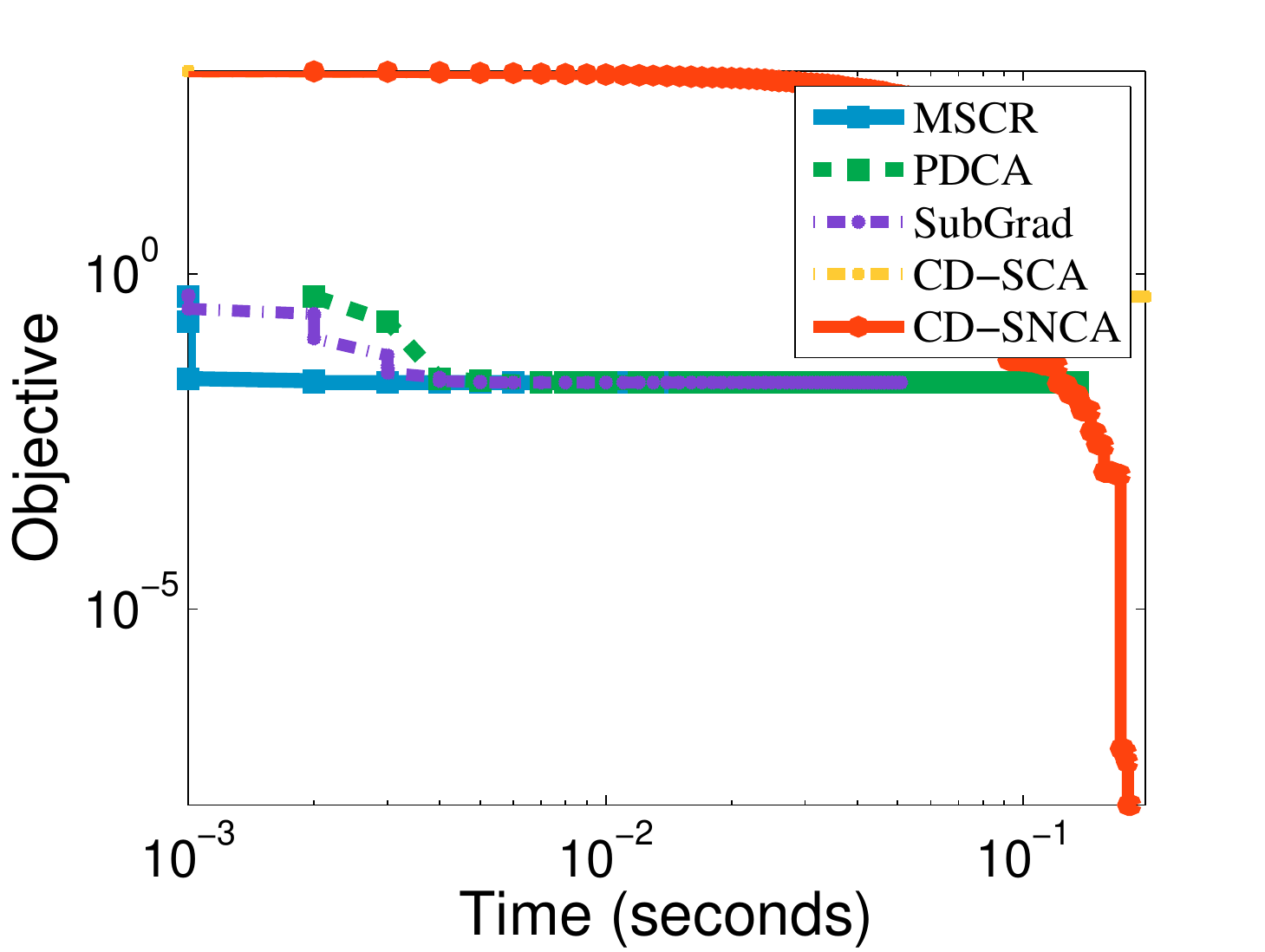}\vspace{-6pt} \caption{\scriptsize randn-1024-256-C}\end{subfigure}~~\begin{subfigure}{0.25\textwidth}\includegraphics[height=\objimghei,width=\textwidth]{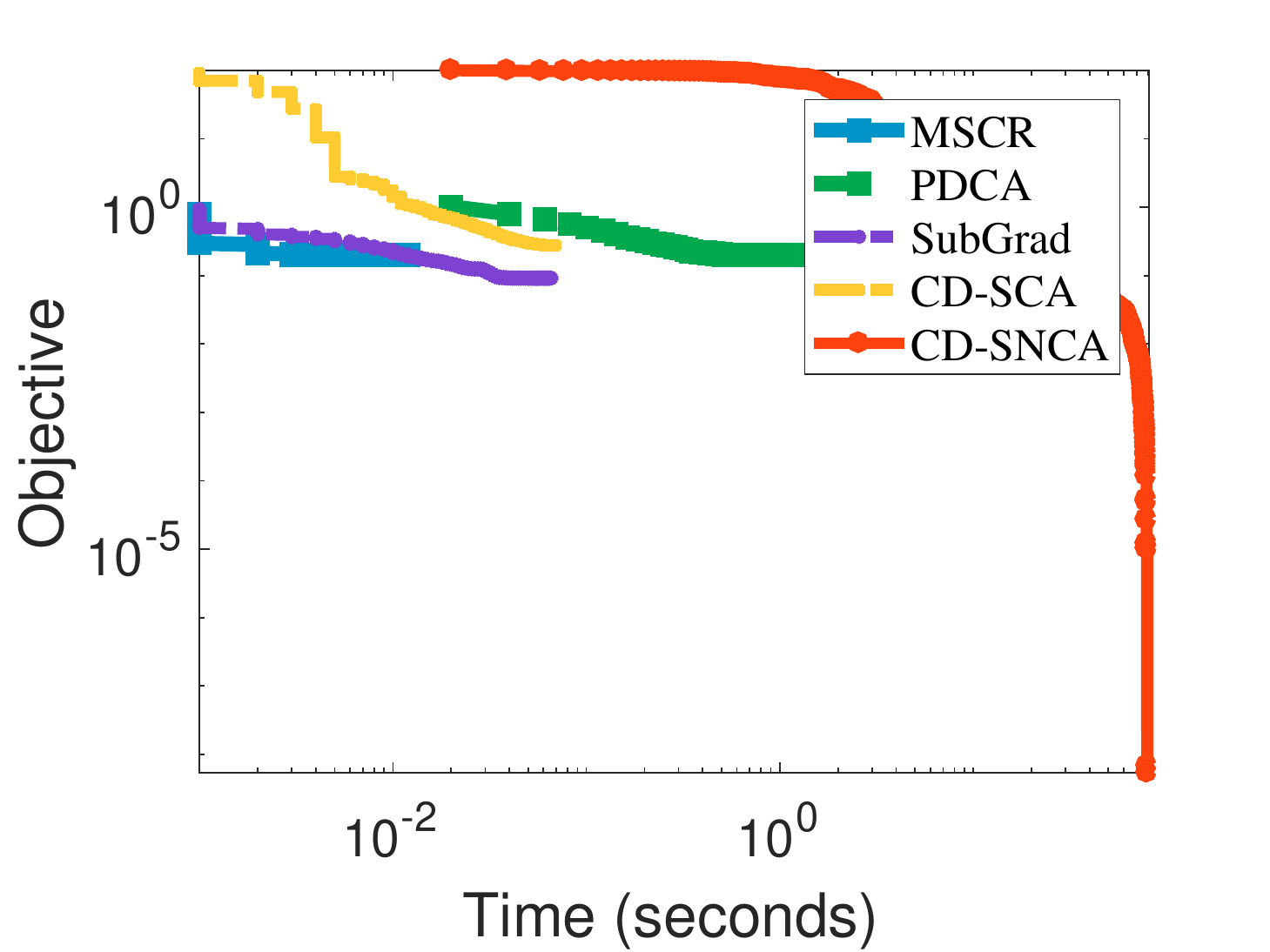}\vspace{-6pt} \caption{\scriptsize randn-2048-256-C  }\end{subfigure}\\

      \begin{subfigure}{0.25\textwidth}\includegraphics[height=\objimghei,width=\textwidth]{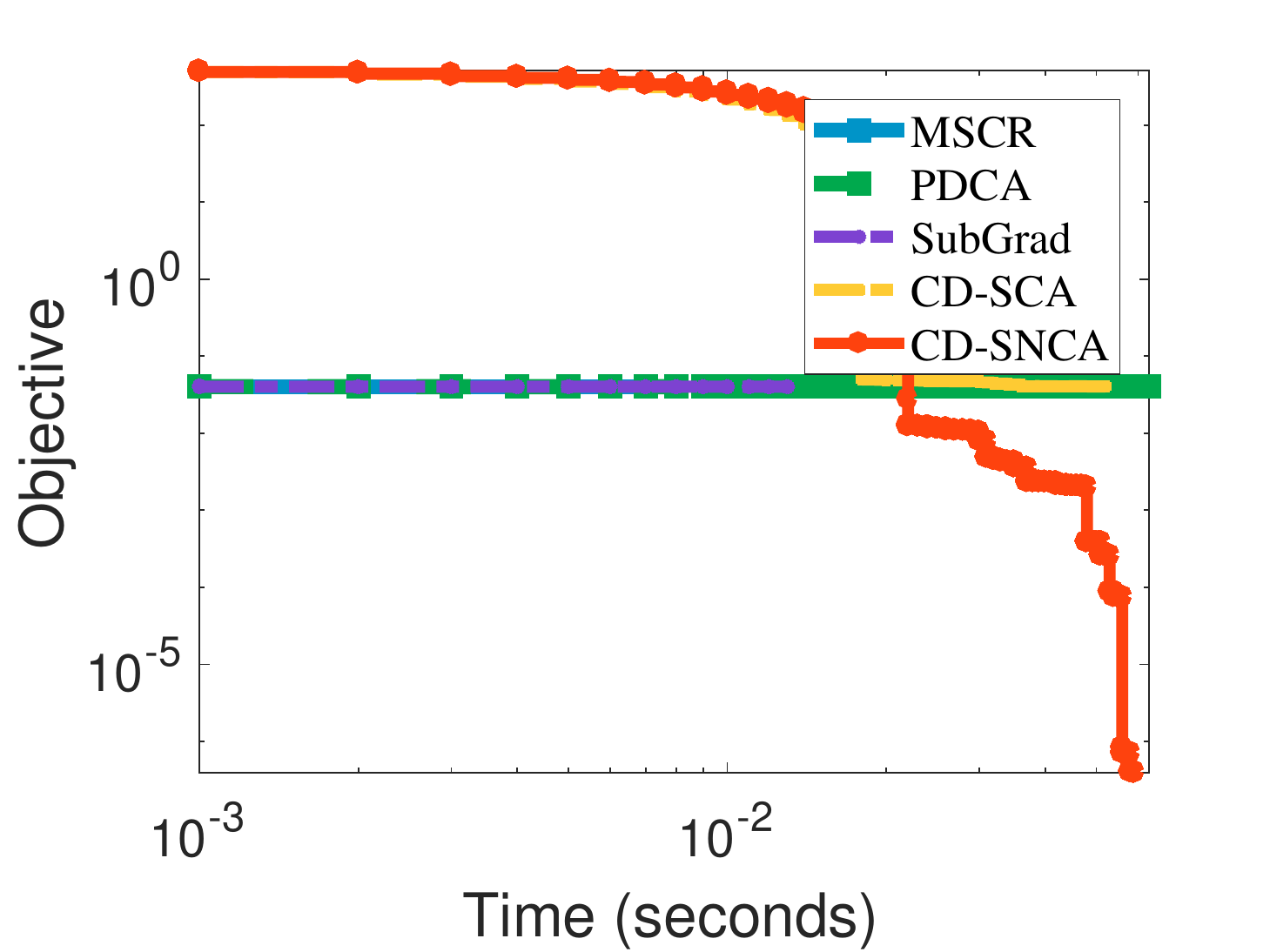}\vspace{-6pt} \caption{\scriptsize e2006-256-1024-C  }\end{subfigure}~~\begin{subfigure}{0.25\textwidth}\includegraphics[height=\objimghei,width=\textwidth]{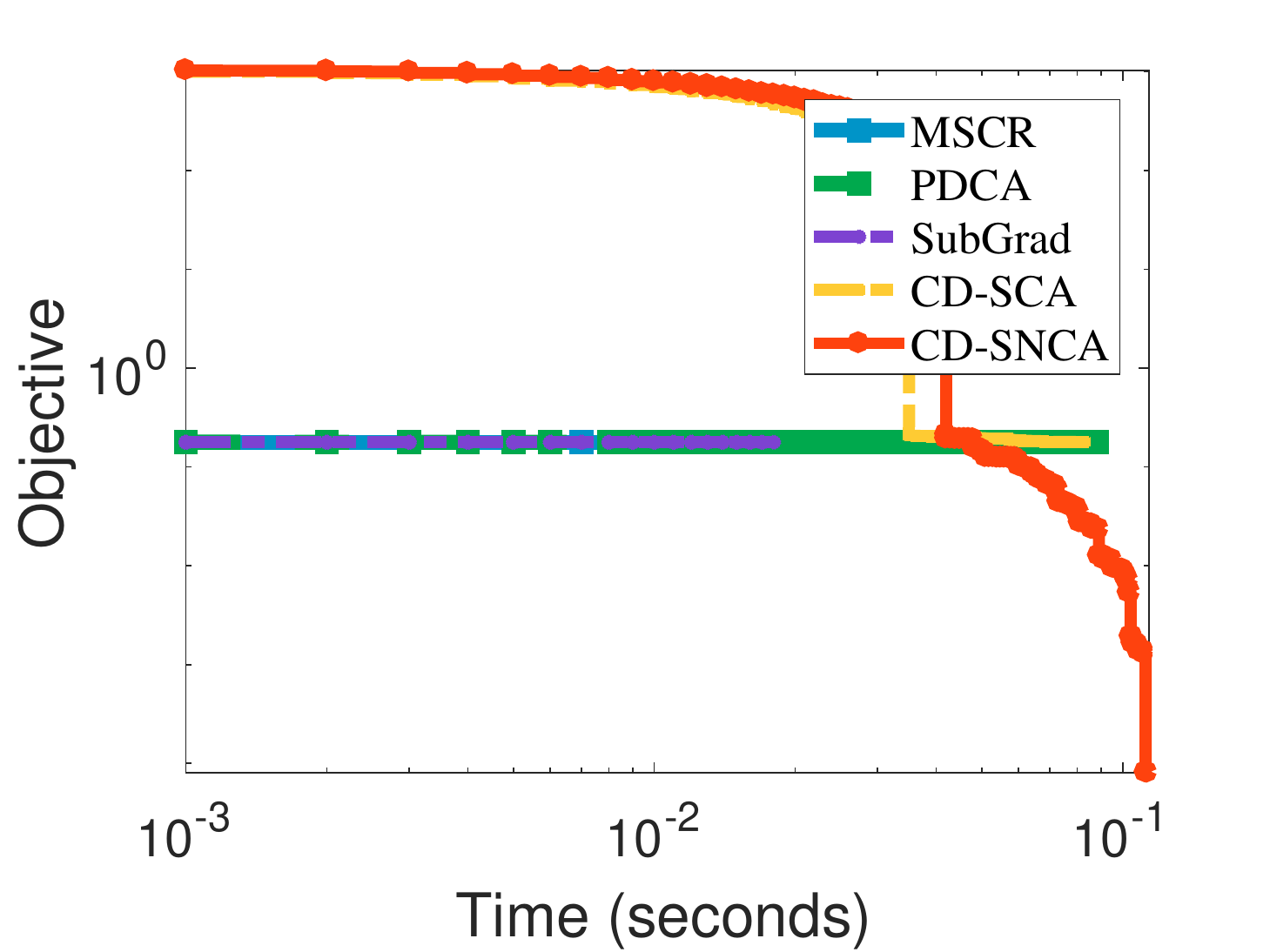}\vspace{-6pt} \caption{ \scriptsize e2006-256-2048-C } \end{subfigure}~~\begin{subfigure}{0.25\textwidth}\includegraphics[height=\objimghei,width=\textwidth]{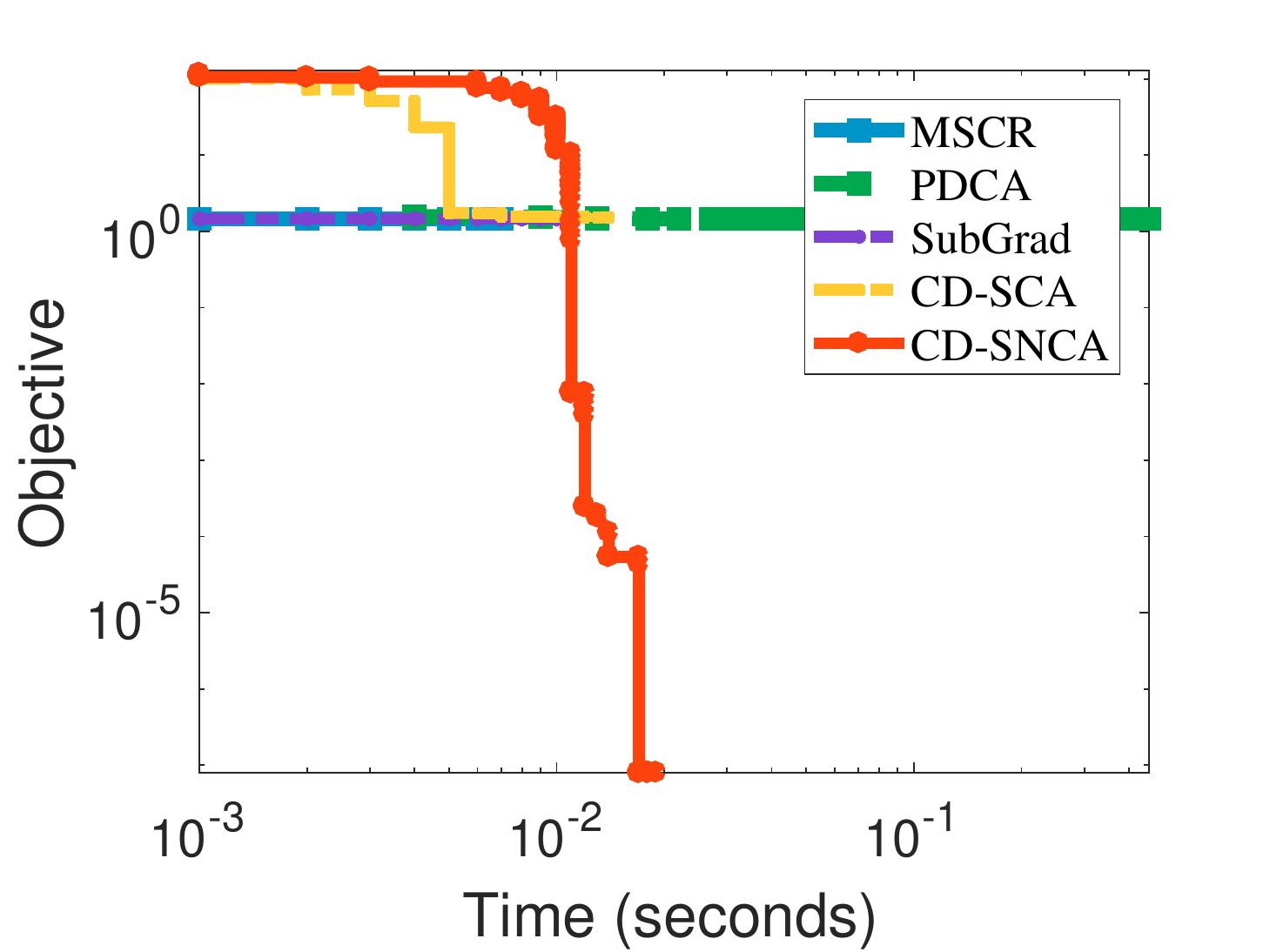}\vspace{-6pt} \caption{\scriptsize e2006-1024-256-C}\end{subfigure}~~\begin{subfigure}{0.25\textwidth}\includegraphics[height=\objimghei,width=\textwidth]{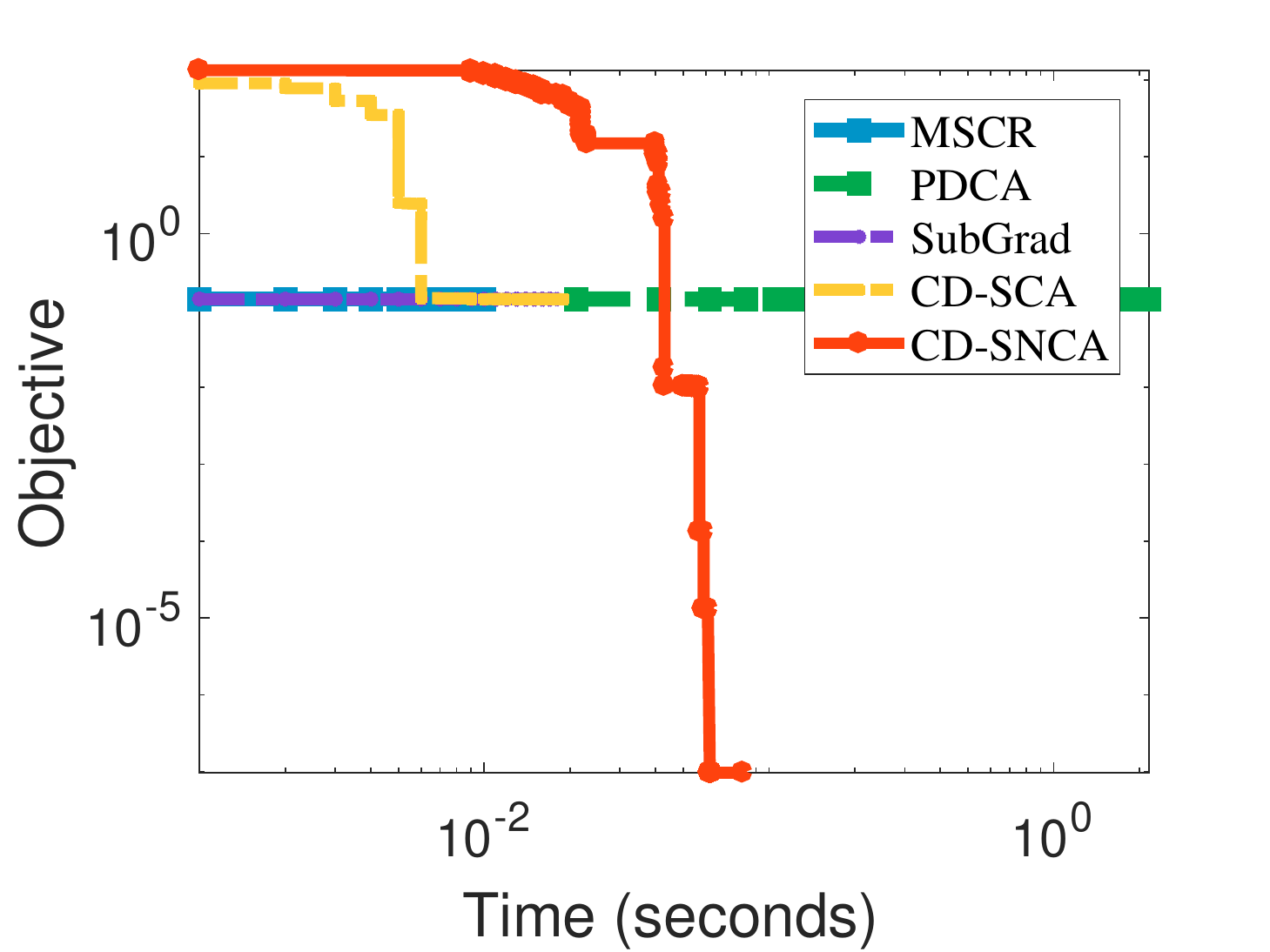}\vspace{-6pt} \caption{\scriptsize e2006-2048-256-C  }\end{subfigure}\\

\centering
\caption{The convergence curve of the compared methods for solving the $\ell_p$ norm generalized eigenvalue problem on different data sets.}

\label{exp:cpu:12}
\end{figure*}

\begin{figure*} [!t]
\centering
      \begin{subfigure}{0.25\textwidth}\includegraphics[height=\objimghei,width=\textwidth]{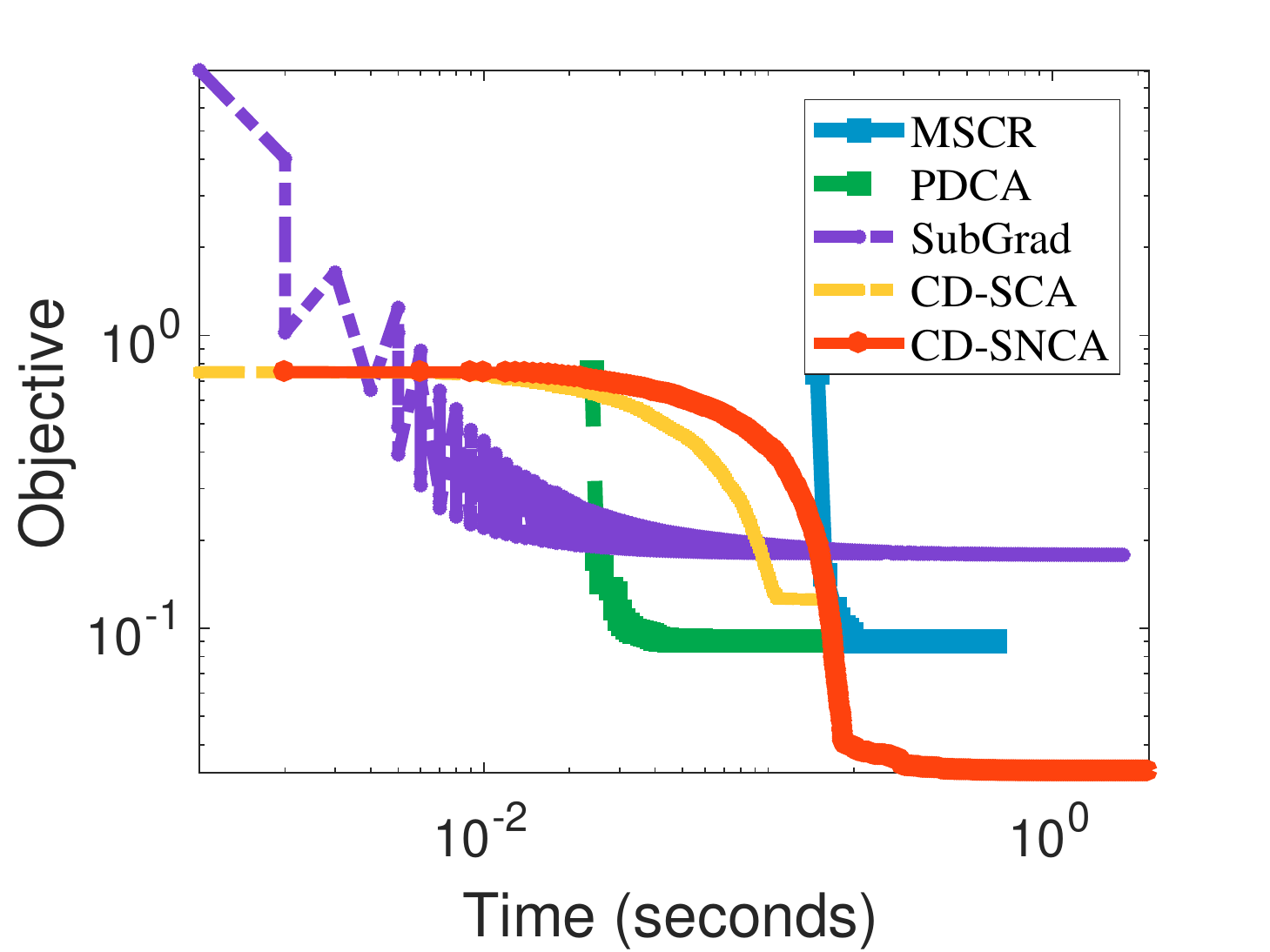}\vspace{-6pt} \caption{\scriptsize randn-256-1024 }\end{subfigure}~~\begin{subfigure}{0.25\textwidth}\includegraphics[height=\objimghei,width=\textwidth]{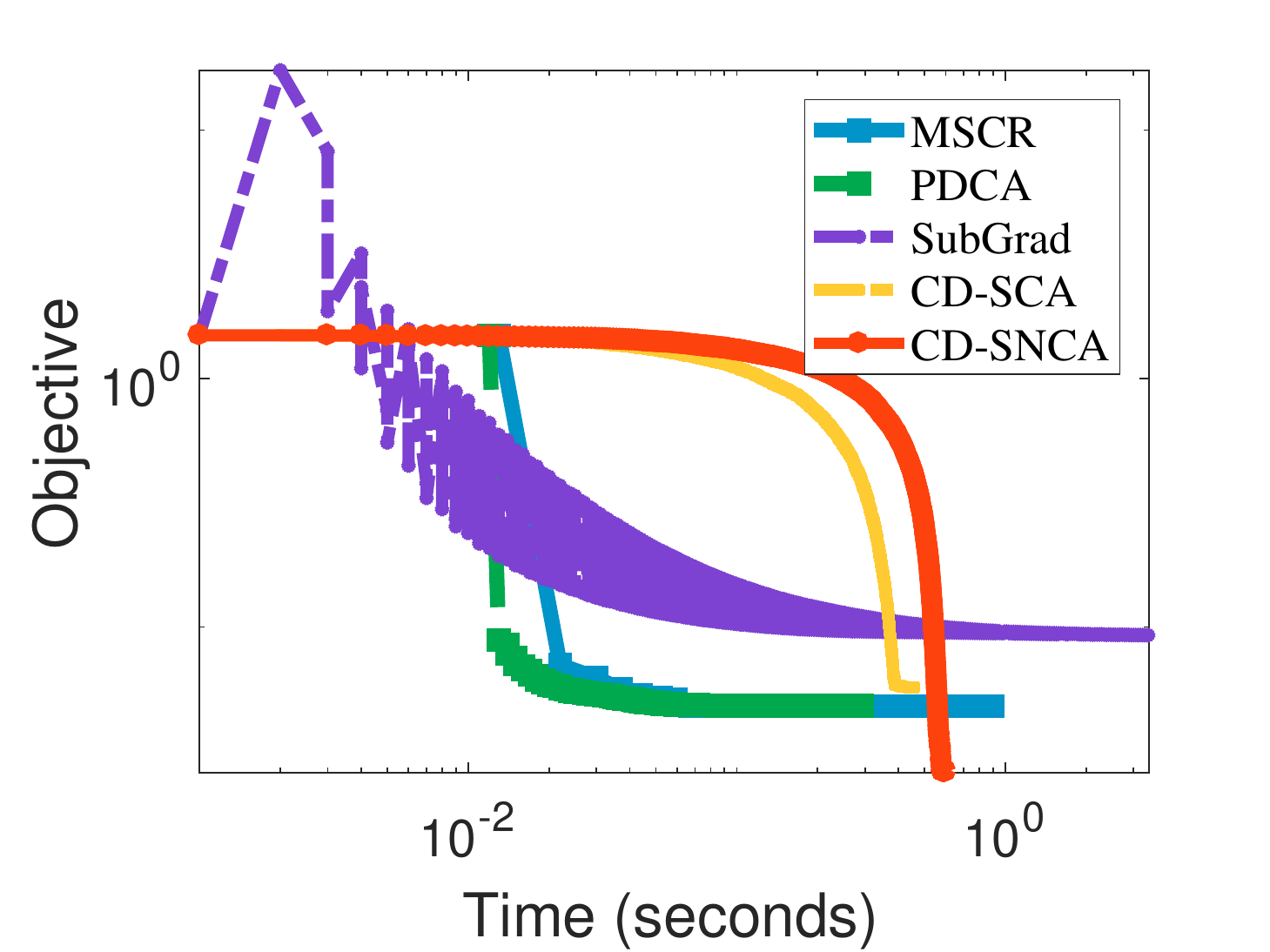}\vspace{-6pt} \caption{ \scriptsize randn-256-2048 } \end{subfigure}~~\begin{subfigure}{0.25\textwidth}\includegraphics[height=\objimghei,width=\textwidth]{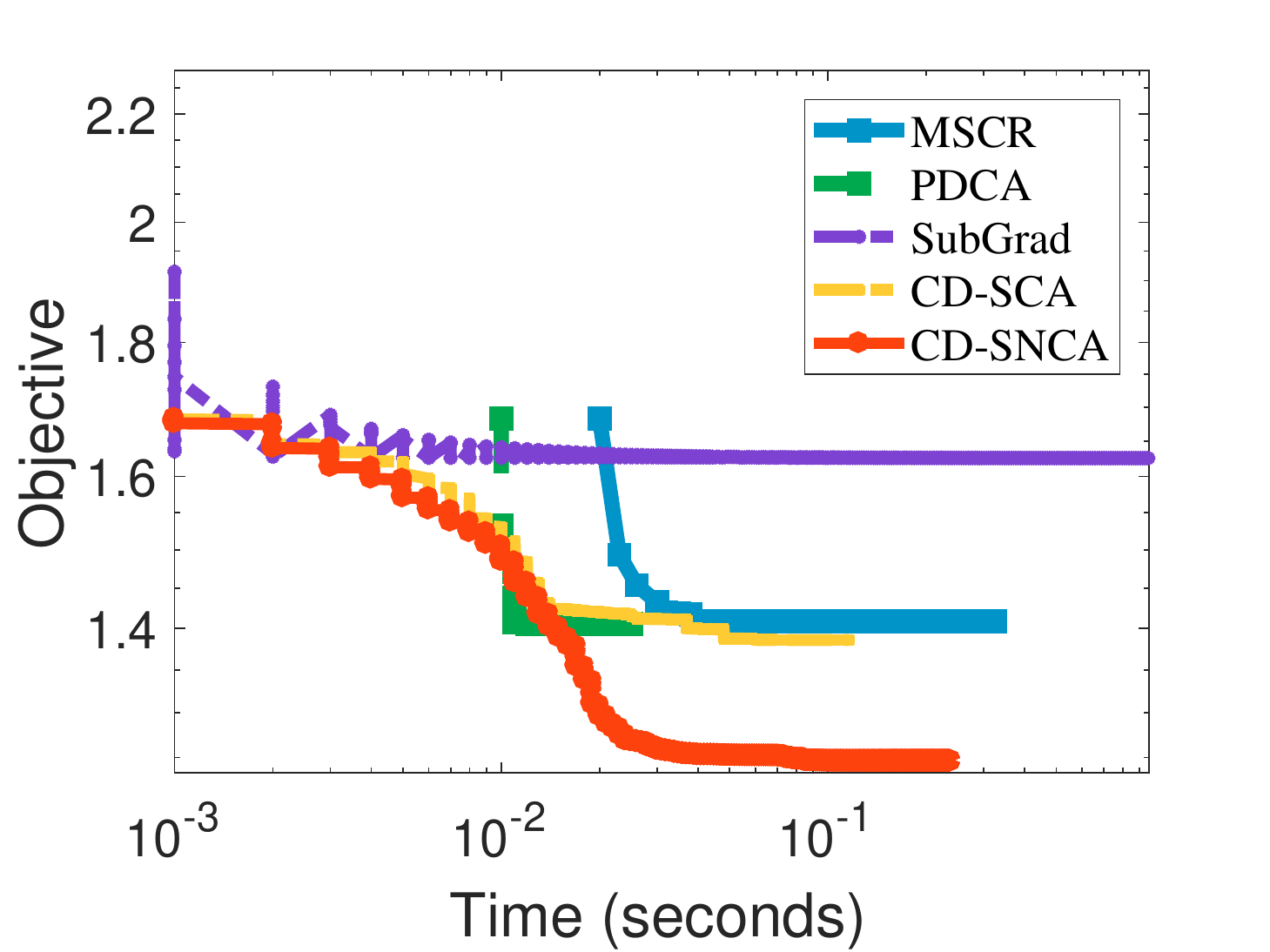}\vspace{-6pt} \caption{\scriptsize randn-1024-256}\end{subfigure}~~\begin{subfigure}{0.25\textwidth}\includegraphics[height=\objimghei,width=\textwidth]{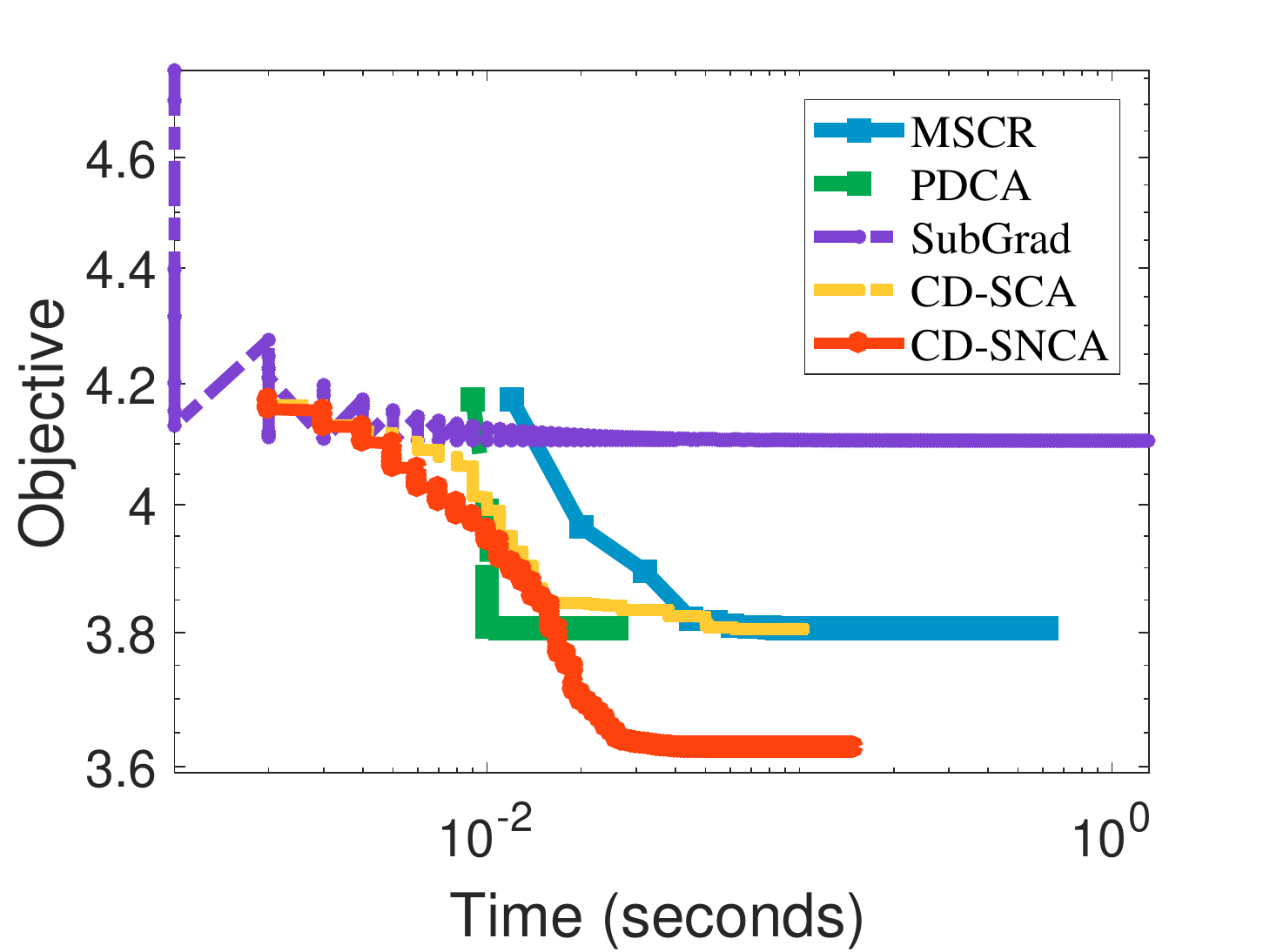}\vspace{-6pt} \caption{\scriptsize randn-2048-256  }\end{subfigure}\\

      \begin{subfigure}{0.25\textwidth}\includegraphics[height=\objimghei,width=\textwidth]{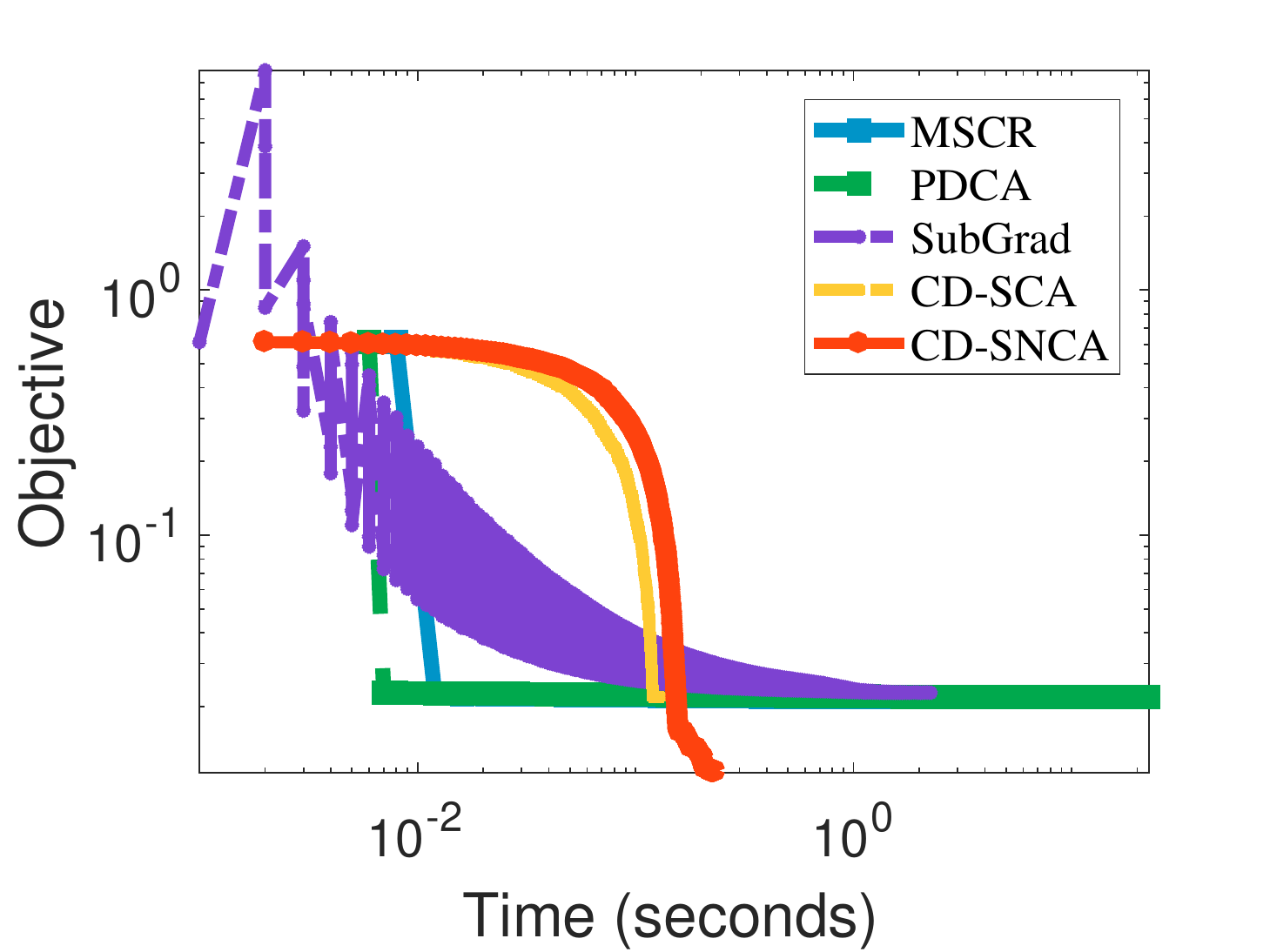}\vspace{-6pt} \caption{\scriptsize e2006-256-1024 }\end{subfigure}~~\begin{subfigure}{0.25\textwidth}\includegraphics[height=\objimghei,width=\textwidth]{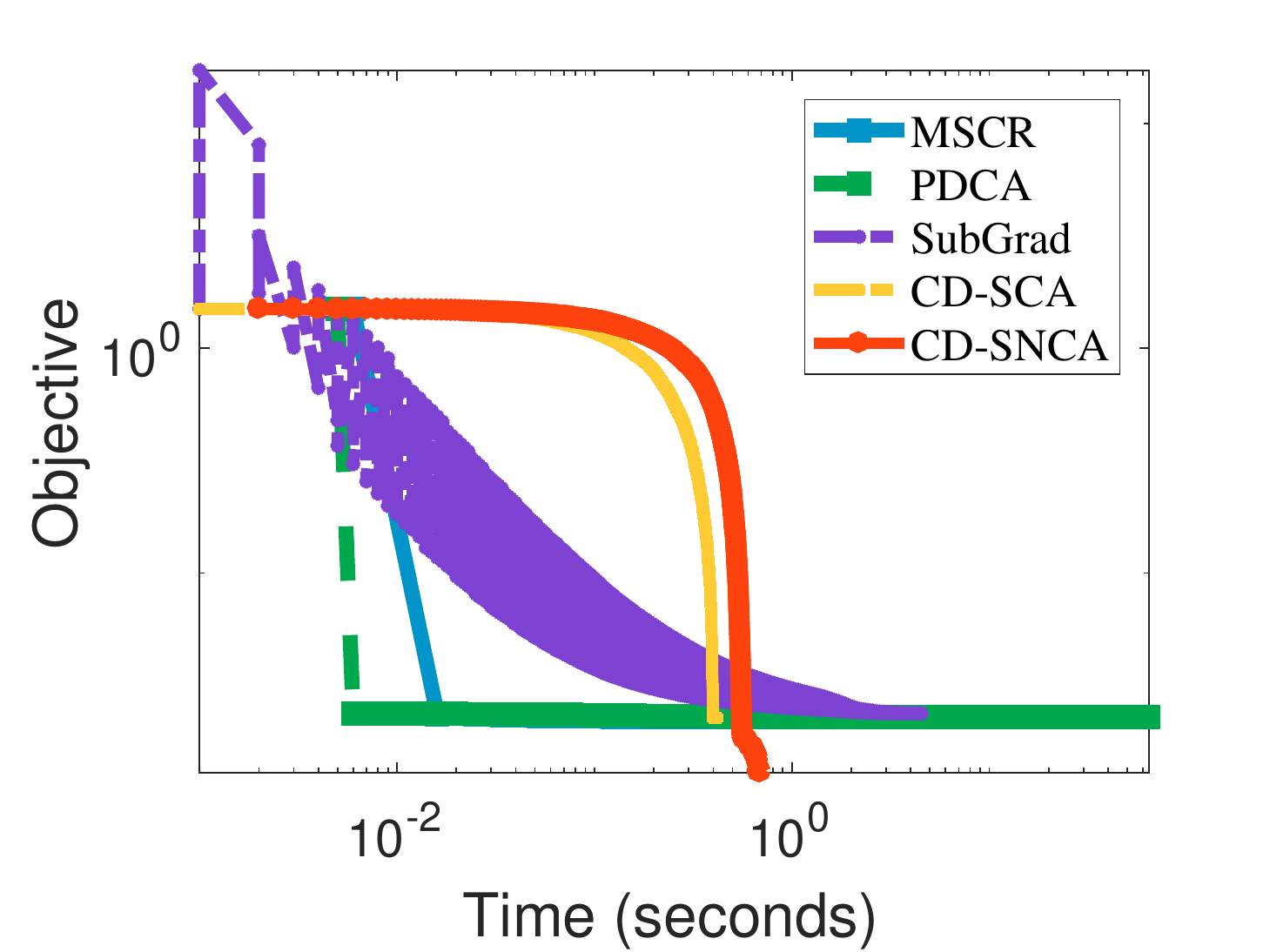}\vspace{-6pt} \caption{ \scriptsize e2006-256-2048 } \end{subfigure}~~\begin{subfigure}{0.25\textwidth}\includegraphics[height=\objimghei,width=\textwidth]{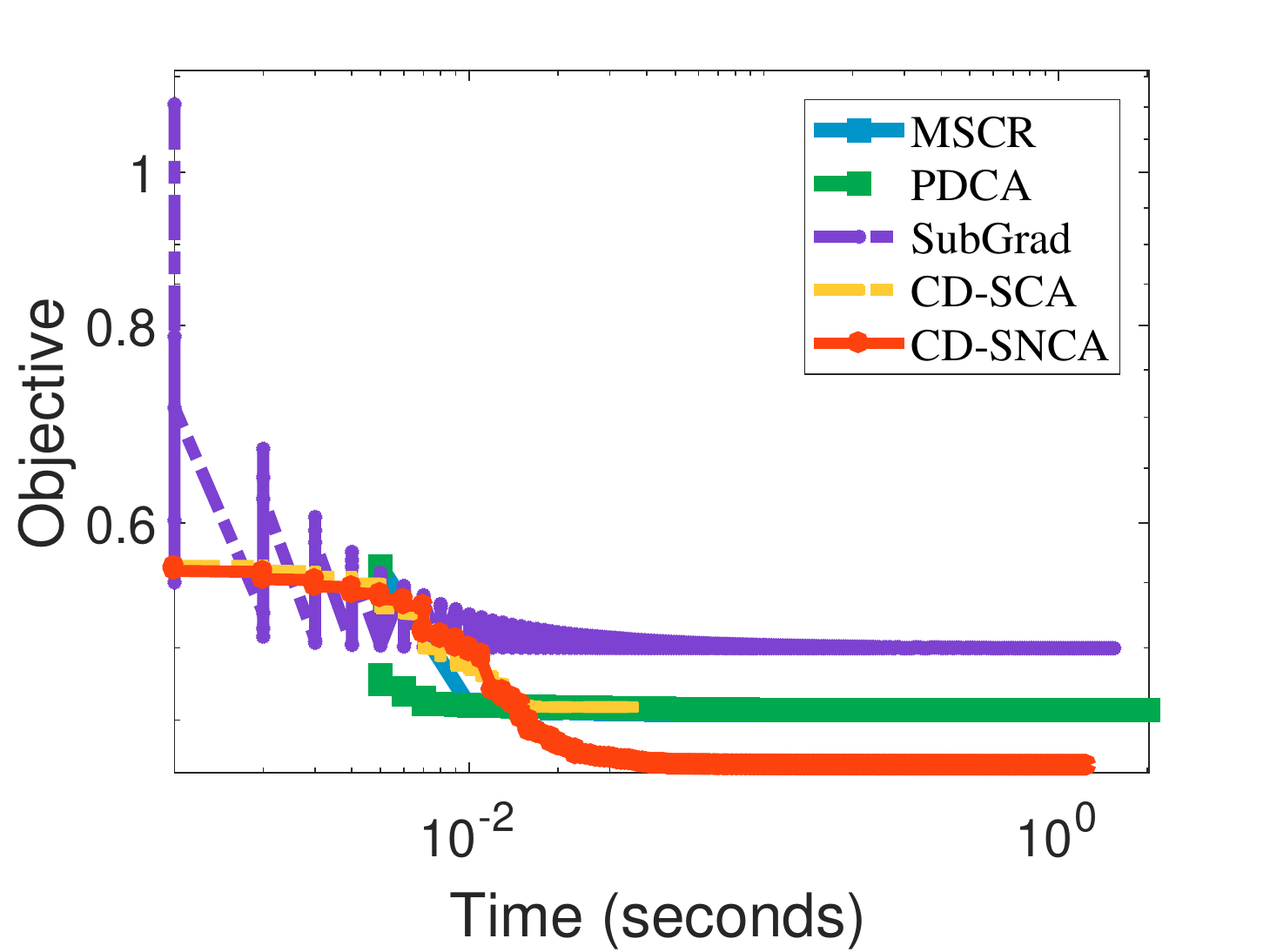}\vspace{-6pt} \caption{\scriptsize e2006-1024-256}\end{subfigure}~~\begin{subfigure}{0.25\textwidth}\includegraphics[height=\objimghei,width=\textwidth]{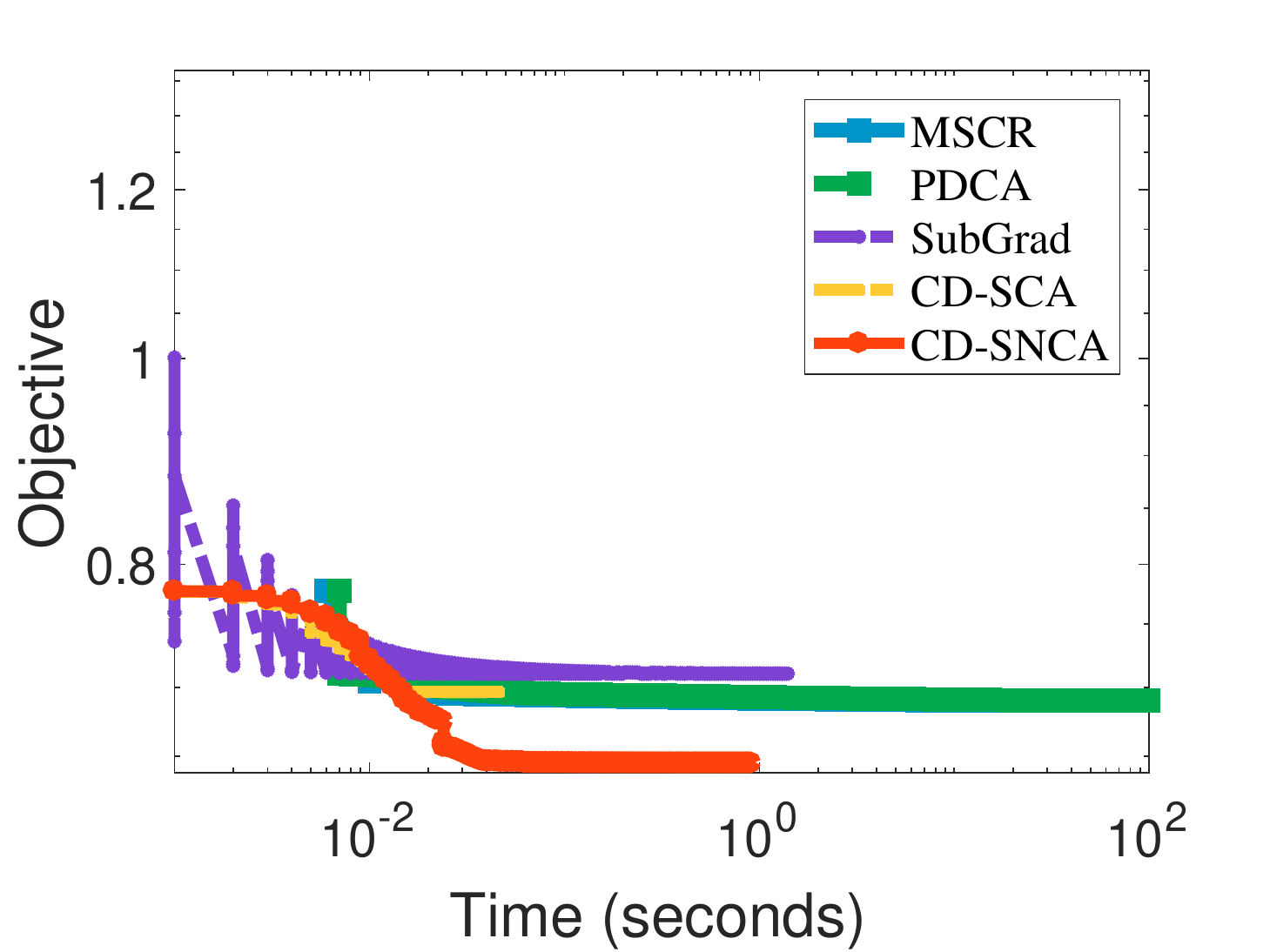}\vspace{-6pt} \caption{\scriptsize e2006-2048-256  }\end{subfigure}\\

      \begin{subfigure}{0.25\textwidth}\includegraphics[height=\objimghei,width=\textwidth]{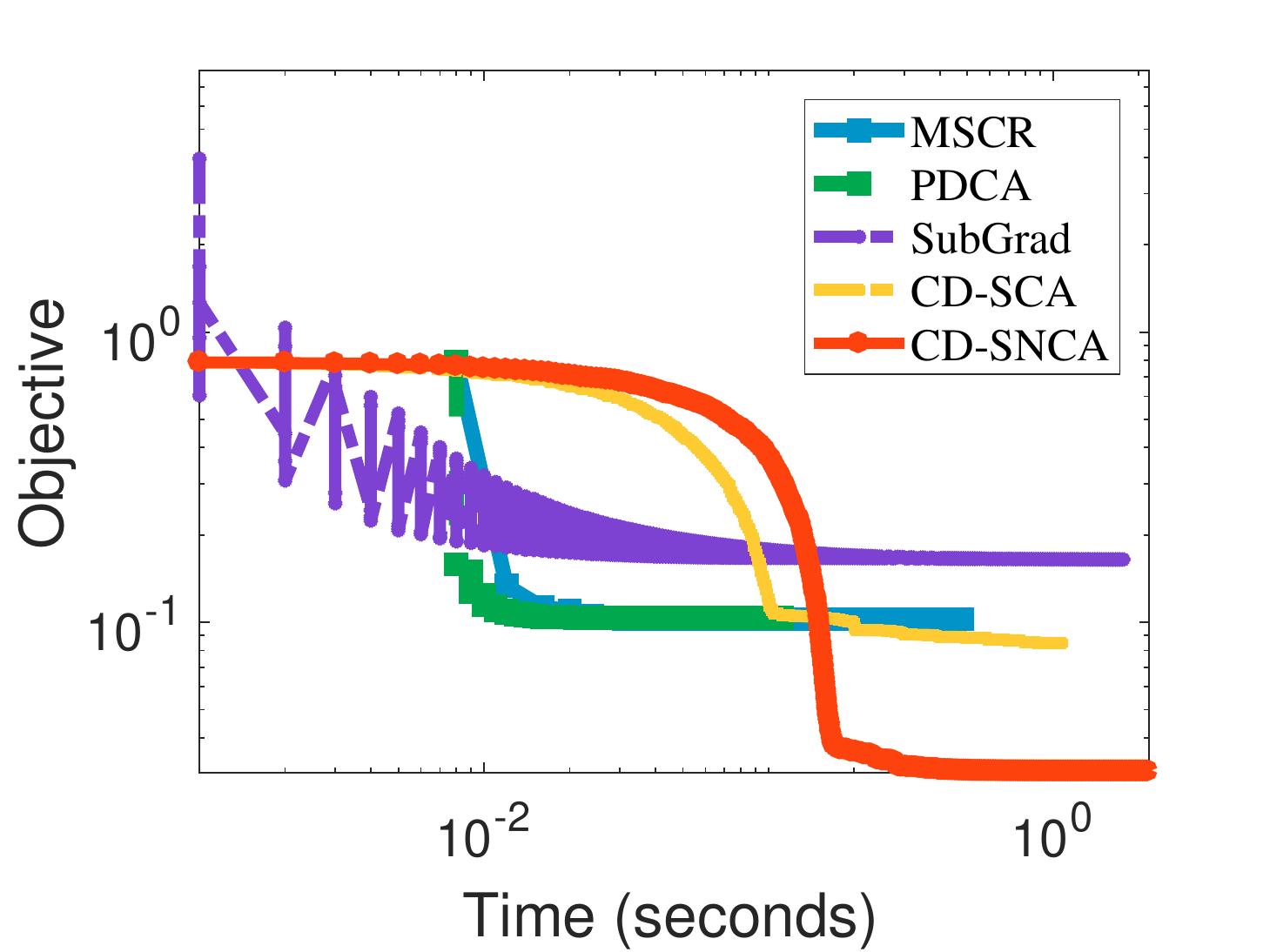}\vspace{-6pt} \caption{\scriptsize randn-256-1024-C  }\end{subfigure}~~\begin{subfigure}{0.25\textwidth}\includegraphics[height=\objimghei,width=\textwidth]{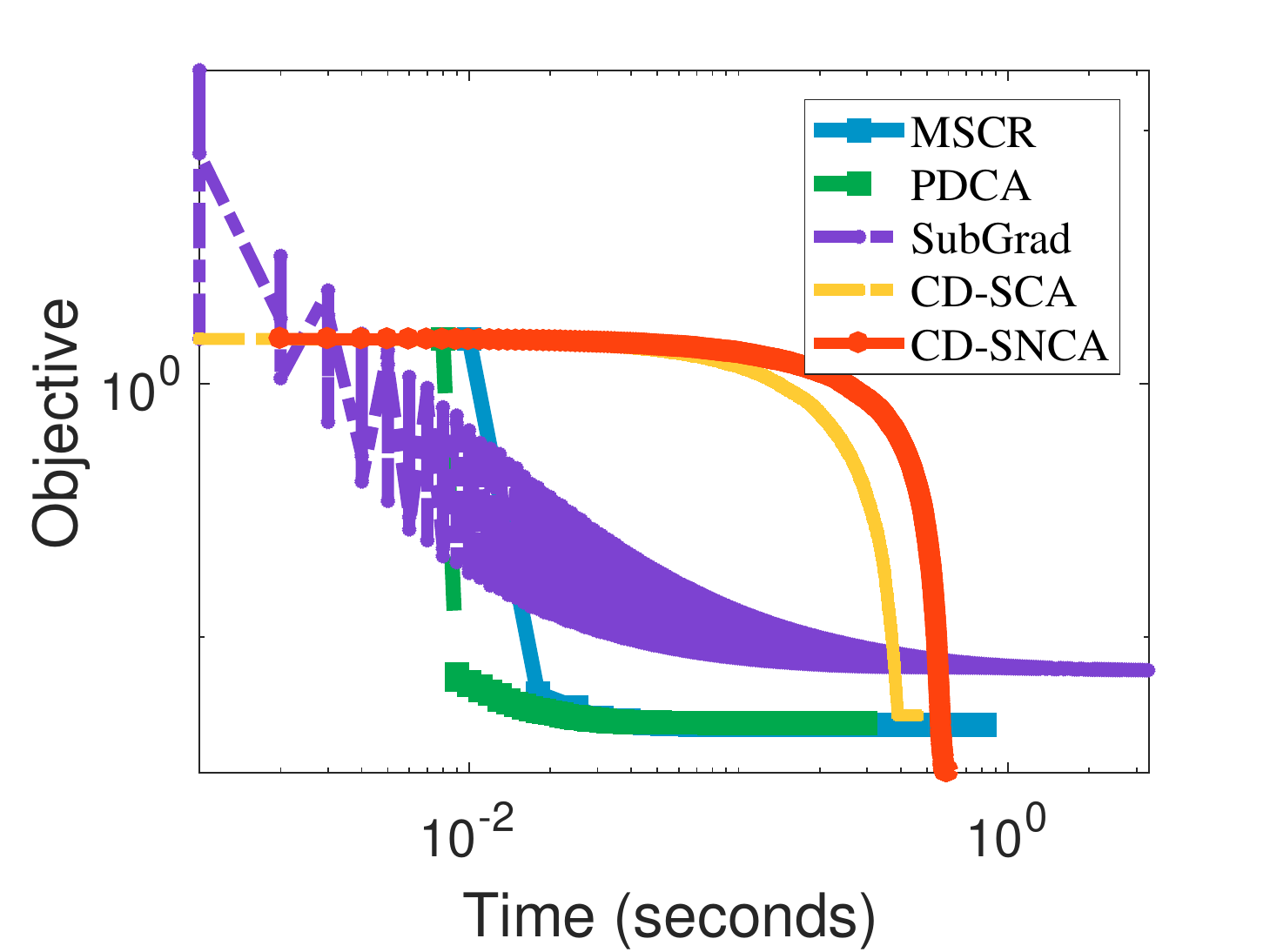}\vspace{-6pt} \caption{ \scriptsize randn-256-2048-C } \end{subfigure}~~\begin{subfigure}{0.25\textwidth}\includegraphics[height=\objimghei,width=\textwidth]{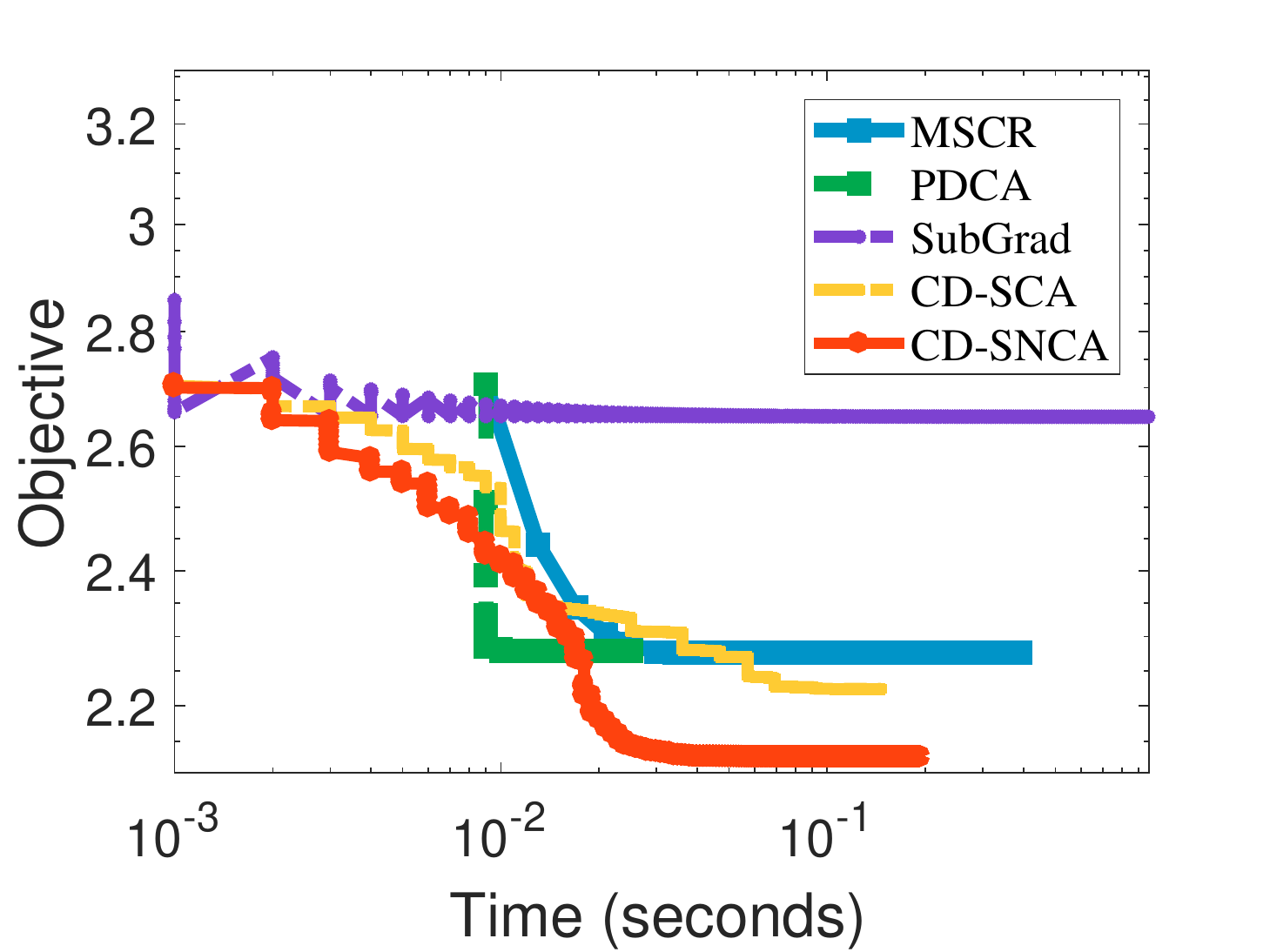}\vspace{-6pt} \caption{\scriptsize randn-1024-256-C}\end{subfigure}~~\begin{subfigure}{0.25\textwidth}\includegraphics[height=\objimghei,width=\textwidth]{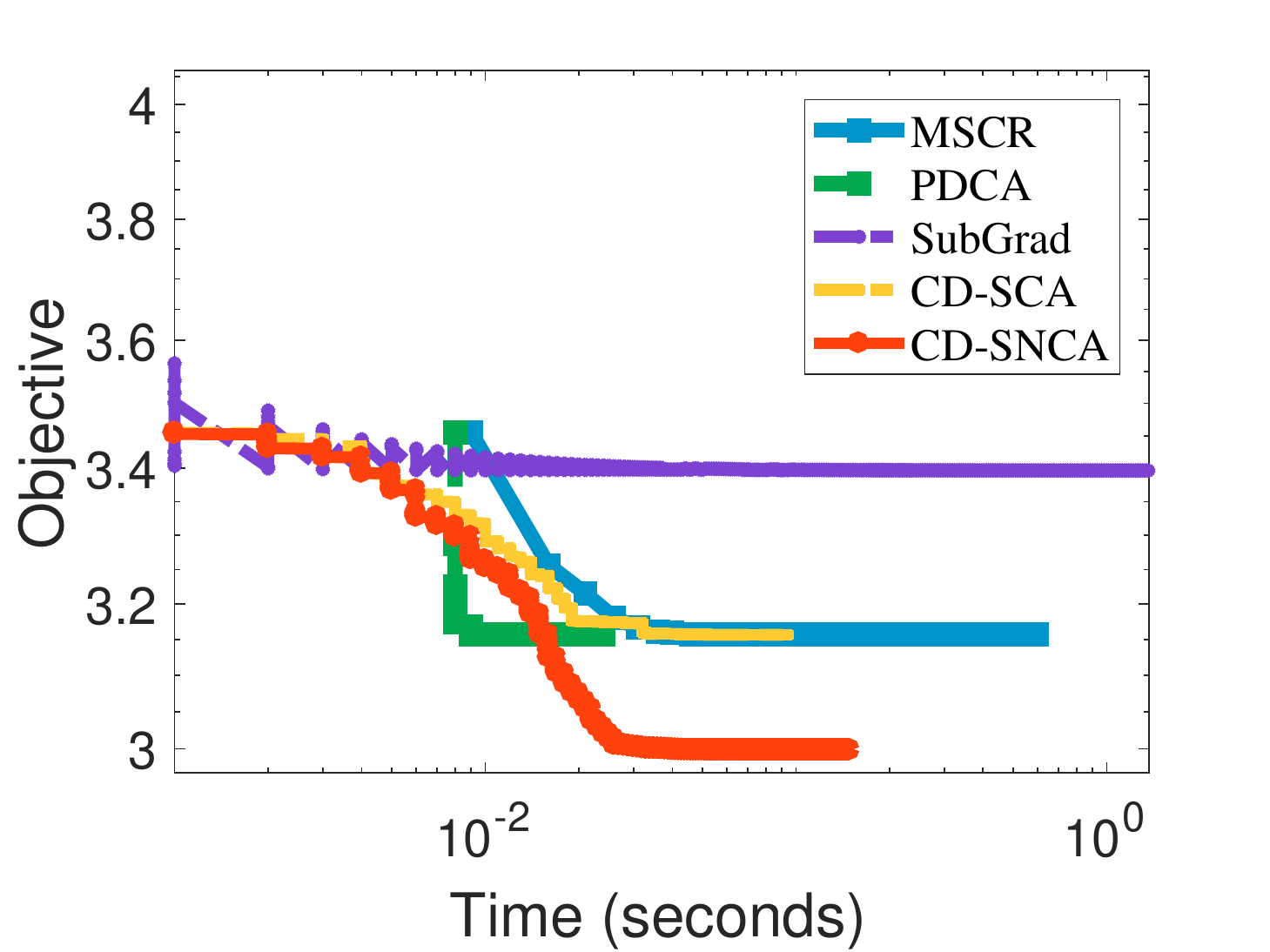}\vspace{-6pt} \caption{\scriptsize randn-2048-256-C  }\end{subfigure}\\

      \begin{subfigure}{0.25\textwidth}\includegraphics[height=\objimghei,width=\textwidth]{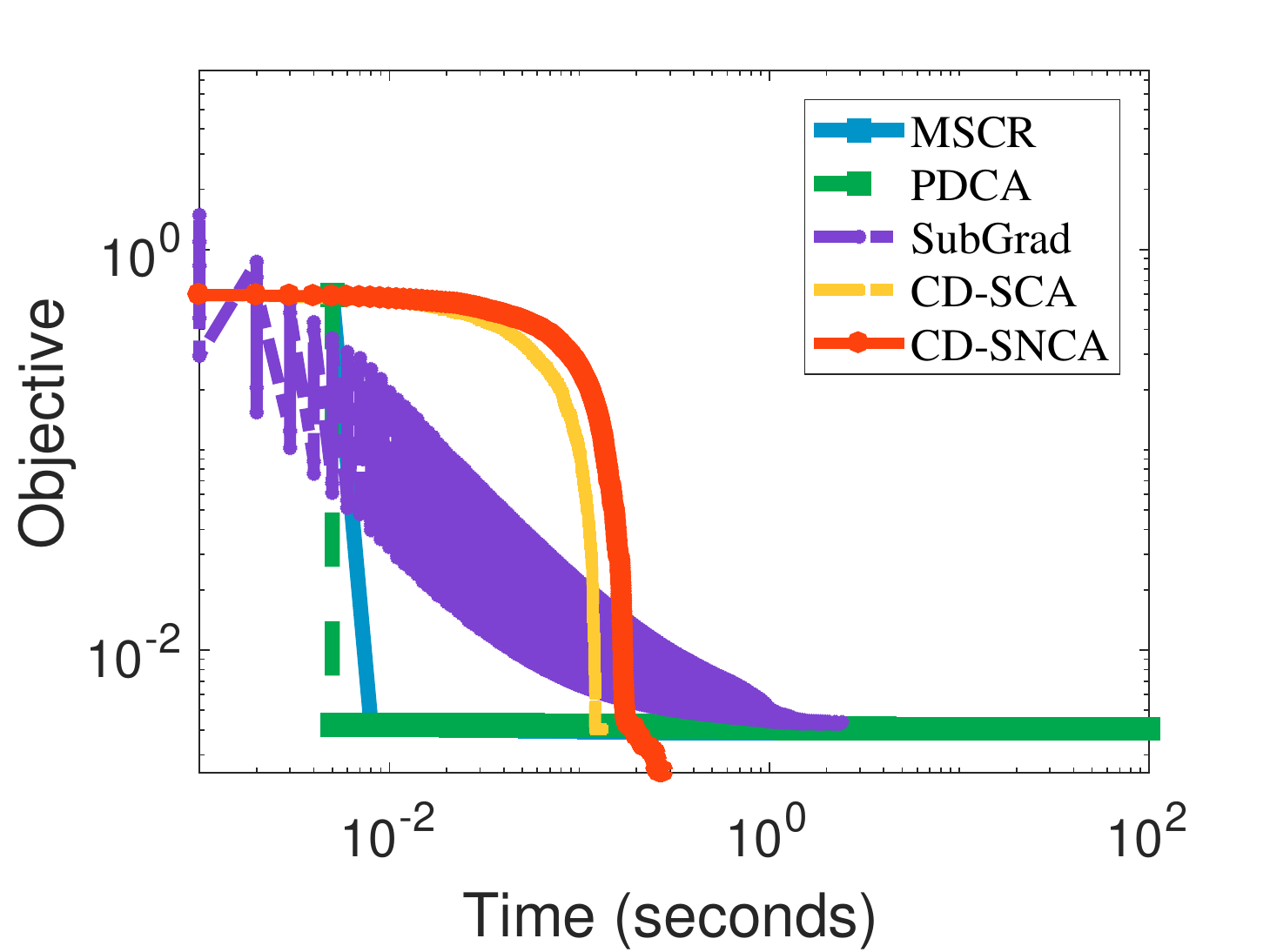}\vspace{-6pt} \caption{\scriptsize e2006-256-1024-C }\end{subfigure}~~\begin{subfigure}{0.25\textwidth}\includegraphics[height=\objimghei,width=\textwidth]{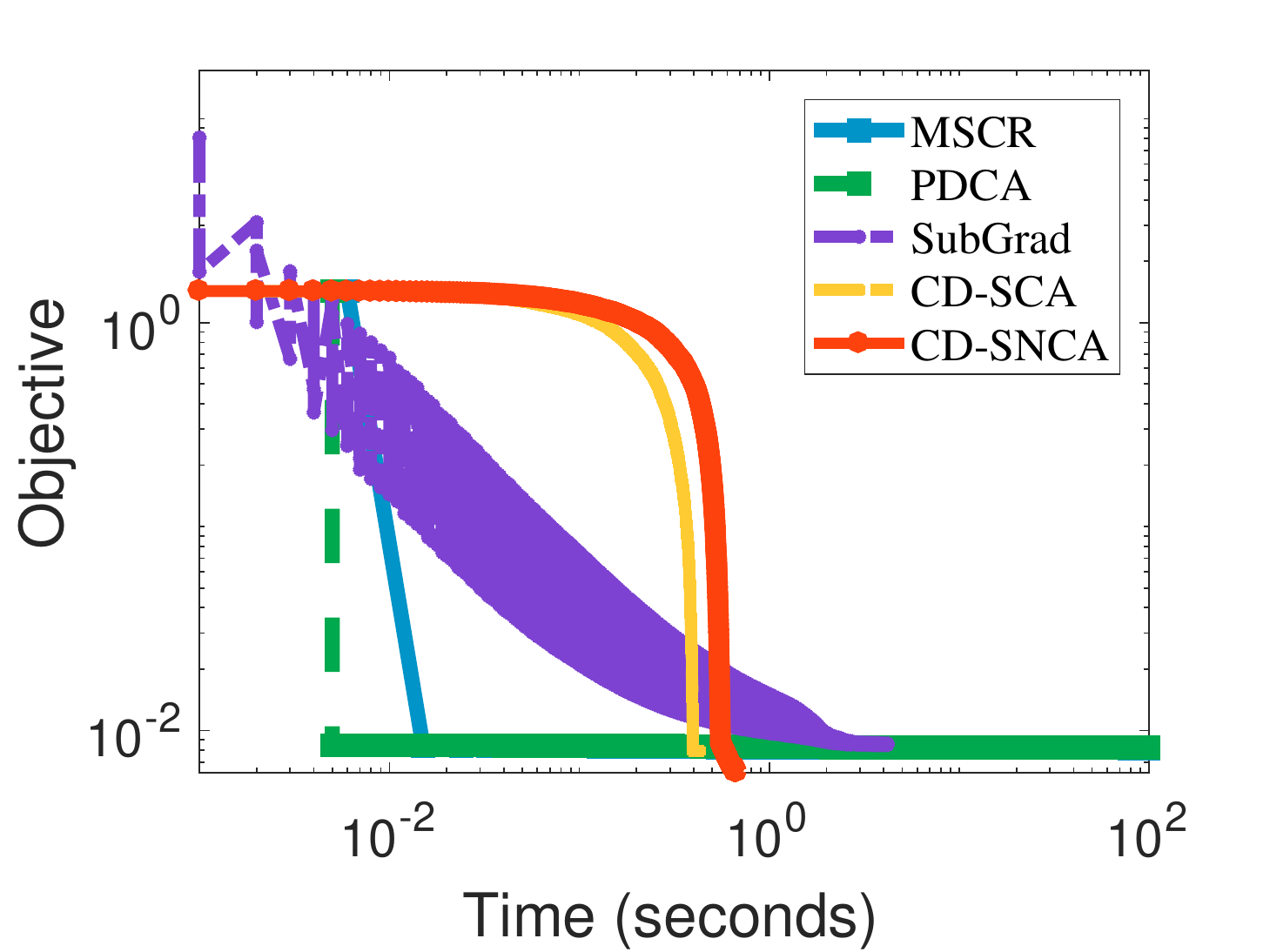}\vspace{-6pt} \caption{ \scriptsize e2006-256-2048-C } \end{subfigure}~~\begin{subfigure}{0.25\textwidth}\includegraphics[height=\objimghei,width=\textwidth]{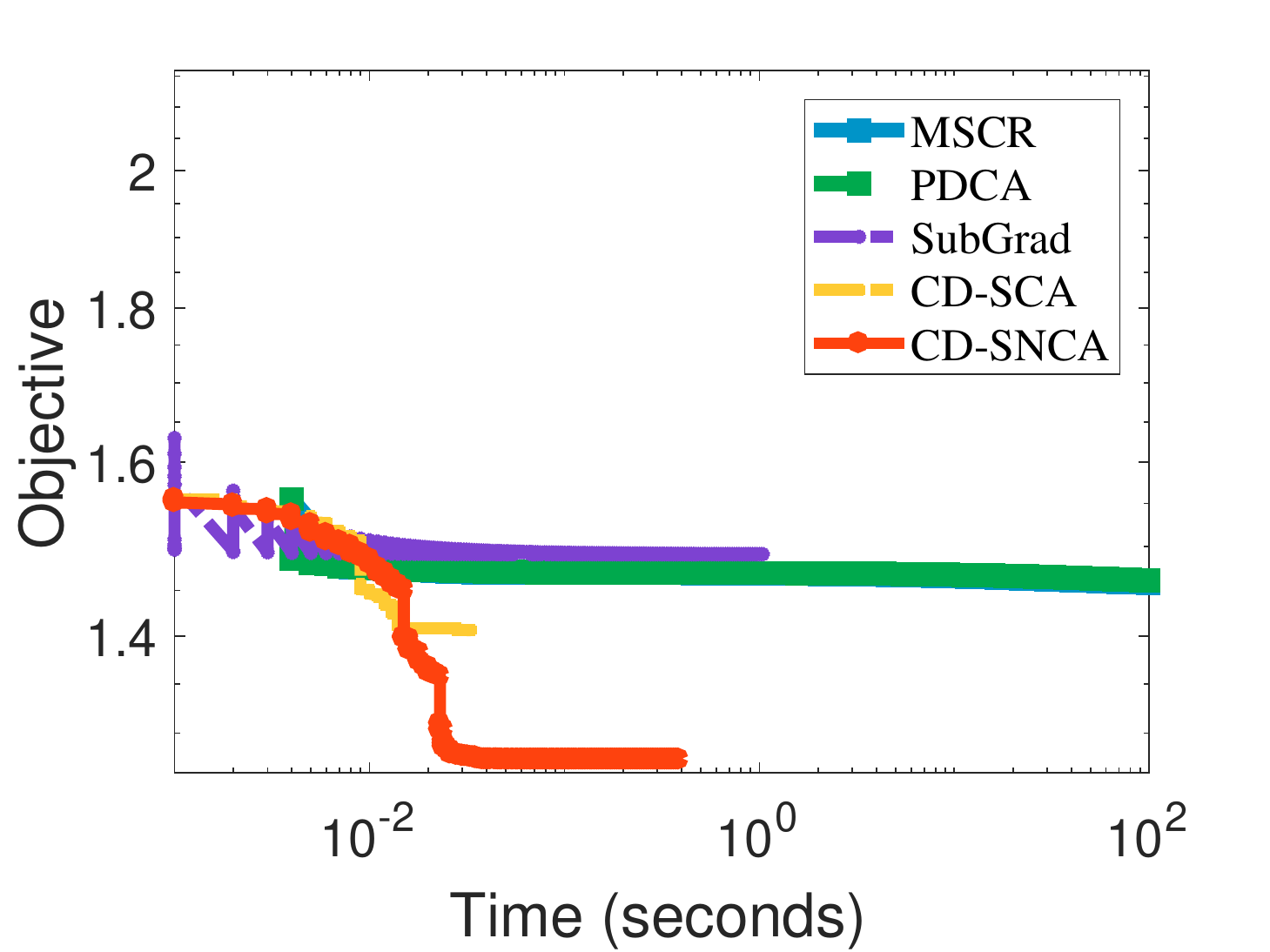}\vspace{-6pt} \caption{\scriptsize e2006-1024-256-C}\end{subfigure}~~\begin{subfigure}{0.25\textwidth}\includegraphics[height=\objimghei,width=\textwidth]{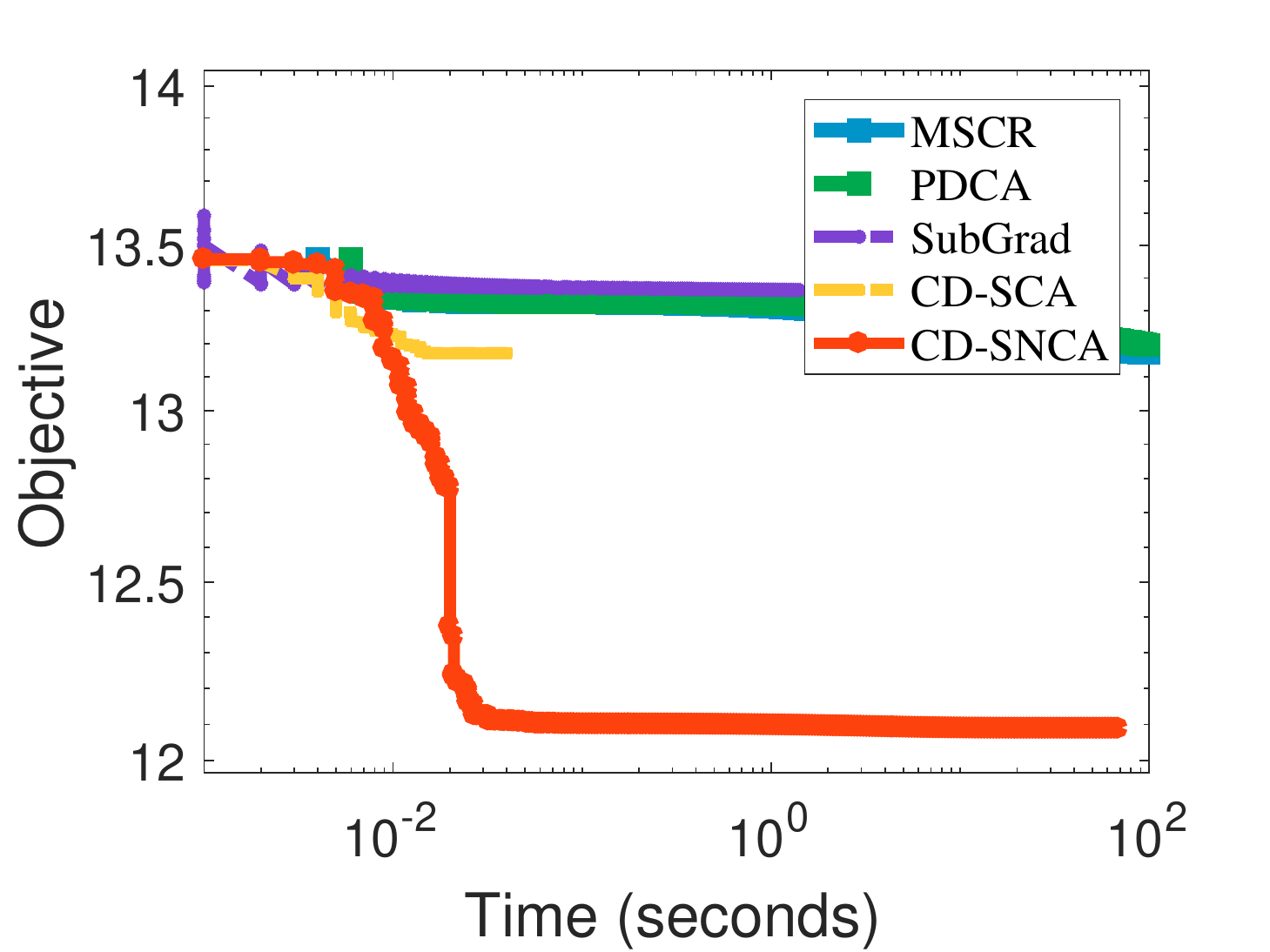}\vspace{-6pt} \caption{\scriptsize e2006-2048-256-C  }\end{subfigure}\\

\centering
\caption{The convergence curve of the compared methods for solving the approximate sparse optimization problem on different data sets.}
\label{exp:cpu:2}

\end{figure*}

\begin{figure*} [!t]
\centering
      \begin{subfigure}{0.25\textwidth}\includegraphics[height=\objimghei,width=\textwidth]{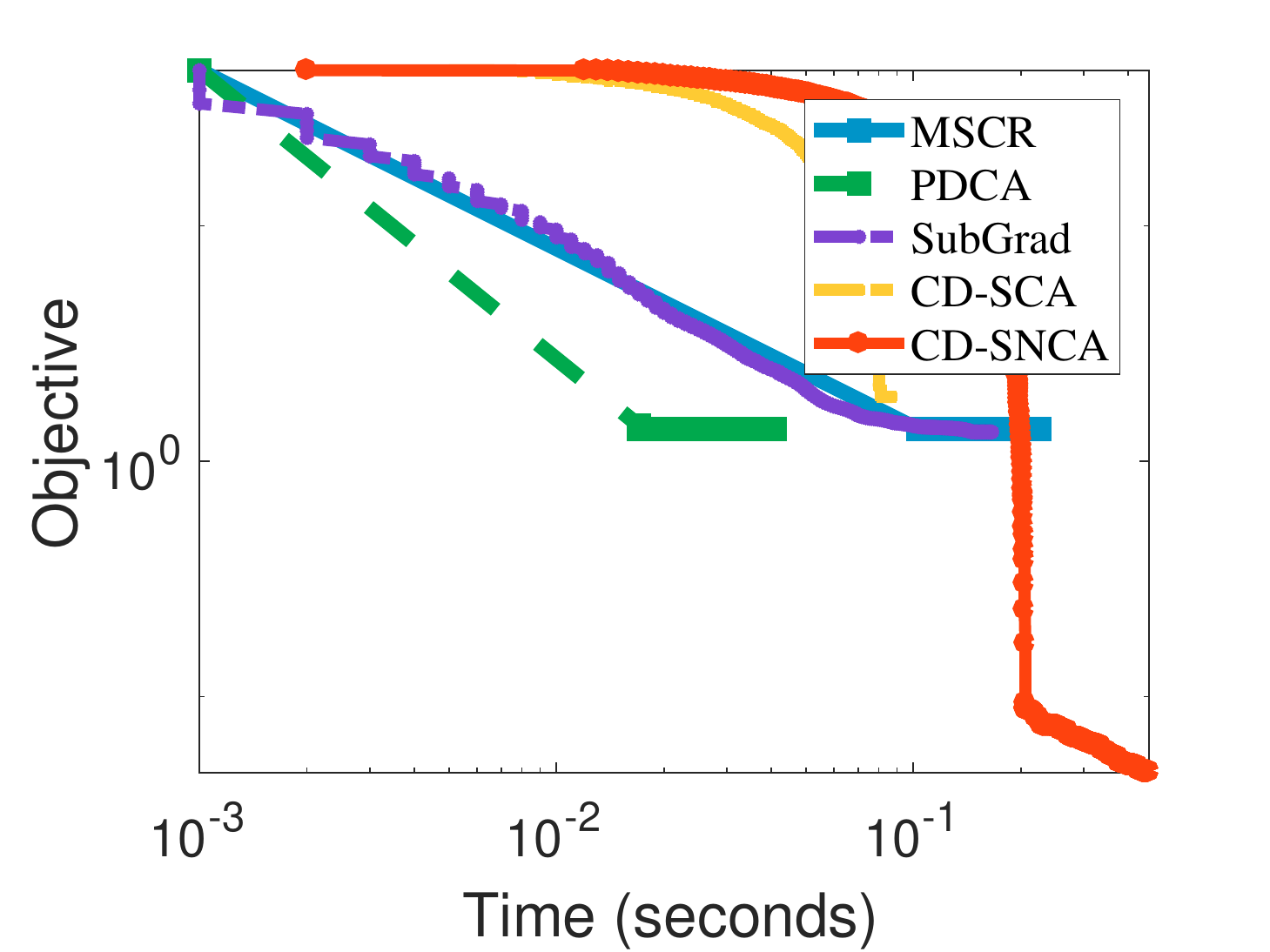}\vspace{-6pt} \caption{\scriptsize randn-256-1024  }\end{subfigure}~~\begin{subfigure}{0.25\textwidth}\includegraphics[height=\objimghei,width=\textwidth]{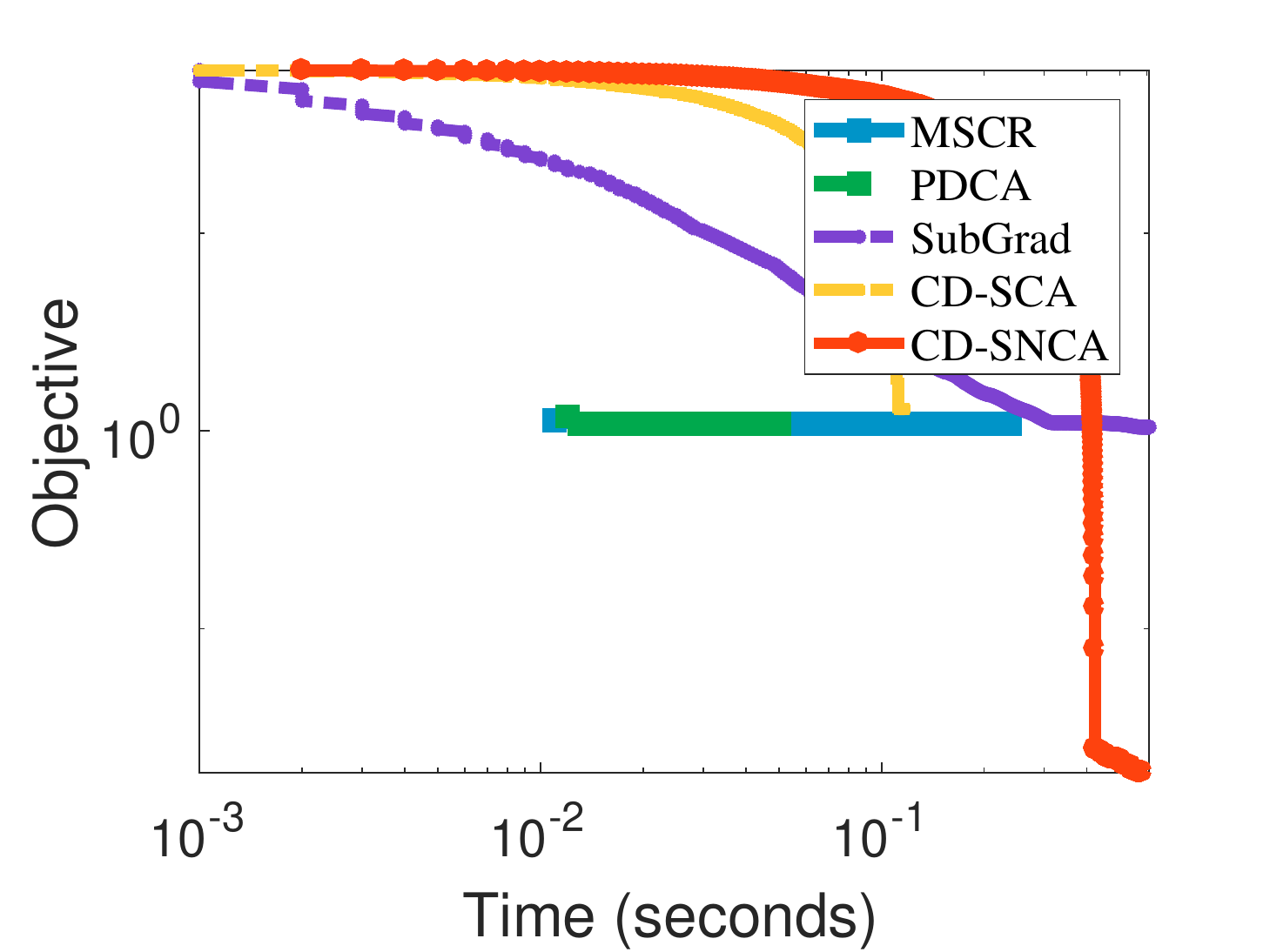}\vspace{-6pt} \caption{ \scriptsize randn-256-2048 } \end{subfigure}~~\begin{subfigure}{0.25\textwidth}\includegraphics[height=\objimghei,width=\textwidth]{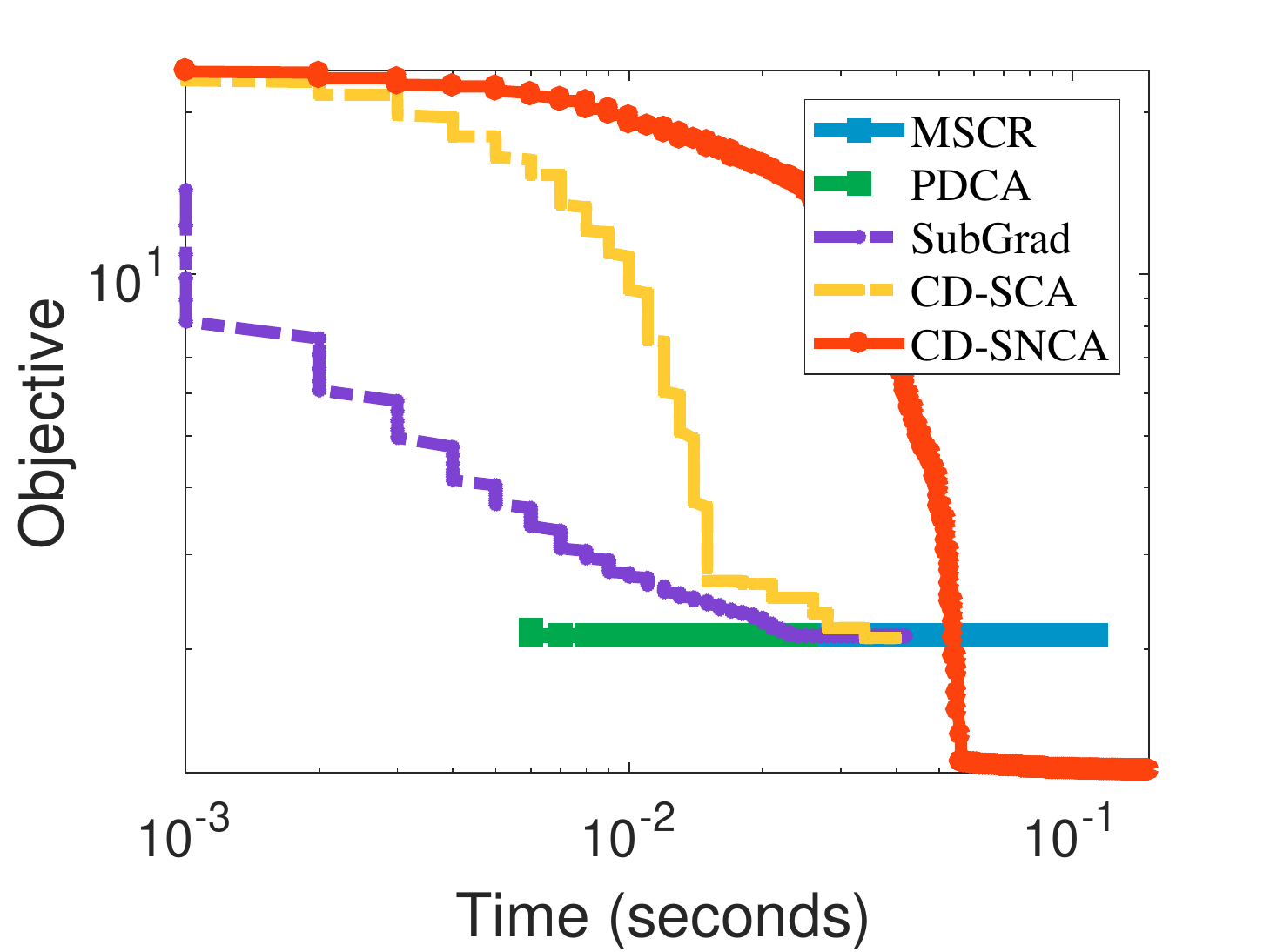}\vspace{-6pt} \caption{\scriptsize randn-1024-256}\end{subfigure}~~\begin{subfigure}{0.25\textwidth}\includegraphics[height=\objimghei,width=\textwidth]{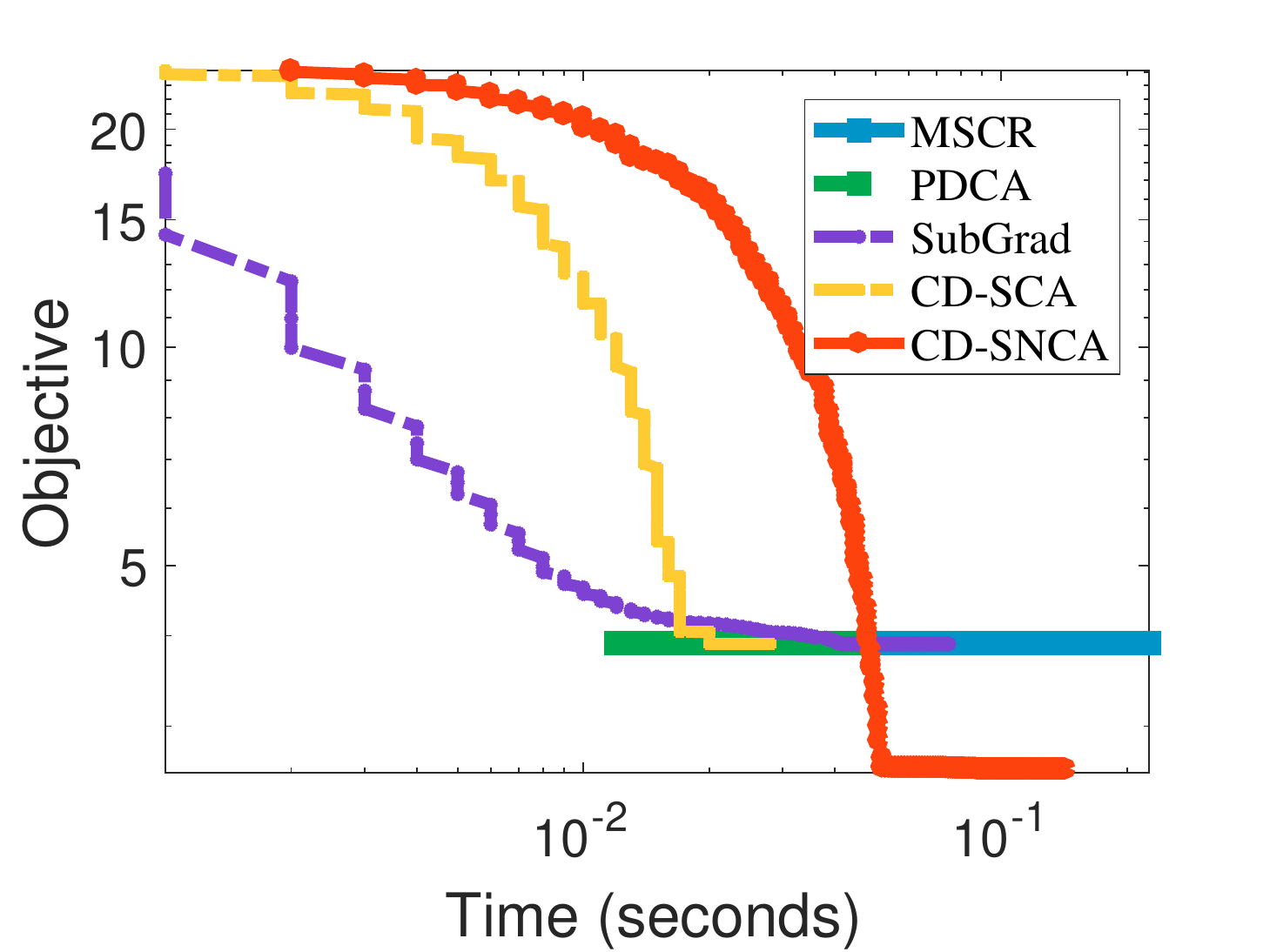}\vspace{-6pt} \caption{\scriptsize randn-2048-256 }\end{subfigure}\\

      \begin{subfigure}{0.25\textwidth}\includegraphics[height=\objimghei,width=\textwidth]{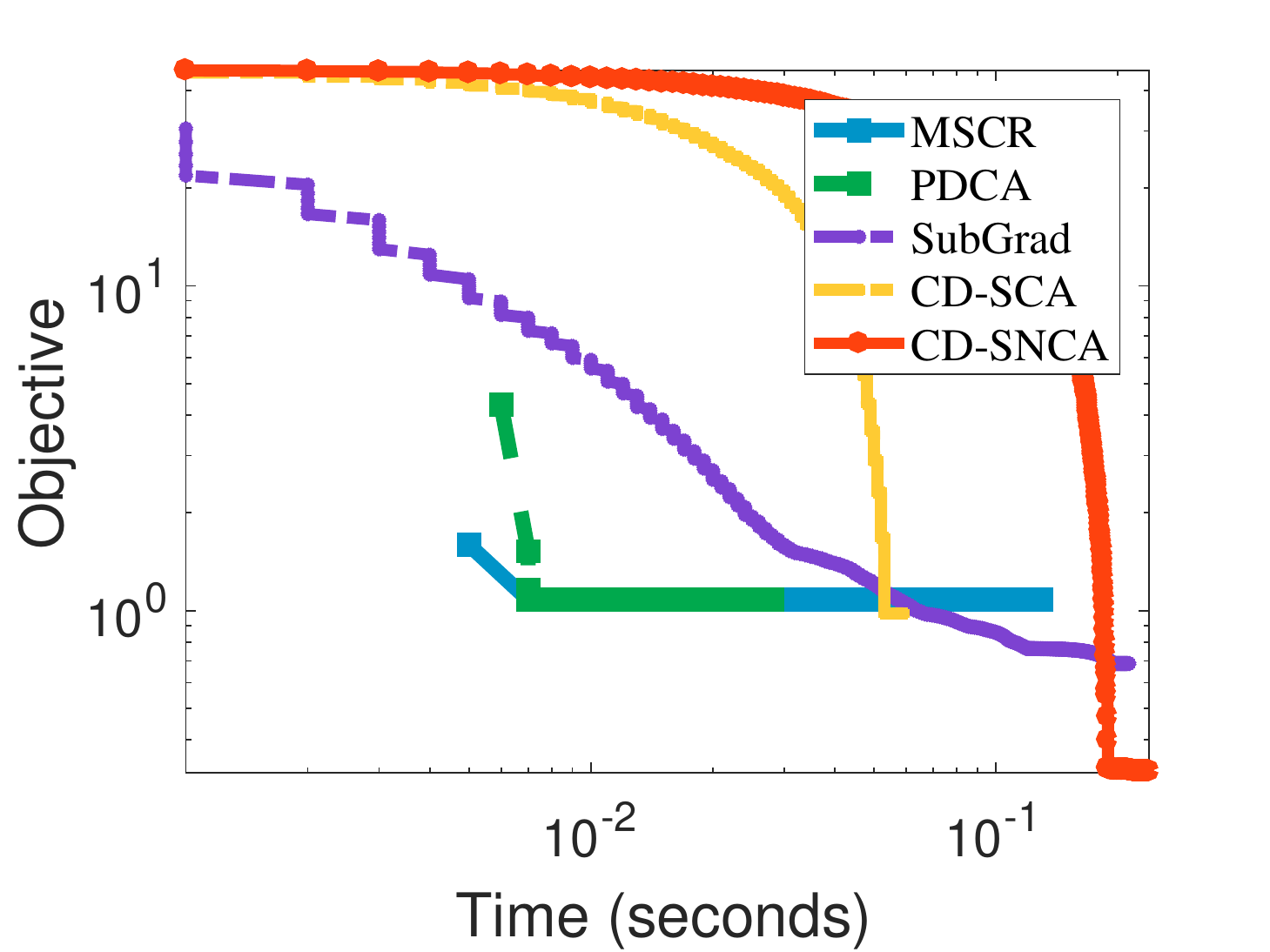}\vspace{-6pt} \caption{\scriptsize e2006-256-1024 }\end{subfigure}~~\begin{subfigure}{0.25\textwidth}\includegraphics[height=\objimghei,width=\textwidth]{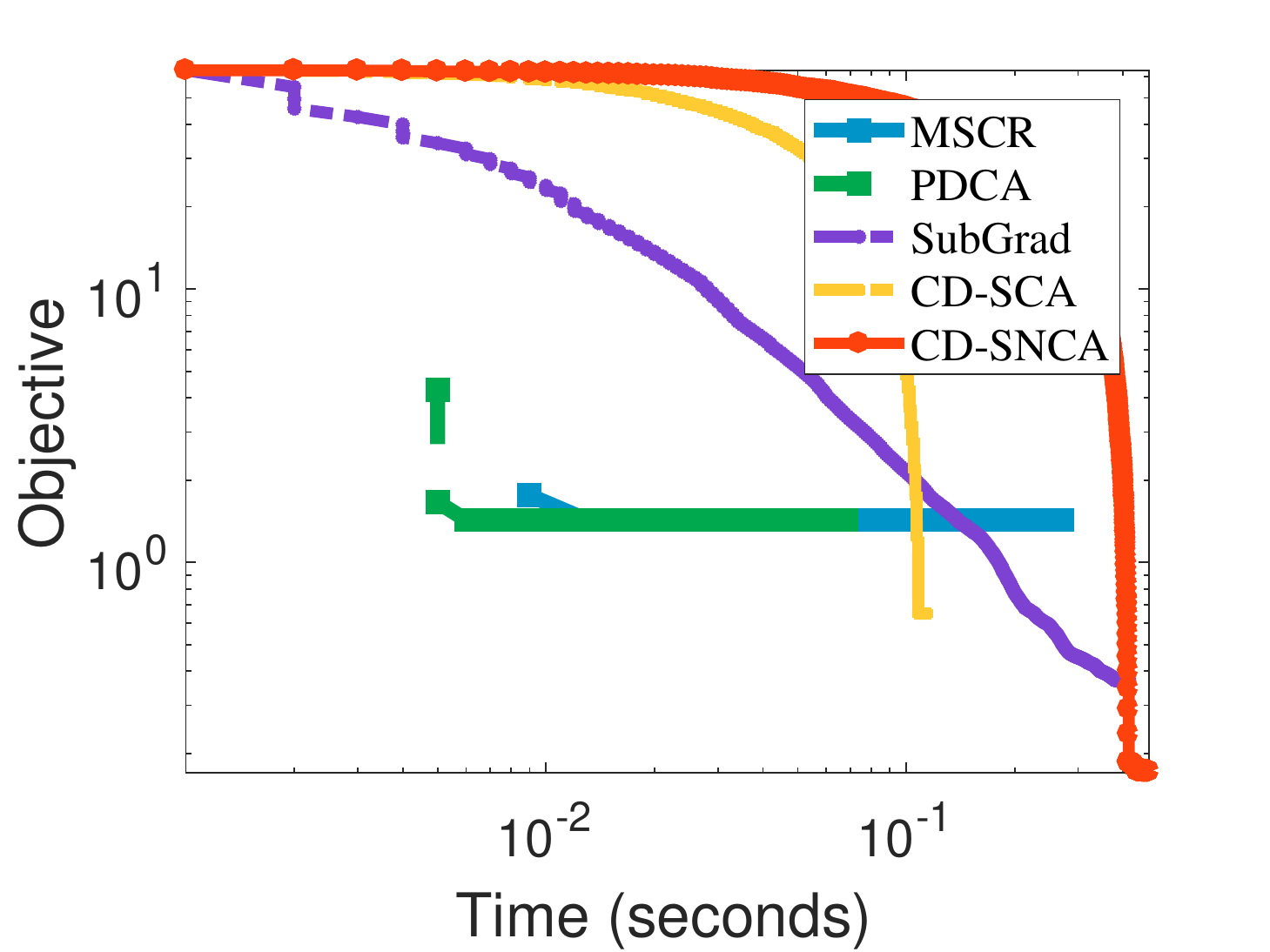}\vspace{-6pt} \caption{ \scriptsize e2006-256-2048 } \end{subfigure}~~\begin{subfigure}{0.25\textwidth}\includegraphics[height=\objimghei,width=\textwidth]{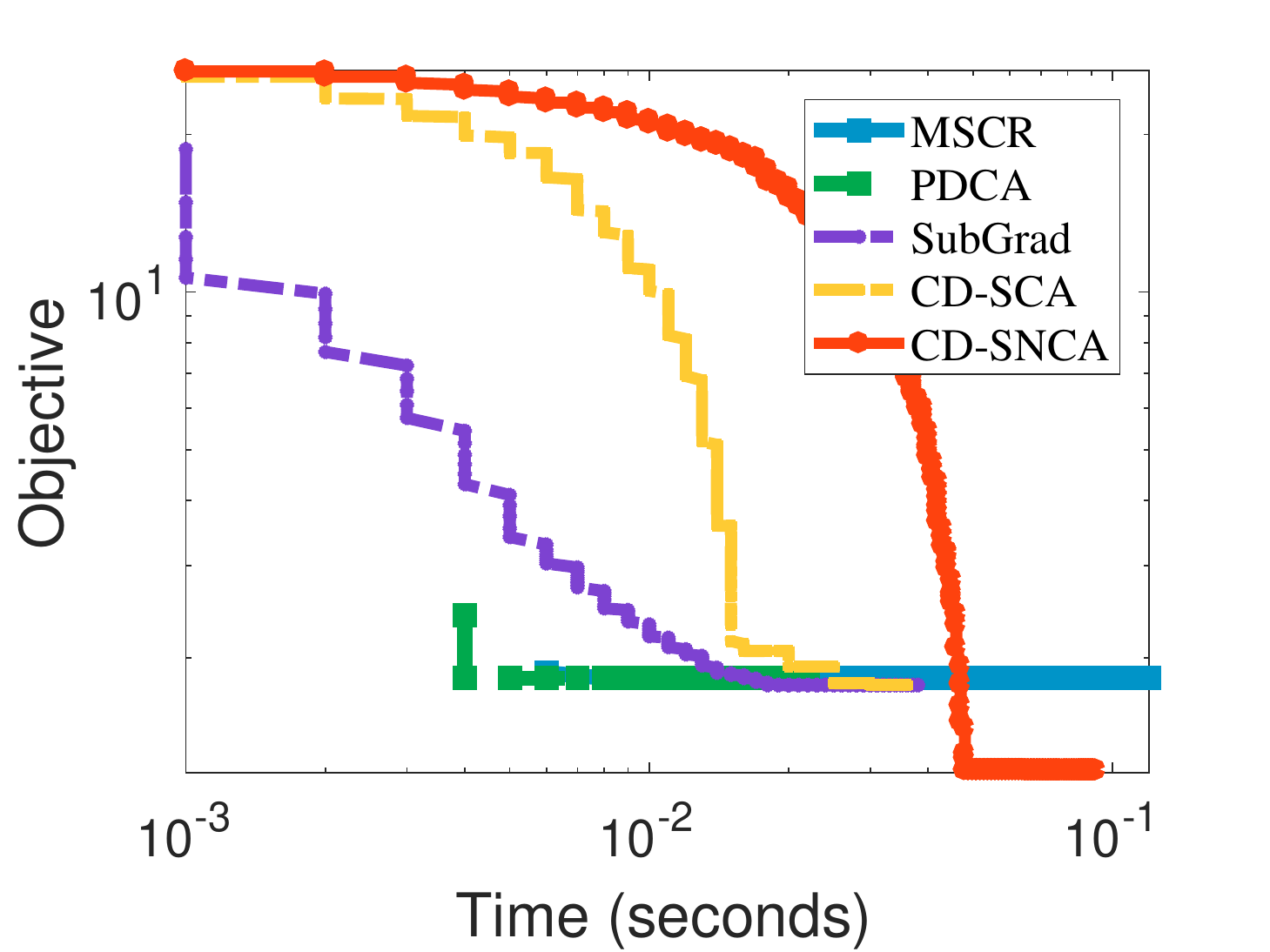}\vspace{-6pt} \caption{\scriptsize e2006-1024-256}\end{subfigure}~~\begin{subfigure}{0.25\textwidth}\includegraphics[height=\objimghei,width=\textwidth]{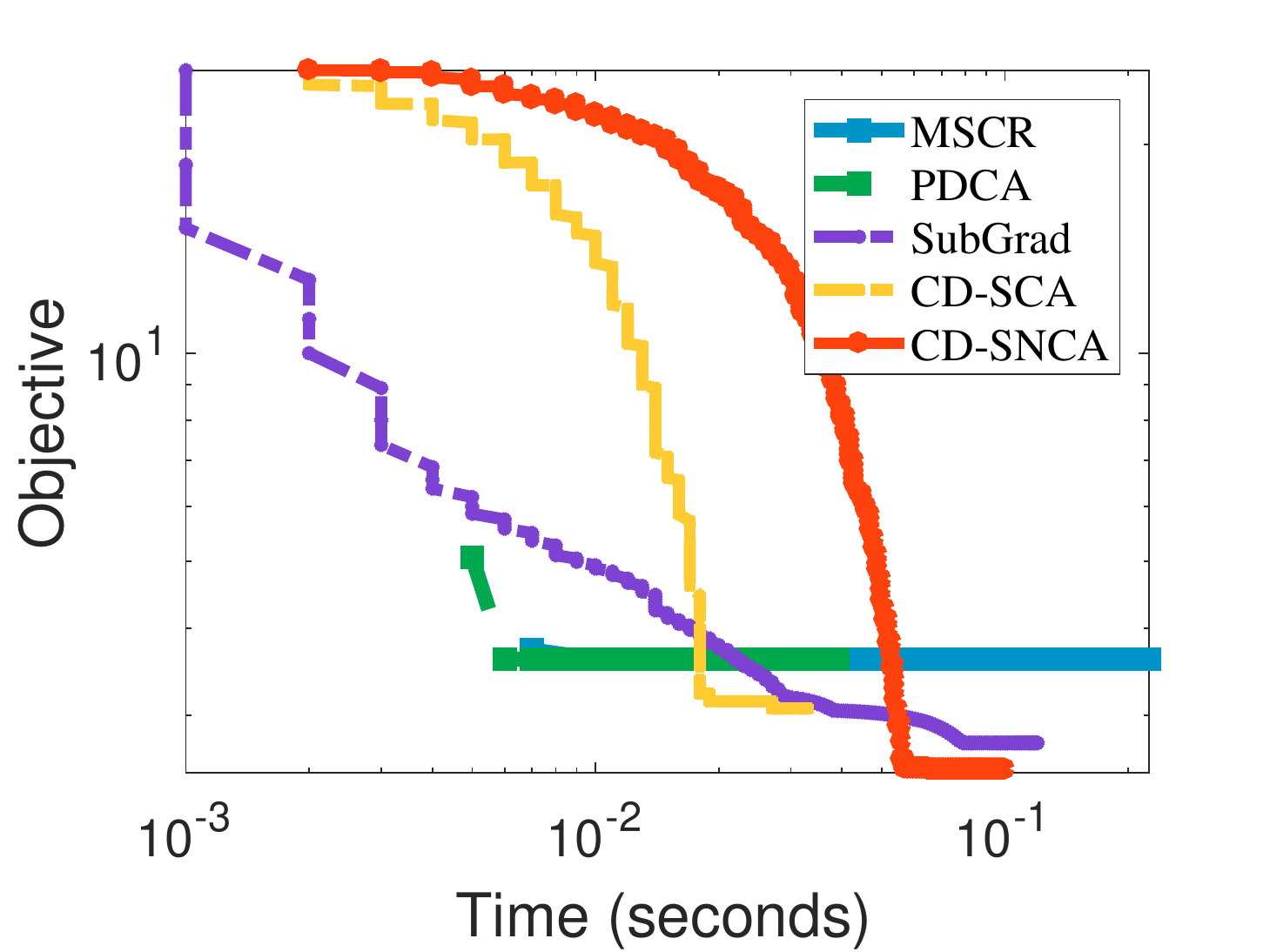}\vspace{-6pt} \caption{\scriptsize e2006-2048-256  }\end{subfigure}\\

      \begin{subfigure}{0.25\textwidth}\includegraphics[height=\objimghei,width=\textwidth]{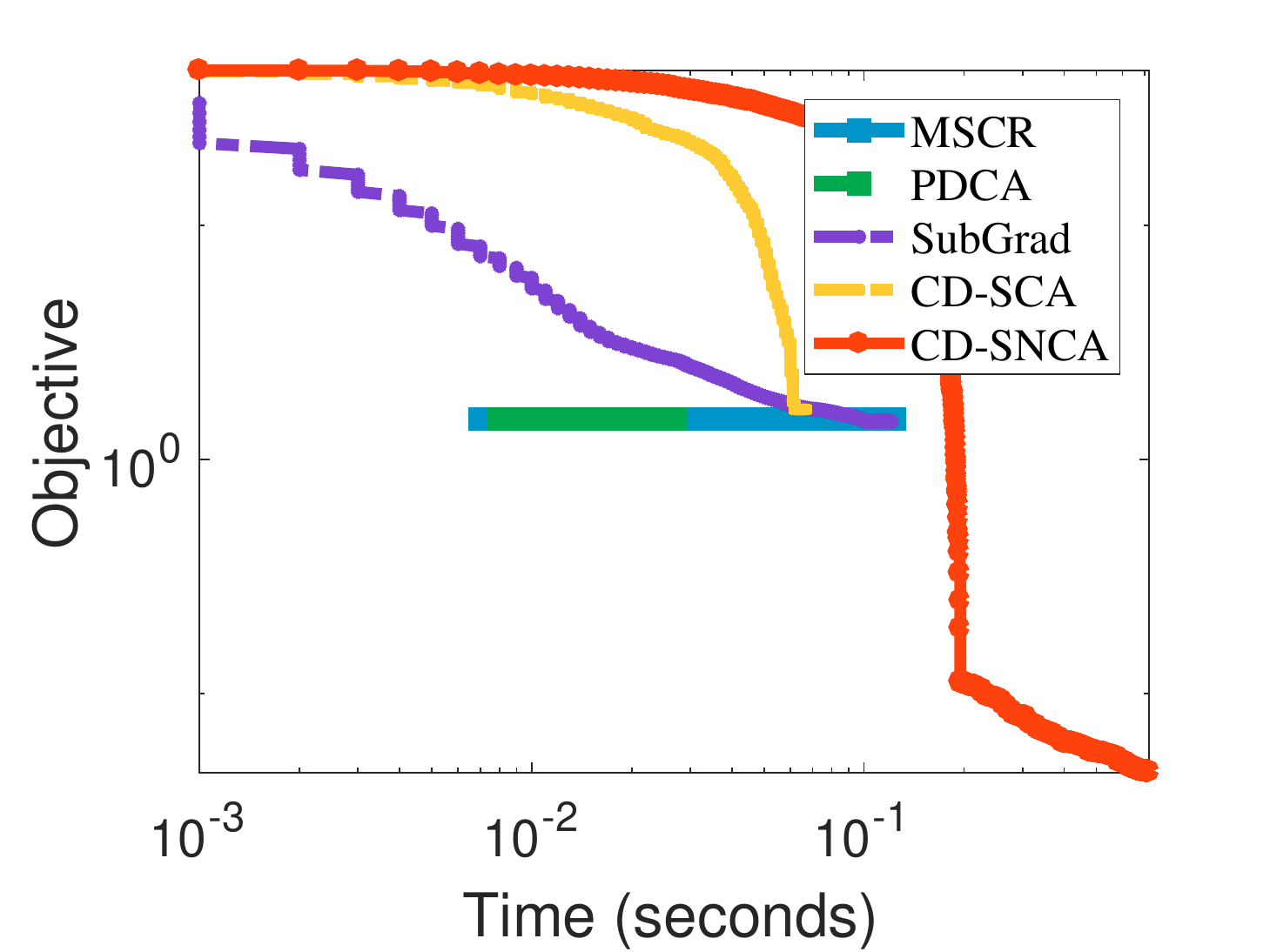}\vspace{-6pt} \caption{\scriptsize randn-256-1024-C }\end{subfigure}~~\begin{subfigure}{0.25\textwidth}\includegraphics[height=\objimghei,width=\textwidth]{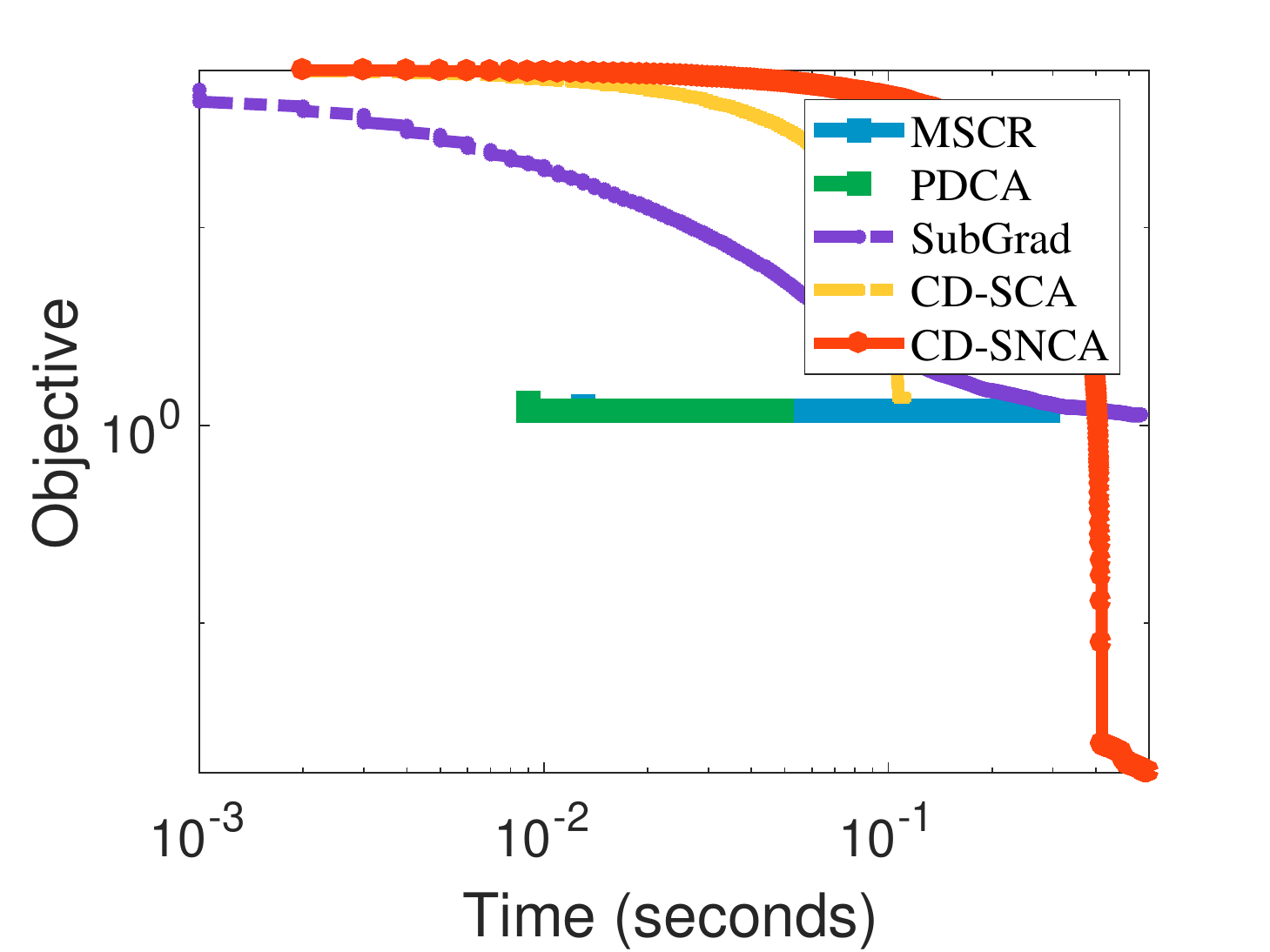}\vspace{-6pt} \caption{ \scriptsize randn-256-2048-C } \end{subfigure}~~\begin{subfigure}{0.25\textwidth}\includegraphics[height=\objimghei,width=\textwidth]{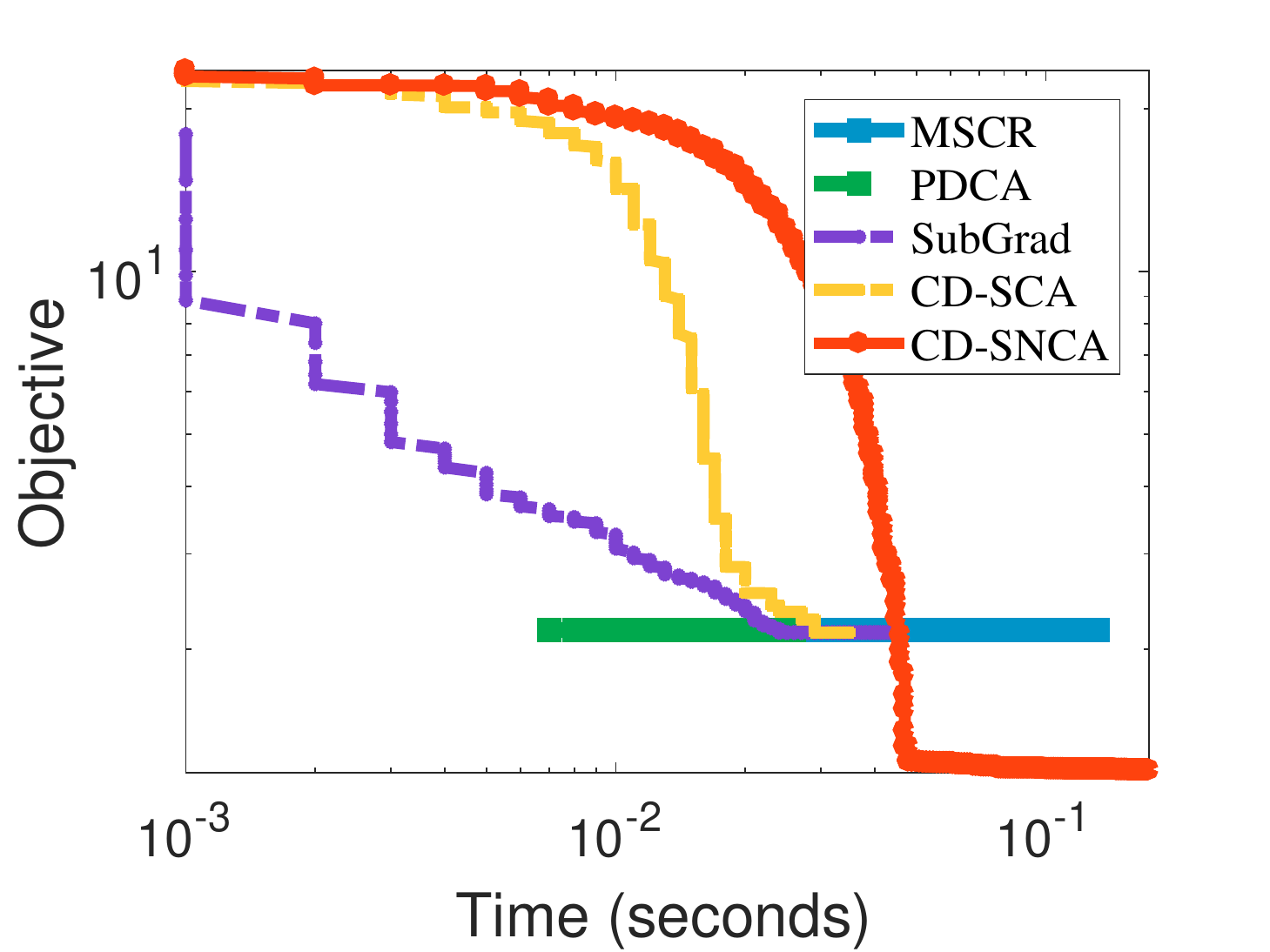}\vspace{-6pt} \caption{\scriptsize randn-1024-256-C}\end{subfigure}~~\begin{subfigure}{0.25\textwidth}\includegraphics[height=\objimghei,width=\textwidth]{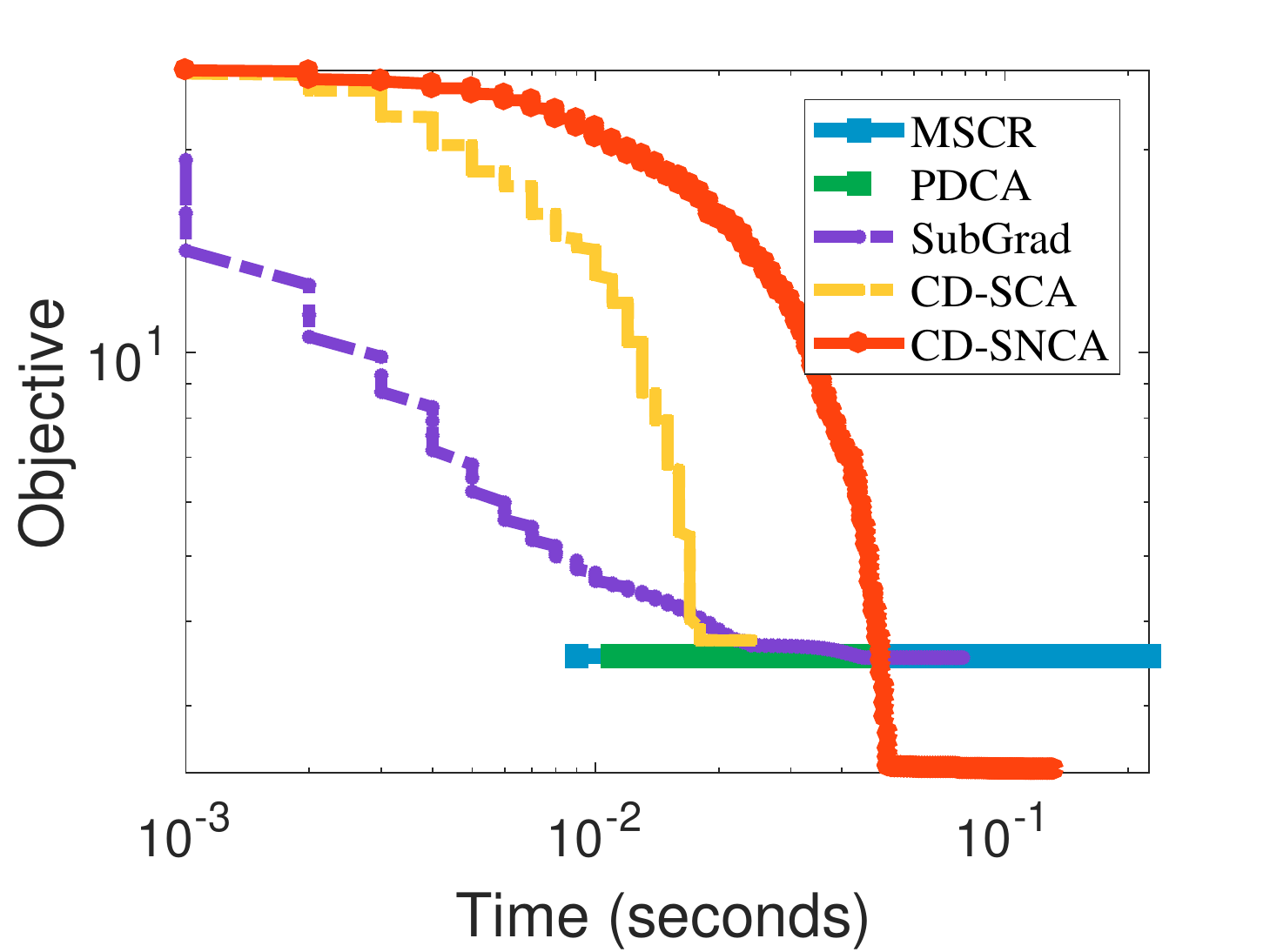}\vspace{-6pt} \caption{\scriptsize randn-2048-256-C }\end{subfigure}\\

      \begin{subfigure}{0.25\textwidth}\includegraphics[height=\objimghei,width=\textwidth]{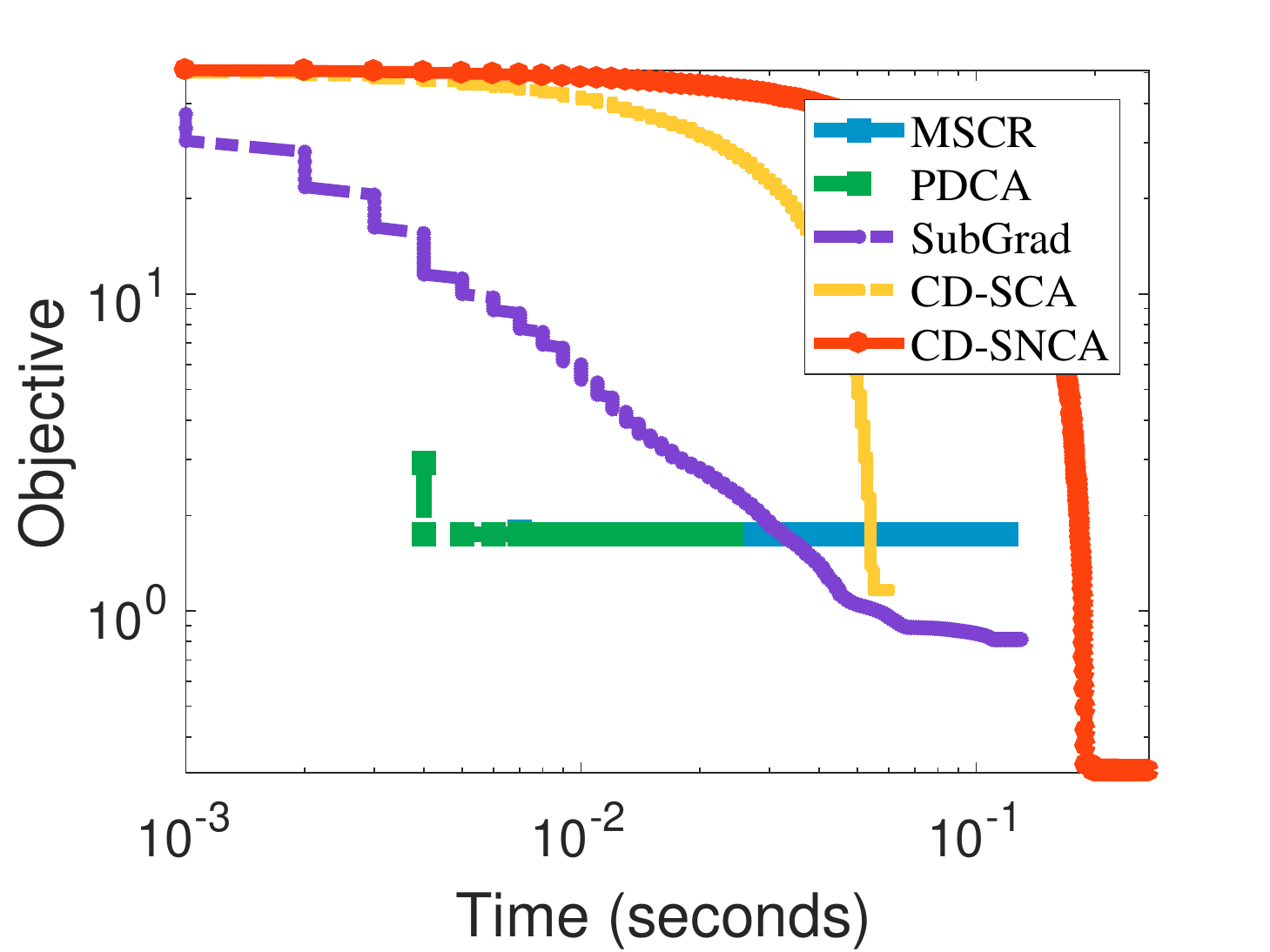}\vspace{-6pt} \caption{\scriptsize e2006-256-1024-C }\end{subfigure}~~\begin{subfigure}{0.25\textwidth}\includegraphics[height=\objimghei,width=\textwidth]{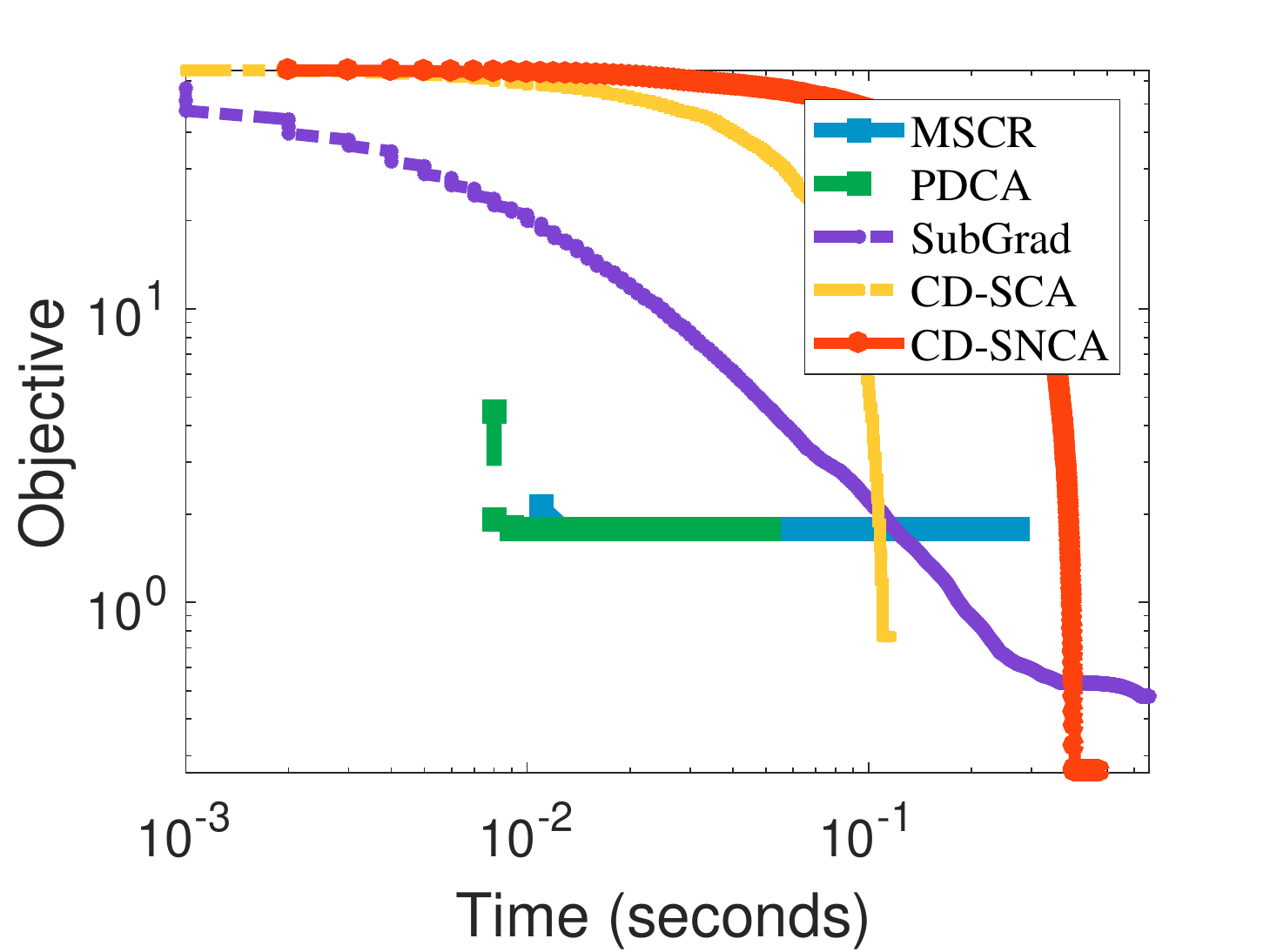}\vspace{-6pt} \caption{ \scriptsize e2006-256-2048-C } \end{subfigure}~~\begin{subfigure}{0.25\textwidth}\includegraphics[height=\objimghei,width=\textwidth]{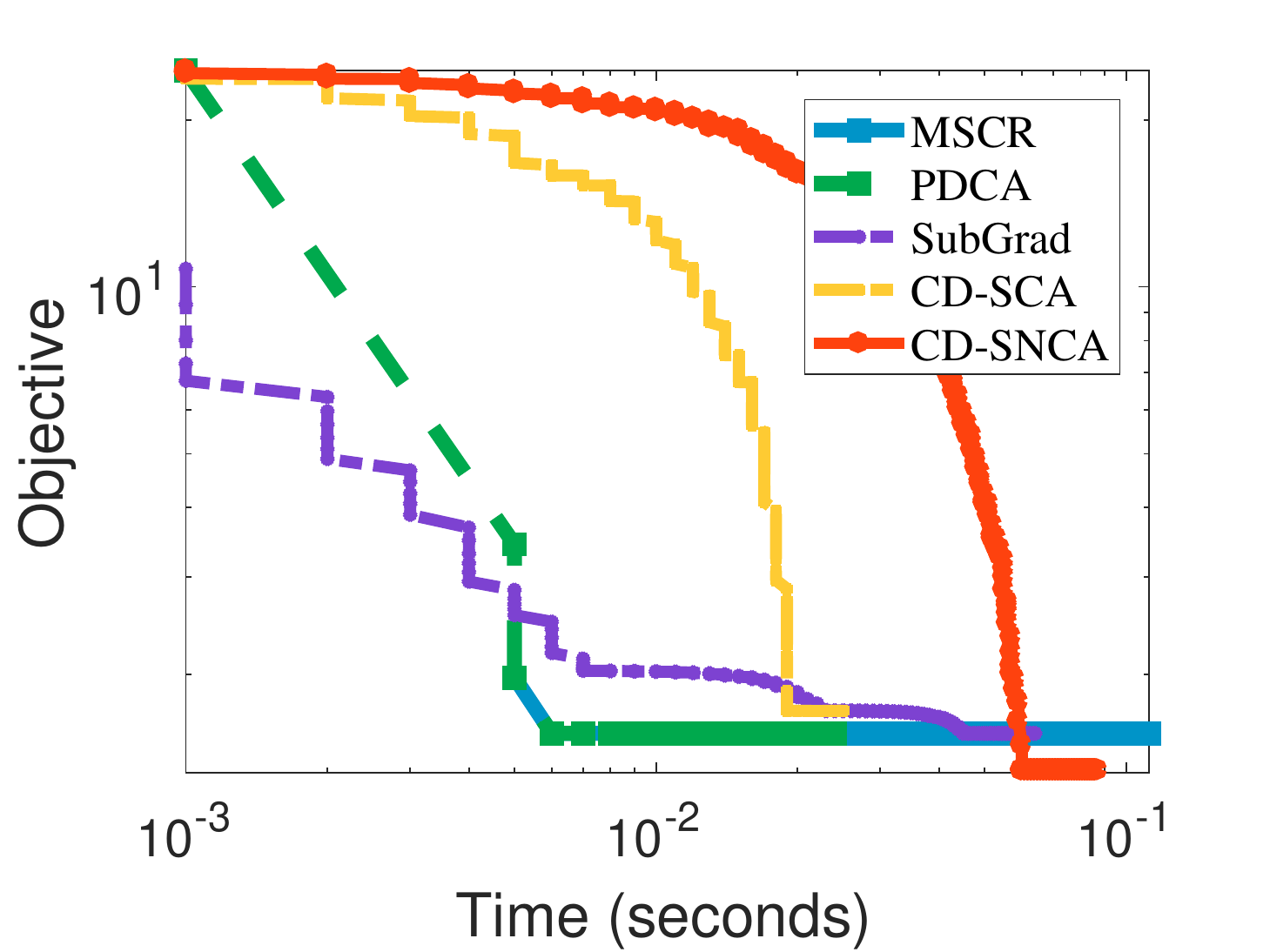}\vspace{-6pt} \caption{\scriptsize e2006-1024-256-C}\end{subfigure}~~\begin{subfigure}{0.25\textwidth}\includegraphics[height=\objimghei,width=\textwidth]{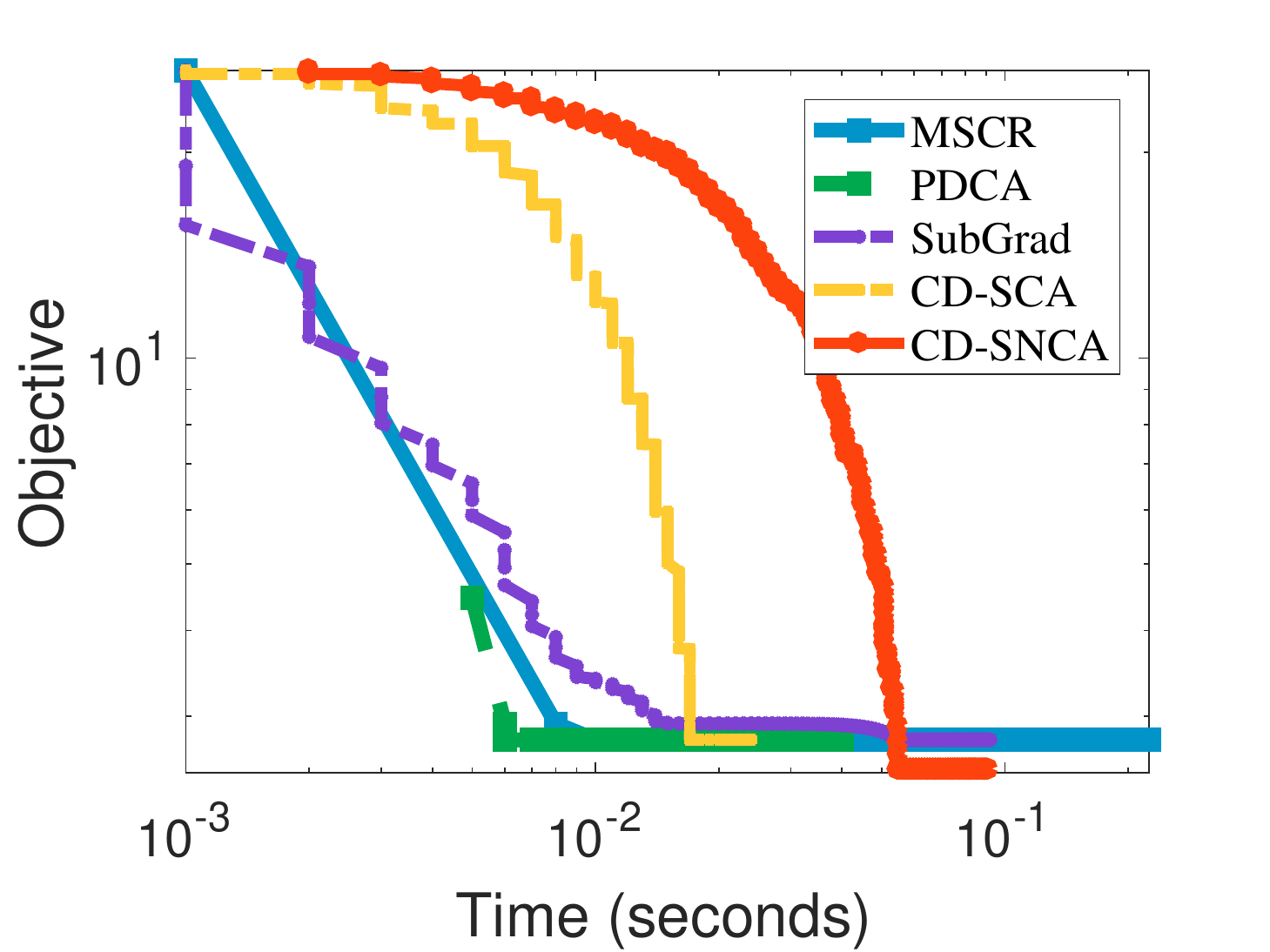}\vspace{-6pt} \caption{\scriptsize e2006-2048-256-C }\end{subfigure}\\

\centering
\caption{The convergence curve of the compared methods for solving the approximate binary optimization problem on different data sets.}
\label{exp:cpu:3}
\end{figure*}

\begin{figure*} [!t]
\centering
      \begin{subfigure}{0.25\textwidth}\includegraphics[height=\objimghei,width=\textwidth]{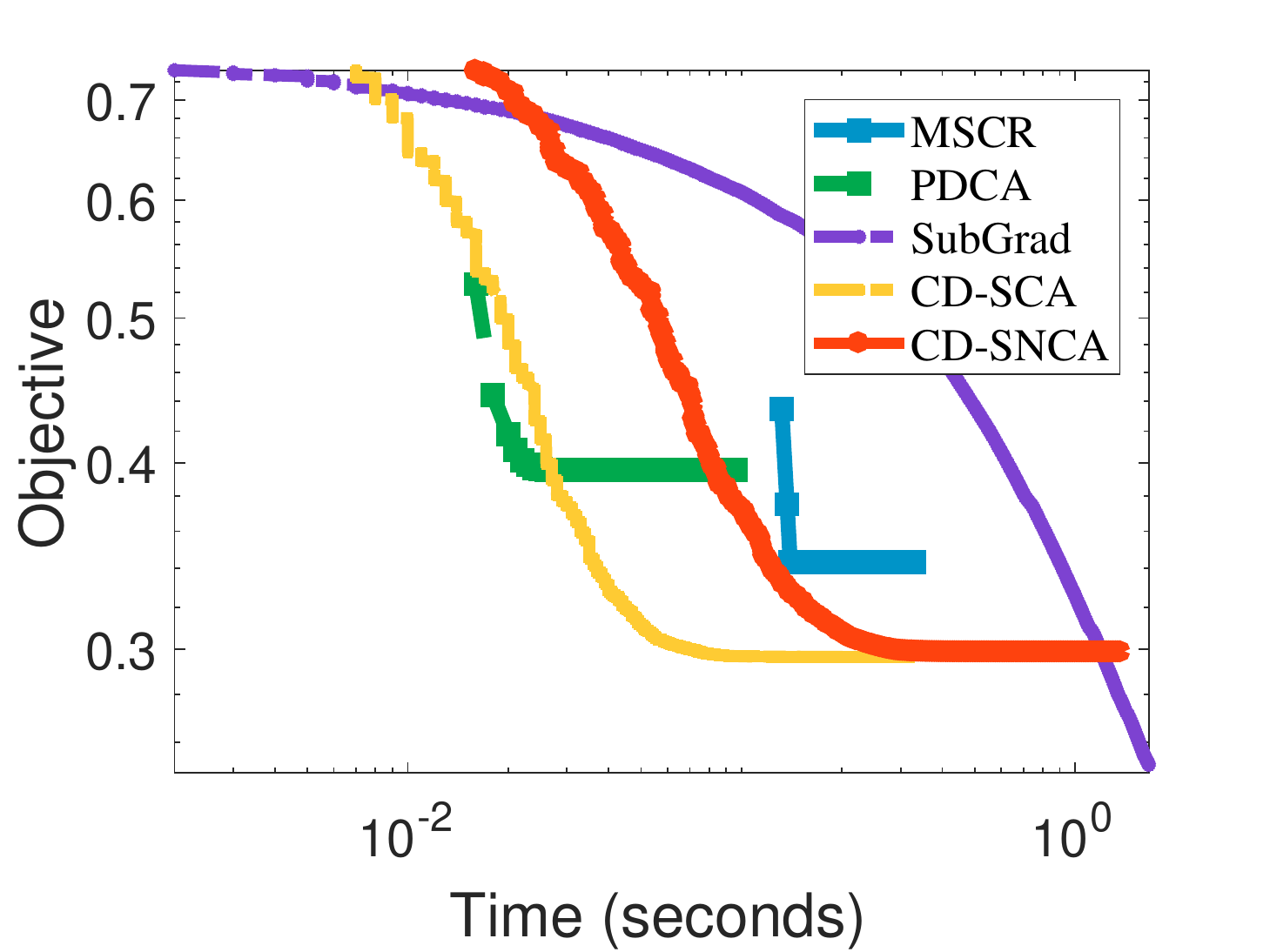}\vspace{-6pt} \caption{\scriptsize randn-256-1024 }\end{subfigure}~~\begin{subfigure}{0.25\textwidth}\includegraphics[height=\objimghei,width=\textwidth]{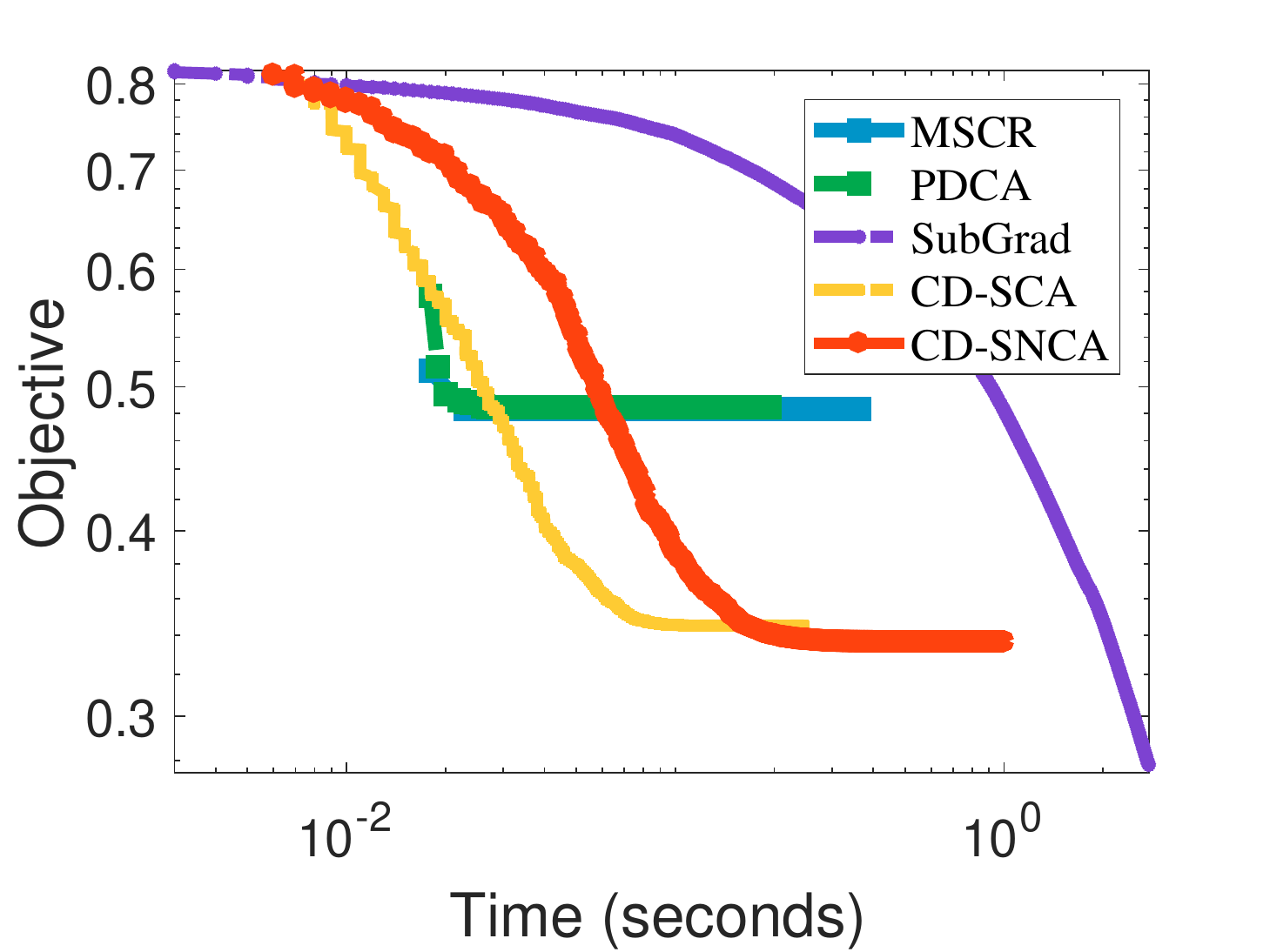}\vspace{-6pt} \caption{ \scriptsize randn-256-2048 } \end{subfigure}~~\begin{subfigure}{0.25\textwidth}\includegraphics[height=\objimghei,width=\textwidth]{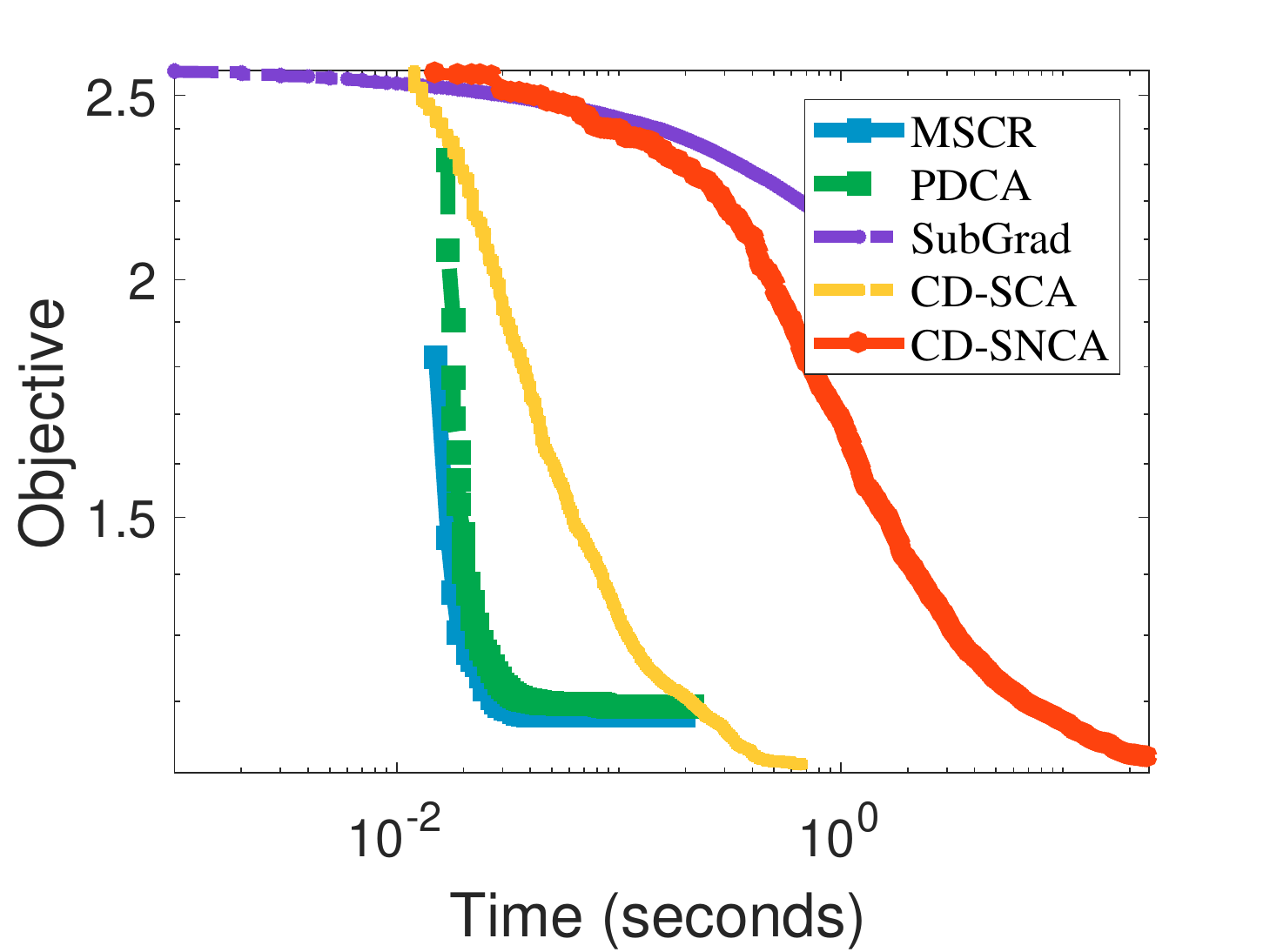}\vspace{-6pt} \caption{\scriptsize randn-1024-256}\end{subfigure}~~\begin{subfigure}{0.25\textwidth}\includegraphics[height=\objimghei,width=\textwidth]{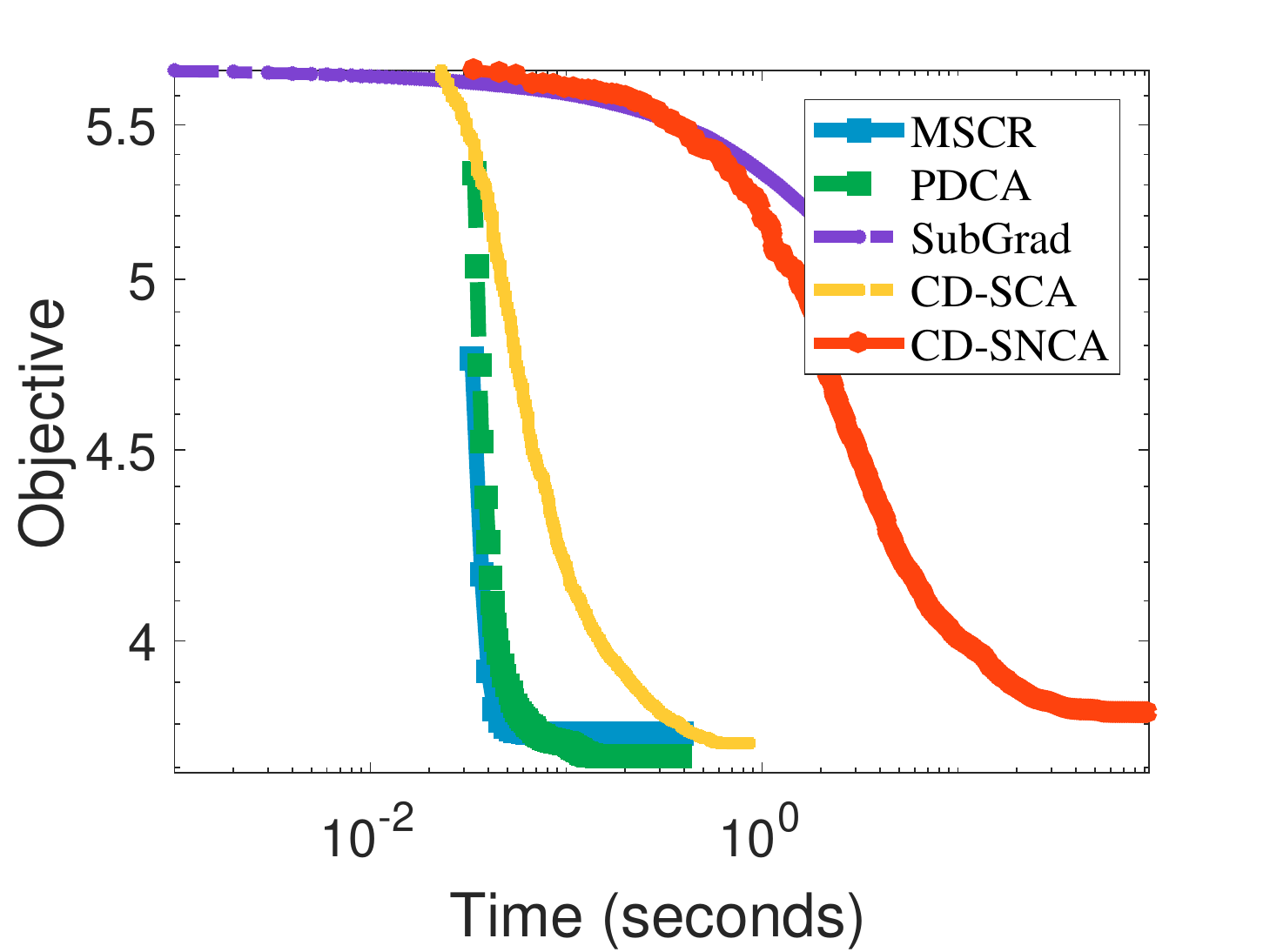}\vspace{-6pt} \caption{\scriptsize randn-2048-256  }\end{subfigure}\\

      \begin{subfigure}{0.25\textwidth}\includegraphics[height=\objimghei,width=\textwidth]{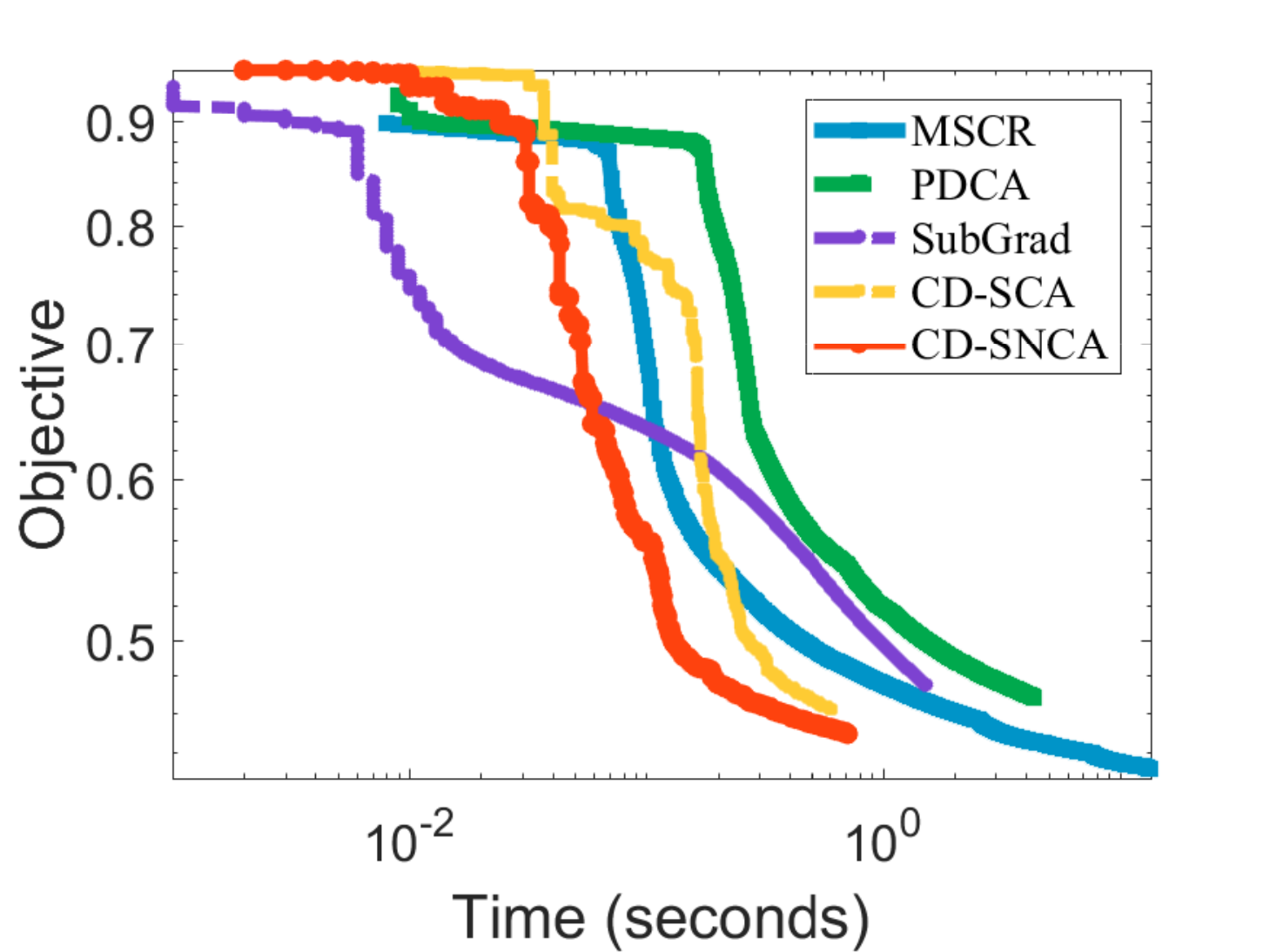}\vspace{-6pt} \caption{\scriptsize e2006-256-1024 }\end{subfigure}~~\begin{subfigure}{0.25\textwidth}\includegraphics[height=\objimghei,width=\textwidth]{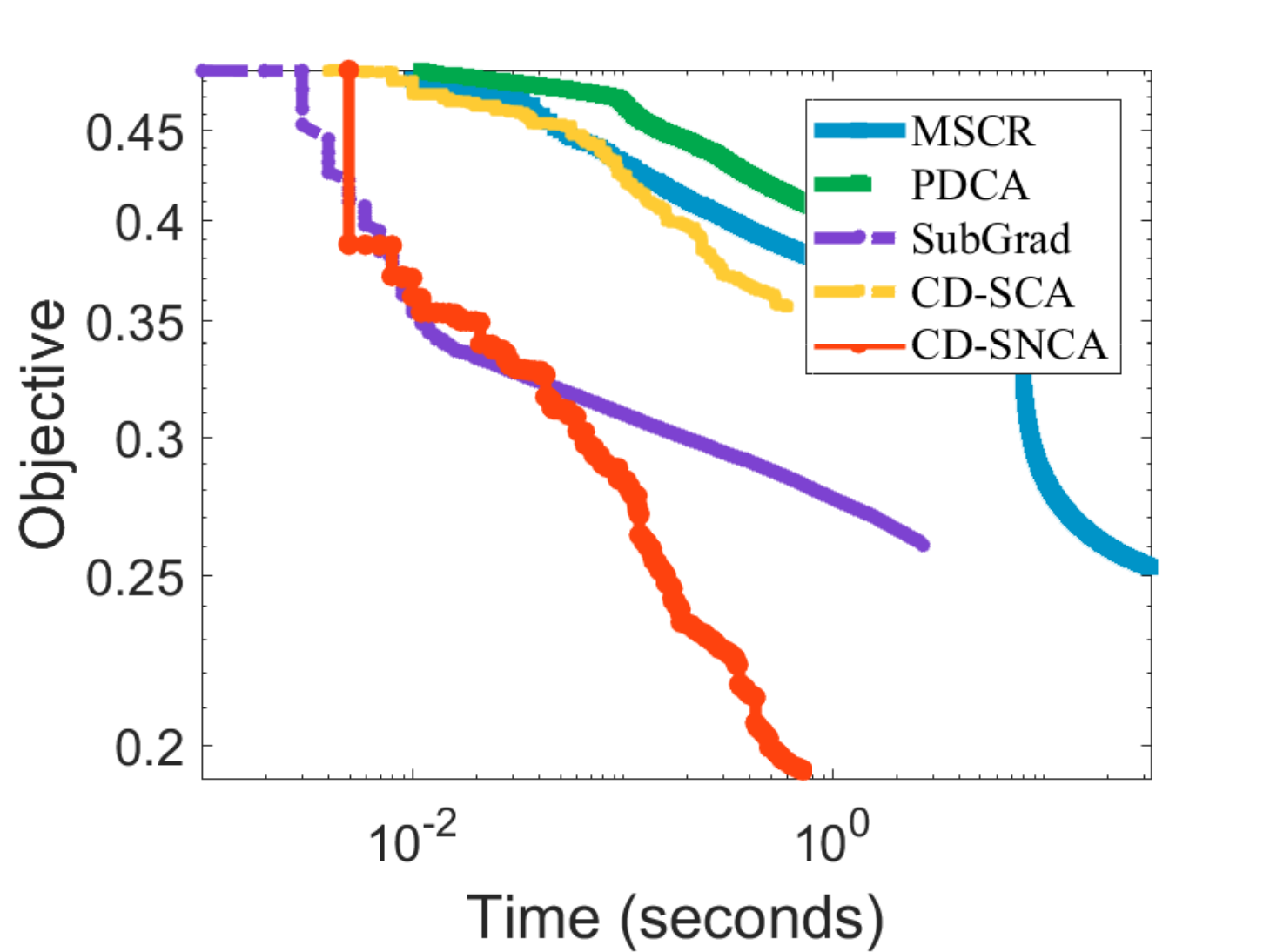}\vspace{-6pt} \caption{ \scriptsize e2006-256-2048 } \end{subfigure}~~\begin{subfigure}{0.25\textwidth}\includegraphics[height=\objimghei,width=\textwidth]{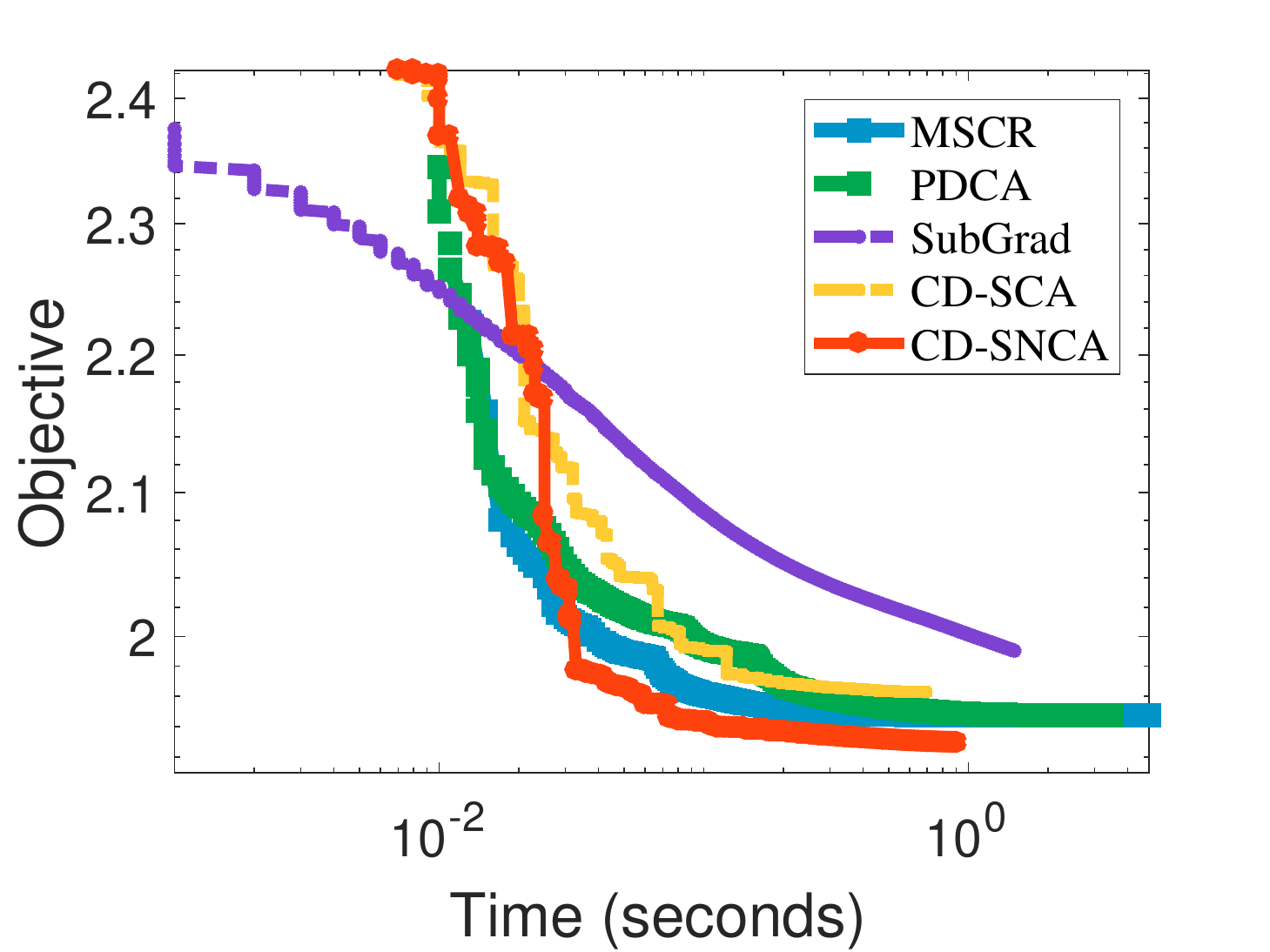}\vspace{-6pt} \caption{\scriptsize e2006-1024-256}\end{subfigure}~~\begin{subfigure}{0.25\textwidth}\includegraphics[height=\objimghei,width=\textwidth]{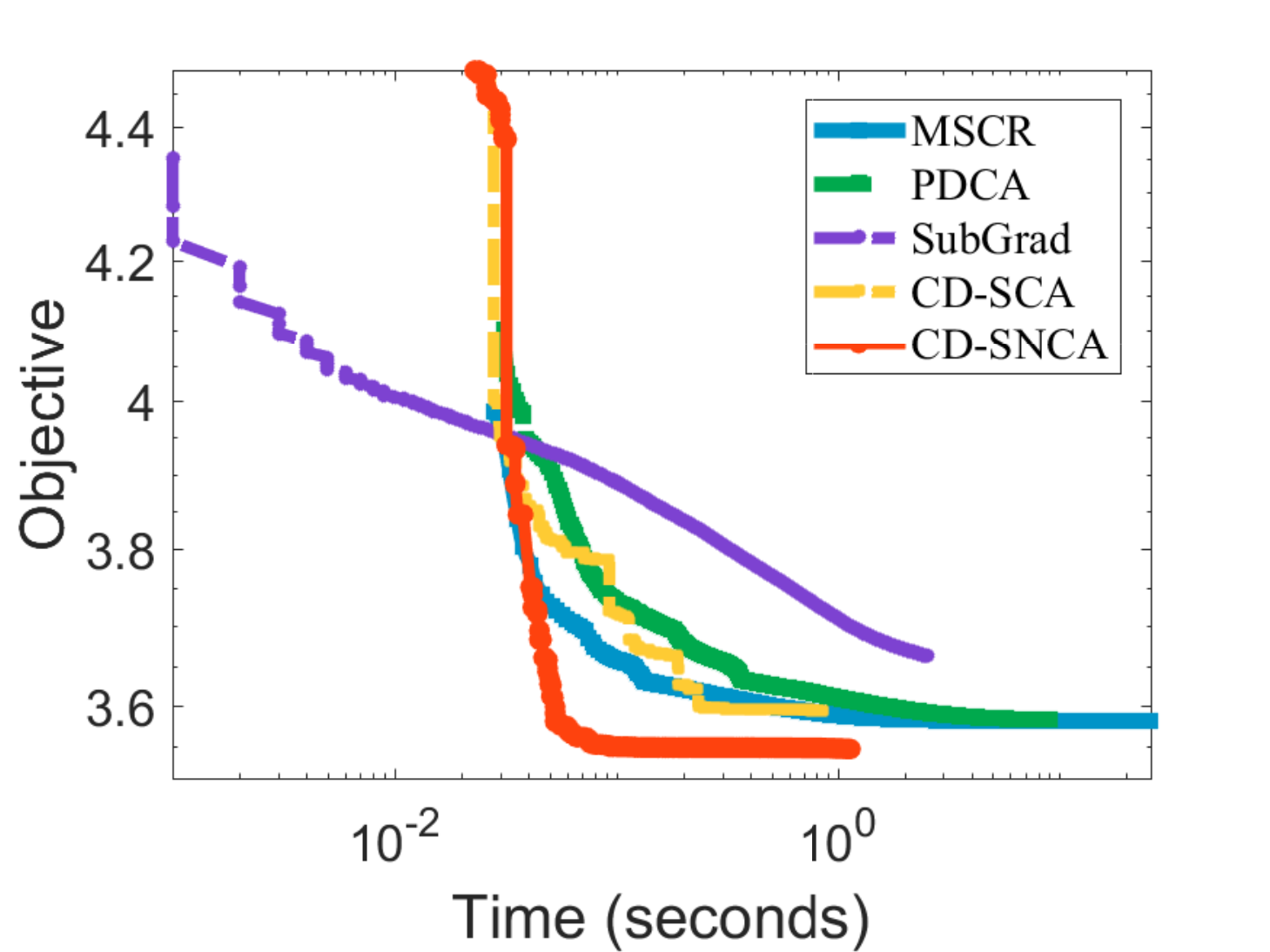}\vspace{-6pt} \caption{\scriptsize e2006-2048-256 }\end{subfigure}\\

      \begin{subfigure}{0.25\textwidth}\includegraphics[height=\objimghei,width=\textwidth]{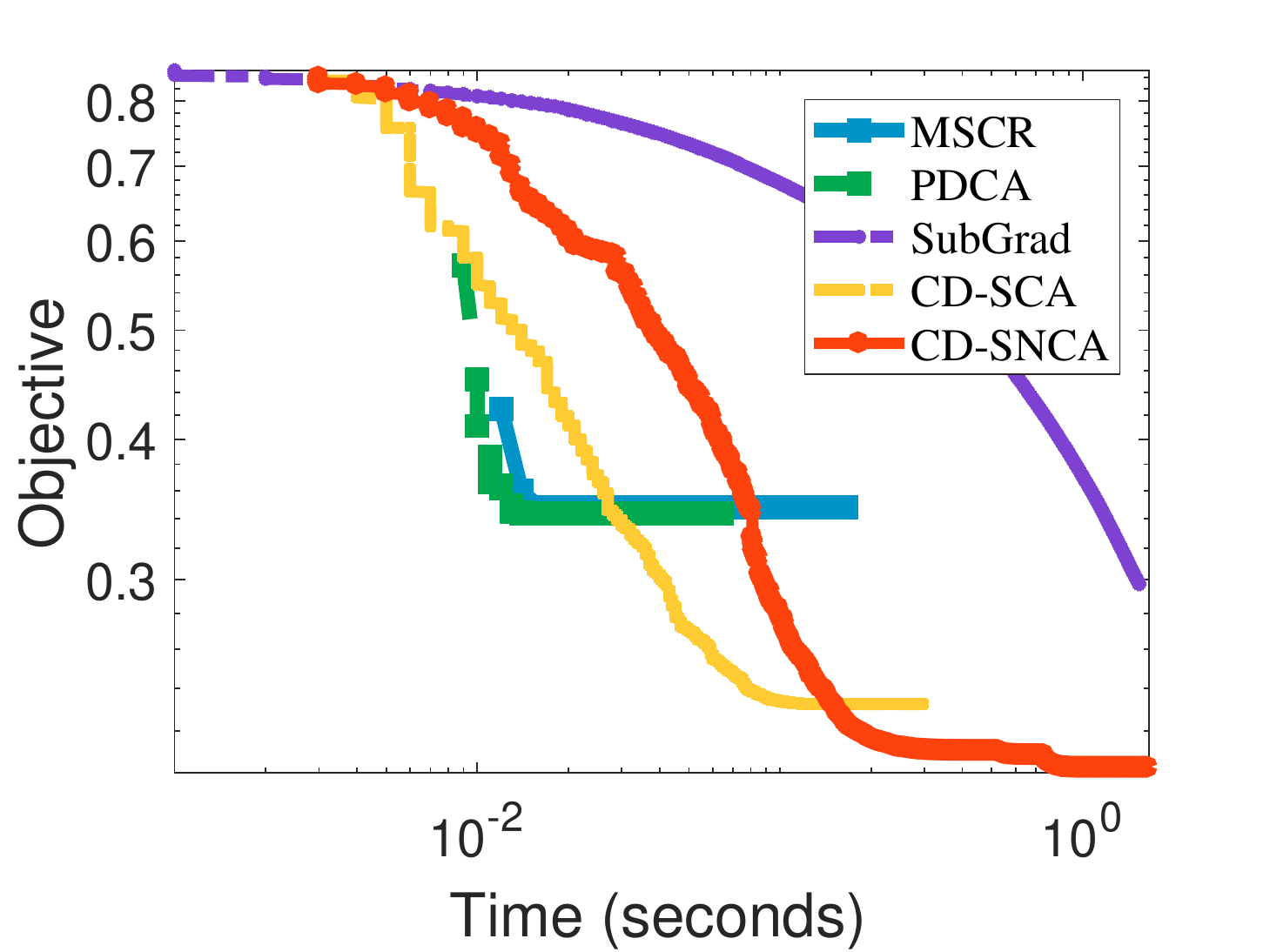}\vspace{-6pt} \caption{\scriptsize randn-256-1024-C }\end{subfigure}~~\begin{subfigure}{0.25\textwidth}\includegraphics[height=\objimghei,width=\textwidth]{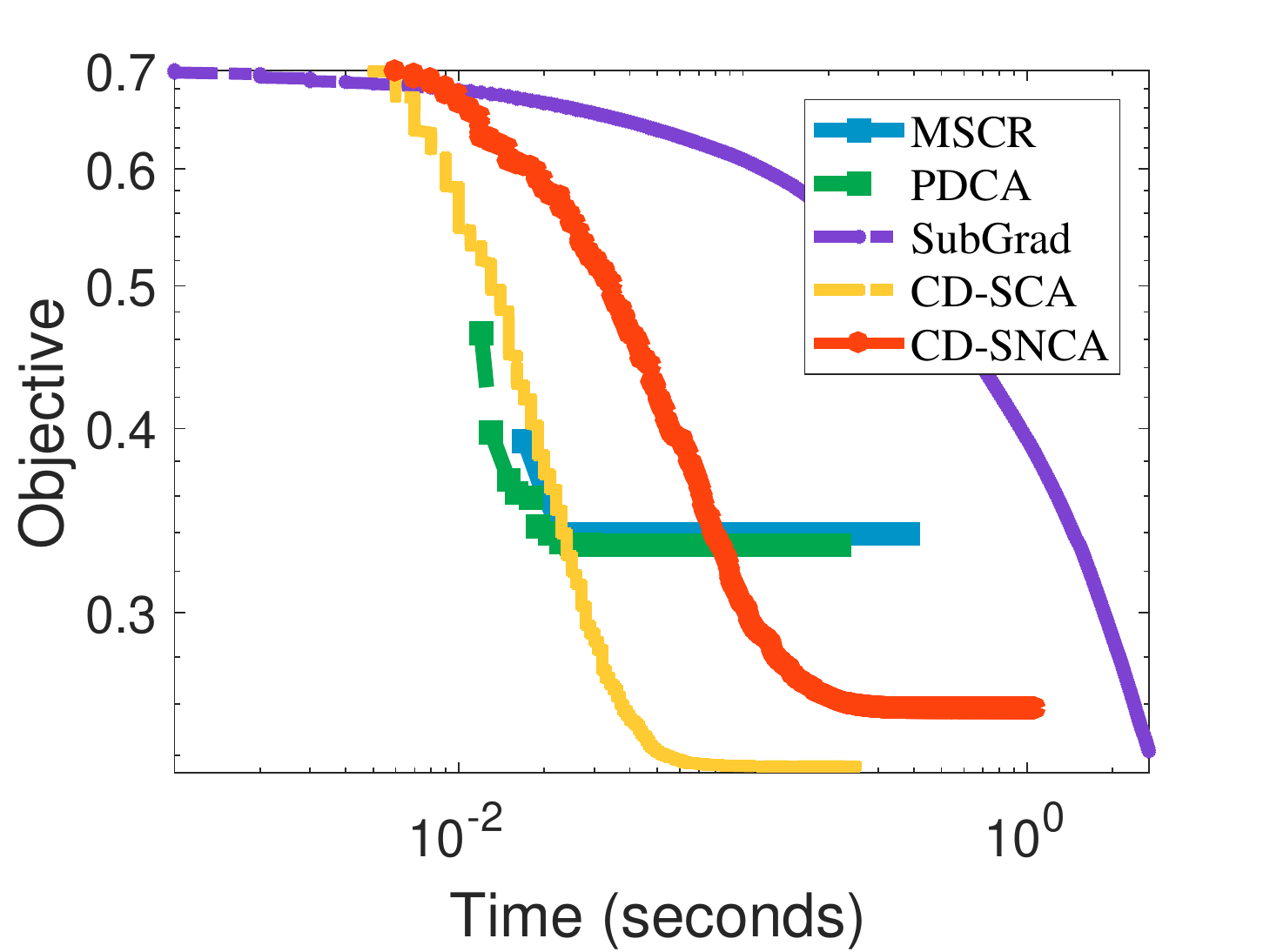}\vspace{-6pt} \caption{ \scriptsize randn-256-2048-C } \end{subfigure}~~\begin{subfigure}{0.25\textwidth}\includegraphics[height=\objimghei,width=\textwidth]{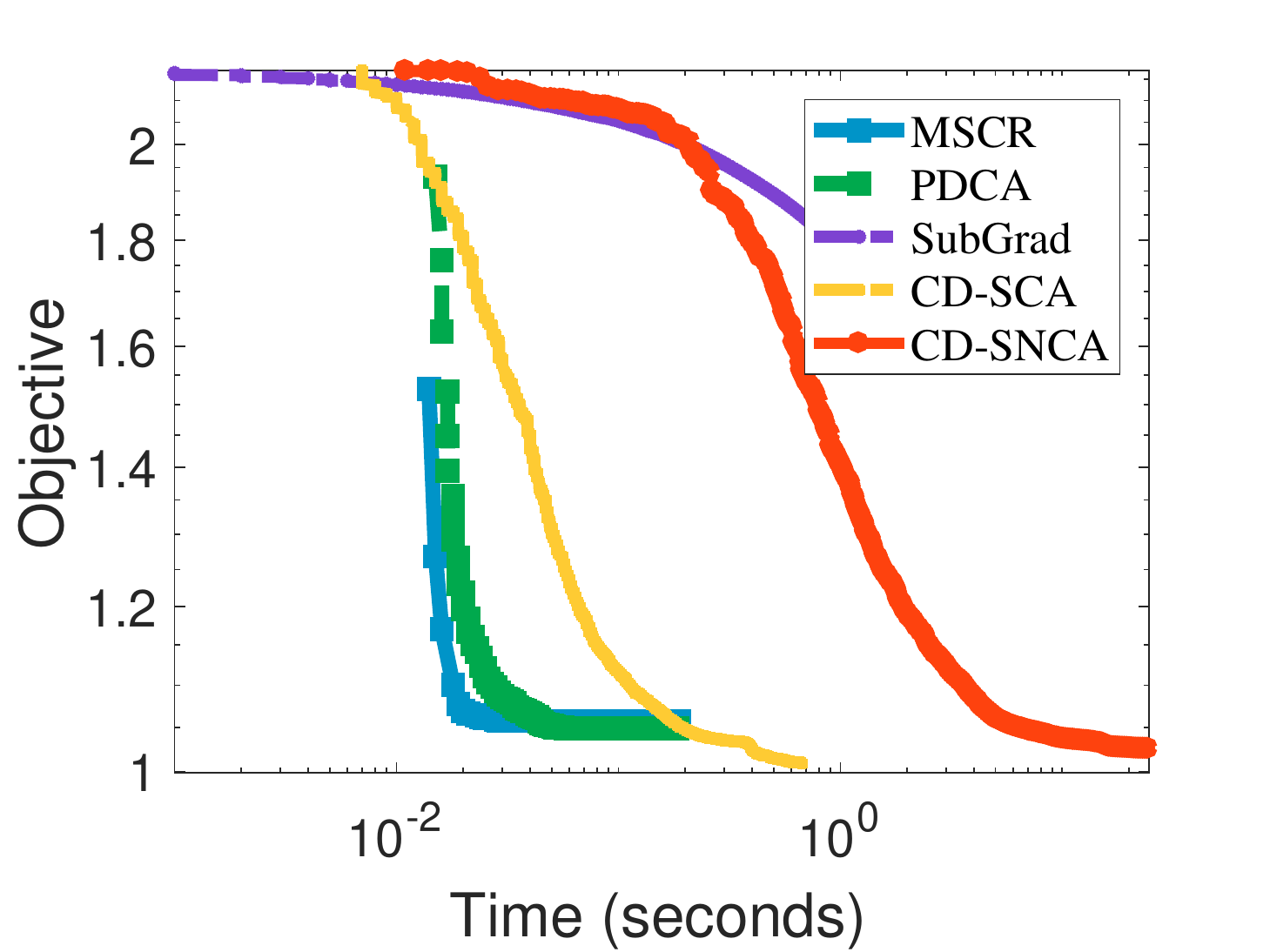}\vspace{-6pt} \caption{\scriptsize randn-1024-256-C }\end{subfigure}~~\begin{subfigure}{0.25\textwidth}\includegraphics[height=\objimghei,width=\textwidth]{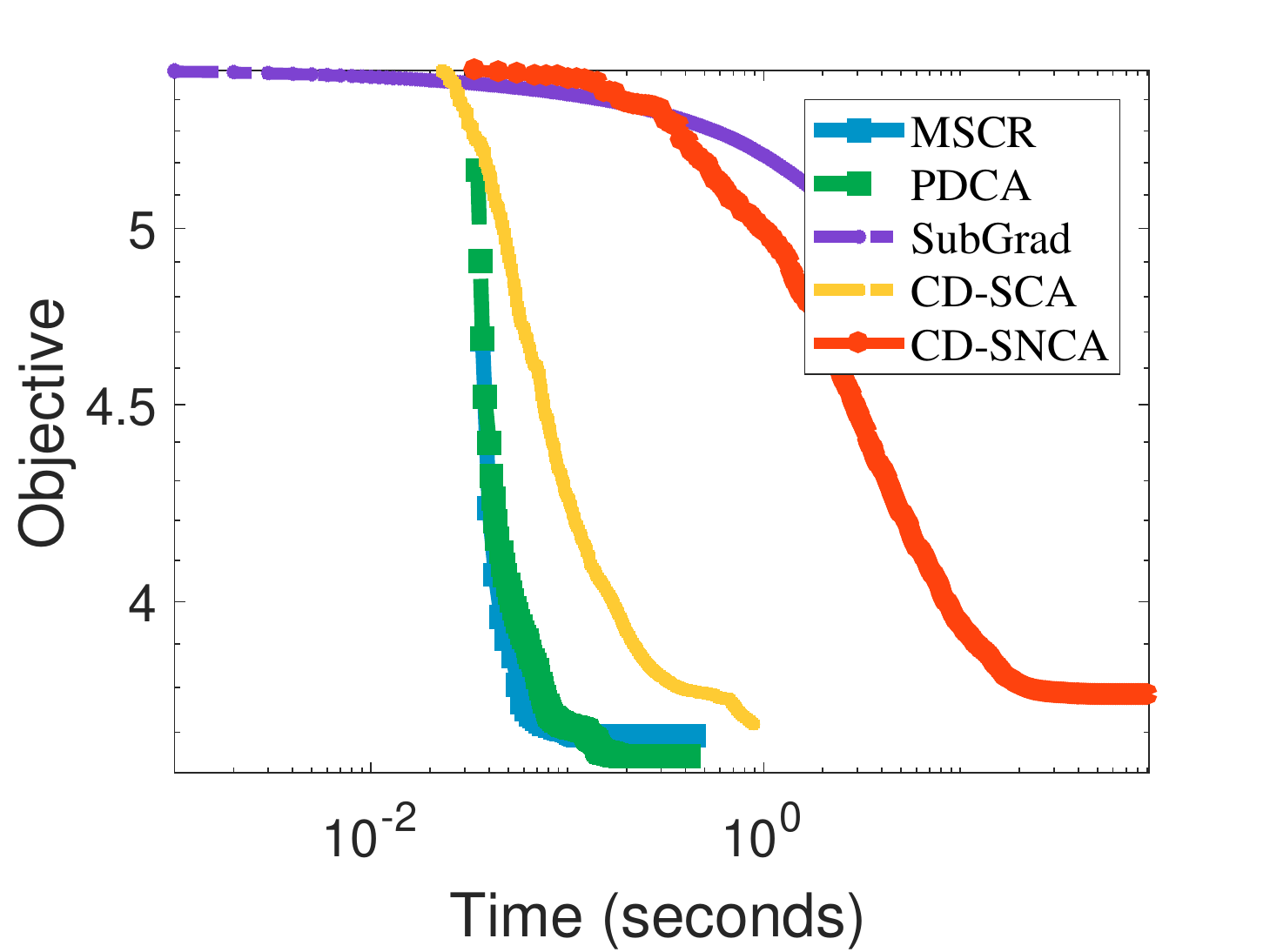}\vspace{-6pt} \caption{\scriptsize randn-2048-256-C  }\end{subfigure}\\

      \begin{subfigure}{0.25\textwidth}\includegraphics[height=\objimghei,width=\textwidth]{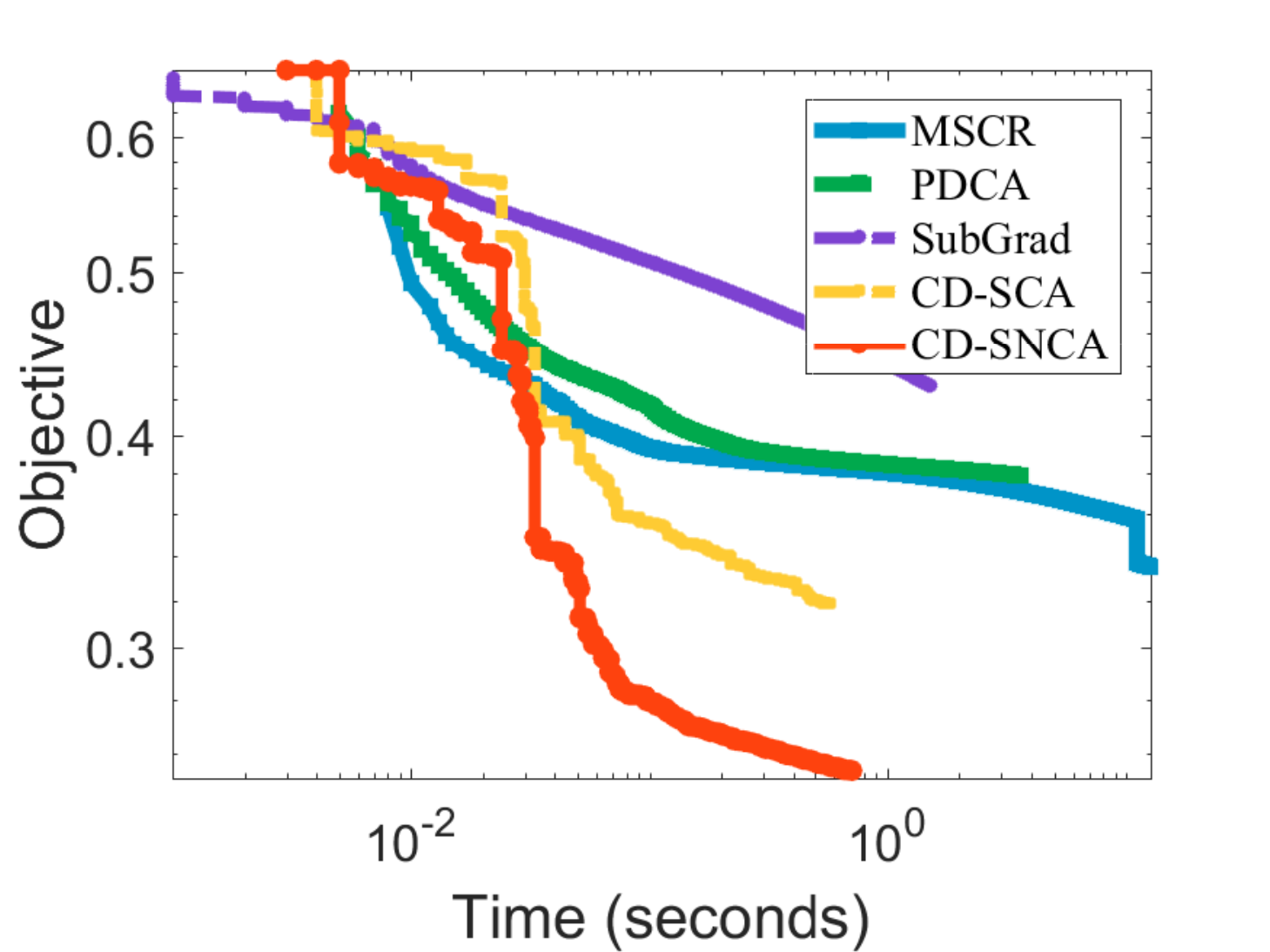}\vspace{-6pt} \caption{\scriptsize e2006-256-1024-C }\end{subfigure}~~\begin{subfigure}{0.25\textwidth}\includegraphics[height=\objimghei,width=\textwidth]{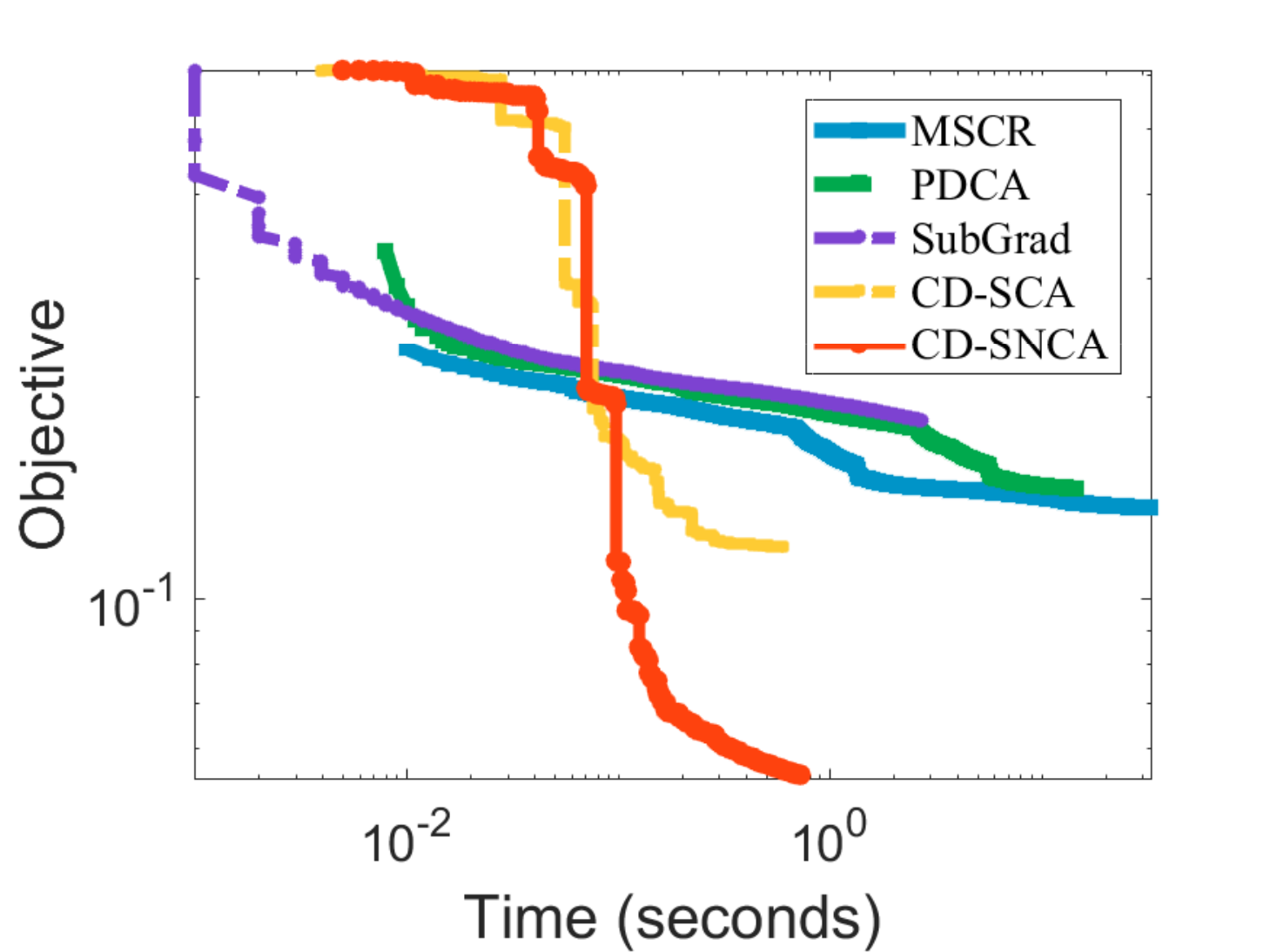}\vspace{-6pt} \caption{ \scriptsize e2006-256-2048-C } \end{subfigure}~~\begin{subfigure}{0.25\textwidth}\includegraphics[height=\objimghei,width=\textwidth]{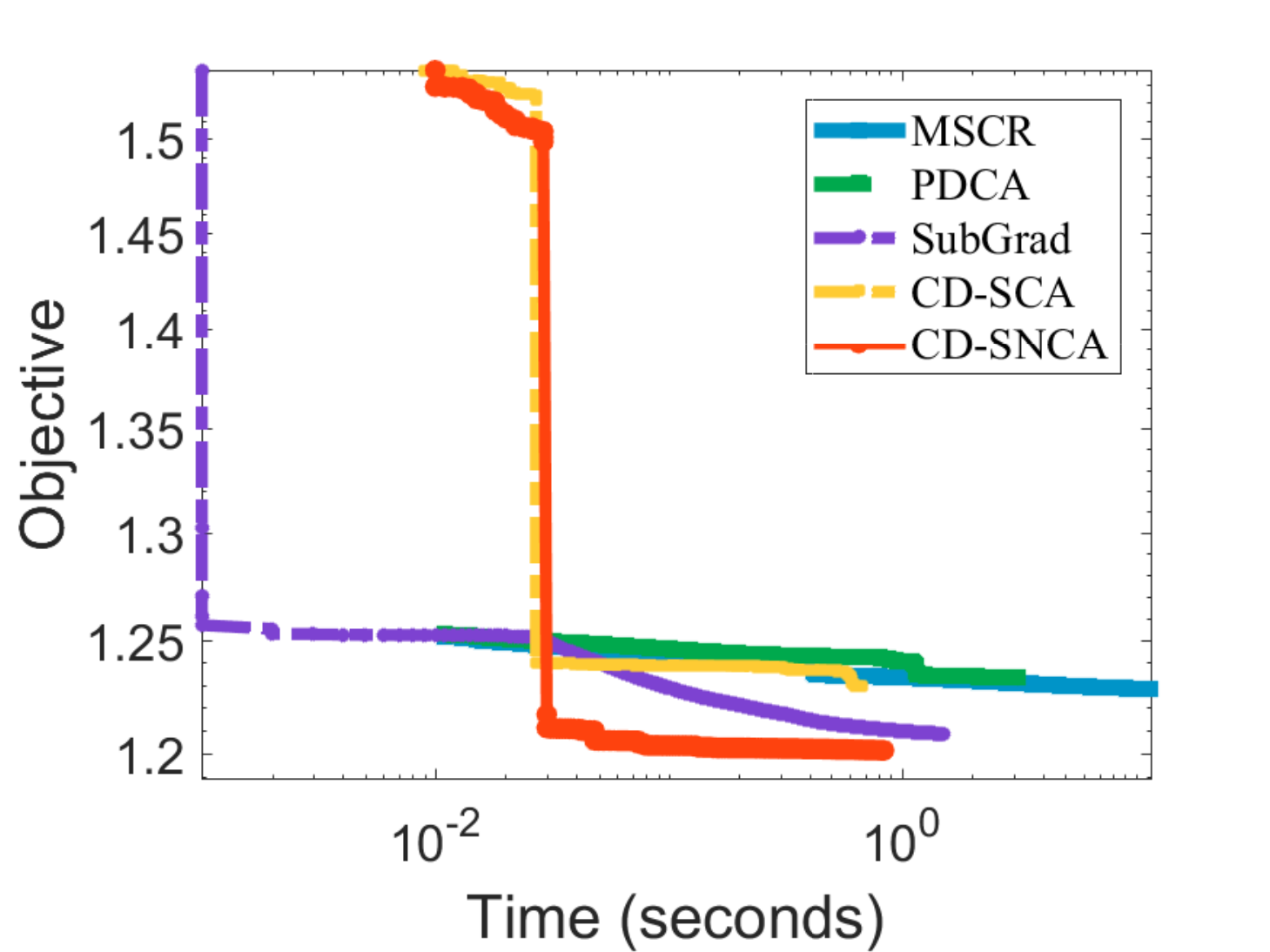}\vspace{-6pt} \caption{\scriptsize e2006-1024-256-C }\end{subfigure}~~\begin{subfigure}{0.25\textwidth}\includegraphics[height=\objimghei,width=\textwidth]{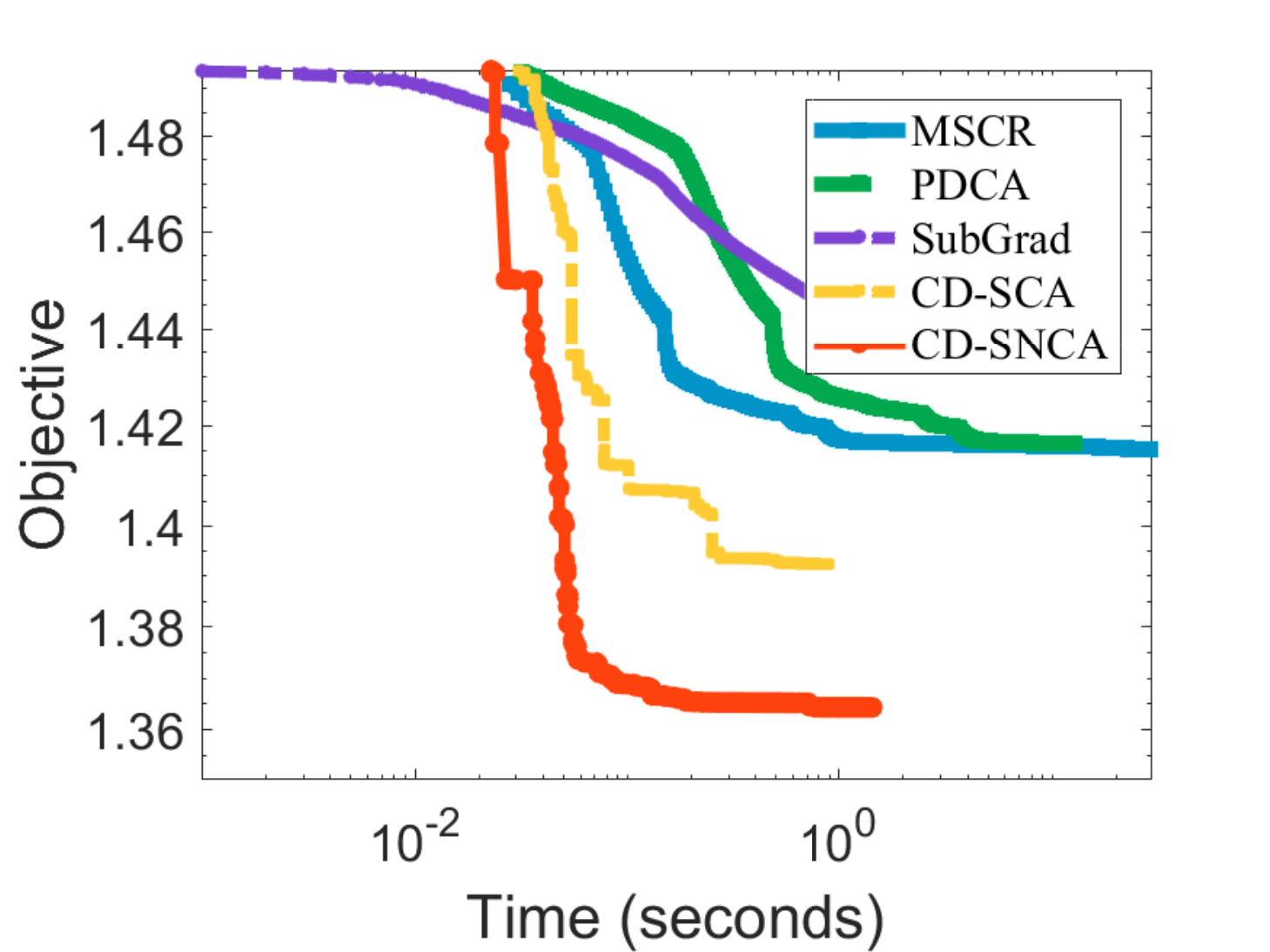}\vspace{-6pt} \caption{\scriptsize e2006-2048-256-C  }\end{subfigure}\\

\centering
\caption{The convergence curve of the compared methods for solving the generalized linear regression problem on different data sets.}
\label{exp:cpu:4}
\end{figure*}

Figure \ref{exp:cpu:12}, Figure \ref{exp:cpu:2}, Figure \ref{exp:cpu:3}, and Figure \ref{exp:cpu:4} show the convergence curve of the compared methods for solving the $\ell_p$ norm generalized eigenvalue problem, the approximate sparse optimization problem, the approximate binary optimization problem, and the generalized linear regression problem, respectively. We conclude that \textit{\textbf{CD-SNCA}} at least achieves comparable efficiency, if it is not faster than the compared methods. However, it generally achieves lower objective values than the other methods.


%
%
%
%

%
%
%

%
%

\section{Discussions}
\label{sect:disc:ext}

This section presents some discussions of our method. First, we discuss the equivalent reformulations for the $\ell_p$ norm generalized eigenvalue problem (see Section \ref{sect:disc:ref}). Second, we use several examples to explain the optimality hierarchy between the optimality conditions (see Section \ref{sect:opt:cond}). Third, we provide a local analysis of CD method for solving the classical PCA problem (see Section \ref{sect:pca}).

\subsection{Equivalent Reformulations for the $\ell_p$ Norm Generalized Eigenvalue Problem}
\label{sect:disc:ref}
First of all, using classical Lagigian dual theory, Problem (\ref{eq:main}) is equvilent to the following optimization problem.
\beq
\min_{\x\in\mathbb{R}^n}~  f(\x) + h(\x),~s.t.~ g(\x)\geq \lambda, \nn\\
\min_{\x\in\mathbb{R}^n}~ h(\x) - g(\x),~s.t.~ f(\x)\leq \delta, \nn
\eeq
\noi for some suitable $\lambda$ and $\delta$. For some special problems where the parameters $\lambda$ and $\delta$ that are not important, the two constrained problems above can be solved by finding the solution to Problem (\ref{eq:main}).

We pay special attention to the following problems with $\Q\succ \mathbf{0}$:
\beq
\min_{\x}~\F_1(\x) \triangleq  \tfrac{\alpha}{2}\x^T\Q\x - \|\A\x\|_p \label{eq:lp:1}\\
\min_{\x}~\F_2(\x) \triangleq - \|\A\x\|_p,~ s.t.~\x^T\Q\x \leq 1 \label{eq:lp:2}\\
\min_{\x}~\F_3(\x) \triangleq\tfrac{1}{2}\x^T\Q\x,~ s.t.~\|\A\x\|_p \geq 1 .\label{eq:lp:3}
\eeq

The following proposition establish the relations between Problem (\ref{eq:lp:1}), Problem (\ref{eq:lp:2}), and Problem (\ref{eq:lp:3}).

\begin{proposition}  \label{prop:equ}

We have the following results.

\noi \textbf{(a)} If $\bar{\x}$ is an optimal solution to $(\ref{eq:lp:1})$, then ${ \pm \bar{\x}}  ({{\bar{\x}^T\Q\bar{\x}}})^{-\frac{1}{2}}$ and $\frac{\pm \bar{\x}}{\|\A\bar{\x}\|_p}$ are respectively optimal solutions to $(\ref{eq:lp:2})$ and $(\ref{eq:lp:3})$.

\noi\textbf{(b)} If $\bar{\y}$ is an optimal solution to $(\ref{eq:lp:2})$, then $\frac{\pm \|\A\bar{\y}\|_p \cdot \bar{\y}}{\alpha \bar{\y}^T \Q\bar{\y}}$ and $\frac{\pm \bar{\y}}{\|\A\bar{\y}\|_p}$ are respectively optimal solutions to $(\ref{eq:lp:1})$ and $(\ref{eq:lp:3})$.

\noi\textbf{(c)} If $\bar{\z}$ is an optimal solution to $(\ref{eq:lp:3})$, then $ \frac{ \pm \bar{\z} \| \A\bar{\z}\|_p }{\alpha \bar{\z}^T \Q\bar{\z}}$ and ${\pm \bar{\z}} ({{\bar{\z}^T\Q\bar{\z}}})^{-\frac{1}{2}}$ are respectively optimal solutions to $(\ref{eq:lp:1})$ and $(\ref{eq:lp:2})$.

\begin{proof}
Using the Lagrangian dual, we introduce a multiplier $\theta_1>0$ for the constraint $\x^T\Q\x \leq 1$ as in Problem (\ref{eq:lp:2}) and a multiplier $\theta_2>0$ for the constraint $-\|\A\x\|_p \leq -1$ as in Problem (\ref{eq:lp:3}), leading to the following optimization problems:
\beq
&\min_{\x}~ \tfrac{\theta_1}{2} \x^T\Q\x - \|\A\x\|_p&,  \nn\\
&\min_{\x}~ \tfrac{1}{2}\x^T\Q\x  - \theta_2 \|\A\x\|_p \Leftrightarrow \min_{\x}~ \tfrac{1}{2 \theta_2}\x^T\Q\x  -  \|\A\x\|_p&. \nn
\eeq
\noi These two problems have the same form as Problem (\ref{eq:lp:1}). It is not hard to notice that the gradient of $\F_1(\cdot)$ can be computed as:
\beq
\nabla \F_1(\x) = \alpha \Q\x - \|\A\x\|_p^{1-p} \A^T ( \text{sign}(\A\x) \odot |\A\x|^{p-1}). \nn
\eeq
By the first-order optimality condition, we have:
\beq
\x = \frac{1}{\alpha}\Q^{-1}\left(\|\A\x\|_p^{1-p} \A^T ( \text{sign}(\A\x) \odot |\A\x|^{p-1})\right) .\nn
\eeq

\noi Therefore, the optimal solution for Problem (\ref{eq:lp:1}), Problem (\ref{eq:lp:2}), and Problem (\ref{eq:lp:3}) only differ by a scale factor.

\noi \textbf{(a)} Since $\bar{\x}$ is the optimal solution to $(\ref{eq:lp:1})$, the optimal solution to Problem $(\ref{eq:lp:2})$ and $(\ref{eq:lp:3})$ can be respectively computed as $\alpha_2 \bar{\x}$ and $\alpha_3 \bar{\x}$ with
\beq
&&\alpha_2 = \arg \min_{\alpha}~F_2(\alpha \bar{\x}),~s.t.~(\alpha \bar{\x})^T\Q(\alpha \bar{\x}) \leq 1 \nn\\
&&\alpha_3 = \arg \min_{\alpha}~F_3(\alpha \bar{\x}),~s.t.~\|\alpha \bar{\x}\|_p \geq 1.\nn
\eeq
\noi After some preliminary calculations, we have: $\alpha_2 = \pm 1 / \sqrt{\bar{\x}^T\Q\bar{\x}}$ and $\alpha_3 = \pm 1 / \|\A\bar{\x}\|_p$.

\noi \textbf{(b)} Since $\bar{\y}$ is an optimal solution to $(\ref{eq:lp:2})$, the optimal solution to Problem $(\ref{eq:lp:1})$ and Problem $(\ref{eq:lp:3})$ can be respectively computed as $\alpha_1 \bar{\y}$ and $\alpha_3 \bar{\y}$ with
\beq
&&\alpha_1 = \arg \min_{\alpha}~F_1(\alpha \bar{\y}) \nn\\
&&\alpha_3 = \arg \min_{\alpha}~F_3(\alpha \bar{\y}),~s.t.~\|\alpha \bar{\y}\|_p \geq 1.\nn
\eeq
\noi Therefore, $\alpha_1 = \pm \frac{ \| \A\bar{\y}\|_p }{\alpha \bar{\y}^T \Q\bar{\y}} $ and $\alpha_3 = \pm 1 / \|\A\bar{\y}\|_p$.

\noi \textbf{(c)} Since $\bar{\z}$ is an optimal solution to $(\ref{eq:lp:3})$, the optimal solution to Problem $(\ref{eq:lp:1})$ and Problem $(\ref{eq:lp:2})$ can be respectively computed as $\alpha_1 \bar{\z}$ and $\alpha_2 \bar{\z}$ with
\beq
&&\alpha_1 = \arg \min_{\alpha}~F_1(\alpha \bar{\z}) \nn\\
&&\alpha_2 = \arg \min_{\alpha}~F_2(\alpha \bar{\z}),~s.t.~(\alpha \bar{\z})^T\Q(\alpha \bar{\z})  \leq 1. \nn
\eeq
\noi Therefore, $\alpha_1 = \pm \frac{ \| \A\bar{\z}\|_p }{\alpha \bar{\z}^T \Q\bar{\z}} $ and $\alpha_2 = \pm 1 / \sqrt{\bar{\z}^T\Q\bar{\z}}$.

\end{proof}

\end{proposition}

\begin{table}[!h]
\begin{center}
\scalebox{0.82}{\begin{tabular}{|c|c|c|c|c|c|c|c|}
\hline
$\y$ & $\x$  & Function Value  & Critical Point  & CWS Point  \\
\hline
$[1;1;1]$      &  $[  1.75;0;-1]$      &   -6.625    & \cone{Yes} &    No           \\
$[1;1;[-1,1]]$      &   NA    &  NA   & No &       No     \\
$[1;1;-1]$ &  $[-0.25;   -2;   -1]$    &  -8.125     & No &    No          \\
$[1;[-1,1];1]$      &    NA   &   NA    & No &    No           \\
$[1;[-1,1];[-1,1]]$      &   NA    &  NA   & No &       No     \\
$[1;[-1,1];-1]$ &     NA  &  NA    & No &    No          \\
$[1;-1;1]$ &   $ [0.25;-2;-3] $  &   -4.1250       & No &    No          \\
$[1;-1;[-1,1]]$ &   $[-0.3333;0.2667;-0.1333]$         &  -1.9956       & No &    No          \\
$[1;-1;-1]$ &   $[-1.75;-4;-3]$    & -16.1250      & No &    No          \\
$[[-1,1];1;1]$      &  NA      &  NA    & No &    No           \\
$[[-1,1];1;[-1,1]]$      &  NA     & NA   & No &       No     \\
$[[-1,1]1;-1]$ &    $[0;-2;-2]$  &   -6.0000     & No &    No          \\
$[[-1,1];[-1,1];1]$      &   NA    &   NA   & No &    No           \\
$[[-1,1];[-1,1];[-1,1]]$      & $[0;0;0]$      &  0  & \cone{Yes} &       No     \\
$[[-1,1];[-1,1];-1]$ &     $[0;0;0]$  &  0    & \cone{Yes} &    No          \\
$[[-1,1];-1;1]$ &    NA   &  NA     & No &    No          \\
$[[-1,1];-1;[-1,1]]$ &    $[0;0;0]$   &  0     & \cone{Yes} &    No          \\
$[[-1,1];-1;-1]$ &    $[0;0;0]$   &  0     & \cone{Yes} &    No          \\
$[-1;1;1]$      &  $[1.25;0;-3]$      &   -7.6250    & \cone{Yes} &    No           \\
$[-1;1;[-1,1]]$      &   NA    &  NA   & No &       No     \\
$[-1;1;-1]$ & $[-0.75;-2;-3]$     &  -12.1250     & No &    No          \\
$[-1;[-1,1];1]$      &   NA    &   NA   & No &    No           \\
$[-1;[-1,1];[-1,1]]$      &  $[0;0;0]$     & 0   & \cone{Yes} &       No     \\
$[-1;[-1,1];-1]$ &    $[0;0;0]$   &  0    & \cone{Yes} &    No          \\
$[-1;-1;1]$ &     $[-0.25;-2;-5]$  &   -6.6250     & No &    No          \\
$[-1;-1;[-1,1]]$ &  $[0;0;0]$     &  0    &\cone{Yes} &    No          \\
$[-1;-1;-1]$ &   $[-2.25;-4;-5]$    &  -18.625     & \cone{Yes} &    \cone{Yes}          \\
\hline
\end{tabular}}
\caption{Solutions satisfying optimality conditions for Problem (\ref{running:1}).  } 
\label{tab:example:ranning:1}
\end{center}
\end{table}
\subsection{Examples for Optimality Hierarchy between the Optimality Conditions} \label{sect:opt:cond}

We show some examples to explain the optimality hierarchy between the optimality conditions. Since the condition of directional point is difficult to verify, we only focus on the condition of critical point and coordinate-wise stationary point in the sequel.

\noi $\bullet$ \textbf{The First Running Example}. We consider the following problem:
\beq
\min_{\x}~ \frac{1}{2}\x^T \Q\x + \la  \x,\p \ra -  \|\A\x\|_1 \label{running:1}
\eeq
\noi with using the following parameters:
\beq
\Q = \left(
            \begin{array}{ccc}
              4 & 0 & 0 \\
              0 & 2 & -1 \\
              0 & -1 & 1\\
            \end{array}
          \right),~\p = \left(
            \begin{array}{c}
              1  \\
              1 \\
              1\\
            \end{array}
          \right),~\A = \left(
            \begin{array}{ccc}
              1 & -1 & 1 \\
              3 & 1 & 0 \\
              4 & 2 & -1\\
            \end{array}
          \right).\nn
\eeq
\noi First, using the Legendre-Fenchel transform, we can rewrite Problem (\ref{running:1}) as the following optimization probelm:
\beq
\min_{\x,~\y}~  \frac{1}{2}\x^T \Q\x + \la  \x,\p \ra -  \la \A\x,~\y \ra,~-1 \leq \y \leq \textbf{1}. \nn
\eeq

\noi Second, we have the following first-order optimality condition for Problem (\ref{running:1}):
\beq \label{running:1:first:order:optimality}
\begin{split}
(\Q\x + \p -  \text{sign}(\A\x))_J = \mathbf{0},~J\triangleq \{j~|~(\A\x)_j\neq0\} \\
-\mathbf{1} \leq (\Q\x + \p)_I \leq  \mathbf{1},~I \triangleq \{i~|~(\A\x)_i=0\}.
\end{split}
\eeq
\noi

\noi Third, we notice the following relations between $\A\x$ and $\y$:
\beq
(\A\x)_i>0 ~ &\Rightarrow& ~\y_i=1 \nn\\
(\A\x)_i<0 ~ &\Rightarrow& ~\y_i=-1 \nn\\
(\A\x)_i=0 ~ &\Rightarrow& ~\y_i \in [-1,1]. \nn
\eeq

\noi We enumerate all possible solutions for $\y$ (as shown in the first column of Table \ref{tab:example:ranning:1}), and then compute the solution satisfying the first-order optimality condition using (\ref{running:1:first:order:optimality}). Table  \ref{tab:example:ranning:1} shows the solutions satisfying optimality conditions for Problem (\ref{running:1}). The condition of the Coordinate-wise Stationary (CWS) point might be a much stronger condition than the condition of critical point.


\begin{table}[!h]
\begin{center}
\scalebox{0.82}{\begin{tabular}{|c|c|c|c|c|c|c|c|}
\hline
$(\lambda_i,~\u_i)$ & $\x$  & Function Value  & Critical Point  & CWS Point  \\
\hline
$( 0.5468,~ [-0.2934 ,   0.8139 ,   0.5015])$      & $\pm[ -0.2169  ,  0.6019    ,0.3709]$      &   -5.7418    & \cone{Yes} &    No           \\
$( 7.8324 ,~ [0.1733 ,  -0.4707 ,   0.8651])$      & $\pm[0.4850 ,  -1.3172   , 2.4212]$      &   -82.2404   & \cone{Yes} &    No           \\
$( 33.6207,~ [-0.9402  , -0.3407 ,   0.0030])$      & $\pm[-5.4514  ,-1.9755  ,  0.0172]$      &   -353.0178    & \cone{Yes} &    \cone{Yes}           \\
      & $[0  , 0  ,  0]$      &   0   & \cone{Yes} &    No           \\
\hline
\end{tabular}}
\caption{Solutions satisfying optimality conditions for Problem (\ref{running:2}).}
\label{tab:example:ranning:2}
\end{center}
\end{table}
\noi $\bullet$ \textbf{The Second Running Example}. We consider the following example:
\beq
\min_{\x}~ \frac{1}{2}\x^T\x -  \|\A\x\|_2 \label{running:2}
\eeq
\noi with using the following parameter:
\beq
\A = \left(
            \begin{array}{ccc}
              1 & -1 & 1 \\
              2 & 0 & 2 \\
              3 & 1 & 0 \\
              4 & 2 & -1 \\
            \end{array}
          \right).\nn
\eeq
\noi Using the first-order optimality condition, one can show that the basic stationary points are $\{\mathbf{0}\}  \cup \{\pm \sqrt{\lambda_i} \u_i\}$, where $(\lambda_i, \u_i)$ are the eigenvalue pairs of the matrix $\A^T\A$. Table \ref{tab:example:ranning:2} shows the solutions satisfying optimality conditions for Problem (\ref{running:2}). There exists two coordinate-wise stationary points. Therefore, coordinate-wise-stationary might be a much stronger condition than criticality.


\begin{table}[!h]
\begin{center}
\scalebox{0.82}{\begin{tabular}{|c|c|c|c|c|c|c|c|}
\hline
 $\y$  & $\x$  & Function Value  & Critical Point  & CWS Point  \\
\hline
$[1;  0 ;0 ;0]$       & $[1 ;   -1   ;  1]$      &    -2.5000    & \cone{Yes} &    No           \\
$[0; 1 ;0; 0]$       & $[2 ;    0  ;   2]$      &   -4.0000    & \cone{Yes} &    No           \\
$[0; 0 ;1 ;0]$       & $[3  ;   1 ;    0]$      &    -9.0000     & \cone{Yes} &    No           \\
$[0; 0 ;0;1]$       & $[4   ;  2    ;-1]$      &   -10.5000    & \cone{Yes} &    \cone{Yes}           \\
$[-1; 0 ;0;0]$      & $[-1   ;  1  ;  -1]$      &   -2.5000    & \cone{Yes} &    No           \\
$[0 ;-1 ;0; 0]$      & $[-2    ; 0 ;   -2]$      &  -4.0000    & \cone{Yes} &    No           \\
$[0 ;0 ;-1 ;0]$      & $[-3  ;  -1;     0]$      &   -9.0000    & \cone{Yes} &    No           \\
$[0 ;0; 0;-1]$      & $[-4   ; -2;     1]$      &  -10.5000    & \cone{Yes} &    \cone{Yes}           \\
\hline
\end{tabular}}
\caption{Solutions satisfying optimality conditions for Problem (\ref{running:3}).}

\label{tab:example:ranning:3}
\end{center}
\end{table}
\noi $\bullet$ \textbf{The Third Running Example}. We consider the following example:
\beq
\min_{\x}~ \frac{1}{2}\x^T\x -  \|\A\x\|_{\infty}  \label{running:3}
\eeq
\noi with using the same value for $\A$ as in Problem (\ref{running:2}). It is equivalent to the following problem:
\beq
\min_{\x,\y}~ \frac{1}{2}\x^T\x -  \la \A\x,~\y \ra,~\|\y\|_{1} \leq 1. \nn
\eeq


\noi We enumerate some possible solutions for $\y$, and then compute the solution satisfying the first-order optimality condition via the optimality of $\x$ that: $\x = \A^T\y$. Table \ref{tab:example:ranning:3} shows the solutions satisfying optimality conditions for Problem (\ref{running:3}). These results consolidate our previous conclusions.


%
%

\subsection{A Local Analysis of CD method for the PCA Problem }\label{sect:pca}

This section provides a local analysis of Algorithm \ref{algo:main} when it is applied to solve the classical PCA problem. We first rewrite the classical PCA problem as an unconstraint smooth optimization problem and then prove that it is smooth and strongly convex in the neighborhood of the global optimal solution. Finally, the local linear convergence rate of the CD method directly follows due to Theorem 1 in \cite{lu2015complexity}.

Given a covariance matrix $\C \in \mathbb{R}^{n\times n}$ with $\C\succcurlyeq \mathbf{0}$, PCA problem can be formulated as:
\beq
\max_{\v}~\v^T\C\v,~s.t.~\|\v\|=1 \nn
\eeq
\noi Using Proposition \ref{prop:equ}, we have the following equivalent problem:
\beq
 \min_{\x}~ \F(\x) = \frac{\alpha}{2}\|\x\|_2^2-\sqrt{\x^T\C\x}. \label{eq:pca111}
\eeq
\noi for any given constant $\alpha>0$. In what follows, consider a positive semidefinite matrix $\C \in \mathbb{R}^{n\times n}$ (which is not necessarily low-rank) with eigenvalue decomposition $\C= \sum_{i=1}^n \lam_i \u_i \u_i^T=\U^T\text{diag}(\lam)\textbf{U}$. We assume that: $\lam_1 \geq \lam_2 \geq ... \geq \lam_n \geq 0$.

\begin{lemma} \label{lemma:vip:PCA}
The following result holds iff $i=1$:
\beq
\O \triangleq  \I - \tfrac{1}{\lam_i}\C+  \u_i \u_i^T \succ \mathbf{0} \label{eq:vip:PCA}
\eeq

\begin{proof}
When $i=1$, it clearly holds that: $\lam_i \I - \C \succeq \mathbf{0}$. Combining with the fact that $\lam_1 \u_1 \u_1^T \succ \mathbf{0}$, we have $\lam_i \I - \C + \lam_i \u_i \u_i^T \succ \mathbf{0}$ for $i=1$.

We now prove that the inequality in (\ref{eq:vip:PCA}) may fail to hold for $i\geq 2$. We have the following results:
\beq
\u_1^T\O \u_1 &=& \lam_i \|\u_1\|_2^2 - \u_1^T\C\u_1 + \lam_i (\u_1^T\u_i)^2 \nn\\
&\overset{(a)}{=}& \lam_i  - \u_1^T\C\u_1 + 0 \nn\\
&\overset{(b)}{=}& \lam_i  - \lam_1 + 0 ~\overset{(c)}{\leq}~0 \nn
\eeq
\noi step $(a)$ uses the fact that $\|\u_i\|=1$ for all $i$, and $\u_i^T\u_j=0$ for all $i\neq j$; step $(b)$ uses the fact that $\C\u_1 = \lam_1 \u_1~\Rightarrow~\u_1^T\C\u_1 = \lam_1$; step $(c)$ uses the fact that $\lam_1 \geq \lam_2 \geq ... \geq \lam_n \geq 0$. Therefore, the matrix $\textbf{O}$ is not positive definite.

We hereby finish the proof of this lemma.

\end{proof}

\end{lemma}

\begin{theorem} \label{pca:critical:points}
We have the following results:

\textbf{(a)} The set of critical points of Problem (\ref{eq:pca111}) are $\{\{\mathbf{0}\}  \cup \{\pm \frac{\sqrt{\lam_k} }{\alpha}\u_k:k=1,...,n\}\}$.

\textbf{(b)} Problem (\ref{eq:pca111}) has at most two local minima $\{\pm \frac{\sqrt{\lam_1} }{\alpha} \u_1\}$ which are the global optima with $\F(\bar{\x}) = - \frac{\lam_1}{2\alpha}$.


\begin{proof}
The first-order and second-order gradient of $\F(\x)$ can be computed as:
\beq
&\nabla \F(\x) = \alpha \x - \frac{\C\x}{\sqrt{\x^T\C \x}} \label{eq:pca:grad} \\
& \nabla^2 \F(\x) \overset{}{=} \frac{ \alpha \sqrt{\x^T\C\x} \cdot \textbf{I} - \C + \alpha^2 \x\x^T }{\sqrt{\x^T\C\x}}  \label{eq:pca:hess}
\eeq

\noi (\textbf{a}) Setting $\nabla \F(\x)=0$ for (\ref{eq:pca:grad}), we have:
\beq
\alpha  \sqrt{\x^T\C \x} \cdot \x =  \C\x \label{eq:first:order:opt}
\eeq
Therefore, we conclude that $\{      \frac{\sqrt{\lam_k} }{\alpha} \u_k,~k=1,...,n\}$ are feasible solutions to (\ref{eq:first:order:opt}). Taking into account that the objective function is symmetric and the trivial solution $\{\mathbf{0}\}$, we finish the proof of the first part of this lemma.

\noi (\textbf{b}) For any nontrivial (nonzero) critical point $\z_i = \pm \frac{\sqrt{\lam_i} }{\alpha} \u_i$, we have the following results from (\ref{eq:pca:hess}):
\beq
&&\nabla^2 \F(\z_i) \nn\\
&=~& \alpha \cdot  \frac{  \sqrt{(\sqrt{\lam_i} \u_i)^T\C(\sqrt{\lam_i} \u_i)} \cdot \textbf{I} - \C +  \lam_i \u_i \u_i^T }{\sqrt{(\sqrt{\lam_i} \u_i)^T\C(\sqrt{\lam_i} \u_i)}} \nn\\
&=~& \alpha \cdot  \frac{ \sqrt{\lam_i} \sqrt{\u_i^T\C\u_i} \cdot \textbf{I} - \C +  \lam_i \u_i \u_i^T }{ \sqrt{\lam_i} \cdot \sqrt{\u_i^T\C\u_i}} \nn\\
&\overset{(a)}{=}~& \alpha \cdot  \frac{ \lam_i \cdot \textbf{I} - \C +  \lam_i \u_i \u_i^T }{ \lam_i } = \alpha \cdot (\textbf{I} - \tfrac{1}{\lam_i}\C +   \u_i \u_i^T) \nn
\eeq
\noi where step $(a)$ uses the fact that $\C\u_i = \lam_i \u_i~\Rightarrow~\u_i^T\C\u_i = \lam_i$. Using Lemma \ref{lemma:vip:PCA}, we conclude that $\nabla^2 F(\z_i) \succ \textbf{0}$ holds only when $i=1$. Therefore, the global optimal solution can be computed as $\bar{\x} = \pm \frac{\sqrt{\lam_1} }{\alpha} \u_1$. Finally, using the fact that $\C\u_1 = \lam_1 \u_1$, we have the following results:
\beq
\F(\bar{\x}) &=& \frac{\alpha}{2}\|\bar{\x}\|_2^2-\sqrt{\bar{\x}^T\C\bar{\x}}\nn\\
&=& \frac{\alpha}{2} \frac{\lam_1}{\alpha^2} \u_1^T \u_1 - \frac{\sqrt{\lam_1}}{\alpha} \sqrt{\u_1^T\C\u_1} \nn\\
&=&  \frac{\alpha}{2} \frac{\lam_1}{\alpha^2} \u_1^T \u_1 - \frac{\sqrt{\lam_1}}{\alpha} \sqrt{\u_1^T \lam_1 \u_1} \nn\\
&=&  -\frac{\lam_1}{2\alpha}. \nn
\eeq

\end{proof}

\end{theorem}

\noi To simplify our discussions, we only focus on $\alpha=1$ for Problem (\ref{eq:pca111}) in the sequel.

The following lemma is useful in our analysis.

\begin{lemma} \label{lemma:cubic:111}
Assume that $0<\varphi<1$. We define
\beq
\K(\varphi) \triangleq  -1 + {\sqrt[{3}]{{\frac {1}{1-\varphi}}+    \bar{\tau}  }}+{\sqrt[{3}]{\frac {1}{1-\varphi}- \bar{\tau}}},~\bar{\tau} = {\sqrt {{\frac {1}{(1-\varphi)^2}}+ \frac{1}{27 (1-\varphi)^3 } }}
\eeq
We have the following result:
\beq
0 \leq \vartheta \leq \K(\varphi) ~~~~ \Rightarrow ~~~~1-\vartheta - (1+\vartheta)^3 (1-\varphi) \geq 0.\nn
\eeq

\begin{proof}

We focus on the following equation:
\beq
1-\vartheta - (1+\vartheta)^3 (1-\varphi) = 0 \nn
\eeq
\noi Dividing both sides by $1-\varphi$, we have the following equivalent equation:
\beq
 (\underbrace{1+\vartheta}_{x})^3 + \underbrace{\frac{1}{(1-\varphi)}}_{p} \cdot (\underbrace{1+\vartheta}_{x}) + \underbrace{\frac{-2}{(1-\varphi)}}_{q}  = 0 \label{eq:relation}
\eeq
\noi Solving the depressed cubic equation $x^3 + p x   + q = 0$ above using the well-known Cardano's formula \footnote{\url{https://en.wikipedia.org/wiki/Cubic_equation}}, we obtain the unique root $\bar{x}$ with
\beq
\bar{x} = {\displaystyle {\sqrt[{3}]{-{\frac {q}{2}}+    \bar{\tau}  }}+{\sqrt[{3}]{-{\frac {q}{2}}- \bar{\tau}}}},~\bar{\tau} = {\sqrt {{\frac {q^{2}}{4}}+{\frac {p^{3}}{27}}}} \nn
\eeq
\noi Using the relations $p=\frac{1}{1-\varphi}$ and $q=-\frac{2}{1-\varphi}$ as in (\ref{eq:relation}), we have:
\beq
\bar{x} =  {\sqrt[{3}]{-{\frac {(-\frac{2}{1-\varphi})}{2}}+    \bar{\tau}  }}+{\sqrt[{3}]{-{\frac {(-\frac{2}{1-\varphi})}{2}}- \bar{\tau}}},~\bar{\tau} = {\sqrt {{\frac {(-\frac{2}{1-\varphi})^{2}}{4}}+{\frac { (\frac{1}{1-\varphi})^{3}}{27}}}} \nn
\eeq

\noi Therefore, $\bar{\vartheta} = \bar{x}-1 = \K(\varphi)$ is the unique root for $1-\vartheta - (1+\vartheta)^3 (1-\varphi) = 0$. Hereby, we finish the proof of this lemma.


\end{proof}

\end{lemma}

\begin{theorem} \label{the:local:smooth:convex}
We define $\delta \triangleq 1-\frac{\lam_2}{\lam_1},\xi \triangleq \frac{\lam_1}{6}\left( -1 - \frac{3}{\sqrt{\lam_1}} + \sqrt{ (1 + \frac{3}{\sqrt{\lam_1}} )^2 + \frac{12}{\lam_1} \delta}\right)$. Assume that $0<\delta<1$. When $\x$ is sufficiently close to the global optimal solution $\bar{\x}$ such that $\|\x - \bar{\x}  \| \leq \varpi$ with $\varpi <  \bar{\varpi} \triangleq \min\{ \sqrt{\lam_1} \mathcal{K}(\frac{{\lam}_2}{{\lam}_1}),   \xi  \} $, we have:

\textbf{(a)} $\sqrt{\lam_1}-\varpi   \leq \|\x\| \leq  \sqrt{\lam_1} + \varpi$.

\textbf{(b)} $\lam_1 - \varpi  \sqrt{\lam_1} \leq \sqrt{\x^T\C\x} \leq  \lam_1 + \varpi \sqrt{\lam_1}$.

\textbf{(c)} $\lam_1 \u_1 \u_1^T + \rho \textbf{I} \succeq \x\x^T \succeq \lam_1 \u_1 \u_1^T - \rho \textbf{I}$ with $\rho \triangleq    3  \varpi^2  + 2\varpi\sqrt{\lam_1} $.

\textbf{(d)} $ \tau \textbf{I} \succeq \nabla^2 \F(\x)\succeq \sigma \textbf{I}$ with $\sigma \triangleq 1 - \frac{\lam_2}{\lam_1}  - \varpi(1+\frac{3}{\sqrt{\lam_1}})  -  \frac{ 3   \varpi^2     }{\lam_1} >0$ and $\tau \triangleq 1 +  \frac{ \lam_1^2 (\sqrt{\lam_1} + \varpi)^2 }{ (\lam_1 - \varpi  \sqrt{\lam_1} )^3} $.

\begin{proof}

\noi \textbf{(a)} We have the following results: $\|\x\| \leq \|\x - \bar{\x}\| + \|\bar{\x}\| = \|\x - \bar{\x}\| + \sqrt{\lam_1}  \leq   \varpi  + \sqrt{\lam_1}$. Moreover, we have: $\|\x\| \geq \|\bar{\x}\| - \|\x - \bar{\x}\| =  \sqrt{\lam_1} - \varpi $.

\noi \textbf{(b)} We note that the matrix norm is defined as: $\|\x\|_{\C}  \triangleq \sqrt{\x^T\C\x}$. It satisfies the triangle inequality since $\C\succeq \mathbf{0}$. We notice that $\|\bar{\x}\| = \sqrt{\lam_1}$, $\|\bar{\x}\|_{\C} = \lam_1$. We define $\Delta \triangleq \x - \bar{\x}$. We have the following results: $\|\x\|_{\C}\leq \|\x - \bar{\x}\|_{\C} + \|\bar{\x}\|_{\C} = \|\x - \bar{\x}\|_{\C} + \lam_1
\leq   \varpi \sqrt{\lam_1}  + \lam_1$. Moreover, we have: $\|\x\|_{\C} \geq \|\bar{\x}\|_{\C} - \|\x - \bar{\x}\|_{\C} =  \lam_1 - \varpi \sqrt{\lam_1}$.


\noi \textbf{(c)} We have the following inequalities:
\beq
\x\x^T &=& \lam_1 \u_1 \u_1^T + \Delta \x^T + \x\Delta^T - \Delta \Delta^T \nn\\
&\succeq& \lam_1 \u_1 \u_1^T - \|\Delta\|\|\x\|- \|\x\|\|\Delta\| - \|\Delta\|\|\Delta\| \nn\\
&\succeq& \lam_1 \u_1 \u_1^T -  \varpi (\sqrt{\lam_1} + \varpi)- (\sqrt{\lam_1} + \varpi)\varpi - \varpi^2 \nn\\
&=& \lam_1 \u_1 \u_1^T - \rho \I \nn
\eeq
\noi Using similar strategies, we have: $\x\x^T  \preceq (\lam_1 \u_1 \u_1^T + \|\Delta\| \|\x\| + \|\x\|\|\Delta\| +\|\Delta\|\|\Delta\| ) \I$.


\noi \textbf{(d)} First, we prove that $\nabla^2 \F(\x)\succeq \sigma \textbf{I}$. We define $\gamma \triangleq \lam_1 - \varpi \sqrt{\lam_1}$ and $\eta \triangleq (\lam_1 + \varpi \sqrt{\lam_1})^3$. Noticing $\varpi < \sqrt{\lam_1} \K( \frac{\lam_2}{\lam_1})$, we invoke Lemma \ref{lemma:cubic:111} with $\vartheta  \triangleq \frac{\varpi}{\sqrt{\lam_1}}$ and $\varphi \triangleq \frac{\lam_2}{\lam_1}$ and obtain:
\beq
&&1 - \vartheta        \geq  (1 + \vartheta)^3 (1-{\lam}_2 /\lam_1) \nn\\
&\Leftrightarrow& \frac{1 - \vartheta}{(1-\vartheta)(1 + \vartheta)^3}        \geq  \frac{(1 + \vartheta)^3 (1-{\lam}_2 /\lam_1)}{ (1-\vartheta)(1 + \vartheta)^3 } \nn\\
&\Leftrightarrow & -  \frac{ 1 }{1 - \vartheta}   +  \frac{  1  }{(1 + \vartheta)^3}      \geq - \frac{ {\lam}_2 /\lam_1 }{1 - \vartheta} \nn\\
&\Leftrightarrow & -  \frac{ {\lam}_1 }{\lam_1 - \varpi \sqrt{\lam_1}}   +  \frac{  \lam_1^3  }{(\lam_1 + \varpi \sqrt{\lam_1})^3}      \geq - \frac{ {\lam}_2}{\lam_1 - \varpi \sqrt{\lam_1}}\nn\\
&\Leftrightarrow &   \frac{  \lam_1^3  }{\eta}  -  \frac{ {\lam}_1 }{\gamma}       \geq - \frac{ {\lam}_2}{\gamma} \label{eq:pca:2222}
\eeq

We now prove that $\varpi \leq 1-\frac{\lam_2}{\lam_1}$. We have the following inequalities: $1 < c \Rightarrow  \frac{12}{\lam_1} \delta < (\frac{6}{\lam_1}\delta)^2 +   \frac{12 c}{\lam_1}\delta \Rightarrow  c^2 + \frac{12}{\lam_1} \delta < c^2 + (\frac{6}{\lam_1}\delta)^2 +   \frac{12 c}{\lam_1}\delta \Rightarrow  c^2 + \frac{12}{\lam_1} \delta < (c + \frac{6}{\lam_1}\delta)^2  \Rightarrow \sqrt{ c^2 + \frac{12}{\lam_1} \delta} < c + \frac{6}{\lam_1}\delta \Rightarrow  \left( -c + \sqrt{ c^2 + \frac{12}{\lam_1} \delta}\right) < \frac{6}{\lam_1}\delta  \Rightarrow \frac{\lam_1}{6}\left( -c + \sqrt{ c^2 + \frac{12}{\lam_1} \delta}\right) < \delta$. Applying this inequality with $c \triangleq 1 + \frac{3}{\sqrt{\lam_1}}$, we have:
\beq
\delta > \frac{\lam_1}{6}\left( -1 - \frac{3}{\sqrt{\lam_1}} + \sqrt{ (1 + \frac{3}{\sqrt{\lam_1}} )^2 + \frac{12}{\lam_1} \delta}\right)   \overset{(a)}{=}  \xi  \overset{(b)}{>}  \varpi \nn
\eeq
\noi step $(a)$ uses the definition of $\xi$; step $(b)$ uses $\varpi <  \bar{\varpi} \triangleq \min\{ \sqrt{\lam_1} \mathcal{K}(\frac{{\lam}_2}{{\lam}_1}), \xi  \} $. Using the definition of $\delta \triangleq 1-\frac{\lam_2}{\lam_1}$, we  conclude that:
\beq \label{eq:imply:varpi}
\varpi \leq \delta = 1-\frac{\lam_2}{\lam_1}.
\eeq

We naturally obtain the following inequalities:
\beq \label{eq:jjjjjjjj}
\nabla^2 \F(\x) &\overset{(a)}{=}~&  \textbf{I} - \frac{\C}{ \sqrt{\x^T\C\x}  } + \frac{ \C\x\x^T \C}{ \x^T\C\x \cdot \sqrt{\x^T\C\x}} \nn\\
& \overset{(b)}{\succeq} ~& \I-  \frac{\C}{\gamma}   +  \frac{\C \x\x^T\C}{\eta}  \nn\\
& \overset{(c)}{\succeq} ~& \I -  \frac{\C}{\gamma}   +  \frac{\C [ \lam_1 \u_1 \u_1^T - \rho \textbf{I} ]\C}{\eta}  \nn\\
&\overset{(d)}{=} ~&\I -  \frac{\C}{\gamma}   +  \frac{  \lam_1^3 \u_1 \u_1^T }{\eta}  -   \frac{ \rho \C \C}{\eta}  \nn\\
&  \overset{(e)}{\succeq} ~&   (1-\frac{\rho \lam_1^2}{\eta}) \I  -  \frac{ \C}{\gamma}   +  \frac{  \lam_1^3 \u_1 \u_1^T }{\eta}    \nn\\
&  \overset{(f)}{\succeq} ~&   (1-\frac{\rho \lam_1^2}{\eta}) \I   -  \frac{\sum_{i=2}^n {\lam}_2\u_i \u_i^T}{\gamma}    +  \left( \frac{  \lam_1^3 {\u}_1 {\u}_1^T }{\eta}     -  \frac{ {\lam}_1 \u_1 \u_1^T}{\gamma}  \right)  \nn\\
&  \overset{(g)}{\succeq} ~&   (1-\frac{\rho \lam_1^2}{\eta}) \I   -  \frac{\sum_{i=2}^n {\lam}_2\u_i \u_i^T}{\gamma}    - \frac{\lam_2}{\gamma} \cdot \left( \u_1 \u_1^T \right)  \nn\\
&  \overset{(h)}{=} ~&  (1    -   \frac{ \rho{\lam}^2_1}{\eta}  -  \frac{ {\lam}_2}{\gamma}  ) \cdot \I  \nn \\
&  \overset{(i)}{=} ~&  (1  -  \frac{ {\lam}_2}{\lam_1 - \varpi \lam_1}    -   \frac{ 3  (  \varpi^2  + \varpi\sqrt{\lam_1} ) {\lam}^2_1}{(\lam_1 + \varpi \lam_1)^3} ) \cdot \I \nn\\
& \overset{(j)}{\succeq} & (1  -  \frac{ \lam_2   }{\lam_1(1 - \varpi)}    -   \frac{ 3  (  \varpi^2  + \varpi\sqrt{\lam_1} )     }{\lam_1}) \cdot \I \nn\\
& \overset{(k)}{\succeq} & (1 - \frac{\lam_2}{\lam_1}  - \varpi  -  \frac{ 3  (  \varpi^2  + \varpi\sqrt{\lam_1} )     }{\lam_1} ) \cdot \I \nn\\
& \overset{}{=} &(1 - \frac{\lam_2}{\lam_1}  - \varpi(1+\frac{3}{\sqrt{\lam_1}})  -  \frac{ 3   \varpi^2     }{\lam_1}   ) \cdot \I \triangleq \sigma \cdot \I
\eeq
\noi where step $(a)$ uses (\ref{eq:pca:hess}) with $\alpha=1$; step $(b)$ uses $\sqrt{\x^T\C\x} \geq \lam_1 - \omega\lam_1$; step $(c)$ uses $\x\x^T \succeq \lam_1 \u_1 \u_1^T - 3  (  \varpi^2  + \varpi\sqrt{\lam_1} ) \I$; step $(d)$ uses the fact that $\C\u\u_1^T \C= \lam_1^2 \u_1\u_1^T$; step $(e)$ uses $-\C\C \succeq -\lam_1^2 \I$; step $(f)$ uses $\C= \sum_{i=1}^n \lam_i\u_i \u_i^T$ and $-\lam_i \geq -\lam_2,~\forall i\geq 2$; step $(g)$ uses the conclusion as in (\ref{eq:pca:2222}); step $(h)$ uses $\I = \sum_{i=1}^n \u_i \u_i^T$; step $(i)$ uses the definition of $\rho,~\eta$ and $\gamma$; step $(j)$ uses the fact that $-\frac{1}{ ( \lam_1 + \varpi  \lam_1)^3} \geq -\frac{1}{ \lam_1^3}$; step $(k)$ uses the result in (\ref{eq:imply:varpi}) and the follow derivations:
\beq
&& \varpi \leq 1-\frac{\lam_2}{\lam_1}  \nn \\
&\Leftrightarrow&  0 \leq \varpi(\lam_1 - \lam_2 - \lam_1 \varpi) \nn\\
&\Leftrightarrow & \frac{\lam_2}{\lam_1(1-\varpi)}   \leq  \frac{\lam_2}{\lam_1} +  \varpi \nn\\
&\Leftrightarrow & -\frac{\lam_2}{\lam_1(1-\varpi)}   \geq  -\frac{\lam_2}{\lam_1} -  \varpi \nn
\eeq

\noi Finally, solving the quadratic equation $ \sigma \triangleq 1 - \frac{\lam_2}{\lam_1}  - \varpi(1+\frac{3}{\sqrt{\lam_1}})  -  \frac{ 3   \varpi^2     }{\lam_1}=0$ yields the positive root $\xi$. We conclude that when $\varpi<\xi$, we have $\sigma>0$.

\noi We now prove that $\nabla^2 \F(\x) \preccurlyeq \tau \textbf{I}$. We have the following results:
\beq
\nabla^2 \F(\x) &\overset{(a)}{=}~&  \textbf{I} - \frac{\C}{ \sqrt{\x^T\C\x}  } + \frac{ \C\x\x^T \C}{ \x^T\C\x \cdot \sqrt{\x^T\C\x}} \nn\\
& \overset{(b)}{\preceq } ~& \I + \mathbf{0} + \frac{ \lam_1^2 \|\x\|_2^2 }{ \|\x\|_{\C}^3   }  \cdot \I \nn\\
& \overset{(c)}{\preceq } ~& (1 +  \frac{ \lam_1^2 (\sqrt{\lam_1} + \varpi)^2 }{ (\lam_1 - \varpi  \sqrt{\lam_1} )^3} )\I \triangleq \tau \I \nn
\eeq
\noi where step $(a)$ uses (\ref{eq:pca:hess}) with $\alpha=1$; step $(b)$ uses $-\C \preceq \mathbf{0}$ and $\C\C \preceq \lam_1^2\I$; step $(c)$ uses $\|\x\| \leq  \sqrt{\lam_1} + \varpi$ and $\|\x\|_{\C} \leq  \lam_1 + \varpi \sqrt{\lam_1}$.

We hereby finish the proof of this theorem.

\end{proof}

\end{theorem}

\noi \textbf{Remarks}. \textit{\textbf{(i)}} The assumption $0<\delta \triangleq 1-\frac{\lam_2}{\lam_1}<1$ is equivalent to $\lam_1>\lam_2>0$, which is mild. \textit{\textbf{(ii)}} Problem (\ref{eq:pca111}) boils down to a smooth and strongly convex optimization problem under some conditions.


\textit{\textbf{CD-SNCA}} with $\theta=0$ essentially reduces to a standard CD method applied to solve a strongly convex smooth problem. Using Theorem 1 of \cite{lu2015complexity}, one can establishes its linear convergence rate.

\begin{theorem} \label{the:pca:local:conv:rate}
\textbf{(Convergence Rate of \textit{\textbf{CD-SNCA}} for the PCA Problem}). We assume that the random-coordinate selection rule is used. Assume that $\|\x^t - \bar{\x}\|\leq \bar{\varpi}$ that $\F(\cdot)$ is $\sigma$-strongly convex and $\tau$-smooth. Here the parameters $\bar{\varpi}, \sigma$ and $\tau$ are define in Theorem \ref{the:local:smooth:convex}. We define $r_t^2 \triangleq  \frac{(1+\sigma)\tau }{2}\|\x-\bar{\x}\|_2^2$ and $\beta \triangleq \frac{2 \sigma}{1+\sigma}$. We have:
\beq
\E_{\xi^{t-1}}[ r_t^2]  \leq (1-\frac{\beta}{n})^{t+1} \left(  r_0^2 + \F(\x^0) - \F(\bar{\x}) \right) \nn
\eeq

\end{theorem}

\noi \textbf{Remarks}. Note that Theorem \ref{the:pca:local:conv:rate} does not rely on neither the \textit{globally bounded nonconvexity assumption} nor the \textit{Luo-Tseng error bound assumption} of $\F(\cdot)$.

\end{document}